\documentclass{article}
\usepackage{amssymb,amsmath,multicol,multirow,amsthm}
\usepackage{eepic}
\usepackage{epsfig,picinpar}

\allinethickness{0.4mm}

\theoremstyle{plain}
\newtheorem{prop}{Proposition}
\newtheorem{proposition}{Proposition}
\newtheorem{lemma}{Lemma}
\newtheorem{theorem}{Theorem}

\theoremstyle{definition}

\newtheorem{definition}{Definition}

\newcommand{\Aut}{\mathop{\rm Aut}\nolimits}

\newcommand{\A}{{\mathcal A}}

\newcommand{\Z}{\mathbb Z}
\newcommand\C{{\mathbb C}}
\newcommand\B{{\mathcal{B}}}
\newcommand{\G}{\mathcal{G}}
\newcommand\s{{\sigma}}

\newcommand{\St}{\mathop{\rm Stab}\nolimits}
\newcommand{\Stg}{\mathop{\rm Stab}\nolimits_G}
\newcommand{\Isom}{\mathop{\rm Isom}\nolimits}
\newcommand\Sym{{\mathsf{Sym}}}
\newcommand{\lims}[1][G]{\mathcal{J}_{#1}}

\title{Classification of groups generated by 3-state automata over a 2-letter
alphabet}

\author{Ievgen Bondarenko, Rostislav Grigorchuk, Rostyslav Kravchenko,\\
Yevgen Muntyan, Volodymyr Nekrashevych,\\
Dmytro Savchuk and Zoran \v{S}uni\'{c}
\thanks{All authors were partially supported by at least one of the NSF grants
DMS-308985, DMS-456185, DMS-600975, and DMS-605019} }

\begin{document}

\maketitle

\section{Introduction}

Automaton groups were formally introduced in the beginning of
1960's~\cite{glushkov:ata,horejs:automata} but it took a while to
realize their importance, utility, and, at the same time,
complexity. Among the publications from the first decade of the
study of automaton groups let us
distinguish~\cite{zarovnyi:group,zarovnyi:wreath} and the book
\cite{gecseg-p:b-automata}.

The first substantial results came only in the 1970's and in the
beginning of the 1980's when it was shown
in~\cite{aleshin:burnside,sushch:burnside,grigorch:burnside,gupta_s:burnside}
that automaton groups provide examples of finitely generated
infinite torsion groups, thus making a contribution to one of the
most famous problems in algebra --- the General Burnside Problem
(more information on all three versions of the problem can be found
in~\cite{adian:b-burnside,golod:icm-burnside,gupta:amm-burnside,kostrikin:around,
zelmanov:icm-restricted,grigorchuk-l:burnside}). The methods used to
study the properties of the examples
from~\cite{aleshin:burnside,sushch:burnside,grigorch:burnside} are
very different. The methods used in~\cite{aleshin:burnside} are
typical for the theory of finite automata (in fact the provided
proof was incorrect; the first correct proof appears
in~\cite{merzlyakov:periodic} as a combination of the results
from~\cite{grigorch:burnside} and~\cite{merzlyakov:periodic}, as
well as in the third edition of the
book~\cite{kargapolov-m:book-with-appendix} and
in~\cite{kudryavtsev-a-p:avtomata-book}). The exposition
in~\cite{sushch:burnside} is based on Kalujnin's tableaux coming
from his theory of iterated wreath products of cyclic groups of
prime order $p$. The approach in~\cite{grigorch:burnside} is based
on the ideas of self-similarity and contraction. These ideas are
apparent both in the proof of the infiniteness and the torsion
property of the group. The self-similarity is apparent from the fact
that the set of all states of the automaton is used as a generating
set for the group (now it is common to call such groups
self-similar). The contraction property here means that the length
of the elements contracts by a factor bounded away from 1 when one
passes to sections. A principal tool introduced in the beginning of
the 1980's was the language of actions on rooted trees suggested by
Gupta and Sidki in~\cite{gupta_s:burnside}, which helped
tremendously in bringing geometric insight to the subject.

A new indication of the importance of automaton groups came when it
was shown that some of them provided the first examples of groups of
intermediate
growth~\cite{grigorch:milnor,grigorch:degrees,grigorch:degrees85}.
This not only answered the question of
J.Milnor~\cite{milnor:problem} about existence of such groups, but
also answered a number of other questions in and around group
theory, including M.~Day's problem~\cite{day:amenable} on existence
of amenable but not elementary amenable groups. Basically, even to
this day, all known examples of groups of intermediate growth and
non-elementary amenable groups are based on automaton groups.

Investigations in the last two
decades~\cite{grigorch:degrees,grigorch:degrees85,gupta_s:burnside,gupta_s:pgroups,
lysionok:presentation,neumann:pride,sidki:subgroups,sidki:presentation,grigorch:hilbert,rozhkov:centralizers,
grigorch:example,gr99:schur,grigorch:jibranch,bartholdi_g:spectrum,bartholdi-g:lie,
grigorch_z:Lamplighter,nekrash:self-similar,grigorchuk-s:hanoi-cr}
show that many automaton groups possess numerous interesting, and
sometimes unusual, properties. This includes just infiniteness (the
groups constructed in~\cite{grigorch:degrees,grigorch:degrees85} as
well as in~\cite{gupta_s:pgroups} answer a question
from~\cite{chandler_m:history} on new examples of infinite groups
with finite quotients), finiteness of width, or more generally
polynomial growth of the dimension of the successive quotients in
the lower central series~\cite{bartholdi-g:lie} (answering a
question of E.~Zelmanov on classification of groups of finite
width), branch
properties~\cite{grigorch:degrees,neumann:pride,grigorch:jibranch}
(answering some questions of S.~Pride and
M.~Edjvet~\cite{pride:largeness,edjvet-p:largeness}), finiteness of
the index of maximal subgroups and presence or absence of the
congruence property~\cite{pervova:dense,pervova:congruence} (related
to topics in pro-finite groups), existence of groups with
exponential but not uniformly exponential
growth~\cite{wilson:nonuniform,wilson:further,bartholdi:nonuniform,nekrashevych:nonuniform}
(providing an answer to a question of M.~Gromov), subgroup
separability and conjugacy separability~\cite{grigorch_w:conjugacy},
further examples of amenable groups but not amenable (or even
sub-exponentially amenable)
groups~\cite{grigorch_z:basilica,bartholdi_v:amenab,grigorchuk-n-s:oberwolfach2},
amenability of groups generated by bounded
automata~\cite{bkn:amenab}, and so on. The word problem can be
solved in contracting self-similar groups by using an extremely
effective \emph{branch
algorithm}~\cite{grigorch:degrees,savchuk:wp}. The conjugacy problem
can also be solved in many
cases~\cite{wilson-z:conjugacy,rozhkov:conjugacy,leonov:conjugacy,grigorch_w:conjugacy}
(in fact we do not know of an example of an automaton group with
unsolvable conjugacy problem). In some instances, it is even known
that the membership problem is
solvable~\cite{grigorch_w:structural}.

In addition to the formulation of many algebraic properties of
groups generated by finite automata, a number of links and
applications were discovered during the last decade. This includes
asymptotic and spectral properties of the Cayley graphs and Schreier
graphs associated to the action on the rooted tree with respect to
the set of generators given by the set of states of the automaton.
For instance, it is shown in~\cite{grigorch_z:Lamplighter} that the
discrete Laplacian on the Cayley graph of the Lamplighter group $\Z
\ltimes (\mathbb{Z}/2\mathbb{Z})^\mathbb{Z}$ has pure point
spectrum. This fact was used to answer a question of M.~Atiyah on
$L^2$-Betti numbers of closed
manifolds~\cite{grigorchuk-lsz:atiyah}. The methods developed in the
study of the spectral properties of Schreier graphs of self-similar
groups can be used to construct Laplacians on fractal sets and to
study their spectral properties
(see~\cite{grigorchuk-n:schur,nekrashevych-t:analysis}.

A new and fruitful direction, bringing further applications of
self-similar groups, was established by the introduction of the
notions of iterated monodromy groups and limit spaces by
V.~Nekrashevych. The theory established a link between contracting
self-similar groups and the geometry of Julia sets of expanding
maps. An example of an application of self-similar groups to
holomorphic dynamics is given by the solution (by L.~Bartholdi and
V.~Nekrashevych in~\cite{bartholdi_n:rabbit}) of the ``twisted
rabbit'' problem of J.~Hubbard. The book~\cite{nekrash:self-similar}
provides a comprehensive introduction to this theory.

In many situations automaton groups serve as renorm groups. For
instance this happens in the study of classical fractals, in the
study of the behavior of dynamical systems~\cite{oliva:phd}, and in
combinatorics --- for example in Hanoi Towers games on $k$ pegs, $k
\geq 3$, as observed by Z.~\v{S}uni\'c
(see~\cite{grigorchuk-s:hanoi-cr}).

There is interest of computer scientists and logicians in automaton
groups, since they may be relevant in the solution of important
complexity problems (see~\cite{rhodes-s:bimachines} for ideas in
this direction). Self-similar groups of intermediate growth are
mentioned by Wolfram in~\cite{wolfram:nks} as examples of ``multiway
systems'' with complex behavior.

Among the major problems in many areas of mathematics are the
classification problems. If the objects are given combinatorially
then it is naturally to try to classify them first by complexity and
then within each complexity class.

A natural complexity parameter in our situation is the pair $(m, n)$
where $m$ is the number of states of the automaton generating the
group and $n$ is the cardinality of the alphabet.

There are 64 invertible 2-state automata acting on a 2-letter
alphabet, but there are only six non-isomorphic $(2, 2)$-automaton
groups, namely, the trivial group, $\Z/2\Z$, $\Z/2\Z\oplus\Z/2\Z$,
$\Z$, the infinite dihedral group $D_\infty$, and the lamplighter
group $\Z \wr \Z/2\Z$~\cite{gns00:automata,grigorch_z:Lamplighter}
(more details are given in Theorem~\ref{thm:class22} below). A
classification of semigroups  generated  by 2-state automata (not
necessary invertible) over a 2-letter alphabet is provided by
I.~Reznikov and V.~Sushchanski{\u\i}~\cite{reznykov_s:growth2x2}.
Some examples from this class, including an automaton generating a
semigroup of intermediate growth, were studied in the subsequent
papers~\cite{reznykov_s:fibonacci,reznykov_s:interm_growth,bartholdi_rs:interm_growth}.

It is not known how many pairwise non-isomorphic groups exists for
any class $(m, n)$ when either $m
> 2$ or $n>2$. Unfortunately, the number of automata that has to be treated grows
super-exponentially with either of the two arguments (there are
$m^{mn}(n!)^m$ invertible $(m,n)$-automata).

Nevertheless, a reasonable task is to consider the problem of
classification for small values of $m$ and $n$ and try to classify
the $(3, 2)$-automaton groups and $(2, 3)$-automaton groups.

Our research group (with some contribution by Y.~Vorobets and
M.~Vorobets) has been working on the problem of classification of
$(3,2)$-automaton groups for the last four yeas and most of the
obtained results are presented in this article. Our research goals
moved in three main directions:

1. Search for new interesting groups and an attempt to use them to
solve known problems. An example of such a group is the Basilica
group (see automaton [852]). It is the first example of an amenable
group (shown in~\cite{bartholdi_v:amenab}) that is not
sub-exponentially amenable group (shown
in~\cite{grigorch_z:basilica}).

2. Recognition of already known groups as self-similar groups, and
use of the self-similar structure in finding new results and
applications for such groups. As examples we can mention the free
group of rank 3 (see automaton [2240]), the free product of three
copies of $\mathbb{Z}/2\mathbb{Z}$ (see automaton [846]),
Baumslag-Solitar groups $BS(1,\pm3)$ (see automata [870] and
[2294]), the Klein bottle group (see automaton [2212]), and the
group of orientation preserving automorphisms of the 2-dimensional
integer lattice (see automaton [2229]).

3. Understanding of typical phenomena that occur for various classes
of automaton groups, formulation and proofs of reasonable
conjectures about the structure of self-similar groups.

The main general results on the class of groups generated by
$(3,2)$-automata are as follows.

\begin{theorem}
There are at most 122 non-isomorphic groups generated by
$(3,2)$-automata.
\end{theorem}

The numbers in brackets in the next two theorems are references to
the numbers of the corresponding automata (more on this encoding
will be said later). Here and thereafter, $C_n$ denotes  the cyclic
group of order $n$.

\begin{theorem}
There are $6$ finite groups in the class: the trivial group $\{1\}$
\textup{[1]}, $C_2$ \textup{[1090]}, $C_2\times C_2$ \textup{[730]},
$D_4$ \textup{[847]}, $C_2\times C_2\times C_2$ \textup{[802]} and
$D_4\times C_2$ \textup{[748]}.
\end{theorem}

\begin{theorem}
There are $6$ abelian groups in the class: the trivial group $\{1\}$
\textup{[1]}, $C_2$ \textup{[1090]}, $C_2\times C_2$ \textup{[730]},
$C_2\times C_2\times C_2$ \textup{[802]}, $\mathbb Z$ \textup{[731]}
and $\mathbb Z^2$ \textup{[771]}.
\end{theorem}

\begin{theorem}\label{thm:freegroup}
The only nonabelian free group in the class is the free group of
rank $3$ generated by the Aleshin-Vorobets-Vorobets automaton
\textup{[2240]}.
\end{theorem}

\begin{theorem}
There are no infinite torsion groups in the class.
\end{theorem}

The short list of general results  does not give full justice to the
work that has been done. Namely, in most individual cases we have
provided a lot of results and detailed information for the group in
question. The variety is rather extreme and it is not surprising at
all that one cannot formulate too many general results.

More work and, likely, some new invariants are required to further
distinguish the 122 groups that are listed in this paper as
potentially non-isomorphic. In some cases one could try to use the
rigidity of actions on rooted trees
(see~\cite{lavrenyuk_n:rigidity}), since in many cases it is easier
to distinguish actions than groups. In the contracting case one
could use, for instance, the geometry of the Schreier graphs and
limit spaces to distinguish the actions.

Next natural step would be to consider the case of $(2,
3)$-automaton groups or 2-generated self-similar groups of binary
tree automorphisms defined by recursions in which every section is
either trivial, a generator, or an inverse of a generator. The cases
$(4, 2)$ and $(5,2)$ also seem to be attractive, as there are many
remarkable groups in these classes.

Another possible direction is to study more carefully only certain
classes of automata (such as the classical linear automata, bounded
and polynomially growing automata in the sense of
Sidki~\cite{sidki:acyclic}, etc.) and the properties of the
corresponding automaton groups.

Many computations used in our work were performed by the package
\texttt{AutomGrp} for \texttt{GAP} system, developed by Y.~Muntyan
and D.~Savchuk~\cite{muntyan_s:automgrp}. The package is not
specific to $(3,2)$-automaton groups (in fact, many functions are
implemented also for groups of tree automorphisms that are not
necessarily generated by automata).

\section{Regular rooted trees, automorphisms, and self-similarity}

Let $X$ be an alphabet on $d$ ($d \geq 2$) letters. Most often we
set $X=\{0,1,\dots,d-1\}$. The set of finite words over $X$, denote
by $X^*$, has the structure of a \emph{regular rooted $d$-ary tree},
which we also denote by $X^*$. The empty word $\emptyset$ is the
\emph{root} of the tree and every vertex $v$ has $d$ children,
namely the words $vx$, for $x$ in $X$. The words of length $n$
constitute \emph{level} $n$ in the tree.

The group of tree automorphisms of $X^*$ is denoted by $\Aut(X^*)$.
Tree automorphisms are precisely the permutations of the vertices
that fix the root and preserve the levels of the tree. Every
automorphism $f$ of $X^*$ can be decomposed as
\begin{equation}\label{decomposition}
  f = \alpha_f(f_0,\dots,f_{d-1})
\end{equation}
where $f_x$, for $x$ in $X$, are automorphisms of $X^*$ and
$\alpha_f$ is a permutation of the set $X$. The permutation
$\alpha_f$ is called the \emph{root permutation} of $f$ and the
automorphisms $f_x$ (denoted also by $f|_x$), $x$ in $X$, are called
\emph{sections} of $f$. The permutation $\alpha_f$ describes the
action of $f$ on the first letter of every word, while the
automorphism $f_x$, for $x$ in $X$, describes the action of $f$ on
the tail of the words in the subtree $xX^*$, consisting of the words
in $X^*$ that start with $x$. Thus the
equality~\eqref{decomposition} describes the action of $f$ through
decomposition into two steps. In the first step the $d$-tuple
$(f_0,\dots,f_{d-1})$ acts on the $d$ subtrees hanging below the
root, and then the permutation $\alpha_f$, permutes these $d$
subtrees. Thus we have
\begin{equation} \label{fxw}
 f(xw) = \alpha_f(x) f_x(w),
\end{equation}
for $x$ in $X$ and $w$ in $X^*$. Second level sections of $f$ are
defined as the sections of the sections of $f$, i.e., $f_{xy}=
(f_x)_y$, for $x,y \in X$, and more generally, for a word $u$ in
$X^*$ and a letter $x$ in $X$ the section of $f$ at $ux$ is defined
as $f_{ux}= (f_u)_x$, while the section of $f$ at the root is $f$
itself.

The group $\Aut(X^*)$ decomposes algebraically as
\begin{equation}
\label{eqn:iter_wreath}
\begin{array}{r}
\Aut(X^*) = \Sym(X)\ltimes\Aut(X^*)^X   = \Sym(X)\wr\Aut(X^*),
\end{array}
\end{equation}
where $\wr$ is the \emph{permutational wreath product} in which the
active group $\Sym(X)$ permutes the coordinates of
$\Aut(X^*)^X=(\Aut(X^*),\dots,\Aut(X^*))$. For arbitrary
automorphisms $f$ and $g$ in $\Aut(X^*)$ we have
\[
 \alpha_f(f_0,\dots,f_{d-1}) \alpha_g(g_0,\dots,g_{d-1}) =
 \alpha_f\alpha_g (f_{g(0)}g_0,\dots,f_{g(d-1)}g_{d-1}).
\]
For future use we note the following formula regarding the sections
of a composition of tree automorphisms. For tree automorphisms $f$
and $g$ and a vertex $u$ in $X^*$,
\begin{equation}\label{chain}
 (fg)_u = f_{g(u)}g_u.
\end{equation}

The group of tree automorphisms $\Aut(X^*)$ is a pro-finite group.
Namely, $\Aut(X^*)$ has the structure of an infinitely iterated
wreath product
\[ \Aut(X^*) = \Sym(X) \wr(\Sym(X)\wr(\Sym(X)\wr\dots)) \]
of the finite group $\Sym(X^*)$ (this follows
from~\eqref{eqn:iter_wreath}). This product is the inverse limit of
the sequence of finitely iterated wreath products of the form
$\Sym(X) \wr(\Sym(X)\wr(\Sym(X)\wr\dots \wr \Sym(X)))$. Every
subgroup of $\Aut(X^*)$ is residually finite. A canonical sequence
of normal subgroups of finite index intersecting trivially is the
sequence of level stabilizers. The $n$-th \emph{level stabilizer} of
a group $G$ of tree automorphisms is the subgroup $\St_{G}(n)$ of
$\Aut(X^*)$ that consists of all tree automorphisms in $G$ that fix
the vertices in the tree $X^*$ up to and including level $n$.

The \emph{boundary} of the tree $X^*$ is the set $X^\omega$ of right
infinite words over $X$ (infinite geodesic rays in $X^*$ connecting
the root to ``infinity''). The boundary has a natural structure of a
metric space in which two infinite words are close if they agree on
long finite prefixes. More precisely, for two distinct rays $\xi$
and $\zeta$, define the distance to be $d(\xi,\zeta) = 1/2^{|\xi
\wedge \zeta|}$, where $|\xi \wedge \zeta|$ denotes the length of
the longest common prefix $\xi \wedge \zeta$ of $\xi$ and $\zeta$.
The induced topology on $X^\omega$ is the Tychonoff product topology
(with $X$ discrete), and $X^\omega$ is a Cantor set. The group of
isometries $\Isom(X^\omega)$ and the group of tree automorphisms
$\Aut(X^*)$ are canonically isomorphic. Namely, the action of the
automorphism group $\Aut(X^*)$ can be extended to an isometric
action on $X^\omega$, simply by declaring that~\eqref{decomposition}
and \eqref{fxw} are valid for right infinite words.

We now turn to the concept of self-similarity. The tree $X^*$ is a
highly self-similar object (the subtree $uX^*$ consisting of words
with prefix $u$ is canonically isomorphic to the whole tree) and we
are interested in groups of tree automorphisms in which this
self-similarity structure is reflected.

\begin{definition}
A group $G$ of tree automorphisms is \emph{self-similar} if, every
section of every automorphism in $G$ is an element of $G$.
\end{definition}

Equivalently, self-similarity can be expressed as follows. A group
$G$ of tree automorphisms is self-similar if, for every $g$ in $G$
and a letter $x$ in $X$, there exists a letter $y$ in $X$ and an
element $h$ in $G$ such that
\[ g(xw) = yh(w), \]
for all words $w$ over $X$.

Self-replicating groups constitute a special class of self-similar
groups. Examples from this class are very common in applications. A
self-similar group $G$ is \emph{self-replicating} if, for every
vertex $u$ in $X^*$, the homomorphism $\varphi_u:\Stg(u) \to G$ from
the stabilizer of the vertex $u$ in $G$ to $G$, given by
$\varphi(g)=g_u$, is surjective.

At the end of the section, let us mention the class of \emph{branch
groups}. Branch groups were introduced~\cite{grigorch:jibranch}
where it is shown that they constitute one of the three classes of
just-infinite groups (infinite groups with no proper, infinite,
homomorphic images). If a class of groups $\mathcal{C}$ is closed
under homomorphic images and if it contains infinite, finitely
generated examples then it contains just-infinite examples (this is
because every infinite, finitely generated group has a just-infinite
image). Such examples are minimal infinite examples in
$\mathcal{C}$. We note that, for example, the group of intermediate
growth constructed in~\cite{grigorch:burnside} is a branch automaton
group that is a just-infinite 2-group. i.e., it is an infinite,
finitely generated, torsion group that has no proper infinite
quotients. The Hanoi Towers group~\cite{grigorchuk-s:standrews} is a
branch group that is not just
infinite~\cite{grigorchuk-n-s:oberwolfach1}. The iterated monodromy
group $IMG(z^2+i)$~\cite{grigorch_ss:img} is a branch groups, while
$\B=IMG(z^2-1)$ is not a branch group, but only weakly branch. More
generally, it is shown in~\cite{bartholdi-n:quadratic1} that the
iterated monodromy groups of post-critically finite quadratic maps
are branch groups in the pre-periodic case and weakly branch groups
in the periodic case (the case refers to the type of post-critical
behavior).

We now define regular (weakly) branch groups. A level transitive
group $G \leq \Aut(X^*)$ of $k$-ary tree automorphisms is a
\emph{regular branch group} over $K$ if $K$ is a normal subgroup of
finite index in $G$ such that $K \times \dots \times K$ is
geometrically contained in $K$. By definition, the subgroup $K$ has
the property that $K \times \dots \times K$ is geometrically
contained in $K$, denoted by $K \times \dots \times K \preceq K$, if
\[ K \times \dots \times K \leq \psi(K \cap \St_G(1)) \]
where $\psi$ is the homomorphism $\psi: \St_G(1) \to \Aut(X^*)
\times \dots \times \Aut(X^*)$ given by $\psi(g) =
(g_0,g_1,\dots,g_{k-1})$. If instead of asking for $K$ to have
finite index in $G$ we only require that $K$ is nontrivial, we say
that $G$ is \emph{regular weakly branch group} over $K$. Note that
if $G$ is level transitive and $K$ is normal in $G$, in order to
show that $G$ is regular (weakly) branch group over $K$, it is
sufficient to show that $K \times 1 \times \dots \times 1 \preceq K$
(i.e. $K \times 1 \times \dots \times 1 \leq \psi(K \cap
\St_G(1))$). More on the class of branch group can be found
in~\cite{grigorch:jibranch} and~\cite{bar_gs:branch}.

\section{Automaton groups}\label{definition}

The full group of tree automorphisms $Aut(X^*)$ is self-similar,
since the section of every tree automorphism is just another tree
automorphism. However, this group is rather large (uncountable). For
various reasons, one may be interested in ways to define (construct)
finitely generated self-similar groups. Automaton groups constitute
a special class of finitely generated self-similar groups. We
provide two ways of thinking about automaton groups. One is through
finite wreath recursions and the other through finite automata.

Every finite system of recursive relations of the form
\begin{equation} \label{s^0}
\left\{
 \begin{array}{ccc}
 s^{(1)} &= & \alpha_1 \left(s^{(1)}_0,s^{(1)}_1,\dots,s^{(1)}_{d-1}\right),  \\
 \dots  \\
 s^{(k)} &= & \alpha_k \left(s^{(k)}_0,s^{(k)}_1,\dots,s^{(k)}_{d-1}\right),
 \end{array}
\right.
\end{equation}
where each symbol $s^{(i)}_j$, $i=1,\dots,k$, $j=0,\dots,d-1$, is a
symbol in the set of symbols $\{s^{(1)}, \dots, s^{(k)}\}$ and
$\alpha_1,\dots,\alpha_k$ are permutations in $\Sym(X)$, has a
unique solution in $\Aut(X^*)$ (in the sense that the above
recursive relations represent the decompositions of the tree
automorphisms $s^{(1)},\dots,s^{(k)}$). Thus, the action of the
automorphism defined by the symbol $s^{(i)}$ is given recursively by
$s^{(i)}(xw) = \alpha_i(x) s^{(i)}_x (w)$.

The group $G$ generated by the automorphisms $s^{(1)}, \dots,
s^{(k)}$ is a finitely generated self-similar group of automorphisms
of $X^*$. This follows since sections of products are products of
sections (see~\eqref{chain}) and all sections of the generators of
$G$ are generators of $G$.

When a self-similar group is defined by a system of the
form~\eqref{s^0}, we say that it is defined by a \emph{wreath
recursion}. We switch now the point of view from wreath recursions
to invertible automata.

\begin{definition} \label{automaton}
A \emph{finite automaton} $A$ is a $4$-tuple $\A=(S,X,\pi,\tau)$
where $S$ is a finite set of \emph{states}, $X$ is a finite
\emph{alphabet} of cardinality $d \geq 2$, $\pi:S \times X \to X$ is
a map, called \emph{output map}, and $\tau:S \times X \to S$ is a
map, called \emph{transition map}. If in addition, for each state
$s$ in $S$, the restriction $\pi_q: X \to X$ given by
$\pi_s(x)=\pi(s,x)$ is a permutation in $\Sym(X)$, the automaton
$\A$ is invertible.
\end{definition}

In fact, we will be only concerned with finite invertible automata
and, in the rest of the text, we will use the word automaton for
such automata.

Each state $s$ of the automaton $\A$ defines a tree automorphism of
$X^*$, which we also denote by $s$. By definition, the root
permutation of the automorphism $s$ (defined by the state $s$) is
the permutation $\pi_s$ and the section of $s$ at $x$ is
$\tau(s,x)$. Therefore
\begin{equation}
\label{eqn:autom_action} s(xw) = \pi_s(x) \tau(s,x)(w)
\end{equation}
for every state $s$ in $S$, letter $x$ in $X$ and word $w$ over $X$.

\begin{definition}
Given an automaton $\A=(S,X,\pi,\tau)$, the group of tree
automorphisms generated by the states of $\A$ is denoted by $G(\A)$
and called the \emph{automaton group} defined by $\A$.
\end{definition}

The generating set $S$ of the automaton group $G(\A)$ generated by
the automaton $\A=(S,X,\pi,\tau)$ is called the \emph{standard}
generating set of $G(\A)$ and plays a distinguished role.

\begin{figure}[!ht]
\begin{center}
\epsfig{file=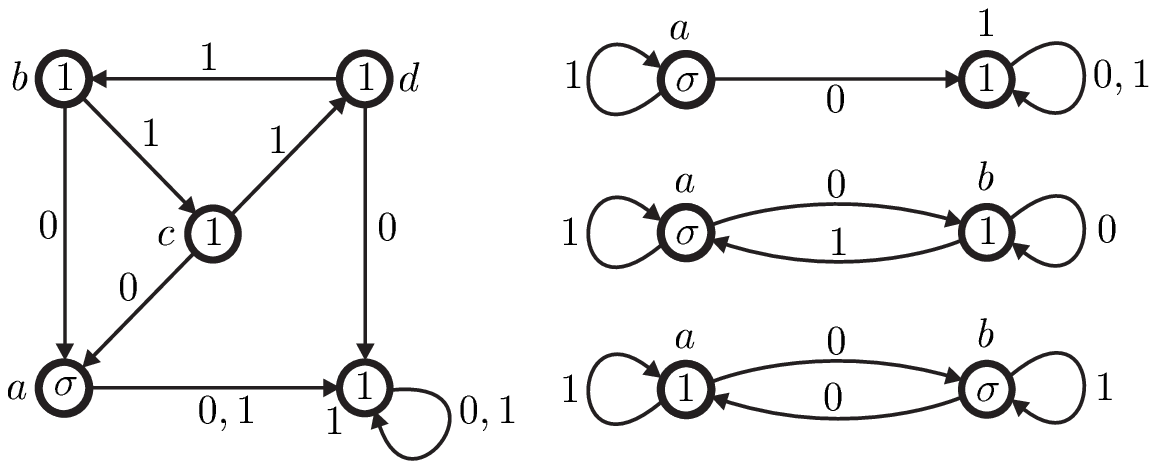,height=110pt}
\end{center}
\caption{An automaton generating $\G$, the binary adding machine,
and two Lamplighter automata}\label{4automata}
\end{figure}

Directed graphs provide convenient representation of automata. The
vertices of the graph, called \emph{Moore diagram} of the automaton
$\A=(S,X,\pi,\tau)$, are the states in $S$. Each state $s$ is
labeled by the root permutation $\alpha_s=\pi_s$ and, for each pair
$(s,x) \in S \times X$, an edge labeled by $x$ connects $s$ to
$s_x=\tau(s,x)$. Several examples are presented in
Figure~\ref{4automata}.
The states of the 5-state automaton in the left half of the figure
generate the group $\G$ of intermediate growth mentioned in the
introduction ($\sigma$ denotes the permutation exchanging $0$ and
$1$, and $1$ denotes the trivial vertex permutation). The top of the
three 2-state automata on the right in Figure~\ref{4automata} is the
so called \emph{binary adding machine}, which generates the infinite
cyclic group $\Z$. The other two automata both generate the
Lamplighter group $L_2=\Z\wr \Z/2\Z = \Z \ltimes \left(\bigoplus
\Z/2\Z \right)^\Z$ (see~\cite{gns00:automata}).

The corresponding wreath recursions for the adding machine and for
the two automata generating the Lamplighter group are given by
\begin{alignat}{6}
  &a = \sigma &&(1,a)  \qquad\qquad &&a = \sigma &&(b,a) \qquad \qquad &&a =  &&(b,a) \notag\\
  &1 =        &&(1,1)               &&b = &&(b,a),                     &&b = \sigma &&(a,b)  \label{lampl}
\end{alignat}
respectively.

The class of \emph{polynomially growing automata} was introduced by
Sidki in~\cite{sidki:acyclic}. Sidki proved in~\cite{sidki:nofree}
that no group generated by such an automaton contains free subgroups
of rank 2. As we already indicated in the introduction, for the
subclass of so called bounded automata the corresponding groups are
amenable~\cite{bkn:amenab}. Recall that an automaton $\A$ is called
\emph{bounded} if, for every state $s$ of $\A$, the function
$f_s(n)$ counting the number of active sections of $s$ at level $n$
is bounded (a state is \emph{active} if its vertex permutation is
nontrivial).

There are other classes of automata (and corresponding automaton
groups) that deserve special attention. We end the section by
mentioning several such classes.

The class of \emph{linear automata} consists of automata in which
both the set of states $S$ and the alphabet $X$ have a structure of
a vector space (over a finite field) and both the output and the
transition function are linear maps (see~\cite{gecseg-p:b-automata}
and~\cite{eilenberg:automata2}).

The class of \emph{bi-invertible automata} consists of automata in
which both the automaton and its dual are invertible. Some of the
automata in our classification are bi-invertible, most notably the
Aleshin-Vorobets-Vorobets automaton [2240] generating the free group
$F_3$ of rank 3 and the Bellaterra automaton [846] generating the
free product $C_2 \ast C_2 \ast C_2$. In fact,  both of these have
even stronger property of being \emph{fully invertible}. Namely, not
only the automaton and its dual are invertible, but also the dual of
the inverse automaton is invertible.

Another important class is the class of automata satisfying the
\emph{open set} condition. Every automaton in this class contains a
\emph{trivial state} (a state defining the trivial tree
automorphism) and this state can be reached from any other state.

One may also study automata that are \emph{strongly connected} (i.e.
automata for which the corresponding Moore diagrams are strongly
connected as directed graphs), automata in which no path contains
more than one active state (such as the automaton defining $\G$ in
Figure~\ref{4automata}), and so on.

\section{Schreier graphs}

Let $G$ be a group generated by a finite set $S$ and let $G$ act on
a set $Y$. We denote by $\Gamma=\Gamma(G,S,Y)$ the \emph{Schreier
graph} of the action of $G$ on $Y$. The vertices of $\Gamma$ are the
elements of $Y$. For every pair $(s,y)$ in $S \times Y$ an edge
labeled by $s$ connects $y$ to $s(y)$. An \emph{orbital Schreier
graph} of the action is the Schreier graph $\Gamma(G,S,y)$ of the
action of $G$ on the $G$-orbit of $y$, for some $y$ in $Y$.

Let $G$ be a group of tree automorphisms of $X^*$ generated by a
finite set $S$. The levels $X^n$, $n \geq 0$, are invariant under
the action of $G$ and we can consider the sequence of finite
Schreier graphs $\Gamma_n(G,S)=\Gamma(G,S,X^n)$, $n \geq 0$. Let
$\xi=x_1x_2x_3\ldots\in X^{\omega}$ be an infinite ray. Then the
pointed Schreier graphs $(\Gamma_n(G,S),x_1x_2\ldots x_n)$ converge
in the local topology (see~\cite{grigorch:degrees}
or~\cite{grigorchuk-z:cortona}) to the pointed orbital Schreier
graph $(\Gamma(G,S,\xi),\xi)$.

Schreier graphs may be sometimes used to compute the spectrum of
some operators related to the group. For a group of tree
automorphisms $G$ generated by a finite symmetric set $S$ there is a
natural unitary representation in the space of bounded linear
operators ${\mathcal H}=B(L_2(X^\omega))$, given by
$\pi_g(f)(x)=f(g^{-1}x)$ (the measure on the boundary $X^\omega$ is
just the product measure associated to the uniform measure on $X$).
Consider the spectrum of the operator
\[ M= \frac{1}{|S|} \sum_{s\in S} \pi_s \]
corresponding to this unitary representation. The spectrum of $M$
for a self-similar group $G$ is approximated by the spectra of the
finite dimensional operators induced by the action of $G$ on the
levels of the tree (see~\cite{bartholdi_g:spectrum}. Denote by
${\mathcal H}_n$ the subspace of ${\mathcal H}=B(L_2(X^\omega))$
spanned by the characteristic functions $f_v$, $v\in X^n$, of the
cylindrical sets corresponding to the $|X|^n$ vertices on level $n$.
The subspace ${\mathcal H}_n$ is invariant under the action of $G$
and ${\mathcal H}_n\subset {\mathcal H}_{n+1}$. Denote by
$\pi^{(n)}_g$ the restriction of $\pi_g$ on ${\mathcal H}_n$. Then,
for $n \geq 0$, the operator
\[
 M_n= \frac{1}{|S|} \sum_{s\in S} \pi_s^{(n)}
 \]
is finite dimensional. Moreover,
\[sp(M)=\overline{\bigcup_{n\geq0}sp(M_n)}, \]
i.e., the spectra of the operators $M_n$ converge to the spectrum of
$M$.

The table of information provided in Section~\ref{app_table}
includes, in each case, the histogram of the spectrum of the
operator $M_9$.

If $P$ is the stabilizer of a point on the boundary $X^\omega$, then
one can consider the quasi-regular representation $\rho_{G/P}$ of
$G$ in $\ell^2(G/P)$.

\begin{theorem}[\cite{bartholdi_g:spectrum}]
If $G$ is amenable or the Schreier graph $G/P$ (the Schreier graph
of the action of $G$ on the cosets of $P$) is amenable then the
spectrum of $M$ and the spectrum of the quasi-regular representation
$\rho_{G/P}$ coincide.
\end{theorem}

In case the parabolic subgroup $P$ is ``small'', the last result may
be used to compute the spectrum of the Markov operator on the Cayley
graph of the group. This approach was successfully used, for
instance, to compute the spectrum of the Lamplighter group
in~\cite{grigorch_z:Lamplighter} (see
also~\cite{kambites-s-s:spectra}).

\section{Contracting groups, limit spaces, and iterated monodromy
groups}\label{s:contracting}

\begin{definition}
A group $G$ generated by an automaton over alphabet $X$ is
\emph{contracting} if there exists a finite subset
$\mathcal{N}\subset G$ such that for every $g\in G$ there exists $n$
(generally depending on $g$) such that section $g_v$ belongs to
$\mathcal{N}$ for all words $v\in X^*$ of length at least $n$. The
smallest set $\mathcal{N}$ with this property is called the
\emph{nucleus} of the group $G$.
\end{definition}

The above definition makes sense for arbitrary self-similar groups
--- not necessarily automaton groups and, moreover, not necessarily
finitely generated groups. In the case of an automaton group the
contracting property may be equivalently stated as follows. An
automaton group $G=G(\A)$ is contracting if there exist constants
$\kappa$, $C$, and $N$, with $0 \leq \kappa < 1$, such that $|g_v|
\leq \kappa |g| + C$, for all vertices $v$ of length at least $N$
and $g \in G$ (the length is measured with respect to the standard
generating set $S$ consisting of the states of $\A$). The
contraction property is a key ingredient in many inductive arguments
and algorithms involving the decomposition
$g=\alpha_g(g_0,\dots,g_{d-1})$. Indeed, the contraction property
implies that, for all sufficiently long elements $g$, all sections
of $g$ at vertices on level at least $N$ are strictly shorter than
$g$.

Contracting groups have rich geometric structure. Each contracting
group is the iterated monodromy group of its \emph{limit dynamical
system}. This system is an (orbispace) self-covering of the
\emph{limit space} of the group. The limit space is a limit of the
graphs of the action of $G$ on the levels $X^n$ of the tree $X^*$
and is defined in the following way.

\begin{definition}
Let $G$ be a contracting group over $X$. Denote by $X^{-\omega}$ the
space of all left-infinite sequences $\ldots x_2x_1$ of elements of
$X$ with the direct product (Tykhonoff) topology. We say that two
sequences $\ldots x_2x_1$ and $\ldots y_2y_1$ are
\emph{asymptotically equivalent} if there exists a sequences $g_k\in
G$, assuming a finite set of values, and such that
\[g_k(x_k\ldots x_1)=y_k\ldots y_1\]
for all $k\ge 1$. The quotient of the space $X^{-\omega}$ by this
equivalence relation is called the \emph{limit space} of $G$.
\end{definition}

The following proposition, proved in~\cite{nekrash:self-similar}
(Proposition~3.6.4) is a convenient way to compute the asymptotic
equivalence.

\begin{proposition}
\label{pr:generatequiv} Let a contracting group $G$ be generated by
a finite automaton $A$. Then the asymptotic equivalence is the
equivalence relation generated by the set of pairs $(\ldots x_2x_1,
\ldots y_2y_1)$ for which there exists a sequence $g_k$ of states of
$A$ such that $g_k(x_k)=y_k$ and $g_k|_{x_k}=g_{k-1}$.
\end{proposition}

The limit dynamical system is the map induced by the shift $\ldots
x_2x_1\mapsto \ldots x_3x_2$. The limit space is a compact
metrizable topological space of finite topological dimension
(see~\cite{nekrash:self-similar}, Theorem~3.6.3). If the group is
self-replicating, then the limit space is locally connected and path
connected.

The main tool of finding the limit space of a contracting group is
realization of the group as the iterated monodromy group of an
expanding partial orbispace self-covering. An exposition of the
theory of such self-coverings is given
in~\cite{nekrash:self-similar}. In particular, if $G$ is the
iterated monodromy group of a post-critically finite complex
rational function, then the limit space of $G$ is homeomorphic to
the Julia set of the function (see Theorems~5.5.3 and~6.4.4
of~\cite{nekrash:self-similar}).

The limit space does not change when we pass from $X$ to $X^n$ in
the natural way (we will change then the limit dynamical system to
its $n$th iterate). It also does not change if we post-conjugate the
wreath recursion by an element of the wreath product $Symm(X)\ltimes
G^X$, i.e., conjugate the group $G$ by an element of the form
$\gamma=\pi(g_0\gamma, g_1\gamma)$, where $\pi\in Symm(X)$ and $g_0,
g_1\in G$.

The limit space can be also visualized using its subdivision into
\emph{tiles}. This method is especially effective, when the group is
generated by bounded automata.

\begin{definition}
\label{def:til} Let $G$ be a contracting group. A \emph{tile}
$\mathcal{T}_G$ of $G$ is the quotient of the space $X^{-\omega}$ by
the equivalence relation, which identifies two sequences $\ldots
x_2x_1$ and $\ldots y_2y_1$ if there exists a sequence $g_k\in G$
assuming a finite number of values and such that
\[g_k(x_k\ldots x_1)=y_k\ldots y_1,\quad g_k|_{x_k\ldots x_1}=1\]
for all $k$.
\end{definition}

Again, an analog of Proposition~\ref{pr:generatequiv} is true: the
equivalence relation from Definition~\ref{def:til} is generated by
the identifications $\ldots x_2x_1=\ldots y_2y_1$ of sequences for
which there exists a sequence $g_k, k=0, 1, 2,\ldots$ of elements of
the nucleus such that $g_k(x_k)=y_k$, $g_k|_{x_k}=g_{k-1}$ and
$g_0=1$.

Suppose that $G$ satisfies the \emph{open set condition}, i.e., the
trivial state can be reached from any other state of the generating
automaton. Then the \emph{boundary} of the tile $\mathcal{T}_G$ is
the image in $\mathcal{T}_G$ of the set of sequences $\ldots x_2x_1$
such that there exists a sequence $g_k\in G$ assuming a finite
number of values and such that $g_k|_{x_k\ldots x_1}\ne 1$. If $G$
is generated by a finite symmetric set $S$, then it is sufficient to
look for the sequence $g_k$ inside $S$.

The limit space of $G$ is obtained from the tile by some
identifications of the points of the boundary. If the group $G$ is
generated by bounded automata, then its boundary consists of a
finite number of points and it is not hard to identify them (i.e.,
to identify the sequences encoding them).

For $v\in X^n$ denote by $\mathcal{T}_Gv$ the image of the
cylindrical set $X^{-\omega}v$ in $\mathcal{T}_G$. It is easy to see
that the map $\ldots x_2x_1\mapsto \ldots x_2x_1v$ induces a
homeomorphism of $\mathcal{T}_G$ with $\mathcal{T}_Gv$ and that
\[\mathcal{T}_G=\bigcup_{v\in X^n}\mathcal{T}_Gv.\]
It is proved in~\cite{nekrash:self-similar} that two pieces
$\mathcal{T}_Gv_1$ and $\mathcal{T}_Gv_2$ intersect if and only if
$g(v_1)=v_2$ for an element $g$ of the nucleus of $G$ and that they
intersect only along images of the boundary of $\mathcal{T}_G$.

This suggests the following procedure of visualizing the limit space
in the case of bounded automata. Identify the points of the boundary
of the tile. We get a finite list $B$ of points, represented by a
finite list $W$ of infinite sequences (some points may be
represented by several sequences). Draw the tile as a graph with
$|B|$ ``boundary points'' (vertices) and identify the boundary
points with the points of $B$ labeled by sequences $W$. Take now
$|X|$ copies of this tile, corresponding to different letters of
$X$. Append the corresponding letters $x\in X$ to the ends of the
labels $w\in W$ of the boundary points of each of the copy of the
tile. Some of the obtained labels will be related by the equivalence
relation of Definition~\ref{def:til}, i.e., represent the same
points of the tile $\mathcal{T}_G$. Glue the corresponding points
together. Some of the obtained labels will belong to $W$. These
points will be the new boundary points. In this way we get a new
graph with labeled boundary points. Repeat now the procedure several
times, rescaling the graph in such a way that the original first
order graphs become small. We will get in this way a graph
resembling the tile $\mathcal{T}_G$ (see Chapter V
in~\cite{bond:PHD_USA} for more details). Making the necessary
identifications of its boundary we get an approximation of the limit
space of $G$. More details on this inductive approximation procedure
can be found in~\cite{nekrash:self-similar} Section~3.10.

The limit space of a finitely generated contracting self-similar
group $G$ can also be viewed as a hyperbolic boundary in the
following way. For a given finite generating system $S$ of $G$
define the \emph{self-similarity graph} $\Sigma(G,S)$ as the graph
with set of vertices $X^{*}$ in which two vertices $v_1,v_2\in
X^{*}$ are connected by an edge if and only if either $v_i=xv_j$,
for some $x\in X$ (vertical edges), or $s(v_i)=v_j$ for some $s\in
S$ (horizontal edges). In case of a contracting group, the
self-similarity graph $\Sigma(G,S)$ is Gromov-hyperbolic and its
hyperbolic boundary is homeomorphic to the limit space $\lims$.

The iterated monodromy group (IMG) construction is dual to the limit
space construction. It may be defined for partial self-coverings of
orbispaces, but we will only provide the definition in case of
topological spaces, since we do not need the more general
construction in this text (all iterated monodromy groups that appear
later are related to partial self-coverings of the Riemann sphere).

Let $\mathcal{M}$ be a path connected and locally path connected
topological space and let $\mathcal{M}_1$ be an open path connected
subset of $\mathcal{M}$. Let $f:\mathcal{M}_1\rightarrow\mathcal{M}$
be a $d$-fold covering. Denote by $f^n$ the $n$-fold iteration of
the map $f$. Then $f^n:\mathcal{M}_n\rightarrow\mathcal{M}$, where
$\mathcal{M}_n = f^{-n}(\mathcal{M})$, is a $d^n$-fold covering.

Fix a base point $t\in\mathcal{M}$ and let $T_t$ be the disjoint
union of the sets $f^{-n}(t), n\geq 0$ (formally speaking, these
sets may not be disjoint in $\mathcal{M}$). The set of pre-images
$T_t$ has a natural structure of a rooted $d$-ary tree. The base
point $t$ is the root, the vertices in $f^{-n}$ constitute level $n$
and every vertex $z$ in $f^{-n}(t)$ is connected by an edge to
$f(z)$ in $f^{-n+1}(t)$, for $n \geq 1$. The fundamental group
$\pi_1(\mathcal{M},t)$ acts naturally, through the monodromy action,
on every level $f^{-n}(t)$ and, in fact, acts by automorphisms on
$T_t$.

\begin{definition}
The \emph{iterated monodromy group} $IMG(f)$ of the covering $f$ is
the quotient of the fundamental group $\pi_1(\mathcal{M},t)$ by the
kernel of its action on the tree of pre-images $T_t$.
\end{definition}

\section{Classification guide}\label{app_notat}

Every 3-state automaton $\A$ with set of states
$S=\{\textbf0,\textbf1,\textbf2\}$ acting on the 2-letter alphabet
$X=\{0,1\}$ is assigned a unique number as follows. Given the wreath
recursion
\[
 \left\{
\begin{array}{l}
  \textbf0=\sigma^{a_{11}}(a_{12},a_{13}), \\
  \textbf1=\sigma^{a_{21}}(a_{22},a_{23}),  \\
  \textbf2=\sigma^{a_{31}}(a_{32},a_{33}),
\end{array}
\right.
\]
defining the automaton $\A$, where $a_{ij}\in
\{\textbf0,\textbf1,\textbf2\}$ for $j\neq1$ and $a_{i1}\in
\{0,1\}$, $i=1,2,3$, assign the number
\[
\begin{array}{l}
\mathop{\rm Number}\nolimits(\mathcal A)= \\
\qquad\qquad a_{12}+3a_{13}+9a_{22}+27a_{23}+81a_{32}+\\
 \qquad\qquad 243a_{33}+729(a_{11}+2a_{21}+4a_{31})+1
\end{array}
\]
to $\A$. With this agreement every $(3,2)$-automaton is assigned a
unique number in the range from $1$ to $5832$. The numbering of the
automata is induced by the lexicographic ordering of all automata in
the class. Each of the automata numbered $1$ through $729$ generates
the trivial group, since all vertex permutations are trivial in this
case. Each of the automata numbered $5104$ through $5832$ generates
the cyclic group $C_2$ of order $2$, since both states represent the
automorphism that acts by changing all letters in every word over
$X$. Therefore the nontrivial part of the classification is
concerned with the automata numbered by $730$ through $5103$.

Denote by $\mathcal A_n$ the automaton numbered by $n$ and by $G_n$
the corresponding group of tree automorphisms. Sometimes we may use
just the number to refer to the corresponding automaton or group.

The following three operations on automata do not change the
isomorphism class of the group generated by the corresponding
automaton (and do not change the action on the tree in essential
way):
\begin{enumerate}
\item[(i)] passing to inverses of all generators,
\item[(ii)] permuting the states of the automaton,
\item[(iii)] permuting the alphabet letters.
\end{enumerate}

\begin{definition}
Two automata $\mathcal A$ and $\mathcal B$ that can be obtained from
one another by using a composition of the operations ($i$)--($iii$),
are called \emph{symmetric}.
\end{definition}

For instance, the two automata in the lower right part of
Figure~\ref{4automata} are symmetric. The wreath recursion for the
automaton obtained by permuting both the names of the states and the
alphabet letters of the first of these two automata is
\begin{alignat*}{2}
  &a =        &&(b,a)  \\
  &b = \sigma &&(b,a)
\end{alignat*}
and this wreath recursion describes exactly the inverses of the tree
automorphism defining the second of the two automata.

Additional identifications can be made after automata minimization
is applied.

\begin{definition}
If the minimization of an automaton $\mathcal A$ is symmetric to the
minimization of an automaton $\mathcal B$, we say that the automata
$\mathcal A$ and $\mathcal B$ are \emph{minimally symmetric} and
write $\mathcal A\sim\mathcal B$.
\end{definition}

There are $194$ classes of $(3,2)$-automata that are pairwise not
minimally symmetric. Of these, $10$ are minimally symmetric to
automata with fewer than $3$ states and, as such, are subject of
Theorem~\ref{thm:class22} (\cite{gns00:automata}, see below).

At present, it is known that there are no more than $122$
non-isomorphic $(3,2)$-automaton groups. Some information on these
groups is given in Section~\ref{app_table}.

The proofs of some particular properties of the $194$ classes of
non-equivalent automata (and in particular, all known isomorphisms)
can be found in Section~\ref{app_proofs}. The few general results
that hold in the whole class were already mentioned in the
introduction.

The table in Section~\ref{app_corrtable} may be used to determine
the equivalence and the group isomorphism class for each automaton.
Every class is numbered by the smallest number of an automaton in
the class. For instance, an entry such as $x\sim y\cong z$ means
that the automata with the smallest number in the equivalence and
the (known) isomorphism class of $x$ are $y$ and $z$, respectively.
While the equivalence classes are easy to determine the isomorphism
class is not. Therefore, there may still be some additional
isomorphisms between some of the classes (which would eventually
cause changes in the $z$ numbers and consolidation of some of the
current isomorphism classes).

If one is interested in some particular $(3,2)$-automaton $\A$, we
recommend the following procedure:
\begin{itemize}
\item
Use the table in Section~\ref{app_corrtable} to find numbers for the
representatives of the equivalence and the isomorphism class of
$\A$. Minimizing the automaton and finding the symmetry is a
straightforward task, which is not presented here.
\item Use Section~\ref{app_table} to find information on the
group generated by $\mathcal A$ (more precisely, the isomorphic
group generated by the chosen representative in the class).
\item Use Section~\ref{app_proofs} to find the proof of the
isomorphism and some known properties.
\end{itemize}

\newpage
\section{Table of equivalence classes (and known isomorphisms)}
\label{app_corrtable} For explanation of the entries see
Section~\ref{app_notat}.

\medskip

\noindent{\footnotesize 1 through $729\sim1\simeq1$,}

\renewcommand{\arraystretch}{0.1}
\begin{multicols}{4} 
{\footnotesize\noindent $730\!\!\sim\!\!730\!\!\cong\!\!730$
$731\!\!\sim\!\!731\!\!\cong\!\!731$
$732\!\!\sim\!\!731\!\!\cong\!\!731$
$733\!\!\sim\!\!731\!\!\cong\!\!731$
$734\!\!\sim\!\!734\!\!\cong\!\!730$
$735\!\!\sim\!\!734\!\!\cong\!\!730$
$736\!\!\sim\!\!731\!\!\cong\!\!731$
$737\!\!\sim\!\!734\!\!\cong\!\!730$
$738\!\!\sim\!\!734\!\!\cong\!\!730$
$739\!\!\sim\!\!739\!\!\cong\!\!739$
$740\!\!\sim\!\!740\!\!\cong\!\!740$
$741\!\!\sim\!\!741\!\!\cong\!\!741$
$742\!\!\sim\!\!740\!\!\cong\!\!740$
$743\!\!\sim\!\!743\!\!\cong\!\!739$
$744\!\!\sim\!\!744\!\!\cong\!\!744$
$745\!\!\sim\!\!741\!\!\cong\!\!741$
$746\!\!\sim\!\!744\!\!\cong\!\!744$
$747\!\!\sim\!\!747\!\!\cong\!\!739$
$748\!\!\sim\!\!748\!\!\cong\!\!748$
$749\!\!\sim\!\!749\!\!\cong\!\!749$
$750\!\!\sim\!\!750\!\!\cong\!\!750$
$751\!\!\sim\!\!749\!\!\cong\!\!749$
$752\!\!\sim\!\!752\!\!\cong\!\!752$
$753\!\!\sim\!\!753\!\!\cong\!\!753$
$754\!\!\sim\!\!750\!\!\cong\!\!750$
$755\!\!\sim\!\!753\!\!\cong\!\!753$
$756\!\!\sim\!\!756\!\!\cong\!\!748$
$757\!\!\sim\!\!739\!\!\cong\!\!739$
$758\!\!\sim\!\!740\!\!\cong\!\!740$
$759\!\!\sim\!\!741\!\!\cong\!\!741$
$760\!\!\sim\!\!740\!\!\cong\!\!740$
$761\!\!\sim\!\!743\!\!\cong\!\!739$
$762\!\!\sim\!\!744\!\!\cong\!\!744$
$763\!\!\sim\!\!741\!\!\cong\!\!741$
$764\!\!\sim\!\!744\!\!\cong\!\!744$
$765\!\!\sim\!\!747\!\!\cong\!\!739$
$766\!\!\sim\!\!766\!\!\cong\!\!730$
$767\!\!\sim\!\!767\!\!\cong\!\!731$
$768\!\!\sim\!\!768\!\!\cong\!\!731$
$769\!\!\sim\!\!767\!\!\cong\!\!731$
$770\!\!\sim\!\!770\!\!\cong\!\!730$
$771\!\!\sim\!\!771\!\!\cong\!\!771$
$772\!\!\sim\!\!768\!\!\cong\!\!731$
$773\!\!\sim\!\!771\!\!\cong\!\!771$
$774\!\!\sim\!\!774\!\!\cong\!\!730$
$775\!\!\sim\!\!775\!\!\cong\!\!775$
$776\!\!\sim\!\!776\!\!\cong\!\!776$
$777\!\!\sim\!\!777\!\!\cong\!\!777$
$778\!\!\sim\!\!776\!\!\cong\!\!776$
$779\!\!\sim\!\!779\!\!\cong\!\!779$
$780\!\!\sim\!\!780\!\!\cong\!\!780$
$781\!\!\sim\!\!777\!\!\cong\!\!777$
$782\!\!\sim\!\!780\!\!\cong\!\!780$
$783\!\!\sim\!\!783\!\!\cong\!\!775$
$784\!\!\sim\!\!748\!\!\cong\!\!748$
$785\!\!\sim\!\!749\!\!\cong\!\!749$
$786\!\!\sim\!\!750\!\!\cong\!\!750$
$787\!\!\sim\!\!749\!\!\cong\!\!749$
$788\!\!\sim\!\!752\!\!\cong\!\!752$
$789\!\!\sim\!\!753\!\!\cong\!\!753$
$790\!\!\sim\!\!750\!\!\cong\!\!750$
$791\!\!\sim\!\!753\!\!\cong\!\!753$
$792\!\!\sim\!\!756\!\!\cong\!\!748$
$793\!\!\sim\!\!775\!\!\cong\!\!775$
$794\!\!\sim\!\!776\!\!\cong\!\!776$
$795\!\!\sim\!\!777\!\!\cong\!\!777$
$796\!\!\sim\!\!776\!\!\cong\!\!776$
$797\!\!\sim\!\!779\!\!\cong\!\!779$
$798\!\!\sim\!\!780\!\!\cong\!\!780$
$799\!\!\sim\!\!777\!\!\cong\!\!777$
$800\!\!\sim\!\!780\!\!\cong\!\!780$
$801\!\!\sim\!\!783\!\!\cong\!\!775$
$802\!\!\sim\!\!802\!\!\cong\!\!802$
$803\!\!\sim\!\!803\!\!\cong\!\!771$
$804\!\!\sim\!\!804\!\!\cong\!\!731$
$805\!\!\sim\!\!803\!\!\cong\!\!771$
$806\!\!\sim\!\!806\!\!\cong\!\!802$
$807\!\!\sim\!\!807\!\!\cong\!\!771$
$808\!\!\sim\!\!804\!\!\cong\!\!731$
$809\!\!\sim\!\!807\!\!\cong\!\!771$
$810\!\!\sim\!\!810\!\!\cong\!\!802$
$811\!\!\sim\!\!748\!\!\cong\!\!748$
$812\!\!\sim\!\!750\!\!\cong\!\!750$
$813\!\!\sim\!\!749\!\!\cong\!\!749$
$814\!\!\sim\!\!750\!\!\cong\!\!750$
$815\!\!\sim\!\!756\!\!\cong\!\!748$
$816\!\!\sim\!\!753\!\!\cong\!\!753$
$817\!\!\sim\!\!749\!\!\cong\!\!749$
$818\!\!\sim\!\!753\!\!\cong\!\!753$
$819\!\!\sim\!\!752\!\!\cong\!\!752$
$820\!\!\sim\!\!820\!\!\cong\!\!820$
$821\!\!\sim\!\!821\!\!\cong\!\!821$
$822\!\!\sim\!\!821\!\!\cong\!\!821$
$823\!\!\sim\!\!821\!\!\cong\!\!821$
$824\!\!\sim\!\!824\!\!\cong\!\!820$
$825\!\!\sim\!\!824\!\!\cong\!\!820$
$826\!\!\sim\!\!821\!\!\cong\!\!821$
$827\!\!\sim\!\!824\!\!\cong\!\!820$
$828\!\!\sim\!\!824\!\!\cong\!\!820$
$829\!\!\sim\!\!820\!\!\cong\!\!820$
$830\!\!\sim\!\!821\!\!\cong\!\!821$
$831\!\!\sim\!\!821\!\!\cong\!\!821$
$832\!\!\sim\!\!821\!\!\cong\!\!821$
$833\!\!\sim\!\!824\!\!\cong\!\!820$
$834\!\!\sim\!\!824\!\!\cong\!\!820$
$835\!\!\sim\!\!821\!\!\cong\!\!821$
$836\!\!\sim\!\!824\!\!\cong\!\!820$
$837\!\!\sim\!\!824\!\!\cong\!\!820$
$838\!\!\sim\!\!838\!\!\cong\!\!838$
$839\!\!\sim\!\!839\!\!\cong\!\!821$
$840\!\!\sim\!\!840\!\!\cong\!\!840$
$841\!\!\sim\!\!839\!\!\cong\!\!821$
$842\!\!\sim\!\!842\!\!\cong\!\!838$
$843\!\!\sim\!\!843\!\!\cong\!\!843$
$844\!\!\sim\!\!840\!\!\cong\!\!840$
$845\!\!\sim\!\!843\!\!\cong\!\!843$
$846\!\!\sim\!\!846\!\!\cong\!\!846$
$847\!\!\sim\!\!847\!\!\cong\!\!847$
$848\!\!\sim\!\!848\!\!\cong\!\!750$
$849\!\!\sim\!\!849\!\!\cong\!\!849$
$850\!\!\sim\!\!848\!\!\cong\!\!750$
$851\!\!\sim\!\!851\!\!\cong\!\!847$
$852\!\!\sim\!\!852\!\!\cong\!\!852$
$853\!\!\sim\!\!849\!\!\cong\!\!849$
$854\!\!\sim\!\!852\!\!\cong\!\!852$
$855\!\!\sim\!\!855\!\!\cong\!\!847$
$856\!\!\sim\!\!856\!\!\cong\!\!856$
$857\!\!\sim\!\!857\!\!\cong\!\!857$
$858\!\!\sim\!\!858\!\!\cong\!\!858$
$859\!\!\sim\!\!857\!\!\cong\!\!857$
$860\!\!\sim\!\!860\!\!\cong\!\!860$
$861\!\!\sim\!\!861\!\!\cong\!\!861$
$862\!\!\sim\!\!858\!\!\cong\!\!858$
$863\!\!\sim\!\!861\!\!\cong\!\!861$
$864\!\!\sim\!\!864\!\!\cong\!\!864$
$865\!\!\sim\!\!865\!\!\cong\!\!820$
$866\!\!\sim\!\!866\!\!\cong\!\!866$
$867\!\!\sim\!\!866\!\!\cong\!\!866$
$868\!\!\sim\!\!866\!\!\cong\!\!866$
$869\!\!\sim\!\!869\!\!\cong\!\!869$
$870\!\!\sim\!\!870\!\!\cong\!\!870$
$871\!\!\sim\!\!866\!\!\cong\!\!866$
$872\!\!\sim\!\!870\!\!\cong\!\!870$
$873\!\!\sim\!\!869\!\!\cong\!\!869$
$874\!\!\sim\!\!874\!\!\cong\!\!874$
$875\!\!\sim\!\!875\!\!\cong\!\!875$
$876\!\!\sim\!\!876\!\!\cong\!\!876$
$877\!\!\sim\!\!875\!\!\cong\!\!875$
$878\!\!\sim\!\!878\!\!\cong\!\!878$
$879\!\!\sim\!\!879\!\!\cong\!\!879$
$880\!\!\sim\!\!876\!\!\cong\!\!876$
$881\!\!\sim\!\!879\!\!\cong\!\!879$
$882\!\!\sim\!\!882\!\!\cong\!\!882$
$883\!\!\sim\!\!883\!\!\cong\!\!883$
$884\!\!\sim\!\!884\!\!\cong\!\!884$
$885\!\!\sim\!\!885\!\!\cong\!\!885$
$886\!\!\sim\!\!884\!\!\cong\!\!884$
$887\!\!\sim\!\!887\!\!\cong\!\!887$
$888\!\!\sim\!\!888\!\!\cong\!\!888$
$889\!\!\sim\!\!885\!\!\cong\!\!885$
$890\!\!\sim\!\!888\!\!\cong\!\!888$
$891\!\!\sim\!\!891\!\!\cong\!\!891$
$892\!\!\sim\!\!739\!\!\cong\!\!739$
$893\!\!\sim\!\!741\!\!\cong\!\!741$
$894\!\!\sim\!\!740\!\!\cong\!\!740$
$895\!\!\sim\!\!741\!\!\cong\!\!741$
$896\!\!\sim\!\!747\!\!\cong\!\!739$
$897\!\!\sim\!\!744\!\!\cong\!\!744$
$898\!\!\sim\!\!740\!\!\cong\!\!740$
$899\!\!\sim\!\!744\!\!\cong\!\!744$
$900\!\!\sim\!\!743\!\!\cong\!\!739$
$901\!\!\sim\!\!820\!\!\cong\!\!820$
$902\!\!\sim\!\!821\!\!\cong\!\!821$
$903\!\!\sim\!\!821\!\!\cong\!\!821$
$904\!\!\sim\!\!821\!\!\cong\!\!821$
$905\!\!\sim\!\!824\!\!\cong\!\!820$
$906\!\!\sim\!\!824\!\!\cong\!\!820$
$907\!\!\sim\!\!821\!\!\cong\!\!821$
$908\!\!\sim\!\!824\!\!\cong\!\!820$
$909\!\!\sim\!\!824\!\!\cong\!\!820$
$910\!\!\sim\!\!820\!\!\cong\!\!820$
$911\!\!\sim\!\!821\!\!\cong\!\!821$
$912\!\!\sim\!\!821\!\!\cong\!\!821$
$913\!\!\sim\!\!821\!\!\cong\!\!821$
$914\!\!\sim\!\!824\!\!\cong\!\!820$
$915\!\!\sim\!\!824\!\!\cong\!\!820$
$916\!\!\sim\!\!821\!\!\cong\!\!821$
$917\!\!\sim\!\!824\!\!\cong\!\!820$
$918\!\!\sim\!\!824\!\!\cong\!\!820$
$919\!\!\sim\!\!919\!\!\cong\!\!820$
$920\!\!\sim\!\!920\!\!\cong\!\!920$
$921\!\!\sim\!\!920\!\!\cong\!\!920$
$922\!\!\sim\!\!920\!\!\cong\!\!920$
$923\!\!\sim\!\!923\!\!\cong\!\!923$
$924\!\!\sim\!\!924\!\!\cong\!\!870$
$925\!\!\sim\!\!920\!\!\cong\!\!920$
$926\!\!\sim\!\!924\!\!\cong\!\!870$
$927\!\!\sim\!\!923\!\!\cong\!\!923$
$928\!\!\sim\!\!928\!\!\cong\!\!820$
$929\!\!\sim\!\!929\!\!\cong\!\!929$
$930\!\!\sim\!\!930\!\!\cong\!\!821$
$931\!\!\sim\!\!929\!\!\cong\!\!929$
$932\!\!\sim\!\!932\!\!\cong\!\!820$
$933\!\!\sim\!\!933\!\!\cong\!\!849$
$934\!\!\sim\!\!930\!\!\cong\!\!821$
$935\!\!\sim\!\!933\!\!\cong\!\!849$
$936\!\!\sim\!\!936\!\!\cong\!\!820$
$937\!\!\sim\!\!937\!\!\cong\!\!937$
$938\!\!\sim\!\!938\!\!\cong\!\!938$
$939\!\!\sim\!\!939\!\!\cong\!\!939$
$940\!\!\sim\!\!938\!\!\cong\!\!938$
$941\!\!\sim\!\!941\!\!\cong\!\!941$
$942\!\!\sim\!\!942\!\!\cong\!\!942$
$943\!\!\sim\!\!939\!\!\cong\!\!939$
$944\!\!\sim\!\!942\!\!\cong\!\!942$
$945\!\!\sim\!\!945\!\!\cong\!\!941$
$946\!\!\sim\!\!838\!\!\cong\!\!838$
$947\!\!\sim\!\!840\!\!\cong\!\!840$
$948\!\!\sim\!\!839\!\!\cong\!\!821$
$949\!\!\sim\!\!840\!\!\cong\!\!840$
$950\!\!\sim\!\!846\!\!\cong\!\!846$
$951\!\!\sim\!\!843\!\!\cong\!\!843$
$952\!\!\sim\!\!839\!\!\cong\!\!821$
$953\!\!\sim\!\!843\!\!\cong\!\!843$
$954\!\!\sim\!\!842\!\!\cong\!\!838$
$955\!\!\sim\!\!955\!\!\cong\!\!937$
$956\!\!\sim\!\!956\!\!\cong\!\!956$
$957\!\!\sim\!\!957\!\!\cong\!\!957$
$958\!\!\sim\!\!956\!\!\cong\!\!956$
$959\!\!\sim\!\!959\!\!\cong\!\!959$
$960\!\!\sim\!\!960\!\!\cong\!\!960$
$961\!\!\sim\!\!957\!\!\cong\!\!957$
$962\!\!\sim\!\!960\!\!\cong\!\!960$
$963\!\!\sim\!\!963\!\!\cong\!\!963$
$964\!\!\sim\!\!964\!\!\cong\!\!739$
$965\!\!\sim\!\!965\!\!\cong\!\!965$
$966\!\!\sim\!\!966\!\!\cong\!\!966$
$967\!\!\sim\!\!965\!\!\cong\!\!965$
$968\!\!\sim\!\!968\!\!\cong\!\!968$
$969\!\!\sim\!\!969\!\!\cong\!\!969$
$970\!\!\sim\!\!966\!\!\cong\!\!966$
$971\!\!\sim\!\!969\!\!\cong\!\!969$
$972\!\!\sim\!\!972\!\!\cong\!\!739$
$973\!\!\sim\!\!748\!\!\cong\!\!748$
$974\!\!\sim\!\!750\!\!\cong\!\!750$
$975\!\!\sim\!\!749\!\!\cong\!\!749$
$976\!\!\sim\!\!750\!\!\cong\!\!750$
$977\!\!\sim\!\!756\!\!\cong\!\!748$
$978\!\!\sim\!\!753\!\!\cong\!\!753$
$979\!\!\sim\!\!749\!\!\cong\!\!749$
$980\!\!\sim\!\!753\!\!\cong\!\!753$
$981\!\!\sim\!\!752\!\!\cong\!\!752$
$982\!\!\sim\!\!838\!\!\cong\!\!838$
$983\!\!\sim\!\!839\!\!\cong\!\!821$
$984\!\!\sim\!\!840\!\!\cong\!\!840$
$985\!\!\sim\!\!839\!\!\cong\!\!821$
$986\!\!\sim\!\!842\!\!\cong\!\!838$
$987\!\!\sim\!\!843\!\!\cong\!\!843$
$988\!\!\sim\!\!840\!\!\cong\!\!840$
$989\!\!\sim\!\!843\!\!\cong\!\!843$
$990\!\!\sim\!\!846\!\!\cong\!\!846$
$991\!\!\sim\!\!865\!\!\cong\!\!820$
$992\!\!\sim\!\!866\!\!\cong\!\!866$
$993\!\!\sim\!\!866\!\!\cong\!\!866$
$994\!\!\sim\!\!866\!\!\cong\!\!866$
$995\!\!\sim\!\!869\!\!\cong\!\!869$
$996\!\!\sim\!\!870\!\!\cong\!\!870$
$997\!\!\sim\!\!866\!\!\cong\!\!866$
$998\!\!\sim\!\!870\!\!\cong\!\!870$
$999\!\!\sim\!\!869\!\!\cong\!\!869$
$1000\!\!\sim\!\!820\!\!\cong\!\!820$
$1001\!\!\sim\!\!821\!\!\cong\!\!821$
$1002\!\!\sim\!\!821\!\!\cong\!\!821$
$1003\!\!\sim\!\!821\!\!\cong\!\!821$
$1004\!\!\sim\!\!824\!\!\cong\!\!820$
$1005\!\!\sim\!\!824\!\!\cong\!\!820$
$1006\!\!\sim\!\!821\!\!\cong\!\!821$
$1007\!\!\sim\!\!824\!\!\cong\!\!820$
$1008\!\!\sim\!\!824\!\!\cong\!\!820$
$1009\!\!\sim\!\!847\!\!\cong\!\!847$
$1010\!\!\sim\!\!848\!\!\cong\!\!750$
$1011\!\!\sim\!\!849\!\!\cong\!\!849$
$1012\!\!\sim\!\!848\!\!\cong\!\!750$
$1013\!\!\sim\!\!851\!\!\cong\!\!847$
$1014\!\!\sim\!\!852\!\!\cong\!\!852$
$1015\!\!\sim\!\!849\!\!\cong\!\!849$
$1016\!\!\sim\!\!852\!\!\cong\!\!852$
$1017\!\!\sim\!\!855\!\!\cong\!\!847$
$1018\!\!\sim\!\!874\!\!\cong\!\!874$
$1019\!\!\sim\!\!875\!\!\cong\!\!875$
$1020\!\!\sim\!\!876\!\!\cong\!\!876$
$1021\!\!\sim\!\!875\!\!\cong\!\!875$
$1022\!\!\sim\!\!878\!\!\cong\!\!878$
$1023\!\!\sim\!\!879\!\!\cong\!\!879$
$1024\!\!\sim\!\!876\!\!\cong\!\!876$
$1025\!\!\sim\!\!879\!\!\cong\!\!879$
$1026\!\!\sim\!\!882\!\!\cong\!\!882$
$1027\!\!\sim\!\!820\!\!\cong\!\!820$
$1028\!\!\sim\!\!821\!\!\cong\!\!821$
$1029\!\!\sim\!\!821\!\!\cong\!\!821$
$1030\!\!\sim\!\!821\!\!\cong\!\!821$
$1031\!\!\sim\!\!824\!\!\cong\!\!820$
$1032\!\!\sim\!\!824\!\!\cong\!\!820$
$1033\!\!\sim\!\!821\!\!\cong\!\!821$
$1034\!\!\sim\!\!824\!\!\cong\!\!820$
$1035\!\!\sim\!\!824\!\!\cong\!\!820$
$1036\!\!\sim\!\!856\!\!\cong\!\!856$
$1037\!\!\sim\!\!857\!\!\cong\!\!857$
$1038\!\!\sim\!\!858\!\!\cong\!\!858$
$1039\!\!\sim\!\!857\!\!\cong\!\!857$
$1040\!\!\sim\!\!860\!\!\cong\!\!860$
$1041\!\!\sim\!\!861\!\!\cong\!\!861$
$1042\!\!\sim\!\!858\!\!\cong\!\!858$
$1043\!\!\sim\!\!861\!\!\cong\!\!861$
$1044\!\!\sim\!\!864\!\!\cong\!\!864$
$1045\!\!\sim\!\!883\!\!\cong\!\!883$
$1046\!\!\sim\!\!884\!\!\cong\!\!884$
$1047\!\!\sim\!\!885\!\!\cong\!\!885$
$1048\!\!\sim\!\!884\!\!\cong\!\!884$
$1049\!\!\sim\!\!887\!\!\cong\!\!887$
$1050\!\!\sim\!\!888\!\!\cong\!\!888$
$1051\!\!\sim\!\!885\!\!\cong\!\!885$
$1052\!\!\sim\!\!888\!\!\cong\!\!888$
$1053\!\!\sim\!\!891\!\!\cong\!\!891$
$1054\!\!\sim\!\!802\!\!\cong\!\!802$
$1055\!\!\sim\!\!804\!\!\cong\!\!731$
$1056\!\!\sim\!\!803\!\!\cong\!\!771$
$1057\!\!\sim\!\!804\!\!\cong\!\!731$
$1058\!\!\sim\!\!810\!\!\cong\!\!802$
$1059\!\!\sim\!\!807\!\!\cong\!\!771$
$1060\!\!\sim\!\!803\!\!\cong\!\!771$
$1061\!\!\sim\!\!807\!\!\cong\!\!771$
$1062\!\!\sim\!\!806\!\!\cong\!\!802$
$1063\!\!\sim\!\!964\!\!\cong\!\!739$
$1064\!\!\sim\!\!966\!\!\cong\!\!966$
$1065\!\!\sim\!\!965\!\!\cong\!\!965$
$1066\!\!\sim\!\!966\!\!\cong\!\!966$
$1067\!\!\sim\!\!972\!\!\cong\!\!739$
$1068\!\!\sim\!\!969\!\!\cong\!\!969$
$1069\!\!\sim\!\!965\!\!\cong\!\!965$
$1070\!\!\sim\!\!969\!\!\cong\!\!969$
$1071\!\!\sim\!\!968\!\!\cong\!\!968$
$1072\!\!\sim\!\!883\!\!\cong\!\!883$
$1073\!\!\sim\!\!885\!\!\cong\!\!885$
$1074\!\!\sim\!\!884\!\!\cong\!\!884$
$1075\!\!\sim\!\!885\!\!\cong\!\!885$
$1076\!\!\sim\!\!891\!\!\cong\!\!891$
$1077\!\!\sim\!\!888\!\!\cong\!\!888$
$1078\!\!\sim\!\!884\!\!\cong\!\!884$
$1079\!\!\sim\!\!888\!\!\cong\!\!888$
$1080\!\!\sim\!\!887\!\!\cong\!\!887$
$1081\!\!\sim\!\!964\!\!\cong\!\!739$
$1082\!\!\sim\!\!966\!\!\cong\!\!966$
$1083\!\!\sim\!\!965\!\!\cong\!\!965$
$1084\!\!\sim\!\!966\!\!\cong\!\!966$
$1085\!\!\sim\!\!972\!\!\cong\!\!739$
$1086\!\!\sim\!\!969\!\!\cong\!\!969$
$1087\!\!\sim\!\!965\!\!\cong\!\!965$
$1088\!\!\sim\!\!969\!\!\cong\!\!969$
$1089\!\!\sim\!\!968\!\!\cong\!\!968$
$1090\!\!\sim\!\!1090\!\!\cong\!\!1090$
$1091\!\!\sim\!\!1091\!\!\cong\!\!731$
$1092\!\!\sim\!\!1091\!\!\cong\!\!731$
$1093\!\!\sim\!\!1091\!\!\cong\!\!731$
$1094\!\!\sim\!\!1094\!\!\cong\!\!1090$
$1095\!\!\sim\!\!1094\!\!\cong\!\!1090$
$1096\!\!\sim\!\!1091\!\!\cong\!\!731$
$1097\!\!\sim\!\!1094\!\!\cong\!\!1090$
$1098\!\!\sim\!\!1094\!\!\cong\!\!1090$
$1099\!\!\sim\!\!1090\!\!\cong\!\!1090$
$1100\!\!\sim\!\!1091\!\!\cong\!\!731$
$1101\!\!\sim\!\!1091\!\!\cong\!\!731$
$1102\!\!\sim\!\!1091\!\!\cong\!\!731$
$1103\!\!\sim\!\!1094\!\!\cong\!\!1090$
$1104\!\!\sim\!\!1094\!\!\cong\!\!1090$
$1105\!\!\sim\!\!1091\!\!\cong\!\!731$
$1106\!\!\sim\!\!1094\!\!\cong\!\!1090$
$1107\!\!\sim\!\!1094\!\!\cong\!\!1090$
$1108\!\!\sim\!\!883\!\!\cong\!\!883$
$1109\!\!\sim\!\!885\!\!\cong\!\!885$
$1110\!\!\sim\!\!884\!\!\cong\!\!884$
$1111\!\!\sim\!\!885\!\!\cong\!\!885$
$1112\!\!\sim\!\!891\!\!\cong\!\!891$
$1113\!\!\sim\!\!888\!\!\cong\!\!888$
$1114\!\!\sim\!\!884\!\!\cong\!\!884$
$1115\!\!\sim\!\!888\!\!\cong\!\!888$
$1116\!\!\sim\!\!887\!\!\cong\!\!887$
$1117\!\!\sim\!\!1090\!\!\cong\!\!1090$
$1118\!\!\sim\!\!1091\!\!\cong\!\!731$
$1119\!\!\sim\!\!1091\!\!\cong\!\!731$
$1120\!\!\sim\!\!1091\!\!\cong\!\!731$
$1121\!\!\sim\!\!1094\!\!\cong\!\!1090$
$1122\!\!\sim\!\!1094\!\!\cong\!\!1090$
$1123\!\!\sim\!\!1091\!\!\cong\!\!731$
$1124\!\!\sim\!\!1094\!\!\cong\!\!1090$
$1125\!\!\sim\!\!1094\!\!\cong\!\!1090$
$1126\!\!\sim\!\!1090\!\!\cong\!\!1090$
$1127\!\!\sim\!\!1091\!\!\cong\!\!731$
$1128\!\!\sim\!\!1091\!\!\cong\!\!731$
$1129\!\!\sim\!\!1091\!\!\cong\!\!731$
$1130\!\!\sim\!\!1094\!\!\cong\!\!1090$
$1131\!\!\sim\!\!1094\!\!\cong\!\!1090$
$1132\!\!\sim\!\!1091\!\!\cong\!\!731$
$1133\!\!\sim\!\!1094\!\!\cong\!\!1090$
$1134\!\!\sim\!\!1094\!\!\cong\!\!1090$
$1135\!\!\sim\!\!775\!\!\cong\!\!775$
$1136\!\!\sim\!\!777\!\!\cong\!\!777$
$1137\!\!\sim\!\!776\!\!\cong\!\!776$
$1138\!\!\sim\!\!777\!\!\cong\!\!777$
$1139\!\!\sim\!\!783\!\!\cong\!\!775$
$1140\!\!\sim\!\!780\!\!\cong\!\!780$
$1141\!\!\sim\!\!776\!\!\cong\!\!776$
$1142\!\!\sim\!\!780\!\!\cong\!\!780$
$1143\!\!\sim\!\!779\!\!\cong\!\!779$
$1144\!\!\sim\!\!955\!\!\cong\!\!937$
$1145\!\!\sim\!\!957\!\!\cong\!\!957$
$1146\!\!\sim\!\!956\!\!\cong\!\!956$
$1147\!\!\sim\!\!957\!\!\cong\!\!957$
$1148\!\!\sim\!\!963\!\!\cong\!\!963$
$1149\!\!\sim\!\!960\!\!\cong\!\!960$
$1150\!\!\sim\!\!956\!\!\cong\!\!956$
$1151\!\!\sim\!\!960\!\!\cong\!\!960$
$1152\!\!\sim\!\!959\!\!\cong\!\!959$
$1153\!\!\sim\!\!874\!\!\cong\!\!874$
$1154\!\!\sim\!\!876\!\!\cong\!\!876$
$1155\!\!\sim\!\!875\!\!\cong\!\!875$
$1156\!\!\sim\!\!876\!\!\cong\!\!876$
$1157\!\!\sim\!\!882\!\!\cong\!\!882$
$1158\!\!\sim\!\!879\!\!\cong\!\!879$
$1159\!\!\sim\!\!875\!\!\cong\!\!875$
$1160\!\!\sim\!\!879\!\!\cong\!\!879$
$1161\!\!\sim\!\!878\!\!\cong\!\!878$
$1162\!\!\sim\!\!937\!\!\cong\!\!937$
$1163\!\!\sim\!\!939\!\!\cong\!\!939$
$1164\!\!\sim\!\!938\!\!\cong\!\!938$
$1165\!\!\sim\!\!939\!\!\cong\!\!939$
$1166\!\!\sim\!\!945\!\!\cong\!\!941$
$1167\!\!\sim\!\!942\!\!\cong\!\!942$
$1168\!\!\sim\!\!938\!\!\cong\!\!938$
$1169\!\!\sim\!\!942\!\!\cong\!\!942$
$1170\!\!\sim\!\!941\!\!\cong\!\!941$
$1171\!\!\sim\!\!1090\!\!\cong\!\!1090$
$1172\!\!\sim\!\!1091\!\!\cong\!\!731$
$1173\!\!\sim\!\!1091\!\!\cong\!\!731$
$1174\!\!\sim\!\!1091\!\!\cong\!\!731$
$1175\!\!\sim\!\!1094\!\!\cong\!\!1090$
$1176\!\!\sim\!\!1094\!\!\cong\!\!1090$
$1177\!\!\sim\!\!1091\!\!\cong\!\!731$
$1178\!\!\sim\!\!1094\!\!\cong\!\!1090$
$1179\!\!\sim\!\!1094\!\!\cong\!\!1090$
$1180\!\!\sim\!\!1090\!\!\cong\!\!1090$
$1181\!\!\sim\!\!1091\!\!\cong\!\!731$
$1182\!\!\sim\!\!1091\!\!\cong\!\!731$
$1183\!\!\sim\!\!1091\!\!\cong\!\!731$
$1184\!\!\sim\!\!1094\!\!\cong\!\!1090$
$1185\!\!\sim\!\!1094\!\!\cong\!\!1090$
$1186\!\!\sim\!\!1091\!\!\cong\!\!731$
$1187\!\!\sim\!\!1094\!\!\cong\!\!1090$
$1188\!\!\sim\!\!1094\!\!\cong\!\!1090$
$1189\!\!\sim\!\!856\!\!\cong\!\!856$
$1190\!\!\sim\!\!858\!\!\cong\!\!858$
$1191\!\!\sim\!\!857\!\!\cong\!\!857$
$1192\!\!\sim\!\!858\!\!\cong\!\!858$
$1193\!\!\sim\!\!864\!\!\cong\!\!864$
$1194\!\!\sim\!\!861\!\!\cong\!\!861$
$1195\!\!\sim\!\!857\!\!\cong\!\!857$
$1196\!\!\sim\!\!861\!\!\cong\!\!861$
$1197\!\!\sim\!\!860\!\!\cong\!\!860$
$1198\!\!\sim\!\!1090\!\!\cong\!\!1090$
$1199\!\!\sim\!\!1091\!\!\cong\!\!731$
$1200\!\!\sim\!\!1091\!\!\cong\!\!731$
$1201\!\!\sim\!\!1091\!\!\cong\!\!731$
$1202\!\!\sim\!\!1094\!\!\cong\!\!1090$
$1203\!\!\sim\!\!1094\!\!\cong\!\!1090$
$1204\!\!\sim\!\!1091\!\!\cong\!\!731$
$1205\!\!\sim\!\!1094\!\!\cong\!\!1090$
$1206\!\!\sim\!\!1094\!\!\cong\!\!1090$
$1207\!\!\sim\!\!1090\!\!\cong\!\!1090$
$1208\!\!\sim\!\!1091\!\!\cong\!\!731$
$1209\!\!\sim\!\!1091\!\!\cong\!\!731$
$1210\!\!\sim\!\!1091\!\!\cong\!\!731$
$1211\!\!\sim\!\!1094\!\!\cong\!\!1090$
$1212\!\!\sim\!\!1094\!\!\cong\!\!1090$
$1213\!\!\sim\!\!1091\!\!\cong\!\!731$
$1214\!\!\sim\!\!1094\!\!\cong\!\!1090$
$1215\!\!\sim\!\!1094\!\!\cong\!\!1090$
$1216\!\!\sim\!\!739\!\!\cong\!\!739$
$1217\!\!\sim\!\!741\!\!\cong\!\!741$
$1218\!\!\sim\!\!740\!\!\cong\!\!740$
$1219\!\!\sim\!\!741\!\!\cong\!\!741$
$1220\!\!\sim\!\!747\!\!\cong\!\!739$
$1221\!\!\sim\!\!744\!\!\cong\!\!744$
$1222\!\!\sim\!\!740\!\!\cong\!\!740$
$1223\!\!\sim\!\!744\!\!\cong\!\!744$
$1224\!\!\sim\!\!743\!\!\cong\!\!739$
$1225\!\!\sim\!\!919\!\!\cong\!\!820$
$1226\!\!\sim\!\!920\!\!\cong\!\!920$
$1227\!\!\sim\!\!920\!\!\cong\!\!920$
$1228\!\!\sim\!\!920\!\!\cong\!\!920$
$1229\!\!\sim\!\!923\!\!\cong\!\!923$
$1230\!\!\sim\!\!924\!\!\cong\!\!870$
$1231\!\!\sim\!\!920\!\!\cong\!\!920$
$1232\!\!\sim\!\!924\!\!\cong\!\!870$
$1233\!\!\sim\!\!923\!\!\cong\!\!923$
$1234\!\!\sim\!\!838\!\!\cong\!\!838$
$1235\!\!\sim\!\!840\!\!\cong\!\!840$
$1236\!\!\sim\!\!839\!\!\cong\!\!821$
$1237\!\!\sim\!\!840\!\!\cong\!\!840$
$1238\!\!\sim\!\!846\!\!\cong\!\!846$
$1239\!\!\sim\!\!843\!\!\cong\!\!843$
$1240\!\!\sim\!\!839\!\!\cong\!\!821$
$1241\!\!\sim\!\!843\!\!\cong\!\!843$
$1242\!\!\sim\!\!842\!\!\cong\!\!838$
$1243\!\!\sim\!\!820\!\!\cong\!\!820$
$1244\!\!\sim\!\!821\!\!\cong\!\!821$
$1245\!\!\sim\!\!821\!\!\cong\!\!821$
$1246\!\!\sim\!\!821\!\!\cong\!\!821$
$1247\!\!\sim\!\!824\!\!\cong\!\!820$
$1248\!\!\sim\!\!824\!\!\cong\!\!820$
$1249\!\!\sim\!\!821\!\!\cong\!\!821$
$1250\!\!\sim\!\!824\!\!\cong\!\!820$
$1251\!\!\sim\!\!824\!\!\cong\!\!820$
$1252\!\!\sim\!\!928\!\!\cong\!\!820$
$1253\!\!\sim\!\!929\!\!\cong\!\!929$
$1254\!\!\sim\!\!930\!\!\cong\!\!821$
$1255\!\!\sim\!\!929\!\!\cong\!\!929$
$1256\!\!\sim\!\!932\!\!\cong\!\!820$
$1257\!\!\sim\!\!933\!\!\cong\!\!849$
$1258\!\!\sim\!\!930\!\!\cong\!\!821$
$1259\!\!\sim\!\!933\!\!\cong\!\!849$
$1260\!\!\sim\!\!936\!\!\cong\!\!820$
$1261\!\!\sim\!\!955\!\!\cong\!\!937$
$1262\!\!\sim\!\!956\!\!\cong\!\!956$
$1263\!\!\sim\!\!957\!\!\cong\!\!957$
$1264\!\!\sim\!\!956\!\!\cong\!\!956$
$1265\!\!\sim\!\!959\!\!\cong\!\!959$
$1266\!\!\sim\!\!960\!\!\cong\!\!960$
$1267\!\!\sim\!\!957\!\!\cong\!\!957$
$1268\!\!\sim\!\!960\!\!\cong\!\!960$
$1269\!\!\sim\!\!963\!\!\cong\!\!963$
$1270\!\!\sim\!\!820\!\!\cong\!\!820$
$1271\!\!\sim\!\!821\!\!\cong\!\!821$
$1272\!\!\sim\!\!821\!\!\cong\!\!821$
$1273\!\!\sim\!\!821\!\!\cong\!\!821$
$1274\!\!\sim\!\!824\!\!\cong\!\!820$
$1275\!\!\sim\!\!824\!\!\cong\!\!820$
$1276\!\!\sim\!\!821\!\!\cong\!\!821$
$1277\!\!\sim\!\!824\!\!\cong\!\!820$
$1278\!\!\sim\!\!824\!\!\cong\!\!820$
$1279\!\!\sim\!\!937\!\!\cong\!\!937$
$1280\!\!\sim\!\!938\!\!\cong\!\!938$
$1281\!\!\sim\!\!939\!\!\cong\!\!939$
$1282\!\!\sim\!\!938\!\!\cong\!\!938$
$1283\!\!\sim\!\!941\!\!\cong\!\!941$
$1284\!\!\sim\!\!942\!\!\cong\!\!942$
$1285\!\!\sim\!\!939\!\!\cong\!\!939$
$1286\!\!\sim\!\!942\!\!\cong\!\!942$
$1287\!\!\sim\!\!945\!\!\cong\!\!941$
$1288\!\!\sim\!\!964\!\!\cong\!\!739$
$1289\!\!\sim\!\!965\!\!\cong\!\!965$
$1290\!\!\sim\!\!966\!\!\cong\!\!966$
$1291\!\!\sim\!\!965\!\!\cong\!\!965$
$1292\!\!\sim\!\!968\!\!\cong\!\!968$
$1293\!\!\sim\!\!969\!\!\cong\!\!969$
$1294\!\!\sim\!\!966\!\!\cong\!\!966$
$1295\!\!\sim\!\!969\!\!\cong\!\!969$
$1296\!\!\sim\!\!972\!\!\cong\!\!739$
$1297\!\!\sim\!\!775\!\!\cong\!\!775$
$1298\!\!\sim\!\!777\!\!\cong\!\!777$
$1299\!\!\sim\!\!776\!\!\cong\!\!776$
$1300\!\!\sim\!\!777\!\!\cong\!\!777$
$1301\!\!\sim\!\!783\!\!\cong\!\!775$
$1302\!\!\sim\!\!780\!\!\cong\!\!780$
$1303\!\!\sim\!\!776\!\!\cong\!\!776$
$1304\!\!\sim\!\!780\!\!\cong\!\!780$
$1305\!\!\sim\!\!779\!\!\cong\!\!779$
$1306\!\!\sim\!\!937\!\!\cong\!\!937$
$1307\!\!\sim\!\!939\!\!\cong\!\!939$
$1308\!\!\sim\!\!938\!\!\cong\!\!938$
$1309\!\!\sim\!\!939\!\!\cong\!\!939$
$1310\!\!\sim\!\!945\!\!\cong\!\!941$
$1311\!\!\sim\!\!942\!\!\cong\!\!942$
$1312\!\!\sim\!\!938\!\!\cong\!\!938$
$1313\!\!\sim\!\!942\!\!\cong\!\!942$
$1314\!\!\sim\!\!941\!\!\cong\!\!941$
$1315\!\!\sim\!\!856\!\!\cong\!\!856$
$1316\!\!\sim\!\!858\!\!\cong\!\!858$
$1317\!\!\sim\!\!857\!\!\cong\!\!857$
$1318\!\!\sim\!\!858\!\!\cong\!\!858$
$1319\!\!\sim\!\!864\!\!\cong\!\!864$
$1320\!\!\sim\!\!861\!\!\cong\!\!861$
$1321\!\!\sim\!\!857\!\!\cong\!\!857$
$1322\!\!\sim\!\!861\!\!\cong\!\!861$
$1323\!\!\sim\!\!860\!\!\cong\!\!860$
$1324\!\!\sim\!\!955\!\!\cong\!\!937$
$1325\!\!\sim\!\!957\!\!\cong\!\!957$
$1326\!\!\sim\!\!956\!\!\cong\!\!956$
$1327\!\!\sim\!\!957\!\!\cong\!\!957$
$1328\!\!\sim\!\!963\!\!\cong\!\!963$
$1329\!\!\sim\!\!960\!\!\cong\!\!960$
$1330\!\!\sim\!\!956\!\!\cong\!\!956$
$1331\!\!\sim\!\!960\!\!\cong\!\!960$
$1332\!\!\sim\!\!959\!\!\cong\!\!959$
$1333\!\!\sim\!\!1090\!\!\cong\!\!1090$
$1334\!\!\sim\!\!1091\!\!\cong\!\!731$
$1335\!\!\sim\!\!1091\!\!\cong\!\!731$
$1336\!\!\sim\!\!1091\!\!\cong\!\!731$
$1337\!\!\sim\!\!1094\!\!\cong\!\!1090$
$1338\!\!\sim\!\!1094\!\!\cong\!\!1090$
$1339\!\!\sim\!\!1091\!\!\cong\!\!731$
$1340\!\!\sim\!\!1094\!\!\cong\!\!1090$
$1341\!\!\sim\!\!1094\!\!\cong\!\!1090$
$1342\!\!\sim\!\!1090\!\!\cong\!\!1090$
$1343\!\!\sim\!\!1091\!\!\cong\!\!731$
$1344\!\!\sim\!\!1091\!\!\cong\!\!731$
$1345\!\!\sim\!\!1091\!\!\cong\!\!731$
$1346\!\!\sim\!\!1094\!\!\cong\!\!1090$
$1347\!\!\sim\!\!1094\!\!\cong\!\!1090$
$1348\!\!\sim\!\!1091\!\!\cong\!\!731$
$1349\!\!\sim\!\!1094\!\!\cong\!\!1090$
$1350\!\!\sim\!\!1094\!\!\cong\!\!1090$
$1351\!\!\sim\!\!874\!\!\cong\!\!874$
$1352\!\!\sim\!\!876\!\!\cong\!\!876$
$1353\!\!\sim\!\!875\!\!\cong\!\!875$
$1354\!\!\sim\!\!876\!\!\cong\!\!876$
$1355\!\!\sim\!\!882\!\!\cong\!\!882$
$1356\!\!\sim\!\!879\!\!\cong\!\!879$
$1357\!\!\sim\!\!875\!\!\cong\!\!875$
$1358\!\!\sim\!\!879\!\!\cong\!\!879$
$1359\!\!\sim\!\!878\!\!\cong\!\!878$
$1360\!\!\sim\!\!1090\!\!\cong\!\!1090$
$1361\!\!\sim\!\!1091\!\!\cong\!\!731$
$1362\!\!\sim\!\!1091\!\!\cong\!\!731$
$1363\!\!\sim\!\!1091\!\!\cong\!\!731$
$1364\!\!\sim\!\!1094\!\!\cong\!\!1090$
$1365\!\!\sim\!\!1094\!\!\cong\!\!1090$
$1366\!\!\sim\!\!1091\!\!\cong\!\!731$
$1367\!\!\sim\!\!1094\!\!\cong\!\!1090$
$1368\!\!\sim\!\!1094\!\!\cong\!\!1090$
$1369\!\!\sim\!\!1090\!\!\cong\!\!1090$
$1370\!\!\sim\!\!1091\!\!\cong\!\!731$
$1371\!\!\sim\!\!1091\!\!\cong\!\!731$
$1372\!\!\sim\!\!1091\!\!\cong\!\!731$
$1373\!\!\sim\!\!1094\!\!\cong\!\!1090$
$1374\!\!\sim\!\!1094\!\!\cong\!\!1090$
$1375\!\!\sim\!\!1091\!\!\cong\!\!731$
$1376\!\!\sim\!\!1094\!\!\cong\!\!1090$
$1377\!\!\sim\!\!1094\!\!\cong\!\!1090$
$1378\!\!\sim\!\!766\!\!\cong\!\!730$
$1379\!\!\sim\!\!768\!\!\cong\!\!731$
$1380\!\!\sim\!\!767\!\!\cong\!\!731$
$1381\!\!\sim\!\!768\!\!\cong\!\!731$
$1382\!\!\sim\!\!774\!\!\cong\!\!730$
$1383\!\!\sim\!\!771\!\!\cong\!\!771$
$1384\!\!\sim\!\!767\!\!\cong\!\!731$
$1385\!\!\sim\!\!771\!\!\cong\!\!771$
$1386\!\!\sim\!\!770\!\!\cong\!\!730$
$1387\!\!\sim\!\!928\!\!\cong\!\!820$
$1388\!\!\sim\!\!930\!\!\cong\!\!821$
$1389\!\!\sim\!\!929\!\!\cong\!\!929$
$1390\!\!\sim\!\!930\!\!\cong\!\!821$
$1391\!\!\sim\!\!936\!\!\cong\!\!820$
$1392\!\!\sim\!\!933\!\!\cong\!\!849$
$1393\!\!\sim\!\!929\!\!\cong\!\!929$
$1394\!\!\sim\!\!933\!\!\cong\!\!849$
$1395\!\!\sim\!\!932\!\!\cong\!\!820$
$1396\!\!\sim\!\!847\!\!\cong\!\!847$
$1397\!\!\sim\!\!849\!\!\cong\!\!849$
$1398\!\!\sim\!\!848\!\!\cong\!\!750$
$1399\!\!\sim\!\!849\!\!\cong\!\!849$
$1400\!\!\sim\!\!855\!\!\cong\!\!847$
$1401\!\!\sim\!\!852\!\!\cong\!\!852$
$1402\!\!\sim\!\!848\!\!\cong\!\!750$
$1403\!\!\sim\!\!852\!\!\cong\!\!852$
$1404\!\!\sim\!\!851\!\!\cong\!\!847$
$1405\!\!\sim\!\!928\!\!\cong\!\!820$
$1406\!\!\sim\!\!930\!\!\cong\!\!821$
$1407\!\!\sim\!\!929\!\!\cong\!\!929$
$1408\!\!\sim\!\!930\!\!\cong\!\!821$
$1409\!\!\sim\!\!936\!\!\cong\!\!820$
$1410\!\!\sim\!\!933\!\!\cong\!\!849$
$1411\!\!\sim\!\!929\!\!\cong\!\!929$
$1412\!\!\sim\!\!933\!\!\cong\!\!849$
$1413\!\!\sim\!\!932\!\!\cong\!\!820$
$1414\!\!\sim\!\!1090\!\!\cong\!\!1090$
$1415\!\!\sim\!\!1091\!\!\cong\!\!731$
$1416\!\!\sim\!\!1091\!\!\cong\!\!731$
$1417\!\!\sim\!\!1091\!\!\cong\!\!731$
$1418\!\!\sim\!\!1094\!\!\cong\!\!1090$
$1419\!\!\sim\!\!1094\!\!\cong\!\!1090$
$1420\!\!\sim\!\!1091\!\!\cong\!\!731$
$1421\!\!\sim\!\!1094\!\!\cong\!\!1090$
$1422\!\!\sim\!\!1094\!\!\cong\!\!1090$
$1423\!\!\sim\!\!1090\!\!\cong\!\!1090$
$1424\!\!\sim\!\!1091\!\!\cong\!\!731$
$1425\!\!\sim\!\!1091\!\!\cong\!\!731$
$1426\!\!\sim\!\!1091\!\!\cong\!\!731$
$1427\!\!\sim\!\!1094\!\!\cong\!\!1090$
$1428\!\!\sim\!\!1094\!\!\cong\!\!1090$
$1429\!\!\sim\!\!1091\!\!\cong\!\!731$
$1430\!\!\sim\!\!1094\!\!\cong\!\!1090$
$1431\!\!\sim\!\!1094\!\!\cong\!\!1090$
$1432\!\!\sim\!\!847\!\!\cong\!\!847$
$1433\!\!\sim\!\!849\!\!\cong\!\!849$
$1434\!\!\sim\!\!848\!\!\cong\!\!750$
$1435\!\!\sim\!\!849\!\!\cong\!\!849$
$1436\!\!\sim\!\!855\!\!\cong\!\!847$
$1437\!\!\sim\!\!852\!\!\cong\!\!852$
$1438\!\!\sim\!\!848\!\!\cong\!\!750$
$1439\!\!\sim\!\!852\!\!\cong\!\!852$
$1440\!\!\sim\!\!851\!\!\cong\!\!847$
$1441\!\!\sim\!\!1090\!\!\cong\!\!1090$
$1442\!\!\sim\!\!1091\!\!\cong\!\!731$
$1443\!\!\sim\!\!1091\!\!\cong\!\!731$
$1444\!\!\sim\!\!1091\!\!\cong\!\!731$
$1445\!\!\sim\!\!1094\!\!\cong\!\!1090$
$1446\!\!\sim\!\!1094\!\!\cong\!\!1090$
$1447\!\!\sim\!\!1091\!\!\cong\!\!731$
$1448\!\!\sim\!\!1094\!\!\cong\!\!1090$
$1449\!\!\sim\!\!1094\!\!\cong\!\!1090$
$1450\!\!\sim\!\!1090\!\!\cong\!\!1090$
$1451\!\!\sim\!\!1091\!\!\cong\!\!731$
$1452\!\!\sim\!\!1091\!\!\cong\!\!731$
$1453\!\!\sim\!\!1091\!\!\cong\!\!731$
$1454\!\!\sim\!\!1094\!\!\cong\!\!1090$
$1455\!\!\sim\!\!1094\!\!\cong\!\!1090$
$1456\!\!\sim\!\!1091\!\!\cong\!\!731$
$1457\!\!\sim\!\!1094\!\!\cong\!\!1090$
$1458\!\!\sim\!\!1094\!\!\cong\!\!1090$
$1459\!\!\sim\!\!1094\!\!\cong\!\!1090$
$1460\!\!\sim\!\!972\!\!\cong\!\!739$
$1461\!\!\sim\!\!1094\!\!\cong\!\!1090$
$1462\!\!\sim\!\!972\!\!\cong\!\!739$
$1463\!\!\sim\!\!810\!\!\cong\!\!802$
$1464\!\!\sim\!\!891\!\!\cong\!\!891$
$1465\!\!\sim\!\!1094\!\!\cong\!\!1090$
$1466\!\!\sim\!\!891\!\!\cong\!\!891$
$1467\!\!\sim\!\!1094\!\!\cong\!\!1090$
$1468\!\!\sim\!\!1091\!\!\cong\!\!731$
$1469\!\!\sim\!\!966\!\!\cong\!\!966$
$1470\!\!\sim\!\!1091\!\!\cong\!\!731$
$1471\!\!\sim\!\!966\!\!\cong\!\!966$
$1472\!\!\sim\!\!804\!\!\cong\!\!731$
$1473\!\!\sim\!\!885\!\!\cong\!\!885$
$1474\!\!\sim\!\!1091\!\!\cong\!\!731$
$1475\!\!\sim\!\!885\!\!\cong\!\!885$
$1476\!\!\sim\!\!1091\!\!\cong\!\!731$
$1477\!\!\sim\!\!1094\!\!\cong\!\!1090$
$1478\!\!\sim\!\!969\!\!\cong\!\!969$
$1479\!\!\sim\!\!1094\!\!\cong\!\!1090$
$1480\!\!\sim\!\!969\!\!\cong\!\!969$
$1481\!\!\sim\!\!807\!\!\cong\!\!771$
$1482\!\!\sim\!\!888\!\!\cong\!\!888$
$1483\!\!\sim\!\!1094\!\!\cong\!\!1090$
$1484\!\!\sim\!\!888\!\!\cong\!\!888$
$1485\!\!\sim\!\!1094\!\!\cong\!\!1090$
$1486\!\!\sim\!\!1091\!\!\cong\!\!731$
$1487\!\!\sim\!\!966\!\!\cong\!\!966$
$1488\!\!\sim\!\!1091\!\!\cong\!\!731$
$1489\!\!\sim\!\!966\!\!\cong\!\!966$
$1490\!\!\sim\!\!804\!\!\cong\!\!731$
$1491\!\!\sim\!\!885\!\!\cong\!\!885$
$1492\!\!\sim\!\!1091\!\!\cong\!\!731$
$1493\!\!\sim\!\!885\!\!\cong\!\!885$
$1494\!\!\sim\!\!1091\!\!\cong\!\!731$
$1495\!\!\sim\!\!1090\!\!\cong\!\!1090$
$1496\!\!\sim\!\!964\!\!\cong\!\!739$
$1497\!\!\sim\!\!1090\!\!\cong\!\!1090$
$1498\!\!\sim\!\!964\!\!\cong\!\!739$
$1499\!\!\sim\!\!802\!\!\cong\!\!802$
$1500\!\!\sim\!\!883\!\!\cong\!\!883$
$1501\!\!\sim\!\!1090\!\!\cong\!\!1090$
$1502\!\!\sim\!\!883\!\!\cong\!\!883$
$1503\!\!\sim\!\!1090\!\!\cong\!\!1090$
$1504\!\!\sim\!\!1091\!\!\cong\!\!731$
$1505\!\!\sim\!\!965\!\!\cong\!\!965$
$1506\!\!\sim\!\!1091\!\!\cong\!\!731$
$1507\!\!\sim\!\!965\!\!\cong\!\!965$
$1508\!\!\sim\!\!803\!\!\cong\!\!771$
$1509\!\!\sim\!\!884\!\!\cong\!\!884$
$1510\!\!\sim\!\!1091\!\!\cong\!\!731$
$1511\!\!\sim\!\!884\!\!\cong\!\!884$
$1512\!\!\sim\!\!1091\!\!\cong\!\!731$
$1513\!\!\sim\!\!1094\!\!\cong\!\!1090$
$1514\!\!\sim\!\!969\!\!\cong\!\!969$
$1515\!\!\sim\!\!1094\!\!\cong\!\!1090$
$1516\!\!\sim\!\!969\!\!\cong\!\!969$
$1517\!\!\sim\!\!807\!\!\cong\!\!771$
$1518\!\!\sim\!\!888\!\!\cong\!\!888$
$1519\!\!\sim\!\!1094\!\!\cong\!\!1090$
$1520\!\!\sim\!\!888\!\!\cong\!\!888$
$1521\!\!\sim\!\!1094\!\!\cong\!\!1090$
$1522\!\!\sim\!\!1091\!\!\cong\!\!731$
$1523\!\!\sim\!\!965\!\!\cong\!\!965$
$1524\!\!\sim\!\!1091\!\!\cong\!\!731$
$1525\!\!\sim\!\!965\!\!\cong\!\!965$
$1526\!\!\sim\!\!803\!\!\cong\!\!771$
$1527\!\!\sim\!\!884\!\!\cong\!\!884$
$1528\!\!\sim\!\!1091\!\!\cong\!\!731$
$1529\!\!\sim\!\!884\!\!\cong\!\!884$
$1530\!\!\sim\!\!1091\!\!\cong\!\!731$
$1531\!\!\sim\!\!1094\!\!\cong\!\!1090$
$1532\!\!\sim\!\!968\!\!\cong\!\!968$
$1533\!\!\sim\!\!1094\!\!\cong\!\!1090$
$1534\!\!\sim\!\!968\!\!\cong\!\!968$
$1535\!\!\sim\!\!806\!\!\cong\!\!802$
$1536\!\!\sim\!\!887\!\!\cong\!\!887$
$1537\!\!\sim\!\!1094\!\!\cong\!\!1090$
$1538\!\!\sim\!\!887\!\!\cong\!\!887$
$1539\!\!\sim\!\!1094\!\!\cong\!\!1090$
$1540\!\!\sim\!\!851\!\!\cong\!\!847$
$1541\!\!\sim\!\!824\!\!\cong\!\!820$
$1542\!\!\sim\!\!878\!\!\cong\!\!878$
$1543\!\!\sim\!\!842\!\!\cong\!\!838$
$1544\!\!\sim\!\!756\!\!\cong\!\!748$
$1545\!\!\sim\!\!869\!\!\cong\!\!869$
$1546\!\!\sim\!\!860\!\!\cong\!\!860$
$1547\!\!\sim\!\!824\!\!\cong\!\!820$
$1548\!\!\sim\!\!887\!\!\cong\!\!887$
$1549\!\!\sim\!\!848\!\!\cong\!\!750$
$1550\!\!\sim\!\!821\!\!\cong\!\!821$
$1551\!\!\sim\!\!875\!\!\cong\!\!875$
$1552\!\!\sim\!\!839\!\!\cong\!\!821$
$1553\!\!\sim\!\!750\!\!\cong\!\!750$
$1554\!\!\sim\!\!866\!\!\cong\!\!866$
$1555\!\!\sim\!\!857\!\!\cong\!\!857$
$1556\!\!\sim\!\!821\!\!\cong\!\!821$
$1557\!\!\sim\!\!884\!\!\cong\!\!884$
$1558\!\!\sim\!\!852\!\!\cong\!\!852$
$1559\!\!\sim\!\!824\!\!\cong\!\!820$
$1560\!\!\sim\!\!879\!\!\cong\!\!879$
$1561\!\!\sim\!\!843\!\!\cong\!\!843$
$1562\!\!\sim\!\!753\!\!\cong\!\!753$
$1563\!\!\sim\!\!870\!\!\cong\!\!870$
$1564\!\!\sim\!\!861\!\!\cong\!\!861$
$1565\!\!\sim\!\!824\!\!\cong\!\!820$
$1566\!\!\sim\!\!888\!\!\cong\!\!888$
$1567\!\!\sim\!\!848\!\!\cong\!\!750$
$1568\!\!\sim\!\!821\!\!\cong\!\!821$
$1569\!\!\sim\!\!875\!\!\cong\!\!875$
$1570\!\!\sim\!\!839\!\!\cong\!\!821$
$1571\!\!\sim\!\!750\!\!\cong\!\!750$
$1572\!\!\sim\!\!866\!\!\cong\!\!866$
$1573\!\!\sim\!\!857\!\!\cong\!\!857$
$1574\!\!\sim\!\!821\!\!\cong\!\!821$
$1575\!\!\sim\!\!884\!\!\cong\!\!884$
$1576\!\!\sim\!\!847\!\!\cong\!\!847$
$1577\!\!\sim\!\!820\!\!\cong\!\!820$
$1578\!\!\sim\!\!874\!\!\cong\!\!874$
$1579\!\!\sim\!\!838\!\!\cong\!\!838$
$1580\!\!\sim\!\!748\!\!\cong\!\!748$
$1581\!\!\sim\!\!865\!\!\cong\!\!820$
$1582\!\!\sim\!\!856\!\!\cong\!\!856$
$1583\!\!\sim\!\!820\!\!\cong\!\!820$
$1584\!\!\sim\!\!883\!\!\cong\!\!883$
$1585\!\!\sim\!\!849\!\!\cong\!\!849$
$1586\!\!\sim\!\!821\!\!\cong\!\!821$
$1587\!\!\sim\!\!876\!\!\cong\!\!876$
$1588\!\!\sim\!\!840\!\!\cong\!\!840$
$1589\!\!\sim\!\!749\!\!\cong\!\!749$
$1590\!\!\sim\!\!866\!\!\cong\!\!866$
$1591\!\!\sim\!\!858\!\!\cong\!\!858$
$1592\!\!\sim\!\!821\!\!\cong\!\!821$
$1593\!\!\sim\!\!885\!\!\cong\!\!885$
$1594\!\!\sim\!\!852\!\!\cong\!\!852$
$1595\!\!\sim\!\!824\!\!\cong\!\!820$
$1596\!\!\sim\!\!879\!\!\cong\!\!879$
$1597\!\!\sim\!\!843\!\!\cong\!\!843$
$1598\!\!\sim\!\!753\!\!\cong\!\!753$
$1599\!\!\sim\!\!870\!\!\cong\!\!870$
$1600\!\!\sim\!\!861\!\!\cong\!\!861$
$1601\!\!\sim\!\!824\!\!\cong\!\!820$
$1602\!\!\sim\!\!888\!\!\cong\!\!888$
$1603\!\!\sim\!\!849\!\!\cong\!\!849$
$1604\!\!\sim\!\!821\!\!\cong\!\!821$
$1605\!\!\sim\!\!876\!\!\cong\!\!876$
$1606\!\!\sim\!\!840\!\!\cong\!\!840$
$1607\!\!\sim\!\!749\!\!\cong\!\!749$
$1608\!\!\sim\!\!866\!\!\cong\!\!866$
$1609\!\!\sim\!\!858\!\!\cong\!\!858$
$1610\!\!\sim\!\!821\!\!\cong\!\!821$
$1611\!\!\sim\!\!885\!\!\cong\!\!885$
$1612\!\!\sim\!\!855\!\!\cong\!\!847$
$1613\!\!\sim\!\!824\!\!\cong\!\!820$
$1614\!\!\sim\!\!882\!\!\cong\!\!882$
$1615\!\!\sim\!\!846\!\!\cong\!\!846$
$1616\!\!\sim\!\!752\!\!\cong\!\!752$
$1617\!\!\sim\!\!869\!\!\cong\!\!869$
$1618\!\!\sim\!\!864\!\!\cong\!\!864$
$1619\!\!\sim\!\!824\!\!\cong\!\!820$
$1620\!\!\sim\!\!891\!\!\cong\!\!891$
$1621\!\!\sim\!\!1094\!\!\cong\!\!1090$
$1622\!\!\sim\!\!945\!\!\cong\!\!941$
$1623\!\!\sim\!\!1094\!\!\cong\!\!1090$
$1624\!\!\sim\!\!963\!\!\cong\!\!963$
$1625\!\!\sim\!\!783\!\!\cong\!\!775$
$1626\!\!\sim\!\!882\!\!\cong\!\!882$
$1627\!\!\sim\!\!1094\!\!\cong\!\!1090$
$1628\!\!\sim\!\!864\!\!\cong\!\!864$
$1629\!\!\sim\!\!1094\!\!\cong\!\!1090$
$1630\!\!\sim\!\!1091\!\!\cong\!\!731$
$1631\!\!\sim\!\!939\!\!\cong\!\!939$
$1632\!\!\sim\!\!1091\!\!\cong\!\!731$
$1633\!\!\sim\!\!957\!\!\cong\!\!957$
$1634\!\!\sim\!\!777\!\!\cong\!\!777$
$1635\!\!\sim\!\!876\!\!\cong\!\!876$
$1636\!\!\sim\!\!1091\!\!\cong\!\!731$
$1637\!\!\sim\!\!858\!\!\cong\!\!858$
$1638\!\!\sim\!\!1091\!\!\cong\!\!731$
$1639\!\!\sim\!\!1094\!\!\cong\!\!1090$
$1640\!\!\sim\!\!942\!\!\cong\!\!942$
$1641\!\!\sim\!\!1094\!\!\cong\!\!1090$
$1642\!\!\sim\!\!960\!\!\cong\!\!960$
$1643\!\!\sim\!\!780\!\!\cong\!\!780$
$1644\!\!\sim\!\!879\!\!\cong\!\!879$
$1645\!\!\sim\!\!1094\!\!\cong\!\!1090$
$1646\!\!\sim\!\!861\!\!\cong\!\!861$
$1647\!\!\sim\!\!1094\!\!\cong\!\!1090$
$1648\!\!\sim\!\!1091\!\!\cong\!\!731$
$1649\!\!\sim\!\!939\!\!\cong\!\!939$
$1650\!\!\sim\!\!1091\!\!\cong\!\!731$
$1651\!\!\sim\!\!957\!\!\cong\!\!957$
$1652\!\!\sim\!\!777\!\!\cong\!\!777$
$1653\!\!\sim\!\!876\!\!\cong\!\!876$
$1654\!\!\sim\!\!1091\!\!\cong\!\!731$
$1655\!\!\sim\!\!858\!\!\cong\!\!858$
$1656\!\!\sim\!\!1091\!\!\cong\!\!731$
$1657\!\!\sim\!\!1090\!\!\cong\!\!1090$
$1658\!\!\sim\!\!937\!\!\cong\!\!937$
$1659\!\!\sim\!\!1090\!\!\cong\!\!1090$
$1660\!\!\sim\!\!955\!\!\cong\!\!937$
$1661\!\!\sim\!\!775\!\!\cong\!\!775$
$1662\!\!\sim\!\!874\!\!\cong\!\!874$
$1663\!\!\sim\!\!1090\!\!\cong\!\!1090$
$1664\!\!\sim\!\!856\!\!\cong\!\!856$
$1665\!\!\sim\!\!1090\!\!\cong\!\!1090$
$1666\!\!\sim\!\!1091\!\!\cong\!\!731$
$1667\!\!\sim\!\!938\!\!\cong\!\!938$
$1668\!\!\sim\!\!1091\!\!\cong\!\!731$
$1669\!\!\sim\!\!956\!\!\cong\!\!956$
$1670\!\!\sim\!\!776\!\!\cong\!\!776$
$1671\!\!\sim\!\!875\!\!\cong\!\!875$
$1672\!\!\sim\!\!1091\!\!\cong\!\!731$
$1673\!\!\sim\!\!857\!\!\cong\!\!857$
$1674\!\!\sim\!\!1091\!\!\cong\!\!731$
$1675\!\!\sim\!\!1094\!\!\cong\!\!1090$
$1676\!\!\sim\!\!942\!\!\cong\!\!942$
$1677\!\!\sim\!\!1094\!\!\cong\!\!1090$
$1678\!\!\sim\!\!960\!\!\cong\!\!960$
$1679\!\!\sim\!\!780\!\!\cong\!\!780$
$1680\!\!\sim\!\!879\!\!\cong\!\!879$
$1681\!\!\sim\!\!1094\!\!\cong\!\!1090$
$1682\!\!\sim\!\!861\!\!\cong\!\!861$
$1683\!\!\sim\!\!1094\!\!\cong\!\!1090$
$1684\!\!\sim\!\!1091\!\!\cong\!\!731$
$1685\!\!\sim\!\!938\!\!\cong\!\!938$
$1686\!\!\sim\!\!1091\!\!\cong\!\!731$
$1687\!\!\sim\!\!956\!\!\cong\!\!956$
$1688\!\!\sim\!\!776\!\!\cong\!\!776$
$1689\!\!\sim\!\!875\!\!\cong\!\!875$
$1690\!\!\sim\!\!1091\!\!\cong\!\!731$
$1691\!\!\sim\!\!857\!\!\cong\!\!857$
$1692\!\!\sim\!\!1091\!\!\cong\!\!731$
$1693\!\!\sim\!\!1094\!\!\cong\!\!1090$
$1694\!\!\sim\!\!941\!\!\cong\!\!941$
$1695\!\!\sim\!\!1094\!\!\cong\!\!1090$
$1696\!\!\sim\!\!959\!\!\cong\!\!959$
$1697\!\!\sim\!\!779\!\!\cong\!\!779$
$1698\!\!\sim\!\!878\!\!\cong\!\!878$
$1699\!\!\sim\!\!1094\!\!\cong\!\!1090$
$1700\!\!\sim\!\!860\!\!\cong\!\!860$
$1701\!\!\sim\!\!1094\!\!\cong\!\!1090$
$1702\!\!\sim\!\!851\!\!\cong\!\!847$
$1703\!\!\sim\!\!842\!\!\cong\!\!838$
$1704\!\!\sim\!\!860\!\!\cong\!\!860$
$1705\!\!\sim\!\!824\!\!\cong\!\!820$
$1706\!\!\sim\!\!756\!\!\cong\!\!748$
$1707\!\!\sim\!\!824\!\!\cong\!\!820$
$1708\!\!\sim\!\!878\!\!\cong\!\!878$
$1709\!\!\sim\!\!869\!\!\cong\!\!869$
$1710\!\!\sim\!\!887\!\!\cong\!\!887$
$1711\!\!\sim\!\!848\!\!\cong\!\!750$
$1712\!\!\sim\!\!839\!\!\cong\!\!821$
$1713\!\!\sim\!\!857\!\!\cong\!\!857$
$1714\!\!\sim\!\!821\!\!\cong\!\!821$
$1715\!\!\sim\!\!750\!\!\cong\!\!750$
$1716\!\!\sim\!\!821\!\!\cong\!\!821$
$1717\!\!\sim\!\!875\!\!\cong\!\!875$
$1718\!\!\sim\!\!866\!\!\cong\!\!866$
$1719\!\!\sim\!\!884\!\!\cong\!\!884$
$1720\!\!\sim\!\!852\!\!\cong\!\!852$
$1721\!\!\sim\!\!843\!\!\cong\!\!843$
$1722\!\!\sim\!\!861\!\!\cong\!\!861$
$1723\!\!\sim\!\!824\!\!\cong\!\!820$
$1724\!\!\sim\!\!753\!\!\cong\!\!753$
$1725\!\!\sim\!\!824\!\!\cong\!\!820$
$1726\!\!\sim\!\!879\!\!\cong\!\!879$
$1727\!\!\sim\!\!870\!\!\cong\!\!870$
$1728\!\!\sim\!\!888\!\!\cong\!\!888$
$1729\!\!\sim\!\!848\!\!\cong\!\!750$
$1730\!\!\sim\!\!839\!\!\cong\!\!821$
$1731\!\!\sim\!\!857\!\!\cong\!\!857$
$1732\!\!\sim\!\!821\!\!\cong\!\!821$
$1733\!\!\sim\!\!750\!\!\cong\!\!750$
$1734\!\!\sim\!\!821\!\!\cong\!\!821$
$1735\!\!\sim\!\!875\!\!\cong\!\!875$
$1736\!\!\sim\!\!866\!\!\cong\!\!866$
$1737\!\!\sim\!\!884\!\!\cong\!\!884$
$1738\!\!\sim\!\!847\!\!\cong\!\!847$
$1739\!\!\sim\!\!838\!\!\cong\!\!838$
$1740\!\!\sim\!\!856\!\!\cong\!\!856$
$1741\!\!\sim\!\!820\!\!\cong\!\!820$
$1742\!\!\sim\!\!748\!\!\cong\!\!748$
$1743\!\!\sim\!\!820\!\!\cong\!\!820$
$1744\!\!\sim\!\!874\!\!\cong\!\!874$
$1745\!\!\sim\!\!865\!\!\cong\!\!820$
$1746\!\!\sim\!\!883\!\!\cong\!\!883$
$1747\!\!\sim\!\!849\!\!\cong\!\!849$
$1748\!\!\sim\!\!840\!\!\cong\!\!840$
$1749\!\!\sim\!\!858\!\!\cong\!\!858$
$1750\!\!\sim\!\!821\!\!\cong\!\!821$
$1751\!\!\sim\!\!749\!\!\cong\!\!749$
$1752\!\!\sim\!\!821\!\!\cong\!\!821$
$1753\!\!\sim\!\!876\!\!\cong\!\!876$
$1754\!\!\sim\!\!866\!\!\cong\!\!866$
$1755\!\!\sim\!\!885\!\!\cong\!\!885$
$1756\!\!\sim\!\!852\!\!\cong\!\!852$
$1757\!\!\sim\!\!843\!\!\cong\!\!843$
$1758\!\!\sim\!\!861\!\!\cong\!\!861$
$1759\!\!\sim\!\!824\!\!\cong\!\!820$
$1760\!\!\sim\!\!753\!\!\cong\!\!753$
$1761\!\!\sim\!\!824\!\!\cong\!\!820$
$1762\!\!\sim\!\!879\!\!\cong\!\!879$
$1763\!\!\sim\!\!870\!\!\cong\!\!870$
$1764\!\!\sim\!\!888\!\!\cong\!\!888$
$1765\!\!\sim\!\!849\!\!\cong\!\!849$
$1766\!\!\sim\!\!840\!\!\cong\!\!840$
$1767\!\!\sim\!\!858\!\!\cong\!\!858$
$1768\!\!\sim\!\!821\!\!\cong\!\!821$
$1769\!\!\sim\!\!749\!\!\cong\!\!749$
$1770\!\!\sim\!\!821\!\!\cong\!\!821$
$1771\!\!\sim\!\!876\!\!\cong\!\!876$
$1772\!\!\sim\!\!866\!\!\cong\!\!866$
$1773\!\!\sim\!\!885\!\!\cong\!\!885$
$1774\!\!\sim\!\!855\!\!\cong\!\!847$
$1775\!\!\sim\!\!846\!\!\cong\!\!846$
$1776\!\!\sim\!\!864\!\!\cong\!\!864$
$1777\!\!\sim\!\!824\!\!\cong\!\!820$
$1778\!\!\sim\!\!752\!\!\cong\!\!752$
$1779\!\!\sim\!\!824\!\!\cong\!\!820$
$1780\!\!\sim\!\!882\!\!\cong\!\!882$
$1781\!\!\sim\!\!869\!\!\cong\!\!869$
$1782\!\!\sim\!\!891\!\!\cong\!\!891$
$1783\!\!\sim\!\!770\!\!\cong\!\!730$
$1784\!\!\sim\!\!743\!\!\cong\!\!739$
$1785\!\!\sim\!\!779\!\!\cong\!\!779$
$1786\!\!\sim\!\!743\!\!\cong\!\!739$
$1787\!\!\sim\!\!734\!\!\cong\!\!730$
$1788\!\!\sim\!\!752\!\!\cong\!\!752$
$1789\!\!\sim\!\!779\!\!\cong\!\!779$
$1790\!\!\sim\!\!752\!\!\cong\!\!752$
$1791\!\!\sim\!\!806\!\!\cong\!\!802$
$1792\!\!\sim\!\!767\!\!\cong\!\!731$
$1793\!\!\sim\!\!740\!\!\cong\!\!740$
$1794\!\!\sim\!\!776\!\!\cong\!\!776$
$1795\!\!\sim\!\!740\!\!\cong\!\!740$
$1796\!\!\sim\!\!731\!\!\cong\!\!731$
$1797\!\!\sim\!\!749\!\!\cong\!\!749$
$1798\!\!\sim\!\!776\!\!\cong\!\!776$
$1799\!\!\sim\!\!749\!\!\cong\!\!749$
$1800\!\!\sim\!\!803\!\!\cong\!\!771$
$1801\!\!\sim\!\!771\!\!\cong\!\!771$
$1802\!\!\sim\!\!744\!\!\cong\!\!744$
$1803\!\!\sim\!\!780\!\!\cong\!\!780$
$1804\!\!\sim\!\!744\!\!\cong\!\!744$
$1805\!\!\sim\!\!734\!\!\cong\!\!730$
$1806\!\!\sim\!\!753\!\!\cong\!\!753$
$1807\!\!\sim\!\!780\!\!\cong\!\!780$
$1808\!\!\sim\!\!753\!\!\cong\!\!753$
$1809\!\!\sim\!\!807\!\!\cong\!\!771$
$1810\!\!\sim\!\!767\!\!\cong\!\!731$
$1811\!\!\sim\!\!740\!\!\cong\!\!740$
$1812\!\!\sim\!\!776\!\!\cong\!\!776$
$1813\!\!\sim\!\!740\!\!\cong\!\!740$
$1814\!\!\sim\!\!731\!\!\cong\!\!731$
$1815\!\!\sim\!\!749\!\!\cong\!\!749$
$1816\!\!\sim\!\!776\!\!\cong\!\!776$
$1817\!\!\sim\!\!749\!\!\cong\!\!749$
$1818\!\!\sim\!\!803\!\!\cong\!\!771$
$1819\!\!\sim\!\!766\!\!\cong\!\!730$
$1820\!\!\sim\!\!739\!\!\cong\!\!739$
$1821\!\!\sim\!\!775\!\!\cong\!\!775$
$1822\!\!\sim\!\!739\!\!\cong\!\!739$
$1823\!\!\sim\!\!730\!\!\cong\!\!730$
$1824\!\!\sim\!\!748\!\!\cong\!\!748$
$1825\!\!\sim\!\!775\!\!\cong\!\!775$
$1826\!\!\sim\!\!748\!\!\cong\!\!748$
$1827\!\!\sim\!\!802\!\!\cong\!\!802$
$1828\!\!\sim\!\!768\!\!\cong\!\!731$
$1829\!\!\sim\!\!741\!\!\cong\!\!741$
$1830\!\!\sim\!\!777\!\!\cong\!\!777$
$1831\!\!\sim\!\!741\!\!\cong\!\!741$
$1832\!\!\sim\!\!731\!\!\cong\!\!731$
$1833\!\!\sim\!\!750\!\!\cong\!\!750$
$1834\!\!\sim\!\!777\!\!\cong\!\!777$
$1835\!\!\sim\!\!750\!\!\cong\!\!750$
$1836\!\!\sim\!\!804\!\!\cong\!\!731$
$1837\!\!\sim\!\!771\!\!\cong\!\!771$
$1838\!\!\sim\!\!744\!\!\cong\!\!744$
$1839\!\!\sim\!\!780\!\!\cong\!\!780$
$1840\!\!\sim\!\!744\!\!\cong\!\!744$
$1841\!\!\sim\!\!734\!\!\cong\!\!730$
$1842\!\!\sim\!\!753\!\!\cong\!\!753$
$1843\!\!\sim\!\!780\!\!\cong\!\!780$
$1844\!\!\sim\!\!753\!\!\cong\!\!753$
$1845\!\!\sim\!\!807\!\!\cong\!\!771$
$1846\!\!\sim\!\!768\!\!\cong\!\!731$
$1847\!\!\sim\!\!741\!\!\cong\!\!741$
$1848\!\!\sim\!\!777\!\!\cong\!\!777$
$1849\!\!\sim\!\!741\!\!\cong\!\!741$
$1850\!\!\sim\!\!731\!\!\cong\!\!731$
$1851\!\!\sim\!\!750\!\!\cong\!\!750$
$1852\!\!\sim\!\!777\!\!\cong\!\!777$
$1853\!\!\sim\!\!750\!\!\cong\!\!750$
$1854\!\!\sim\!\!804\!\!\cong\!\!731$
$1855\!\!\sim\!\!774\!\!\cong\!\!730$
$1856\!\!\sim\!\!747\!\!\cong\!\!739$
$1857\!\!\sim\!\!783\!\!\cong\!\!775$
$1858\!\!\sim\!\!747\!\!\cong\!\!739$
$1859\!\!\sim\!\!734\!\!\cong\!\!730$
$1860\!\!\sim\!\!756\!\!\cong\!\!748$
$1861\!\!\sim\!\!783\!\!\cong\!\!775$
$1862\!\!\sim\!\!756\!\!\cong\!\!748$
$1863\!\!\sim\!\!810\!\!\cong\!\!802$
$1864\!\!\sim\!\!932\!\!\cong\!\!820$
$1865\!\!\sim\!\!923\!\!\cong\!\!923$
$1866\!\!\sim\!\!941\!\!\cong\!\!941$
$1867\!\!\sim\!\!824\!\!\cong\!\!820$
$1868\!\!\sim\!\!747\!\!\cong\!\!739$
$1869\!\!\sim\!\!824\!\!\cong\!\!820$
$1870\!\!\sim\!\!959\!\!\cong\!\!959$
$1871\!\!\sim\!\!846\!\!\cong\!\!846$
$1872\!\!\sim\!\!968\!\!\cong\!\!968$
$1873\!\!\sim\!\!929\!\!\cong\!\!929$
$1874\!\!\sim\!\!920\!\!\cong\!\!920$
$1875\!\!\sim\!\!938\!\!\cong\!\!938$
$1876\!\!\sim\!\!821\!\!\cong\!\!821$
$1877\!\!\sim\!\!741\!\!\cong\!\!741$
$1878\!\!\sim\!\!821\!\!\cong\!\!821$
$1879\!\!\sim\!\!956\!\!\cong\!\!956$
$1880\!\!\sim\!\!840\!\!\cong\!\!840$
$1881\!\!\sim\!\!965\!\!\cong\!\!965$
$1882\!\!\sim\!\!933\!\!\cong\!\!849$
$1883\!\!\sim\!\!924\!\!\cong\!\!870$
$1884\!\!\sim\!\!942\!\!\cong\!\!942$
$1885\!\!\sim\!\!824\!\!\cong\!\!820$
$1886\!\!\sim\!\!744\!\!\cong\!\!744$
$1887\!\!\sim\!\!824\!\!\cong\!\!820$
$1888\!\!\sim\!\!960\!\!\cong\!\!960$
$1889\!\!\sim\!\!843\!\!\cong\!\!843$
$1890\!\!\sim\!\!969\!\!\cong\!\!969$
$1891\!\!\sim\!\!929\!\!\cong\!\!929$
$1892\!\!\sim\!\!920\!\!\cong\!\!920$
$1893\!\!\sim\!\!938\!\!\cong\!\!938$
$1894\!\!\sim\!\!821\!\!\cong\!\!821$
$1895\!\!\sim\!\!741\!\!\cong\!\!741$
$1896\!\!\sim\!\!821\!\!\cong\!\!821$
$1897\!\!\sim\!\!956\!\!\cong\!\!956$
$1898\!\!\sim\!\!840\!\!\cong\!\!840$
$1899\!\!\sim\!\!965\!\!\cong\!\!965$
$1900\!\!\sim\!\!928\!\!\cong\!\!820$
$1901\!\!\sim\!\!919\!\!\cong\!\!820$
$1902\!\!\sim\!\!937\!\!\cong\!\!937$
$1903\!\!\sim\!\!820\!\!\cong\!\!820$
$1904\!\!\sim\!\!739\!\!\cong\!\!739$
$1905\!\!\sim\!\!820\!\!\cong\!\!820$
$1906\!\!\sim\!\!955\!\!\cong\!\!937$
$1907\!\!\sim\!\!838\!\!\cong\!\!838$
$1908\!\!\sim\!\!964\!\!\cong\!\!739$
$1909\!\!\sim\!\!930\!\!\cong\!\!821$
$1910\!\!\sim\!\!920\!\!\cong\!\!920$
$1911\!\!\sim\!\!939\!\!\cong\!\!939$
$1912\!\!\sim\!\!821\!\!\cong\!\!821$
$1913\!\!\sim\!\!740\!\!\cong\!\!740$
$1914\!\!\sim\!\!821\!\!\cong\!\!821$
$1915\!\!\sim\!\!957\!\!\cong\!\!957$
$1916\!\!\sim\!\!839\!\!\cong\!\!821$
$1917\!\!\sim\!\!966\!\!\cong\!\!966$
$1918\!\!\sim\!\!933\!\!\cong\!\!849$
$1919\!\!\sim\!\!924\!\!\cong\!\!870$
$1920\!\!\sim\!\!942\!\!\cong\!\!942$
$1921\!\!\sim\!\!824\!\!\cong\!\!820$
$1922\!\!\sim\!\!744\!\!\cong\!\!744$
$1923\!\!\sim\!\!824\!\!\cong\!\!820$
$1924\!\!\sim\!\!960\!\!\cong\!\!960$
$1925\!\!\sim\!\!843\!\!\cong\!\!843$
$1926\!\!\sim\!\!969\!\!\cong\!\!969$
$1927\!\!\sim\!\!930\!\!\cong\!\!821$
$1928\!\!\sim\!\!920\!\!\cong\!\!920$
$1929\!\!\sim\!\!939\!\!\cong\!\!939$
$1930\!\!\sim\!\!821\!\!\cong\!\!821$
$1931\!\!\sim\!\!740\!\!\cong\!\!740$
$1932\!\!\sim\!\!821\!\!\cong\!\!821$
$1933\!\!\sim\!\!957\!\!\cong\!\!957$
$1934\!\!\sim\!\!839\!\!\cong\!\!821$
$1935\!\!\sim\!\!966\!\!\cong\!\!966$
$1936\!\!\sim\!\!936\!\!\cong\!\!820$
$1937\!\!\sim\!\!923\!\!\cong\!\!923$
$1938\!\!\sim\!\!945\!\!\cong\!\!941$
$1939\!\!\sim\!\!824\!\!\cong\!\!820$
$1940\!\!\sim\!\!743\!\!\cong\!\!739$
$1941\!\!\sim\!\!824\!\!\cong\!\!820$
$1942\!\!\sim\!\!963\!\!\cong\!\!963$
$1943\!\!\sim\!\!842\!\!\cong\!\!838$
$1944\!\!\sim\!\!972\!\!\cong\!\!739$
$1945\!\!\sim\!\!1094\!\!\cong\!\!1090$
$1946\!\!\sim\!\!963\!\!\cong\!\!963$
$1947\!\!\sim\!\!1094\!\!\cong\!\!1090$
$1948\!\!\sim\!\!945\!\!\cong\!\!941$
$1949\!\!\sim\!\!783\!\!\cong\!\!775$
$1950\!\!\sim\!\!864\!\!\cong\!\!864$
$1951\!\!\sim\!\!1094\!\!\cong\!\!1090$
$1952\!\!\sim\!\!882\!\!\cong\!\!882$
$1953\!\!\sim\!\!1094\!\!\cong\!\!1090$
$1954\!\!\sim\!\!1091\!\!\cong\!\!731$
$1955\!\!\sim\!\!957\!\!\cong\!\!957$
$1956\!\!\sim\!\!1091\!\!\cong\!\!731$
$1957\!\!\sim\!\!939\!\!\cong\!\!939$
$1958\!\!\sim\!\!777\!\!\cong\!\!777$
$1959\!\!\sim\!\!858\!\!\cong\!\!858$
$1960\!\!\sim\!\!1091\!\!\cong\!\!731$
$1961\!\!\sim\!\!876\!\!\cong\!\!876$
$1962\!\!\sim\!\!1091\!\!\cong\!\!731$
$1963\!\!\sim\!\!1094\!\!\cong\!\!1090$
$1964\!\!\sim\!\!960\!\!\cong\!\!960$
$1965\!\!\sim\!\!1094\!\!\cong\!\!1090$
$1966\!\!\sim\!\!942\!\!\cong\!\!942$
$1967\!\!\sim\!\!780\!\!\cong\!\!780$
$1968\!\!\sim\!\!861\!\!\cong\!\!861$
$1969\!\!\sim\!\!1094\!\!\cong\!\!1090$
$1970\!\!\sim\!\!879\!\!\cong\!\!879$
$1971\!\!\sim\!\!1094\!\!\cong\!\!1090$
$1972\!\!\sim\!\!1091\!\!\cong\!\!731$
$1973\!\!\sim\!\!957\!\!\cong\!\!957$
$1974\!\!\sim\!\!1091\!\!\cong\!\!731$
$1975\!\!\sim\!\!939\!\!\cong\!\!939$
$1976\!\!\sim\!\!777\!\!\cong\!\!777$
$1977\!\!\sim\!\!858\!\!\cong\!\!858$
$1978\!\!\sim\!\!1091\!\!\cong\!\!731$
$1979\!\!\sim\!\!876\!\!\cong\!\!876$
$1980\!\!\sim\!\!1091\!\!\cong\!\!731$
$1981\!\!\sim\!\!1090\!\!\cong\!\!1090$
$1982\!\!\sim\!\!955\!\!\cong\!\!937$
$1983\!\!\sim\!\!1090\!\!\cong\!\!1090$
$1984\!\!\sim\!\!937\!\!\cong\!\!937$
$1985\!\!\sim\!\!775\!\!\cong\!\!775$
$1986\!\!\sim\!\!856\!\!\cong\!\!856$
$1987\!\!\sim\!\!1090\!\!\cong\!\!1090$
$1988\!\!\sim\!\!874\!\!\cong\!\!874$
$1989\!\!\sim\!\!1090\!\!\cong\!\!1090$
$1990\!\!\sim\!\!1091\!\!\cong\!\!731$
$1991\!\!\sim\!\!956\!\!\cong\!\!956$
$1992\!\!\sim\!\!1091\!\!\cong\!\!731$
$1993\!\!\sim\!\!938\!\!\cong\!\!938$
$1994\!\!\sim\!\!776\!\!\cong\!\!776$
$1995\!\!\sim\!\!857\!\!\cong\!\!857$
$1996\!\!\sim\!\!1091\!\!\cong\!\!731$
$1997\!\!\sim\!\!875\!\!\cong\!\!875$
$1998\!\!\sim\!\!1091\!\!\cong\!\!731$
$1999\!\!\sim\!\!1094\!\!\cong\!\!1090$
$2000\!\!\sim\!\!960\!\!\cong\!\!960$
$2001\!\!\sim\!\!1094\!\!\cong\!\!1090$
$2002\!\!\sim\!\!942\!\!\cong\!\!942$
$2003\!\!\sim\!\!780\!\!\cong\!\!780$
$2004\!\!\sim\!\!861\!\!\cong\!\!861$
$2005\!\!\sim\!\!1094\!\!\cong\!\!1090$
$2006\!\!\sim\!\!879\!\!\cong\!\!879$
$2007\!\!\sim\!\!1094\!\!\cong\!\!1090$
$2008\!\!\sim\!\!1091\!\!\cong\!\!731$
$2009\!\!\sim\!\!956\!\!\cong\!\!956$
$2010\!\!\sim\!\!1091\!\!\cong\!\!731$
$2011\!\!\sim\!\!938\!\!\cong\!\!938$
$2012\!\!\sim\!\!776\!\!\cong\!\!776$
$2013\!\!\sim\!\!857\!\!\cong\!\!857$
$2014\!\!\sim\!\!1091\!\!\cong\!\!731$
$2015\!\!\sim\!\!875\!\!\cong\!\!875$
$2016\!\!\sim\!\!1091\!\!\cong\!\!731$
$2017\!\!\sim\!\!1094\!\!\cong\!\!1090$
$2018\!\!\sim\!\!959\!\!\cong\!\!959$
$2019\!\!\sim\!\!1094\!\!\cong\!\!1090$
$2020\!\!\sim\!\!941\!\!\cong\!\!941$
$2021\!\!\sim\!\!779\!\!\cong\!\!779$
$2022\!\!\sim\!\!860\!\!\cong\!\!860$
$2023\!\!\sim\!\!1094\!\!\cong\!\!1090$
$2024\!\!\sim\!\!878\!\!\cong\!\!878$
$2025\!\!\sim\!\!1094\!\!\cong\!\!1090$
$2026\!\!\sim\!\!932\!\!\cong\!\!820$
$2027\!\!\sim\!\!824\!\!\cong\!\!820$
$2028\!\!\sim\!\!959\!\!\cong\!\!959$
$2029\!\!\sim\!\!923\!\!\cong\!\!923$
$2030\!\!\sim\!\!747\!\!\cong\!\!739$
$2031\!\!\sim\!\!846\!\!\cong\!\!846$
$2032\!\!\sim\!\!941\!\!\cong\!\!941$
$2033\!\!\sim\!\!824\!\!\cong\!\!820$
$2034\!\!\sim\!\!968\!\!\cong\!\!968$
$2035\!\!\sim\!\!929\!\!\cong\!\!929$
$2036\!\!\sim\!\!821\!\!\cong\!\!821$
$2037\!\!\sim\!\!956\!\!\cong\!\!956$
$2038\!\!\sim\!\!920\!\!\cong\!\!920$
$2039\!\!\sim\!\!741\!\!\cong\!\!741$
$2040\!\!\sim\!\!840\!\!\cong\!\!840$
$2041\!\!\sim\!\!938\!\!\cong\!\!938$
$2042\!\!\sim\!\!821\!\!\cong\!\!821$
$2043\!\!\sim\!\!965\!\!\cong\!\!965$
$2044\!\!\sim\!\!933\!\!\cong\!\!849$
$2045\!\!\sim\!\!824\!\!\cong\!\!820$
$2046\!\!\sim\!\!960\!\!\cong\!\!960$
$2047\!\!\sim\!\!924\!\!\cong\!\!870$
$2048\!\!\sim\!\!744\!\!\cong\!\!744$
$2049\!\!\sim\!\!843\!\!\cong\!\!843$
$2050\!\!\sim\!\!942\!\!\cong\!\!942$
$2051\!\!\sim\!\!824\!\!\cong\!\!820$
$2052\!\!\sim\!\!969\!\!\cong\!\!969$
$2053\!\!\sim\!\!929\!\!\cong\!\!929$
$2054\!\!\sim\!\!821\!\!\cong\!\!821$
$2055\!\!\sim\!\!956\!\!\cong\!\!956$
$2056\!\!\sim\!\!920\!\!\cong\!\!920$
$2057\!\!\sim\!\!741\!\!\cong\!\!741$
$2058\!\!\sim\!\!840\!\!\cong\!\!840$
$2059\!\!\sim\!\!938\!\!\cong\!\!938$
$2060\!\!\sim\!\!821\!\!\cong\!\!821$
$2061\!\!\sim\!\!965\!\!\cong\!\!965$
$2062\!\!\sim\!\!928\!\!\cong\!\!820$
$2063\!\!\sim\!\!820\!\!\cong\!\!820$
$2064\!\!\sim\!\!955\!\!\cong\!\!937$
$2065\!\!\sim\!\!919\!\!\cong\!\!820$
$2066\!\!\sim\!\!739\!\!\cong\!\!739$
$2067\!\!\sim\!\!838\!\!\cong\!\!838$
$2068\!\!\sim\!\!937\!\!\cong\!\!937$
$2069\!\!\sim\!\!820\!\!\cong\!\!820$
$2070\!\!\sim\!\!964\!\!\cong\!\!739$
$2071\!\!\sim\!\!930\!\!\cong\!\!821$
$2072\!\!\sim\!\!821\!\!\cong\!\!821$
$2073\!\!\sim\!\!957\!\!\cong\!\!957$
$2074\!\!\sim\!\!920\!\!\cong\!\!920$
$2075\!\!\sim\!\!740\!\!\cong\!\!740$
$2076\!\!\sim\!\!839\!\!\cong\!\!821$
$2077\!\!\sim\!\!939\!\!\cong\!\!939$
$2078\!\!\sim\!\!821\!\!\cong\!\!821$
$2079\!\!\sim\!\!966\!\!\cong\!\!966$
$2080\!\!\sim\!\!933\!\!\cong\!\!849$
$2081\!\!\sim\!\!824\!\!\cong\!\!820$
$2082\!\!\sim\!\!960\!\!\cong\!\!960$
$2083\!\!\sim\!\!924\!\!\cong\!\!870$
$2084\!\!\sim\!\!744\!\!\cong\!\!744$
$2085\!\!\sim\!\!843\!\!\cong\!\!843$
$2086\!\!\sim\!\!942\!\!\cong\!\!942$
$2087\!\!\sim\!\!824\!\!\cong\!\!820$
$2088\!\!\sim\!\!969\!\!\cong\!\!969$
$2089\!\!\sim\!\!930\!\!\cong\!\!821$
$2090\!\!\sim\!\!821\!\!\cong\!\!821$
$2091\!\!\sim\!\!957\!\!\cong\!\!957$
$2092\!\!\sim\!\!920\!\!\cong\!\!920$
$2093\!\!\sim\!\!740\!\!\cong\!\!740$
$2094\!\!\sim\!\!839\!\!\cong\!\!821$
$2095\!\!\sim\!\!939\!\!\cong\!\!939$
$2096\!\!\sim\!\!821\!\!\cong\!\!821$
$2097\!\!\sim\!\!966\!\!\cong\!\!966$
$2098\!\!\sim\!\!936\!\!\cong\!\!820$
$2099\!\!\sim\!\!824\!\!\cong\!\!820$
$2100\!\!\sim\!\!963\!\!\cong\!\!963$
$2101\!\!\sim\!\!923\!\!\cong\!\!923$
$2102\!\!\sim\!\!743\!\!\cong\!\!739$
$2103\!\!\sim\!\!842\!\!\cong\!\!838$
$2104\!\!\sim\!\!945\!\!\cong\!\!941$
$2105\!\!\sim\!\!824\!\!\cong\!\!820$
$2106\!\!\sim\!\!972\!\!\cong\!\!739$
$2107\!\!\sim\!\!1094\!\!\cong\!\!1090$
$2108\!\!\sim\!\!936\!\!\cong\!\!820$
$2109\!\!\sim\!\!1094\!\!\cong\!\!1090$
$2110\!\!\sim\!\!936\!\!\cong\!\!820$
$2111\!\!\sim\!\!774\!\!\cong\!\!730$
$2112\!\!\sim\!\!855\!\!\cong\!\!847$
$2113\!\!\sim\!\!1094\!\!\cong\!\!1090$
$2114\!\!\sim\!\!855\!\!\cong\!\!847$
$2115\!\!\sim\!\!1094\!\!\cong\!\!1090$
$2116\!\!\sim\!\!1091\!\!\cong\!\!731$
$2117\!\!\sim\!\!930\!\!\cong\!\!821$
$2118\!\!\sim\!\!1091\!\!\cong\!\!731$
$2119\!\!\sim\!\!930\!\!\cong\!\!821$
$2120\!\!\sim\!\!768\!\!\cong\!\!731$
$2121\!\!\sim\!\!849\!\!\cong\!\!849$
$2122\!\!\sim\!\!1091\!\!\cong\!\!731$
$2123\!\!\sim\!\!849\!\!\cong\!\!849$
$2124\!\!\sim\!\!1091\!\!\cong\!\!731$
$2125\!\!\sim\!\!1094\!\!\cong\!\!1090$
$2126\!\!\sim\!\!933\!\!\cong\!\!849$
$2127\!\!\sim\!\!1094\!\!\cong\!\!1090$
$2128\!\!\sim\!\!933\!\!\cong\!\!849$
$2129\!\!\sim\!\!771\!\!\cong\!\!771$
$2130\!\!\sim\!\!852\!\!\cong\!\!852$
$2131\!\!\sim\!\!1094\!\!\cong\!\!1090$
$2132\!\!\sim\!\!852\!\!\cong\!\!852$
$2133\!\!\sim\!\!1094\!\!\cong\!\!1090$
$2134\!\!\sim\!\!1091\!\!\cong\!\!731$
$2135\!\!\sim\!\!930\!\!\cong\!\!821$
$2136\!\!\sim\!\!1091\!\!\cong\!\!731$
$2137\!\!\sim\!\!930\!\!\cong\!\!821$
$2138\!\!\sim\!\!768\!\!\cong\!\!731$
$2139\!\!\sim\!\!849\!\!\cong\!\!849$
$2140\!\!\sim\!\!1091\!\!\cong\!\!731$
$2141\!\!\sim\!\!849\!\!\cong\!\!849$
$2142\!\!\sim\!\!1091\!\!\cong\!\!731$
$2143\!\!\sim\!\!1090\!\!\cong\!\!1090$
$2144\!\!\sim\!\!928\!\!\cong\!\!820$
$2145\!\!\sim\!\!1090\!\!\cong\!\!1090$
$2146\!\!\sim\!\!928\!\!\cong\!\!820$
$2147\!\!\sim\!\!766\!\!\cong\!\!730$
$2148\!\!\sim\!\!847\!\!\cong\!\!847$
$2149\!\!\sim\!\!1090\!\!\cong\!\!1090$
$2150\!\!\sim\!\!847\!\!\cong\!\!847$
$2151\!\!\sim\!\!1090\!\!\cong\!\!1090$
$2152\!\!\sim\!\!1091\!\!\cong\!\!731$
$2153\!\!\sim\!\!929\!\!\cong\!\!929$
$2154\!\!\sim\!\!1091\!\!\cong\!\!731$
$2155\!\!\sim\!\!929\!\!\cong\!\!929$
$2156\!\!\sim\!\!767\!\!\cong\!\!731$
$2157\!\!\sim\!\!848\!\!\cong\!\!750$
$2158\!\!\sim\!\!1091\!\!\cong\!\!731$
$2159\!\!\sim\!\!848\!\!\cong\!\!750$
$2160\!\!\sim\!\!1091\!\!\cong\!\!731$
$2161\!\!\sim\!\!1094\!\!\cong\!\!1090$
$2162\!\!\sim\!\!933\!\!\cong\!\!849$
$2163\!\!\sim\!\!1094\!\!\cong\!\!1090$
$2164\!\!\sim\!\!933\!\!\cong\!\!849$
$2165\!\!\sim\!\!771\!\!\cong\!\!771$
$2166\!\!\sim\!\!852\!\!\cong\!\!852$
$2167\!\!\sim\!\!1094\!\!\cong\!\!1090$
$2168\!\!\sim\!\!852\!\!\cong\!\!852$
$2169\!\!\sim\!\!1094\!\!\cong\!\!1090$
$2170\!\!\sim\!\!1091\!\!\cong\!\!731$
$2171\!\!\sim\!\!929\!\!\cong\!\!929$
$2172\!\!\sim\!\!1091\!\!\cong\!\!731$
$2173\!\!\sim\!\!929\!\!\cong\!\!929$
$2174\!\!\sim\!\!767\!\!\cong\!\!731$
$2175\!\!\sim\!\!848\!\!\cong\!\!750$
$2176\!\!\sim\!\!1091\!\!\cong\!\!731$
$2177\!\!\sim\!\!848\!\!\cong\!\!750$
$2178\!\!\sim\!\!1091\!\!\cong\!\!731$
$2179\!\!\sim\!\!1094\!\!\cong\!\!1090$
$2180\!\!\sim\!\!932\!\!\cong\!\!820$
$2181\!\!\sim\!\!1094\!\!\cong\!\!1090$
$2182\!\!\sim\!\!932\!\!\cong\!\!820$
$2183\!\!\sim\!\!770\!\!\cong\!\!730$
$2184\!\!\sim\!\!851\!\!\cong\!\!847$
$2185\!\!\sim\!\!1094\!\!\cong\!\!1090$
$2186\!\!\sim\!\!851\!\!\cong\!\!847$
$2187\!\!\sim\!\!1094\!\!\cong\!\!1090$
$2188\!\!\sim\!\!730\!\!\cong\!\!730$
$2189\!\!\sim\!\!730\!\!\cong\!\!730$
$2190\!\!\sim\!\!2190\!\!\cong\!\!750$
$2191\!\!\sim\!\!730\!\!\cong\!\!730$
$2192\!\!\sim\!\!730\!\!\cong\!\!730$
$2193\!\!\sim\!\!2193\!\!\cong\!\!2193$
$2194\!\!\sim\!\!2190\!\!\cong\!\!750$
$2195\!\!\sim\!\!2193\!\!\cong\!\!2193$
$2196\!\!\sim\!\!2196\!\!\cong\!\!802$
$2197\!\!\sim\!\!730\!\!\cong\!\!730$
$2198\!\!\sim\!\!730\!\!\cong\!\!730$
$2199\!\!\sim\!\!2199\!\!\cong\!\!2199$
$2200\!\!\sim\!\!730\!\!\cong\!\!730$
$2201\!\!\sim\!\!730\!\!\cong\!\!730$
$2202\!\!\sim\!\!2202\!\!\cong\!\!2202$
$2203\!\!\sim\!\!2203\!\!\cong\!\!2203$
$2204\!\!\sim\!\!2204\!\!\cong\!\!2204$
$2205\!\!\sim\!\!2205\!\!\cong\!\!775$
$2206\!\!\sim\!\!2206\!\!\cong\!\!748$
$2207\!\!\sim\!\!2207\!\!\cong\!\!2207$
$2208\!\!\sim\!\!731\!\!\cong\!\!731$
$2209\!\!\sim\!\!2209\!\!\cong\!\!2209$
$2210\!\!\sim\!\!2210\!\!\cong\!\!2210$
$2211\!\!\sim\!\!731\!\!\cong\!\!731$
$2212\!\!\sim\!\!2212\!\!\cong\!\!2212$
$2213\!\!\sim\!\!2213\!\!\cong\!\!2213$
$2214\!\!\sim\!\!2214\!\!\cong\!\!748$
$2215\!\!\sim\!\!730\!\!\cong\!\!730$
$2216\!\!\sim\!\!730\!\!\cong\!\!730$
$2217\!\!\sim\!\!2203\!\!\cong\!\!2203$
$2218\!\!\sim\!\!730\!\!\cong\!\!730$
$2219\!\!\sim\!\!730\!\!\cong\!\!730$
$2220\!\!\sim\!\!2204\!\!\cong\!\!2204$
$2221\!\!\sim\!\!2199\!\!\cong\!\!2199$
$2222\!\!\sim\!\!2202\!\!\cong\!\!2202$
$2223\!\!\sim\!\!2205\!\!\cong\!\!775$
$2224\!\!\sim\!\!730\!\!\cong\!\!730$
$2225\!\!\sim\!\!730\!\!\cong\!\!730$
$2226\!\!\sim\!\!2226\!\!\cong\!\!820$
$2227\!\!\sim\!\!730\!\!\cong\!\!730$
$2228\!\!\sim\!\!730\!\!\cong\!\!730$
$2229\!\!\sim\!\!2229\!\!\cong\!\!2229$
$2230\!\!\sim\!\!2226\!\!\cong\!\!820$
$2231\!\!\sim\!\!2229\!\!\cong\!\!2229$
$2232\!\!\sim\!\!2232\!\!\cong\!\!730$
$2233\!\!\sim\!\!2233\!\!\cong\!\!2233$
$2234\!\!\sim\!\!2234\!\!\cong\!\!2234$
$2235\!\!\sim\!\!731\!\!\cong\!\!731$
$2236\!\!\sim\!\!2236\!\!\cong\!\!2236$
$2237\!\!\sim\!\!2237\!\!\cong\!\!2237$
$2238\!\!\sim\!\!731\!\!\cong\!\!731$
$2239\!\!\sim\!\!2239\!\!\cong\!\!2239$
$2240\!\!\sim\!\!2240\!\!\cong\!\!2240$
$2241\!\!\sim\!\!2241\!\!\cong\!\!739$
$2242\!\!\sim\!\!2206\!\!\cong\!\!748$
$2243\!\!\sim\!\!2209\!\!\cong\!\!2209$
$2244\!\!\sim\!\!2212\!\!\cong\!\!2212$
$2245\!\!\sim\!\!2207\!\!\cong\!\!2207$
$2246\!\!\sim\!\!2210\!\!\cong\!\!2210$
$2247\!\!\sim\!\!2213\!\!\cong\!\!2213$
$2248\!\!\sim\!\!731\!\!\cong\!\!731$
$2249\!\!\sim\!\!731\!\!\cong\!\!731$
$2250\!\!\sim\!\!2214\!\!\cong\!\!748$
$2251\!\!\sim\!\!2233\!\!\cong\!\!2233$
$2252\!\!\sim\!\!2236\!\!\cong\!\!2236$
$2253\!\!\sim\!\!2239\!\!\cong\!\!2239$
$2254\!\!\sim\!\!2234\!\!\cong\!\!2234$
$2255\!\!\sim\!\!2237\!\!\cong\!\!2237$
$2256\!\!\sim\!\!2240\!\!\cong\!\!2240$
$2257\!\!\sim\!\!731\!\!\cong\!\!731$
$2258\!\!\sim\!\!731\!\!\cong\!\!731$
$2259\!\!\sim\!\!2241\!\!\cong\!\!739$
$2260\!\!\sim\!\!2260\!\!\cong\!\!802$
$2261\!\!\sim\!\!2261\!\!\cong\!\!2261$
$2262\!\!\sim\!\!2262\!\!\cong\!\!750$
$2263\!\!\sim\!\!2261\!\!\cong\!\!2261$
$2264\!\!\sim\!\!2264\!\!\cong\!\!730$
$2265\!\!\sim\!\!2265\!\!\cong\!\!2265$
$2266\!\!\sim\!\!2262\!\!\cong\!\!750$
$2267\!\!\sim\!\!2265\!\!\cong\!\!2265$
$2268\!\!\sim\!\!734\!\!\cong\!\!730$
$2269\!\!\sim\!\!730\!\!\cong\!\!730$
$2270\!\!\sim\!\!730\!\!\cong\!\!730$
$2271\!\!\sim\!\!2271\!\!\cong\!\!2271$
$2272\!\!\sim\!\!730\!\!\cong\!\!730$
$2273\!\!\sim\!\!730\!\!\cong\!\!730$
$2274\!\!\sim\!\!2274\!\!\cong\!\!2274$
$2275\!\!\sim\!\!2271\!\!\cong\!\!2271$
$2276\!\!\sim\!\!2274\!\!\cong\!\!2274$
$2277\!\!\sim\!\!2277\!\!\cong\!\!2277$
$2278\!\!\sim\!\!730\!\!\cong\!\!730$
$2279\!\!\sim\!\!730\!\!\cong\!\!730$
$2280\!\!\sim\!\!2280\!\!\cong\!\!2280$
$2281\!\!\sim\!\!730\!\!\cong\!\!730$
$2282\!\!\sim\!\!730\!\!\cong\!\!730$
$2283\!\!\sim\!\!2283\!\!\cong\!\!2283$
$2284\!\!\sim\!\!2284\!\!\cong\!\!2284$
$2285\!\!\sim\!\!2285\!\!\cong\!\!2285$
$2286\!\!\sim\!\!2286\!\!\cong\!\!2286$
$2287\!\!\sim\!\!2287\!\!\cong\!\!2287$
$2288\!\!\sim\!\!2285\!\!\cong\!\!2285$
$2289\!\!\sim\!\!731\!\!\cong\!\!731$
$2290\!\!\sim\!\!2283\!\!\cong\!\!2283$
$2291\!\!\sim\!\!2274\!\!\cong\!\!2274$
$2292\!\!\sim\!\!731\!\!\cong\!\!731$
$2293\!\!\sim\!\!2293\!\!\cong\!\!2293$
$2294\!\!\sim\!\!2294\!\!\cong\!\!2294$
$2295\!\!\sim\!\!2295\!\!\cong\!\!2295$
$2296\!\!\sim\!\!730\!\!\cong\!\!730$
$2297\!\!\sim\!\!730\!\!\cong\!\!730$
$2298\!\!\sim\!\!2284\!\!\cong\!\!2284$
$2299\!\!\sim\!\!730\!\!\cong\!\!730$
$2300\!\!\sim\!\!730\!\!\cong\!\!730$
$2301\!\!\sim\!\!2285\!\!\cong\!\!2285$
$2302\!\!\sim\!\!2280\!\!\cong\!\!2280$
$2303\!\!\sim\!\!2283\!\!\cong\!\!2283$
$2304\!\!\sim\!\!2286\!\!\cong\!\!2286$
$2305\!\!\sim\!\!730\!\!\cong\!\!730$
$2306\!\!\sim\!\!730\!\!\cong\!\!730$
$2307\!\!\sim\!\!2307\!\!\cong\!\!2307$
$2308\!\!\sim\!\!730\!\!\cong\!\!730$
$2309\!\!\sim\!\!730\!\!\cong\!\!730$
$2310\!\!\sim\!\!2287\!\!\cong\!\!2287$
$2311\!\!\sim\!\!2307\!\!\cong\!\!2307$
$2312\!\!\sim\!\!2287\!\!\cong\!\!2287$
$2313\!\!\sim\!\!2313\!\!\cong\!\!2277$
$2314\!\!\sim\!\!2307\!\!\cong\!\!2307$
$2315\!\!\sim\!\!2284\!\!\cong\!\!2284$
$2316\!\!\sim\!\!731\!\!\cong\!\!731$
$2317\!\!\sim\!\!2280\!\!\cong\!\!2280$
$2318\!\!\sim\!\!2271\!\!\cong\!\!2271$
$2319\!\!\sim\!\!731\!\!\cong\!\!731$
$2320\!\!\sim\!\!2320\!\!\cong\!\!2294$
$2321\!\!\sim\!\!2293\!\!\cong\!\!2293$
$2322\!\!\sim\!\!2322\!\!\cong\!\!2322$
$2323\!\!\sim\!\!2287\!\!\cong\!\!2287$
$2324\!\!\sim\!\!2283\!\!\cong\!\!2283$
$2325\!\!\sim\!\!2293\!\!\cong\!\!2293$
$2326\!\!\sim\!\!2285\!\!\cong\!\!2285$
$2327\!\!\sim\!\!2274\!\!\cong\!\!2274$
$2328\!\!\sim\!\!2294\!\!\cong\!\!2294$
$2329\!\!\sim\!\!731\!\!\cong\!\!731$
$2330\!\!\sim\!\!731\!\!\cong\!\!731$
$2331\!\!\sim\!\!2295\!\!\cong\!\!2295$
$2332\!\!\sim\!\!2307\!\!\cong\!\!2307$
$2333\!\!\sim\!\!2280\!\!\cong\!\!2280$
$2334\!\!\sim\!\!2320\!\!\cong\!\!2294$
$2335\!\!\sim\!\!2284\!\!\cong\!\!2284$
$2336\!\!\sim\!\!2271\!\!\cong\!\!2271$
$2337\!\!\sim\!\!2293\!\!\cong\!\!2293$
$2338\!\!\sim\!\!731\!\!\cong\!\!731$
$2339\!\!\sim\!\!731\!\!\cong\!\!731$
$2340\!\!\sim\!\!2322\!\!\cong\!\!2322$
$2341\!\!\sim\!\!2313\!\!\cong\!\!2277$
$2342\!\!\sim\!\!2286\!\!\cong\!\!2286$
$2343\!\!\sim\!\!2322\!\!\cong\!\!2322$
$2344\!\!\sim\!\!2286\!\!\cong\!\!2286$
$2345\!\!\sim\!\!2277\!\!\cong\!\!2277$
$2346\!\!\sim\!\!2295\!\!\cong\!\!2295$
$2347\!\!\sim\!\!2322\!\!\cong\!\!2322$
$2348\!\!\sim\!\!2295\!\!\cong\!\!2295$
$2349\!\!\sim\!\!734\!\!\cong\!\!730$
$2350\!\!\sim\!\!820\!\!\cong\!\!820$
$2351\!\!\sim\!\!820\!\!\cong\!\!820$
$2352\!\!\sim\!\!2352\!\!\cong\!\!740$
$2353\!\!\sim\!\!820\!\!\cong\!\!820$
$2354\!\!\sim\!\!820\!\!\cong\!\!820$
$2355\!\!\sim\!\!2355\!\!\cong\!\!2355$
$2356\!\!\sim\!\!2352\!\!\cong\!\!740$
$2357\!\!\sim\!\!2355\!\!\cong\!\!2355$
$2358\!\!\sim\!\!2358\!\!\cong\!\!820$
$2359\!\!\sim\!\!820\!\!\cong\!\!820$
$2360\!\!\sim\!\!820\!\!\cong\!\!820$
$2361\!\!\sim\!\!2361\!\!\cong\!\!2361$
$2362\!\!\sim\!\!820\!\!\cong\!\!820$
$2363\!\!\sim\!\!820\!\!\cong\!\!820$
$2364\!\!\sim\!\!2364\!\!\cong\!\!2364$
$2365\!\!\sim\!\!2365\!\!\cong\!\!2365$
$2366\!\!\sim\!\!2366\!\!\cong\!\!2366$
$2367\!\!\sim\!\!2367\!\!\cong\!\!2367$
$2368\!\!\sim\!\!2368\!\!\cong\!\!739$
$2369\!\!\sim\!\!2369\!\!\cong\!\!2369$
$2370\!\!\sim\!\!821\!\!\cong\!\!821$
$2371\!\!\sim\!\!2371\!\!\cong\!\!2371$
$2372\!\!\sim\!\!2372\!\!\cong\!\!2372$
$2373\!\!\sim\!\!821\!\!\cong\!\!821$
$2374\!\!\sim\!\!2374\!\!\cong\!\!821$
$2375\!\!\sim\!\!2375\!\!\cong\!\!2375$
$2376\!\!\sim\!\!2376\!\!\cong\!\!739$
$2377\!\!\sim\!\!820\!\!\cong\!\!820$
$2378\!\!\sim\!\!820\!\!\cong\!\!820$
$2379\!\!\sim\!\!2365\!\!\cong\!\!2365$
$2380\!\!\sim\!\!820\!\!\cong\!\!820$
$2381\!\!\sim\!\!820\!\!\cong\!\!820$
$2382\!\!\sim\!\!2366\!\!\cong\!\!2366$
$2383\!\!\sim\!\!2361\!\!\cong\!\!2361$
$2384\!\!\sim\!\!2364\!\!\cong\!\!2364$
$2385\!\!\sim\!\!2367\!\!\cong\!\!2367$
$2386\!\!\sim\!\!820\!\!\cong\!\!820$
$2387\!\!\sim\!\!820\!\!\cong\!\!820$
$2388\!\!\sim\!\!2388\!\!\cong\!\!821$
$2389\!\!\sim\!\!820\!\!\cong\!\!820$
$2390\!\!\sim\!\!820\!\!\cong\!\!820$
$2391\!\!\sim\!\!2391\!\!\cong\!\!2391$
$2392\!\!\sim\!\!2388\!\!\cong\!\!821$
$2393\!\!\sim\!\!2391\!\!\cong\!\!2391$
$2394\!\!\sim\!\!2394\!\!\cong\!\!820$
$2395\!\!\sim\!\!2395\!\!\cong\!\!2395$
$2396\!\!\sim\!\!2396\!\!\cong\!\!2396$
$2397\!\!\sim\!\!821\!\!\cong\!\!821$
$2398\!\!\sim\!\!2398\!\!\cong\!\!2398$
$2399\!\!\sim\!\!2399\!\!\cong\!\!2399$
$2400\!\!\sim\!\!821\!\!\cong\!\!821$
$2401\!\!\sim\!\!2401\!\!\cong\!\!2401$
$2402\!\!\sim\!\!2402\!\!\cong\!\!2402$
$2403\!\!\sim\!\!2403\!\!\cong\!\!2287$
$2404\!\!\sim\!\!2368\!\!\cong\!\!739$
$2405\!\!\sim\!\!2371\!\!\cong\!\!2371$
$2406\!\!\sim\!\!2374\!\!\cong\!\!821$
$2407\!\!\sim\!\!2369\!\!\cong\!\!2369$
$2408\!\!\sim\!\!2372\!\!\cong\!\!2372$
$2409\!\!\sim\!\!2375\!\!\cong\!\!2375$
$2410\!\!\sim\!\!821\!\!\cong\!\!821$
$2411\!\!\sim\!\!821\!\!\cong\!\!821$
$2412\!\!\sim\!\!2376\!\!\cong\!\!739$
$2413\!\!\sim\!\!2395\!\!\cong\!\!2395$
$2414\!\!\sim\!\!2398\!\!\cong\!\!2398$
$2415\!\!\sim\!\!2401\!\!\cong\!\!2401$
$2416\!\!\sim\!\!2396\!\!\cong\!\!2396$
$2417\!\!\sim\!\!2399\!\!\cong\!\!2399$
$2418\!\!\sim\!\!2402\!\!\cong\!\!2402$
$2419\!\!\sim\!\!821\!\!\cong\!\!821$
$2420\!\!\sim\!\!821\!\!\cong\!\!821$
$2421\!\!\sim\!\!2403\!\!\cong\!\!2287$
$2422\!\!\sim\!\!2422\!\!\cong\!\!820$
$2423\!\!\sim\!\!2423\!\!\cong\!\!2423$
$2424\!\!\sim\!\!2424\!\!\cong\!\!966$
$2425\!\!\sim\!\!2423\!\!\cong\!\!2423$
$2426\!\!\sim\!\!2426\!\!\cong\!\!2277$
$2427\!\!\sim\!\!2427\!\!\cong\!\!2427$
$2428\!\!\sim\!\!2424\!\!\cong\!\!966$
$2429\!\!\sim\!\!2427\!\!\cong\!\!2427$
$2430\!\!\sim\!\!824\!\!\cong\!\!820$
$2431\!\!\sim\!\!730\!\!\cong\!\!730$
$2432\!\!\sim\!\!730\!\!\cong\!\!730$
$2433\!\!\sim\!\!2271\!\!\cong\!\!2271$
$2434\!\!\sim\!\!730\!\!\cong\!\!730$
$2435\!\!\sim\!\!730\!\!\cong\!\!730$
$2436\!\!\sim\!\!2274\!\!\cong\!\!2274$
$2437\!\!\sim\!\!2271\!\!\cong\!\!2271$
$2438\!\!\sim\!\!2274\!\!\cong\!\!2274$
$2439\!\!\sim\!\!2277\!\!\cong\!\!2277$
$2440\!\!\sim\!\!730\!\!\cong\!\!730$
$2441\!\!\sim\!\!730\!\!\cong\!\!730$
$2442\!\!\sim\!\!2280\!\!\cong\!\!2280$
$2443\!\!\sim\!\!730\!\!\cong\!\!730$
$2444\!\!\sim\!\!730\!\!\cong\!\!730$
$2445\!\!\sim\!\!2283\!\!\cong\!\!2283$
$2446\!\!\sim\!\!2284\!\!\cong\!\!2284$
$2447\!\!\sim\!\!2285\!\!\cong\!\!2285$
$2448\!\!\sim\!\!2286\!\!\cong\!\!2286$
$2449\!\!\sim\!\!2287\!\!\cong\!\!2287$
$2450\!\!\sim\!\!2285\!\!\cong\!\!2285$
$2451\!\!\sim\!\!731\!\!\cong\!\!731$
$2452\!\!\sim\!\!2283\!\!\cong\!\!2283$
$2453\!\!\sim\!\!2274\!\!\cong\!\!2274$
$2454\!\!\sim\!\!731\!\!\cong\!\!731$
$2455\!\!\sim\!\!2293\!\!\cong\!\!2293$
$2456\!\!\sim\!\!2294\!\!\cong\!\!2294$
$2457\!\!\sim\!\!2295\!\!\cong\!\!2295$
$2458\!\!\sim\!\!730\!\!\cong\!\!730$
$2459\!\!\sim\!\!730\!\!\cong\!\!730$
$2460\!\!\sim\!\!2284\!\!\cong\!\!2284$
$2461\!\!\sim\!\!730\!\!\cong\!\!730$
$2462\!\!\sim\!\!730\!\!\cong\!\!730$
$2463\!\!\sim\!\!2285\!\!\cong\!\!2285$
$2464\!\!\sim\!\!2280\!\!\cong\!\!2280$
$2465\!\!\sim\!\!2283\!\!\cong\!\!2283$
$2466\!\!\sim\!\!2286\!\!\cong\!\!2286$
$2467\!\!\sim\!\!730\!\!\cong\!\!730$
$2468\!\!\sim\!\!730\!\!\cong\!\!730$
$2469\!\!\sim\!\!2307\!\!\cong\!\!2307$
$2470\!\!\sim\!\!730\!\!\cong\!\!730$
$2471\!\!\sim\!\!730\!\!\cong\!\!730$
$2472\!\!\sim\!\!2287\!\!\cong\!\!2287$
$2473\!\!\sim\!\!2307\!\!\cong\!\!2307$
$2474\!\!\sim\!\!2287\!\!\cong\!\!2287$
$2475\!\!\sim\!\!2313\!\!\cong\!\!2277$
$2476\!\!\sim\!\!2307\!\!\cong\!\!2307$
$2477\!\!\sim\!\!2284\!\!\cong\!\!2284$
$2478\!\!\sim\!\!731\!\!\cong\!\!731$
$2479\!\!\sim\!\!2280\!\!\cong\!\!2280$
$2480\!\!\sim\!\!2271\!\!\cong\!\!2271$
$2481\!\!\sim\!\!731\!\!\cong\!\!731$
$2482\!\!\sim\!\!2320\!\!\cong\!\!2294$
$2483\!\!\sim\!\!2293\!\!\cong\!\!2293$
$2484\!\!\sim\!\!2322\!\!\cong\!\!2322$
$2485\!\!\sim\!\!2287\!\!\cong\!\!2287$
$2486\!\!\sim\!\!2283\!\!\cong\!\!2283$
$2487\!\!\sim\!\!2293\!\!\cong\!\!2293$
$2488\!\!\sim\!\!2285\!\!\cong\!\!2285$
$2489\!\!\sim\!\!2274\!\!\cong\!\!2274$
$2490\!\!\sim\!\!2294\!\!\cong\!\!2294$
$2491\!\!\sim\!\!731\!\!\cong\!\!731$
$2492\!\!\sim\!\!731\!\!\cong\!\!731$
$2493\!\!\sim\!\!2295\!\!\cong\!\!2295$
$2494\!\!\sim\!\!2307\!\!\cong\!\!2307$
$2495\!\!\sim\!\!2280\!\!\cong\!\!2280$
$2496\!\!\sim\!\!2320\!\!\cong\!\!2294$
$2497\!\!\sim\!\!2284\!\!\cong\!\!2284$
$2498\!\!\sim\!\!2271\!\!\cong\!\!2271$
$2499\!\!\sim\!\!2293\!\!\cong\!\!2293$
$2500\!\!\sim\!\!731\!\!\cong\!\!731$
$2501\!\!\sim\!\!731\!\!\cong\!\!731$
$2502\!\!\sim\!\!2322\!\!\cong\!\!2322$
$2503\!\!\sim\!\!2313\!\!\cong\!\!2277$
$2504\!\!\sim\!\!2286\!\!\cong\!\!2286$
$2505\!\!\sim\!\!2322\!\!\cong\!\!2322$
$2506\!\!\sim\!\!2286\!\!\cong\!\!2286$
$2507\!\!\sim\!\!2277\!\!\cong\!\!2277$
$2508\!\!\sim\!\!2295\!\!\cong\!\!2295$
$2509\!\!\sim\!\!2322\!\!\cong\!\!2322$
$2510\!\!\sim\!\!2295\!\!\cong\!\!2295$
$2511\!\!\sim\!\!734\!\!\cong\!\!730$
$2512\!\!\sim\!\!730\!\!\cong\!\!730$
$2513\!\!\sim\!\!730\!\!\cong\!\!730$
$2514\!\!\sim\!\!2237\!\!\cong\!\!2237$
$2515\!\!\sim\!\!730\!\!\cong\!\!730$
$2516\!\!\sim\!\!730\!\!\cong\!\!730$
$2517\!\!\sim\!\!2210\!\!\cong\!\!2210$
$2518\!\!\sim\!\!2237\!\!\cong\!\!2237$
$2519\!\!\sim\!\!2210\!\!\cong\!\!2210$
$2520\!\!\sim\!\!2264\!\!\cong\!\!730$
$2521\!\!\sim\!\!730\!\!\cong\!\!730$
$2522\!\!\sim\!\!730\!\!\cong\!\!730$
$2523\!\!\sim\!\!2236\!\!\cong\!\!2236$
$2524\!\!\sim\!\!730\!\!\cong\!\!730$
$2525\!\!\sim\!\!730\!\!\cong\!\!730$
$2526\!\!\sim\!\!2209\!\!\cong\!\!2209$
$2527\!\!\sim\!\!2234\!\!\cong\!\!2234$
$2528\!\!\sim\!\!2207\!\!\cong\!\!2207$
$2529\!\!\sim\!\!2261\!\!\cong\!\!2261$
$2530\!\!\sim\!\!2229\!\!\cong\!\!2229$
$2531\!\!\sim\!\!2204\!\!\cong\!\!2204$
$2532\!\!\sim\!\!731\!\!\cong\!\!731$
$2533\!\!\sim\!\!2202\!\!\cong\!\!2202$
$2534\!\!\sim\!\!2193\!\!\cong\!\!2193$
$2535\!\!\sim\!\!731\!\!\cong\!\!731$
$2536\!\!\sim\!\!2240\!\!\cong\!\!2240$
$2537\!\!\sim\!\!2213\!\!\cong\!\!2213$
$2538\!\!\sim\!\!2265\!\!\cong\!\!2265$
$2539\!\!\sim\!\!730\!\!\cong\!\!730$
$2540\!\!\sim\!\!730\!\!\cong\!\!730$
$2541\!\!\sim\!\!2234\!\!\cong\!\!2234$
$2542\!\!\sim\!\!730\!\!\cong\!\!730$
$2543\!\!\sim\!\!730\!\!\cong\!\!730$
$2544\!\!\sim\!\!2207\!\!\cong\!\!2207$
$2545\!\!\sim\!\!2236\!\!\cong\!\!2236$
$2546\!\!\sim\!\!2209\!\!\cong\!\!2209$
$2547\!\!\sim\!\!2261\!\!\cong\!\!2261$
$2548\!\!\sim\!\!730\!\!\cong\!\!730$
$2549\!\!\sim\!\!730\!\!\cong\!\!730$
$2550\!\!\sim\!\!2233\!\!\cong\!\!2233$
$2551\!\!\sim\!\!730\!\!\cong\!\!730$
$2552\!\!\sim\!\!730\!\!\cong\!\!730$
$2553\!\!\sim\!\!2206\!\!\cong\!\!748$
$2554\!\!\sim\!\!2233\!\!\cong\!\!2233$
$2555\!\!\sim\!\!2206\!\!\cong\!\!748$
$2556\!\!\sim\!\!2260\!\!\cong\!\!802$
$2557\!\!\sim\!\!2226\!\!\cong\!\!820$
$2558\!\!\sim\!\!2203\!\!\cong\!\!2203$
$2559\!\!\sim\!\!731\!\!\cong\!\!731$
$2560\!\!\sim\!\!2199\!\!\cong\!\!2199$
$2561\!\!\sim\!\!2190\!\!\cong\!\!750$
$2562\!\!\sim\!\!731\!\!\cong\!\!731$
$2563\!\!\sim\!\!2239\!\!\cong\!\!2239$
$2564\!\!\sim\!\!2212\!\!\cong\!\!2212$
$2565\!\!\sim\!\!2262\!\!\cong\!\!750$
$2566\!\!\sim\!\!2229\!\!\cong\!\!2229$
$2567\!\!\sim\!\!2202\!\!\cong\!\!2202$
$2568\!\!\sim\!\!2240\!\!\cong\!\!2240$
$2569\!\!\sim\!\!2204\!\!\cong\!\!2204$
$2570\!\!\sim\!\!2193\!\!\cong\!\!2193$
$2571\!\!\sim\!\!2213\!\!\cong\!\!2213$
$2572\!\!\sim\!\!731\!\!\cong\!\!731$
$2573\!\!\sim\!\!731\!\!\cong\!\!731$
$2574\!\!\sim\!\!2265\!\!\cong\!\!2265$
$2575\!\!\sim\!\!2226\!\!\cong\!\!820$
$2576\!\!\sim\!\!2199\!\!\cong\!\!2199$
$2577\!\!\sim\!\!2239\!\!\cong\!\!2239$
$2578\!\!\sim\!\!2203\!\!\cong\!\!2203$
$2579\!\!\sim\!\!2190\!\!\cong\!\!750$
$2580\!\!\sim\!\!2212\!\!\cong\!\!2212$
$2581\!\!\sim\!\!731\!\!\cong\!\!731$
$2582\!\!\sim\!\!731\!\!\cong\!\!731$
$2583\!\!\sim\!\!2262\!\!\cong\!\!750$
$2584\!\!\sim\!\!2232\!\!\cong\!\!730$
$2585\!\!\sim\!\!2205\!\!\cong\!\!775$
$2586\!\!\sim\!\!2241\!\!\cong\!\!739$
$2587\!\!\sim\!\!2205\!\!\cong\!\!775$
$2588\!\!\sim\!\!2196\!\!\cong\!\!802$
$2589\!\!\sim\!\!2214\!\!\cong\!\!748$
$2590\!\!\sim\!\!2241\!\!\cong\!\!739$
$2591\!\!\sim\!\!2214\!\!\cong\!\!748$
$2592\!\!\sim\!\!734\!\!\cong\!\!730$
$2593\!\!\sim\!\!820\!\!\cong\!\!820$
$2594\!\!\sim\!\!820\!\!\cong\!\!820$
$2595\!\!\sim\!\!2399\!\!\cong\!\!2399$
$2596\!\!\sim\!\!820\!\!\cong\!\!820$
$2597\!\!\sim\!\!820\!\!\cong\!\!820$
$2598\!\!\sim\!\!2372\!\!\cong\!\!2372$
$2599\!\!\sim\!\!2399\!\!\cong\!\!2399$
$2600\!\!\sim\!\!2372\!\!\cong\!\!2372$
$2601\!\!\sim\!\!2426\!\!\cong\!\!2277$
$2602\!\!\sim\!\!820\!\!\cong\!\!820$
$2603\!\!\sim\!\!820\!\!\cong\!\!820$
$2604\!\!\sim\!\!2398\!\!\cong\!\!2398$
$2605\!\!\sim\!\!820\!\!\cong\!\!820$
$2606\!\!\sim\!\!820\!\!\cong\!\!820$
$2607\!\!\sim\!\!2371\!\!\cong\!\!2371$
$2608\!\!\sim\!\!2396\!\!\cong\!\!2396$
$2609\!\!\sim\!\!2369\!\!\cong\!\!2369$
$2610\!\!\sim\!\!2423\!\!\cong\!\!2423$
$2611\!\!\sim\!\!2391\!\!\cong\!\!2391$
$2612\!\!\sim\!\!2366\!\!\cong\!\!2366$
$2613\!\!\sim\!\!821\!\!\cong\!\!821$
$2614\!\!\sim\!\!2364\!\!\cong\!\!2364$
$2615\!\!\sim\!\!2355\!\!\cong\!\!2355$
$2616\!\!\sim\!\!821\!\!\cong\!\!821$
$2617\!\!\sim\!\!2402\!\!\cong\!\!2402$
$2618\!\!\sim\!\!2375\!\!\cong\!\!2375$
$2619\!\!\sim\!\!2427\!\!\cong\!\!2427$
$2620\!\!\sim\!\!820\!\!\cong\!\!820$
$2621\!\!\sim\!\!820\!\!\cong\!\!820$
$2622\!\!\sim\!\!2396\!\!\cong\!\!2396$
$2623\!\!\sim\!\!820\!\!\cong\!\!820$
$2624\!\!\sim\!\!820\!\!\cong\!\!820$
$2625\!\!\sim\!\!2369\!\!\cong\!\!2369$
$2626\!\!\sim\!\!2398\!\!\cong\!\!2398$
$2627\!\!\sim\!\!2371\!\!\cong\!\!2371$
$2628\!\!\sim\!\!2423\!\!\cong\!\!2423$
$2629\!\!\sim\!\!820\!\!\cong\!\!820$
$2630\!\!\sim\!\!820\!\!\cong\!\!820$
$2631\!\!\sim\!\!2395\!\!\cong\!\!2395$
$2632\!\!\sim\!\!820\!\!\cong\!\!820$
$2633\!\!\sim\!\!820\!\!\cong\!\!820$
$2634\!\!\sim\!\!2368\!\!\cong\!\!739$
$2635\!\!\sim\!\!2395\!\!\cong\!\!2395$
$2636\!\!\sim\!\!2368\!\!\cong\!\!739$
$2637\!\!\sim\!\!2422\!\!\cong\!\!820$
$2638\!\!\sim\!\!2388\!\!\cong\!\!821$
$2639\!\!\sim\!\!2365\!\!\cong\!\!2365$
$2640\!\!\sim\!\!821\!\!\cong\!\!821$
$2641\!\!\sim\!\!2361\!\!\cong\!\!2361$
$2642\!\!\sim\!\!2352\!\!\cong\!\!740$
$2643\!\!\sim\!\!821\!\!\cong\!\!821$
$2644\!\!\sim\!\!2401\!\!\cong\!\!2401$
$2645\!\!\sim\!\!2374\!\!\cong\!\!821$
$2646\!\!\sim\!\!2424\!\!\cong\!\!966$
$2647\!\!\sim\!\!2391\!\!\cong\!\!2391$
$2648\!\!\sim\!\!2364\!\!\cong\!\!2364$
$2649\!\!\sim\!\!2402\!\!\cong\!\!2402$
$2650\!\!\sim\!\!2366\!\!\cong\!\!2366$
$2651\!\!\sim\!\!2355\!\!\cong\!\!2355$
$2652\!\!\sim\!\!2375\!\!\cong\!\!2375$
$2653\!\!\sim\!\!821\!\!\cong\!\!821$
$2654\!\!\sim\!\!821\!\!\cong\!\!821$
$2655\!\!\sim\!\!2427\!\!\cong\!\!2427$
$2656\!\!\sim\!\!2388\!\!\cong\!\!821$
$2657\!\!\sim\!\!2361\!\!\cong\!\!2361$
$2658\!\!\sim\!\!2401\!\!\cong\!\!2401$
$2659\!\!\sim\!\!2365\!\!\cong\!\!2365$
$2660\!\!\sim\!\!2352\!\!\cong\!\!740$
$2661\!\!\sim\!\!2374\!\!\cong\!\!821$
$2662\!\!\sim\!\!821\!\!\cong\!\!821$
$2663\!\!\sim\!\!821\!\!\cong\!\!821$
$2664\!\!\sim\!\!2424\!\!\cong\!\!966$
$2665\!\!\sim\!\!2394\!\!\cong\!\!820$
$2666\!\!\sim\!\!2367\!\!\cong\!\!2367$
$2667\!\!\sim\!\!2403\!\!\cong\!\!2287$
$2668\!\!\sim\!\!2367\!\!\cong\!\!2367$
$2669\!\!\sim\!\!2358\!\!\cong\!\!820$
$2670\!\!\sim\!\!2376\!\!\cong\!\!739$
$2671\!\!\sim\!\!2403\!\!\cong\!\!2287$
$2672\!\!\sim\!\!2376\!\!\cong\!\!739$
$2673\!\!\sim\!\!824\!\!\cong\!\!820$
$2674\!\!\sim\!\!820\!\!\cong\!\!820$
$2675\!\!\sim\!\!820\!\!\cong\!\!820$
$2676\!\!\sim\!\!2352\!\!\cong\!\!740$
$2677\!\!\sim\!\!820\!\!\cong\!\!820$
$2678\!\!\sim\!\!820\!\!\cong\!\!820$
$2679\!\!\sim\!\!2355\!\!\cong\!\!2355$
$2680\!\!\sim\!\!2352\!\!\cong\!\!740$
$2681\!\!\sim\!\!2355\!\!\cong\!\!2355$
$2682\!\!\sim\!\!2358\!\!\cong\!\!820$
$2683\!\!\sim\!\!820\!\!\cong\!\!820$
$2684\!\!\sim\!\!820\!\!\cong\!\!820$
$2685\!\!\sim\!\!2361\!\!\cong\!\!2361$
$2686\!\!\sim\!\!820\!\!\cong\!\!820$
$2687\!\!\sim\!\!820\!\!\cong\!\!820$
$2688\!\!\sim\!\!2364\!\!\cong\!\!2364$
$2689\!\!\sim\!\!2365\!\!\cong\!\!2365$
$2690\!\!\sim\!\!2366\!\!\cong\!\!2366$
$2691\!\!\sim\!\!2367\!\!\cong\!\!2367$
$2692\!\!\sim\!\!2368\!\!\cong\!\!739$
$2693\!\!\sim\!\!2369\!\!\cong\!\!2369$
$2694\!\!\sim\!\!821\!\!\cong\!\!821$
$2695\!\!\sim\!\!2371\!\!\cong\!\!2371$
$2696\!\!\sim\!\!2372\!\!\cong\!\!2372$
$2697\!\!\sim\!\!821\!\!\cong\!\!821$
$2698\!\!\sim\!\!2374\!\!\cong\!\!821$
$2699\!\!\sim\!\!2375\!\!\cong\!\!2375$
$2700\!\!\sim\!\!2376\!\!\cong\!\!739$
$2701\!\!\sim\!\!820\!\!\cong\!\!820$
$2702\!\!\sim\!\!820\!\!\cong\!\!820$
$2703\!\!\sim\!\!2365\!\!\cong\!\!2365$
$2704\!\!\sim\!\!820\!\!\cong\!\!820$
$2705\!\!\sim\!\!820\!\!\cong\!\!820$
$2706\!\!\sim\!\!2366\!\!\cong\!\!2366$
$2707\!\!\sim\!\!2361\!\!\cong\!\!2361$
$2708\!\!\sim\!\!2364\!\!\cong\!\!2364$
$2709\!\!\sim\!\!2367\!\!\cong\!\!2367$
$2710\!\!\sim\!\!820\!\!\cong\!\!820$
$2711\!\!\sim\!\!820\!\!\cong\!\!820$
$2712\!\!\sim\!\!2388\!\!\cong\!\!821$
$2713\!\!\sim\!\!820\!\!\cong\!\!820$
$2714\!\!\sim\!\!820\!\!\cong\!\!820$
$2715\!\!\sim\!\!2391\!\!\cong\!\!2391$
$2716\!\!\sim\!\!2388\!\!\cong\!\!821$
$2717\!\!\sim\!\!2391\!\!\cong\!\!2391$
$2718\!\!\sim\!\!2394\!\!\cong\!\!820$
$2719\!\!\sim\!\!2395\!\!\cong\!\!2395$
$2720\!\!\sim\!\!2396\!\!\cong\!\!2396$
$2721\!\!\sim\!\!821\!\!\cong\!\!821$
$2722\!\!\sim\!\!2398\!\!\cong\!\!2398$
$2723\!\!\sim\!\!2399\!\!\cong\!\!2399$
$2724\!\!\sim\!\!821\!\!\cong\!\!821$
$2725\!\!\sim\!\!2401\!\!\cong\!\!2401$
$2726\!\!\sim\!\!2402\!\!\cong\!\!2402$
$2727\!\!\sim\!\!2403\!\!\cong\!\!2287$
$2728\!\!\sim\!\!2368\!\!\cong\!\!739$
$2729\!\!\sim\!\!2371\!\!\cong\!\!2371$
$2730\!\!\sim\!\!2374\!\!\cong\!\!821$
$2731\!\!\sim\!\!2369\!\!\cong\!\!2369$
$2732\!\!\sim\!\!2372\!\!\cong\!\!2372$
$2733\!\!\sim\!\!2375\!\!\cong\!\!2375$
$2734\!\!\sim\!\!821\!\!\cong\!\!821$
$2735\!\!\sim\!\!821\!\!\cong\!\!821$
$2736\!\!\sim\!\!2376\!\!\cong\!\!739$
$2737\!\!\sim\!\!2395\!\!\cong\!\!2395$
$2738\!\!\sim\!\!2398\!\!\cong\!\!2398$
$2739\!\!\sim\!\!2401\!\!\cong\!\!2401$
$2740\!\!\sim\!\!2396\!\!\cong\!\!2396$
$2741\!\!\sim\!\!2399\!\!\cong\!\!2399$
$2742\!\!\sim\!\!2402\!\!\cong\!\!2402$
$2743\!\!\sim\!\!821\!\!\cong\!\!821$
$2744\!\!\sim\!\!821\!\!\cong\!\!821$
$2745\!\!\sim\!\!2403\!\!\cong\!\!2287$
$2746\!\!\sim\!\!2422\!\!\cong\!\!820$
$2747\!\!\sim\!\!2423\!\!\cong\!\!2423$
$2748\!\!\sim\!\!2424\!\!\cong\!\!966$
$2749\!\!\sim\!\!2423\!\!\cong\!\!2423$
$2750\!\!\sim\!\!2426\!\!\cong\!\!2277$
$2751\!\!\sim\!\!2427\!\!\cong\!\!2427$
$2752\!\!\sim\!\!2424\!\!\cong\!\!966$
$2753\!\!\sim\!\!2427\!\!\cong\!\!2427$
$2754\!\!\sim\!\!824\!\!\cong\!\!820$
$2755\!\!\sim\!\!820\!\!\cong\!\!820$
$2756\!\!\sim\!\!820\!\!\cong\!\!820$
$2757\!\!\sim\!\!2399\!\!\cong\!\!2399$
$2758\!\!\sim\!\!820\!\!\cong\!\!820$
$2759\!\!\sim\!\!820\!\!\cong\!\!820$
$2760\!\!\sim\!\!2372\!\!\cong\!\!2372$
$2761\!\!\sim\!\!2399\!\!\cong\!\!2399$
$2762\!\!\sim\!\!2372\!\!\cong\!\!2372$
$2763\!\!\sim\!\!2426\!\!\cong\!\!2277$
$2764\!\!\sim\!\!820\!\!\cong\!\!820$
$2765\!\!\sim\!\!820\!\!\cong\!\!820$
$2766\!\!\sim\!\!2398\!\!\cong\!\!2398$
$2767\!\!\sim\!\!820\!\!\cong\!\!820$
$2768\!\!\sim\!\!820\!\!\cong\!\!820$
$2769\!\!\sim\!\!2371\!\!\cong\!\!2371$
$2770\!\!\sim\!\!2396\!\!\cong\!\!2396$
$2771\!\!\sim\!\!2369\!\!\cong\!\!2369$
$2772\!\!\sim\!\!2423\!\!\cong\!\!2423$
$2773\!\!\sim\!\!2391\!\!\cong\!\!2391$
$2774\!\!\sim\!\!2366\!\!\cong\!\!2366$
$2775\!\!\sim\!\!821\!\!\cong\!\!821$
$2776\!\!\sim\!\!2364\!\!\cong\!\!2364$
$2777\!\!\sim\!\!2355\!\!\cong\!\!2355$
$2778\!\!\sim\!\!821\!\!\cong\!\!821$
$2779\!\!\sim\!\!2402\!\!\cong\!\!2402$
$2780\!\!\sim\!\!2375\!\!\cong\!\!2375$
$2781\!\!\sim\!\!2427\!\!\cong\!\!2427$
$2782\!\!\sim\!\!820\!\!\cong\!\!820$
$2783\!\!\sim\!\!820\!\!\cong\!\!820$
$2784\!\!\sim\!\!2396\!\!\cong\!\!2396$
$2785\!\!\sim\!\!820\!\!\cong\!\!820$
$2786\!\!\sim\!\!820\!\!\cong\!\!820$
$2787\!\!\sim\!\!2369\!\!\cong\!\!2369$
$2788\!\!\sim\!\!2398\!\!\cong\!\!2398$
$2789\!\!\sim\!\!2371\!\!\cong\!\!2371$
$2790\!\!\sim\!\!2423\!\!\cong\!\!2423$
$2791\!\!\sim\!\!820\!\!\cong\!\!820$
$2792\!\!\sim\!\!820\!\!\cong\!\!820$
$2793\!\!\sim\!\!2395\!\!\cong\!\!2395$
$2794\!\!\sim\!\!820\!\!\cong\!\!820$
$2795\!\!\sim\!\!820\!\!\cong\!\!820$
$2796\!\!\sim\!\!2368\!\!\cong\!\!739$
$2797\!\!\sim\!\!2395\!\!\cong\!\!2395$
$2798\!\!\sim\!\!2368\!\!\cong\!\!739$
$2799\!\!\sim\!\!2422\!\!\cong\!\!820$
$2800\!\!\sim\!\!2388\!\!\cong\!\!821$
$2801\!\!\sim\!\!2365\!\!\cong\!\!2365$
$2802\!\!\sim\!\!821\!\!\cong\!\!821$
$2803\!\!\sim\!\!2361\!\!\cong\!\!2361$
$2804\!\!\sim\!\!2352\!\!\cong\!\!740$
$2805\!\!\sim\!\!821\!\!\cong\!\!821$
$2806\!\!\sim\!\!2401\!\!\cong\!\!2401$
$2807\!\!\sim\!\!2374\!\!\cong\!\!821$
$2808\!\!\sim\!\!2424\!\!\cong\!\!966$
$2809\!\!\sim\!\!2391\!\!\cong\!\!2391$
$2810\!\!\sim\!\!2364\!\!\cong\!\!2364$
$2811\!\!\sim\!\!2402\!\!\cong\!\!2402$
$2812\!\!\sim\!\!2366\!\!\cong\!\!2366$
$2813\!\!\sim\!\!2355\!\!\cong\!\!2355$
$2814\!\!\sim\!\!2375\!\!\cong\!\!2375$
$2815\!\!\sim\!\!821\!\!\cong\!\!821$
$2816\!\!\sim\!\!821\!\!\cong\!\!821$
$2817\!\!\sim\!\!2427\!\!\cong\!\!2427$
$2818\!\!\sim\!\!2388\!\!\cong\!\!821$
$2819\!\!\sim\!\!2361\!\!\cong\!\!2361$
$2820\!\!\sim\!\!2401\!\!\cong\!\!2401$
$2821\!\!\sim\!\!2365\!\!\cong\!\!2365$
$2822\!\!\sim\!\!2352\!\!\cong\!\!740$
$2823\!\!\sim\!\!2374\!\!\cong\!\!821$
$2824\!\!\sim\!\!821\!\!\cong\!\!821$
$2825\!\!\sim\!\!821\!\!\cong\!\!821$
$2826\!\!\sim\!\!2424\!\!\cong\!\!966$
$2827\!\!\sim\!\!2394\!\!\cong\!\!820$
$2828\!\!\sim\!\!2367\!\!\cong\!\!2367$
$2829\!\!\sim\!\!2403\!\!\cong\!\!2287$
$2830\!\!\sim\!\!2367\!\!\cong\!\!2367$
$2831\!\!\sim\!\!2358\!\!\cong\!\!820$
$2832\!\!\sim\!\!2376\!\!\cong\!\!739$
$2833\!\!\sim\!\!2403\!\!\cong\!\!2287$
$2834\!\!\sim\!\!2376\!\!\cong\!\!739$
$2835\!\!\sim\!\!824\!\!\cong\!\!820$
$2836\!\!\sim\!\!1090\!\!\cong\!\!1090$
$2837\!\!\sim\!\!1090\!\!\cong\!\!1090$
$2838\!\!\sim\!\!2838\!\!\cong\!\!750$
$2839\!\!\sim\!\!1090\!\!\cong\!\!1090$
$2840\!\!\sim\!\!1090\!\!\cong\!\!1090$
$2841\!\!\sim\!\!2841\!\!\cong\!\!2841$
$2842\!\!\sim\!\!2838\!\!\cong\!\!750$
$2843\!\!\sim\!\!2841\!\!\cong\!\!2841$
$2844\!\!\sim\!\!2844\!\!\cong\!\!730$
$2845\!\!\sim\!\!1090\!\!\cong\!\!1090$
$2846\!\!\sim\!\!1090\!\!\cong\!\!1090$
$2847\!\!\sim\!\!2847\!\!\cong\!\!929$
$2848\!\!\sim\!\!1090\!\!\cong\!\!1090$
$2849\!\!\sim\!\!1090\!\!\cong\!\!1090$
$2850\!\!\sim\!\!2850\!\!\cong\!\!2850$
$2851\!\!\sim\!\!2851\!\!\cong\!\!929$
$2852\!\!\sim\!\!2852\!\!\cong\!\!849$
$2853\!\!\sim\!\!2853\!\!\cong\!\!2853$
$2854\!\!\sim\!\!2854\!\!\cong\!\!847$
$2855\!\!\sim\!\!2852\!\!\cong\!\!849$
$2856\!\!\sim\!\!1091\!\!\cong\!\!731$
$2857\!\!\sim\!\!2850\!\!\cong\!\!2850$
$2858\!\!\sim\!\!2841\!\!\cong\!\!2841$
$2859\!\!\sim\!\!1091\!\!\cong\!\!731$
$2860\!\!\sim\!\!2860\!\!\cong\!\!2212$
$2861\!\!\sim\!\!2861\!\!\cong\!\!731$
$2862\!\!\sim\!\!2862\!\!\cong\!\!847$
$2863\!\!\sim\!\!1090\!\!\cong\!\!1090$
$2864\!\!\sim\!\!1090\!\!\cong\!\!1090$
$2865\!\!\sim\!\!2851\!\!\cong\!\!929$
$2866\!\!\sim\!\!1090\!\!\cong\!\!1090$
$2867\!\!\sim\!\!1090\!\!\cong\!\!1090$
$2868\!\!\sim\!\!2852\!\!\cong\!\!849$
$2869\!\!\sim\!\!2847\!\!\cong\!\!929$
$2870\!\!\sim\!\!2850\!\!\cong\!\!2850$
$2871\!\!\sim\!\!2853\!\!\cong\!\!2853$
$2872\!\!\sim\!\!1090\!\!\cong\!\!1090$
$2873\!\!\sim\!\!1090\!\!\cong\!\!1090$
$2874\!\!\sim\!\!2874\!\!\cong\!\!820$
$2875\!\!\sim\!\!1090\!\!\cong\!\!1090$
$2876\!\!\sim\!\!1090\!\!\cong\!\!1090$
$2877\!\!\sim\!\!2854\!\!\cong\!\!847$
$2878\!\!\sim\!\!2874\!\!\cong\!\!820$
$2879\!\!\sim\!\!2854\!\!\cong\!\!847$
$2880\!\!\sim\!\!2880\!\!\cong\!\!730$
$2881\!\!\sim\!\!2874\!\!\cong\!\!820$
$2882\!\!\sim\!\!2851\!\!\cong\!\!929$
$2883\!\!\sim\!\!1091\!\!\cong\!\!731$
$2884\!\!\sim\!\!2847\!\!\cong\!\!929$
$2885\!\!\sim\!\!2838\!\!\cong\!\!750$
$2886\!\!\sim\!\!1091\!\!\cong\!\!731$
$2887\!\!\sim\!\!2887\!\!\cong\!\!731$
$2888\!\!\sim\!\!2860\!\!\cong\!\!2212$
$2889\!\!\sim\!\!2889\!\!\cong\!\!750$
$2890\!\!\sim\!\!2854\!\!\cong\!\!847$
$2891\!\!\sim\!\!2850\!\!\cong\!\!2850$
$2892\!\!\sim\!\!2860\!\!\cong\!\!2212$
$2893\!\!\sim\!\!2852\!\!\cong\!\!849$
$2894\!\!\sim\!\!2841\!\!\cong\!\!2841$
$2895\!\!\sim\!\!2861\!\!\cong\!\!731$
$2896\!\!\sim\!\!1091\!\!\cong\!\!731$
$2897\!\!\sim\!\!1091\!\!\cong\!\!731$
$2898\!\!\sim\!\!2862\!\!\cong\!\!847$
$2899\!\!\sim\!\!2874\!\!\cong\!\!820$
$2900\!\!\sim\!\!2847\!\!\cong\!\!929$
$2901\!\!\sim\!\!2887\!\!\cong\!\!731$
$2902\!\!\sim\!\!2851\!\!\cong\!\!929$
$2903\!\!\sim\!\!2838\!\!\cong\!\!750$
$2904\!\!\sim\!\!2860\!\!\cong\!\!2212$
$2905\!\!\sim\!\!1091\!\!\cong\!\!731$
$2906\!\!\sim\!\!1091\!\!\cong\!\!731$
$2907\!\!\sim\!\!2889\!\!\cong\!\!750$
$2908\!\!\sim\!\!2880\!\!\cong\!\!730$
$2909\!\!\sim\!\!2853\!\!\cong\!\!2853$
$2910\!\!\sim\!\!2889\!\!\cong\!\!750$
$2911\!\!\sim\!\!2853\!\!\cong\!\!2853$
$2912\!\!\sim\!\!2844\!\!\cong\!\!730$
$2913\!\!\sim\!\!2862\!\!\cong\!\!847$
$2914\!\!\sim\!\!2889\!\!\cong\!\!750$
$2915\!\!\sim\!\!2862\!\!\cong\!\!847$
$2916\!\!\sim\!\!1094\!\!\cong\!\!1090$
$2917\!\!\sim\!\!1094\!\!\cong\!\!1090$
$2918\!\!\sim\!\!1094\!\!\cong\!\!1090$
$2919\!\!\sim\!\!972\!\!\cong\!\!739$
$2920\!\!\sim\!\!1094\!\!\cong\!\!1090$
$2921\!\!\sim\!\!1094\!\!\cong\!\!1090$
$2922\!\!\sim\!\!891\!\!\cong\!\!891$
$2923\!\!\sim\!\!972\!\!\cong\!\!739$
$2924\!\!\sim\!\!891\!\!\cong\!\!891$
$2925\!\!\sim\!\!810\!\!\cong\!\!802$
$2926\!\!\sim\!\!1094\!\!\cong\!\!1090$
$2927\!\!\sim\!\!1094\!\!\cong\!\!1090$
$2928\!\!\sim\!\!945\!\!\cong\!\!941$
$2929\!\!\sim\!\!1094\!\!\cong\!\!1090$
$2930\!\!\sim\!\!1094\!\!\cong\!\!1090$
$2931\!\!\sim\!\!864\!\!\cong\!\!864$
$2932\!\!\sim\!\!963\!\!\cong\!\!963$
$2933\!\!\sim\!\!882\!\!\cong\!\!882$
$2934\!\!\sim\!\!783\!\!\cong\!\!775$
$2935\!\!\sim\!\!851\!\!\cong\!\!847$
$2936\!\!\sim\!\!878\!\!\cong\!\!878$
$2937\!\!\sim\!\!824\!\!\cong\!\!820$
$2938\!\!\sim\!\!860\!\!\cong\!\!860$
$2939\!\!\sim\!\!887\!\!\cong\!\!887$
$2940\!\!\sim\!\!824\!\!\cong\!\!820$
$2941\!\!\sim\!\!842\!\!\cong\!\!838$
$2942\!\!\sim\!\!869\!\!\cong\!\!869$
$2943\!\!\sim\!\!756\!\!\cong\!\!748$
$2944\!\!\sim\!\!1094\!\!\cong\!\!1090$
$2945\!\!\sim\!\!1094\!\!\cong\!\!1090$
$2946\!\!\sim\!\!963\!\!\cong\!\!963$
$2947\!\!\sim\!\!1094\!\!\cong\!\!1090$
$2948\!\!\sim\!\!1094\!\!\cong\!\!1090$
$2949\!\!\sim\!\!882\!\!\cong\!\!882$
$2950\!\!\sim\!\!945\!\!\cong\!\!941$
$2951\!\!\sim\!\!864\!\!\cong\!\!864$
$2952\!\!\sim\!\!783\!\!\cong\!\!775$
$2953\!\!\sim\!\!1094\!\!\cong\!\!1090$
$2954\!\!\sim\!\!1094\!\!\cong\!\!1090$
$2955\!\!\sim\!\!936\!\!\cong\!\!820$
$2956\!\!\sim\!\!1094\!\!\cong\!\!1090$
$2957\!\!\sim\!\!1094\!\!\cong\!\!1090$
$2958\!\!\sim\!\!855\!\!\cong\!\!847$
$2959\!\!\sim\!\!936\!\!\cong\!\!820$
$2960\!\!\sim\!\!855\!\!\cong\!\!847$
$2961\!\!\sim\!\!774\!\!\cong\!\!730$
$2962\!\!\sim\!\!932\!\!\cong\!\!820$
$2963\!\!\sim\!\!959\!\!\cong\!\!959$
$2964\!\!\sim\!\!824\!\!\cong\!\!820$
$2965\!\!\sim\!\!941\!\!\cong\!\!941$
$2966\!\!\sim\!\!968\!\!\cong\!\!968$
$2967\!\!\sim\!\!824\!\!\cong\!\!820$
$2968\!\!\sim\!\!923\!\!\cong\!\!923$
$2969\!\!\sim\!\!846\!\!\cong\!\!846$
$2970\!\!\sim\!\!747\!\!\cong\!\!739$
$2971\!\!\sim\!\!851\!\!\cong\!\!847$
$2972\!\!\sim\!\!860\!\!\cong\!\!860$
$2973\!\!\sim\!\!842\!\!\cong\!\!838$
$2974\!\!\sim\!\!878\!\!\cong\!\!878$
$2975\!\!\sim\!\!887\!\!\cong\!\!887$
$2976\!\!\sim\!\!869\!\!\cong\!\!869$
$2977\!\!\sim\!\!824\!\!\cong\!\!820$
$2978\!\!\sim\!\!824\!\!\cong\!\!820$
$2979\!\!\sim\!\!756\!\!\cong\!\!748$
$2980\!\!\sim\!\!932\!\!\cong\!\!820$
$2981\!\!\sim\!\!941\!\!\cong\!\!941$
$2982\!\!\sim\!\!923\!\!\cong\!\!923$
$2983\!\!\sim\!\!959\!\!\cong\!\!959$
$2984\!\!\sim\!\!968\!\!\cong\!\!968$
$2985\!\!\sim\!\!846\!\!\cong\!\!846$
$2986\!\!\sim\!\!824\!\!\cong\!\!820$
$2987\!\!\sim\!\!824\!\!\cong\!\!820$
$2988\!\!\sim\!\!747\!\!\cong\!\!739$
$2989\!\!\sim\!\!770\!\!\cong\!\!730$
$2990\!\!\sim\!\!779\!\!\cong\!\!779$
$2991\!\!\sim\!\!743\!\!\cong\!\!739$
$2992\!\!\sim\!\!779\!\!\cong\!\!779$
$2993\!\!\sim\!\!806\!\!\cong\!\!802$
$2994\!\!\sim\!\!752\!\!\cong\!\!752$
$2995\!\!\sim\!\!743\!\!\cong\!\!739$
$2996\!\!\sim\!\!752\!\!\cong\!\!752$
$2997\!\!\sim\!\!734\!\!\cong\!\!730$
$2998\!\!\sim\!\!1094\!\!\cong\!\!1090$
$2999\!\!\sim\!\!1094\!\!\cong\!\!1090$
$3000\!\!\sim\!\!969\!\!\cong\!\!969$
$3001\!\!\sim\!\!1094\!\!\cong\!\!1090$
$3002\!\!\sim\!\!1094\!\!\cong\!\!1090$
$3003\!\!\sim\!\!888\!\!\cong\!\!888$
$3004\!\!\sim\!\!969\!\!\cong\!\!969$
$3005\!\!\sim\!\!888\!\!\cong\!\!888$
$3006\!\!\sim\!\!807\!\!\cong\!\!771$
$3007\!\!\sim\!\!1094\!\!\cong\!\!1090$
$3008\!\!\sim\!\!1094\!\!\cong\!\!1090$
$3009\!\!\sim\!\!942\!\!\cong\!\!942$
$3010\!\!\sim\!\!1094\!\!\cong\!\!1090$
$3011\!\!\sim\!\!1094\!\!\cong\!\!1090$
$3012\!\!\sim\!\!861\!\!\cong\!\!861$
$3013\!\!\sim\!\!960\!\!\cong\!\!960$
$3014\!\!\sim\!\!879\!\!\cong\!\!879$
$3015\!\!\sim\!\!780\!\!\cong\!\!780$
$3016\!\!\sim\!\!852\!\!\cong\!\!852$
$3017\!\!\sim\!\!879\!\!\cong\!\!879$
$3018\!\!\sim\!\!824\!\!\cong\!\!820$
$3019\!\!\sim\!\!861\!\!\cong\!\!861$
$3020\!\!\sim\!\!888\!\!\cong\!\!888$
$3021\!\!\sim\!\!824\!\!\cong\!\!820$
$3022\!\!\sim\!\!843\!\!\cong\!\!843$
$3023\!\!\sim\!\!870\!\!\cong\!\!870$
$3024\!\!\sim\!\!753\!\!\cong\!\!753$
$3025\!\!\sim\!\!1094\!\!\cong\!\!1090$
$3026\!\!\sim\!\!1094\!\!\cong\!\!1090$
$3027\!\!\sim\!\!960\!\!\cong\!\!960$
$3028\!\!\sim\!\!1094\!\!\cong\!\!1090$
$3029\!\!\sim\!\!1094\!\!\cong\!\!1090$
$3030\!\!\sim\!\!879\!\!\cong\!\!879$
$3031\!\!\sim\!\!942\!\!\cong\!\!942$
$3032\!\!\sim\!\!861\!\!\cong\!\!861$
$3033\!\!\sim\!\!780\!\!\cong\!\!780$
$3034\!\!\sim\!\!1094\!\!\cong\!\!1090$
$3035\!\!\sim\!\!1094\!\!\cong\!\!1090$
$3036\!\!\sim\!\!933\!\!\cong\!\!849$
$3037\!\!\sim\!\!1094\!\!\cong\!\!1090$
$3038\!\!\sim\!\!1094\!\!\cong\!\!1090$
$3039\!\!\sim\!\!852\!\!\cong\!\!852$
$3040\!\!\sim\!\!933\!\!\cong\!\!849$
$3041\!\!\sim\!\!852\!\!\cong\!\!852$
$3042\!\!\sim\!\!771\!\!\cong\!\!771$
$3043\!\!\sim\!\!933\!\!\cong\!\!849$
$3044\!\!\sim\!\!960\!\!\cong\!\!960$
$3045\!\!\sim\!\!824\!\!\cong\!\!820$
$3046\!\!\sim\!\!942\!\!\cong\!\!942$
$3047\!\!\sim\!\!969\!\!\cong\!\!969$
$3048\!\!\sim\!\!824\!\!\cong\!\!820$
$3049\!\!\sim\!\!924\!\!\cong\!\!870$
$3050\!\!\sim\!\!843\!\!\cong\!\!843$
$3051\!\!\sim\!\!744\!\!\cong\!\!744$
$3052\!\!\sim\!\!852\!\!\cong\!\!852$
$3053\!\!\sim\!\!861\!\!\cong\!\!861$
$3054\!\!\sim\!\!843\!\!\cong\!\!843$
$3055\!\!\sim\!\!879\!\!\cong\!\!879$
$3056\!\!\sim\!\!888\!\!\cong\!\!888$
$3057\!\!\sim\!\!870\!\!\cong\!\!870$
$3058\!\!\sim\!\!824\!\!\cong\!\!820$
$3059\!\!\sim\!\!824\!\!\cong\!\!820$
$3060\!\!\sim\!\!753\!\!\cong\!\!753$
$3061\!\!\sim\!\!933\!\!\cong\!\!849$
$3062\!\!\sim\!\!942\!\!\cong\!\!942$
$3063\!\!\sim\!\!924\!\!\cong\!\!870$
$3064\!\!\sim\!\!960\!\!\cong\!\!960$
$3065\!\!\sim\!\!969\!\!\cong\!\!969$
$3066\!\!\sim\!\!843\!\!\cong\!\!843$
$3067\!\!\sim\!\!824\!\!\cong\!\!820$
$3068\!\!\sim\!\!824\!\!\cong\!\!820$
$3069\!\!\sim\!\!744\!\!\cong\!\!744$
$3070\!\!\sim\!\!771\!\!\cong\!\!771$
$3071\!\!\sim\!\!780\!\!\cong\!\!780$
$3072\!\!\sim\!\!744\!\!\cong\!\!744$
$3073\!\!\sim\!\!780\!\!\cong\!\!780$
$3074\!\!\sim\!\!807\!\!\cong\!\!771$
$3075\!\!\sim\!\!753\!\!\cong\!\!753$
$3076\!\!\sim\!\!744\!\!\cong\!\!744$
$3077\!\!\sim\!\!753\!\!\cong\!\!753$
$3078\!\!\sim\!\!734\!\!\cong\!\!730$
$3079\!\!\sim\!\!1091\!\!\cong\!\!731$
$3080\!\!\sim\!\!1091\!\!\cong\!\!731$
$3081\!\!\sim\!\!966\!\!\cong\!\!966$
$3082\!\!\sim\!\!1091\!\!\cong\!\!731$
$3083\!\!\sim\!\!1091\!\!\cong\!\!731$
$3084\!\!\sim\!\!885\!\!\cong\!\!885$
$3085\!\!\sim\!\!966\!\!\cong\!\!966$
$3086\!\!\sim\!\!885\!\!\cong\!\!885$
$3087\!\!\sim\!\!804\!\!\cong\!\!731$
$3088\!\!\sim\!\!1091\!\!\cong\!\!731$
$3089\!\!\sim\!\!1091\!\!\cong\!\!731$
$3090\!\!\sim\!\!939\!\!\cong\!\!939$
$3091\!\!\sim\!\!1091\!\!\cong\!\!731$
$3092\!\!\sim\!\!1091\!\!\cong\!\!731$
$3093\!\!\sim\!\!858\!\!\cong\!\!858$
$3094\!\!\sim\!\!957\!\!\cong\!\!957$
$3095\!\!\sim\!\!876\!\!\cong\!\!876$
$3096\!\!\sim\!\!777\!\!\cong\!\!777$
$3097\!\!\sim\!\!848\!\!\cong\!\!750$
$3098\!\!\sim\!\!875\!\!\cong\!\!875$
$3099\!\!\sim\!\!821\!\!\cong\!\!821$
$3100\!\!\sim\!\!857\!\!\cong\!\!857$
$3101\!\!\sim\!\!884\!\!\cong\!\!884$
$3102\!\!\sim\!\!821\!\!\cong\!\!821$
$3103\!\!\sim\!\!839\!\!\cong\!\!821$
$3104\!\!\sim\!\!866\!\!\cong\!\!866$
$3105\!\!\sim\!\!750\!\!\cong\!\!750$
$3106\!\!\sim\!\!1091\!\!\cong\!\!731$
$3107\!\!\sim\!\!1091\!\!\cong\!\!731$
$3108\!\!\sim\!\!957\!\!\cong\!\!957$
$3109\!\!\sim\!\!1091\!\!\cong\!\!731$
$3110\!\!\sim\!\!1091\!\!\cong\!\!731$
$3111\!\!\sim\!\!876\!\!\cong\!\!876$
$3112\!\!\sim\!\!939\!\!\cong\!\!939$
$3113\!\!\sim\!\!858\!\!\cong\!\!858$
$3114\!\!\sim\!\!777\!\!\cong\!\!777$
$3115\!\!\sim\!\!1091\!\!\cong\!\!731$
$3116\!\!\sim\!\!1091\!\!\cong\!\!731$
$3117\!\!\sim\!\!930\!\!\cong\!\!821$
$3118\!\!\sim\!\!1091\!\!\cong\!\!731$
$3119\!\!\sim\!\!1091\!\!\cong\!\!731$
$3120\!\!\sim\!\!849\!\!\cong\!\!849$
$3121\!\!\sim\!\!930\!\!\cong\!\!821$
$3122\!\!\sim\!\!849\!\!\cong\!\!849$
$3123\!\!\sim\!\!768\!\!\cong\!\!731$
$3124\!\!\sim\!\!929\!\!\cong\!\!929$
$3125\!\!\sim\!\!956\!\!\cong\!\!956$
$3126\!\!\sim\!\!821\!\!\cong\!\!821$
$3127\!\!\sim\!\!938\!\!\cong\!\!938$
$3128\!\!\sim\!\!965\!\!\cong\!\!965$
$3129\!\!\sim\!\!821\!\!\cong\!\!821$
$3130\!\!\sim\!\!920\!\!\cong\!\!920$
$3131\!\!\sim\!\!840\!\!\cong\!\!840$
$3132\!\!\sim\!\!741\!\!\cong\!\!741$
$3133\!\!\sim\!\!848\!\!\cong\!\!750$
$3134\!\!\sim\!\!857\!\!\cong\!\!857$
$3135\!\!\sim\!\!839\!\!\cong\!\!821$
$3136\!\!\sim\!\!875\!\!\cong\!\!875$
$3137\!\!\sim\!\!884\!\!\cong\!\!884$
$3138\!\!\sim\!\!866\!\!\cong\!\!866$
$3139\!\!\sim\!\!821\!\!\cong\!\!821$
$3140\!\!\sim\!\!821\!\!\cong\!\!821$
$3141\!\!\sim\!\!750\!\!\cong\!\!750$
$3142\!\!\sim\!\!929\!\!\cong\!\!929$
$3143\!\!\sim\!\!938\!\!\cong\!\!938$
$3144\!\!\sim\!\!920\!\!\cong\!\!920$
$3145\!\!\sim\!\!956\!\!\cong\!\!956$
$3146\!\!\sim\!\!965\!\!\cong\!\!965$
$3147\!\!\sim\!\!840\!\!\cong\!\!840$
$3148\!\!\sim\!\!821\!\!\cong\!\!821$
$3149\!\!\sim\!\!821\!\!\cong\!\!821$
$3150\!\!\sim\!\!741\!\!\cong\!\!741$
$3151\!\!\sim\!\!767\!\!\cong\!\!731$
$3152\!\!\sim\!\!776\!\!\cong\!\!776$
$3153\!\!\sim\!\!740\!\!\cong\!\!740$
$3154\!\!\sim\!\!776\!\!\cong\!\!776$
$3155\!\!\sim\!\!803\!\!\cong\!\!771$
$3156\!\!\sim\!\!749\!\!\cong\!\!749$
$3157\!\!\sim\!\!740\!\!\cong\!\!740$
$3158\!\!\sim\!\!749\!\!\cong\!\!749$
$3159\!\!\sim\!\!731\!\!\cong\!\!731$
$3160\!\!\sim\!\!1094\!\!\cong\!\!1090$
$3161\!\!\sim\!\!1094\!\!\cong\!\!1090$
$3162\!\!\sim\!\!969\!\!\cong\!\!969$
$3163\!\!\sim\!\!1094\!\!\cong\!\!1090$
$3164\!\!\sim\!\!1094\!\!\cong\!\!1090$
$3165\!\!\sim\!\!888\!\!\cong\!\!888$
$3166\!\!\sim\!\!969\!\!\cong\!\!969$
$3167\!\!\sim\!\!888\!\!\cong\!\!888$
$3168\!\!\sim\!\!807\!\!\cong\!\!771$
$3169\!\!\sim\!\!1094\!\!\cong\!\!1090$
$3170\!\!\sim\!\!1094\!\!\cong\!\!1090$
$3171\!\!\sim\!\!942\!\!\cong\!\!942$
$3172\!\!\sim\!\!1094\!\!\cong\!\!1090$
$3173\!\!\sim\!\!1094\!\!\cong\!\!1090$
$3174\!\!\sim\!\!861\!\!\cong\!\!861$
$3175\!\!\sim\!\!960\!\!\cong\!\!960$
$3176\!\!\sim\!\!879\!\!\cong\!\!879$
$3177\!\!\sim\!\!780\!\!\cong\!\!780$
$3178\!\!\sim\!\!852\!\!\cong\!\!852$
$3179\!\!\sim\!\!879\!\!\cong\!\!879$
$3180\!\!\sim\!\!824\!\!\cong\!\!820$
$3181\!\!\sim\!\!861\!\!\cong\!\!861$
$3182\!\!\sim\!\!888\!\!\cong\!\!888$
$3183\!\!\sim\!\!824\!\!\cong\!\!820$
$3184\!\!\sim\!\!843\!\!\cong\!\!843$
$3185\!\!\sim\!\!870\!\!\cong\!\!870$
$3186\!\!\sim\!\!753\!\!\cong\!\!753$
$3187\!\!\sim\!\!1094\!\!\cong\!\!1090$
$3188\!\!\sim\!\!1094\!\!\cong\!\!1090$
$3189\!\!\sim\!\!960\!\!\cong\!\!960$
$3190\!\!\sim\!\!1094\!\!\cong\!\!1090$
$3191\!\!\sim\!\!1094\!\!\cong\!\!1090$
$3192\!\!\sim\!\!879\!\!\cong\!\!879$
$3193\!\!\sim\!\!942\!\!\cong\!\!942$
$3194\!\!\sim\!\!861\!\!\cong\!\!861$
$3195\!\!\sim\!\!780\!\!\cong\!\!780$
$3196\!\!\sim\!\!1094\!\!\cong\!\!1090$
$3197\!\!\sim\!\!1094\!\!\cong\!\!1090$
$3198\!\!\sim\!\!933\!\!\cong\!\!849$
$3199\!\!\sim\!\!1094\!\!\cong\!\!1090$
$3200\!\!\sim\!\!1094\!\!\cong\!\!1090$
$3201\!\!\sim\!\!852\!\!\cong\!\!852$
$3202\!\!\sim\!\!933\!\!\cong\!\!849$
$3203\!\!\sim\!\!852\!\!\cong\!\!852$
$3204\!\!\sim\!\!771\!\!\cong\!\!771$
$3205\!\!\sim\!\!933\!\!\cong\!\!849$
$3206\!\!\sim\!\!960\!\!\cong\!\!960$
$3207\!\!\sim\!\!824\!\!\cong\!\!820$
$3208\!\!\sim\!\!942\!\!\cong\!\!942$
$3209\!\!\sim\!\!969\!\!\cong\!\!969$
$3210\!\!\sim\!\!824\!\!\cong\!\!820$
$3211\!\!\sim\!\!924\!\!\cong\!\!870$
$3212\!\!\sim\!\!843\!\!\cong\!\!843$
$3213\!\!\sim\!\!744\!\!\cong\!\!744$
$3214\!\!\sim\!\!852\!\!\cong\!\!852$
$3215\!\!\sim\!\!861\!\!\cong\!\!861$
$3216\!\!\sim\!\!843\!\!\cong\!\!843$
$3217\!\!\sim\!\!879\!\!\cong\!\!879$
$3218\!\!\sim\!\!888\!\!\cong\!\!888$
$3219\!\!\sim\!\!870\!\!\cong\!\!870$
$3220\!\!\sim\!\!824\!\!\cong\!\!820$
$3221\!\!\sim\!\!824\!\!\cong\!\!820$
$3222\!\!\sim\!\!753\!\!\cong\!\!753$
$3223\!\!\sim\!\!933\!\!\cong\!\!849$
$3224\!\!\sim\!\!942\!\!\cong\!\!942$
$3225\!\!\sim\!\!924\!\!\cong\!\!870$
$3226\!\!\sim\!\!960\!\!\cong\!\!960$
$3227\!\!\sim\!\!969\!\!\cong\!\!969$
$3228\!\!\sim\!\!843\!\!\cong\!\!843$
$3229\!\!\sim\!\!824\!\!\cong\!\!820$
$3230\!\!\sim\!\!824\!\!\cong\!\!820$
$3231\!\!\sim\!\!744\!\!\cong\!\!744$
$3232\!\!\sim\!\!771\!\!\cong\!\!771$
$3233\!\!\sim\!\!780\!\!\cong\!\!780$
$3234\!\!\sim\!\!744\!\!\cong\!\!744$
$3235\!\!\sim\!\!780\!\!\cong\!\!780$
$3236\!\!\sim\!\!807\!\!\cong\!\!771$
$3237\!\!\sim\!\!753\!\!\cong\!\!753$
$3238\!\!\sim\!\!744\!\!\cong\!\!744$
$3239\!\!\sim\!\!753\!\!\cong\!\!753$
$3240\!\!\sim\!\!734\!\!\cong\!\!730$
$3241\!\!\sim\!\!1094\!\!\cong\!\!1090$
$3242\!\!\sim\!\!1094\!\!\cong\!\!1090$
$3243\!\!\sim\!\!968\!\!\cong\!\!968$
$3244\!\!\sim\!\!1094\!\!\cong\!\!1090$
$3245\!\!\sim\!\!1094\!\!\cong\!\!1090$
$3246\!\!\sim\!\!887\!\!\cong\!\!887$
$3247\!\!\sim\!\!968\!\!\cong\!\!968$
$3248\!\!\sim\!\!887\!\!\cong\!\!887$
$3249\!\!\sim\!\!806\!\!\cong\!\!802$
$3250\!\!\sim\!\!1094\!\!\cong\!\!1090$
$3251\!\!\sim\!\!1094\!\!\cong\!\!1090$
$3252\!\!\sim\!\!941\!\!\cong\!\!941$
$3253\!\!\sim\!\!1094\!\!\cong\!\!1090$
$3254\!\!\sim\!\!1094\!\!\cong\!\!1090$
$3255\!\!\sim\!\!860\!\!\cong\!\!860$
$3256\!\!\sim\!\!959\!\!\cong\!\!959$
$3257\!\!\sim\!\!878\!\!\cong\!\!878$
$3258\!\!\sim\!\!779\!\!\cong\!\!779$
$3259\!\!\sim\!\!855\!\!\cong\!\!847$
$3260\!\!\sim\!\!882\!\!\cong\!\!882$
$3261\!\!\sim\!\!824\!\!\cong\!\!820$
$3262\!\!\sim\!\!864\!\!\cong\!\!864$
$3263\!\!\sim\!\!891\!\!\cong\!\!891$
$3264\!\!\sim\!\!824\!\!\cong\!\!820$
$3265\!\!\sim\!\!846\!\!\cong\!\!846$
$3266\!\!\sim\!\!869\!\!\cong\!\!869$
$3267\!\!\sim\!\!752\!\!\cong\!\!752$
$3268\!\!\sim\!\!1094\!\!\cong\!\!1090$
$3269\!\!\sim\!\!1094\!\!\cong\!\!1090$
$3270\!\!\sim\!\!959\!\!\cong\!\!959$
$3271\!\!\sim\!\!1094\!\!\cong\!\!1090$
$3272\!\!\sim\!\!1094\!\!\cong\!\!1090$
$3273\!\!\sim\!\!878\!\!\cong\!\!878$
$3274\!\!\sim\!\!941\!\!\cong\!\!941$
$3275\!\!\sim\!\!860\!\!\cong\!\!860$
$3276\!\!\sim\!\!779\!\!\cong\!\!779$
$3277\!\!\sim\!\!1094\!\!\cong\!\!1090$
$3278\!\!\sim\!\!1094\!\!\cong\!\!1090$
$3279\!\!\sim\!\!932\!\!\cong\!\!820$
$3280\!\!\sim\!\!1094\!\!\cong\!\!1090$
$3281\!\!\sim\!\!1094\!\!\cong\!\!1090$
$3282\!\!\sim\!\!851\!\!\cong\!\!847$
$3283\!\!\sim\!\!932\!\!\cong\!\!820$
$3284\!\!\sim\!\!851\!\!\cong\!\!847$
$3285\!\!\sim\!\!770\!\!\cong\!\!730$
$3286\!\!\sim\!\!936\!\!\cong\!\!820$
$3287\!\!\sim\!\!963\!\!\cong\!\!963$
$3288\!\!\sim\!\!824\!\!\cong\!\!820$
$3289\!\!\sim\!\!945\!\!\cong\!\!941$
$3290\!\!\sim\!\!972\!\!\cong\!\!739$
$3291\!\!\sim\!\!824\!\!\cong\!\!820$
$3292\!\!\sim\!\!923\!\!\cong\!\!923$
$3293\!\!\sim\!\!842\!\!\cong\!\!838$
$3294\!\!\sim\!\!743\!\!\cong\!\!739$
$3295\!\!\sim\!\!855\!\!\cong\!\!847$
$3296\!\!\sim\!\!864\!\!\cong\!\!864$
$3297\!\!\sim\!\!846\!\!\cong\!\!846$
$3298\!\!\sim\!\!882\!\!\cong\!\!882$
$3299\!\!\sim\!\!891\!\!\cong\!\!891$
$3300\!\!\sim\!\!869\!\!\cong\!\!869$
$3301\!\!\sim\!\!824\!\!\cong\!\!820$
$3302\!\!\sim\!\!824\!\!\cong\!\!820$
$3303\!\!\sim\!\!752\!\!\cong\!\!752$
$3304\!\!\sim\!\!936\!\!\cong\!\!820$
$3305\!\!\sim\!\!945\!\!\cong\!\!941$
$3306\!\!\sim\!\!923\!\!\cong\!\!923$
$3307\!\!\sim\!\!963\!\!\cong\!\!963$
$3308\!\!\sim\!\!972\!\!\cong\!\!739$
$3309\!\!\sim\!\!842\!\!\cong\!\!838$
$3310\!\!\sim\!\!824\!\!\cong\!\!820$
$3311\!\!\sim\!\!824\!\!\cong\!\!820$
$3312\!\!\sim\!\!743\!\!\cong\!\!739$
$3313\!\!\sim\!\!774\!\!\cong\!\!730$
$3314\!\!\sim\!\!783\!\!\cong\!\!775$
$3315\!\!\sim\!\!747\!\!\cong\!\!739$
$3316\!\!\sim\!\!783\!\!\cong\!\!775$
$3317\!\!\sim\!\!810\!\!\cong\!\!802$
$3318\!\!\sim\!\!756\!\!\cong\!\!748$
$3319\!\!\sim\!\!747\!\!\cong\!\!739$
$3320\!\!\sim\!\!756\!\!\cong\!\!748$
$3321\!\!\sim\!\!734\!\!\cong\!\!730$
$3322\!\!\sim\!\!1091\!\!\cong\!\!731$
$3323\!\!\sim\!\!1091\!\!\cong\!\!731$
$3324\!\!\sim\!\!965\!\!\cong\!\!965$
$3325\!\!\sim\!\!1091\!\!\cong\!\!731$
$3326\!\!\sim\!\!1091\!\!\cong\!\!731$
$3327\!\!\sim\!\!884\!\!\cong\!\!884$
$3328\!\!\sim\!\!965\!\!\cong\!\!965$
$3329\!\!\sim\!\!884\!\!\cong\!\!884$
$3330\!\!\sim\!\!803\!\!\cong\!\!771$
$3331\!\!\sim\!\!1091\!\!\cong\!\!731$
$3332\!\!\sim\!\!1091\!\!\cong\!\!731$
$3333\!\!\sim\!\!938\!\!\cong\!\!938$
$3334\!\!\sim\!\!1091\!\!\cong\!\!731$
$3335\!\!\sim\!\!1091\!\!\cong\!\!731$
$3336\!\!\sim\!\!857\!\!\cong\!\!857$
$3337\!\!\sim\!\!956\!\!\cong\!\!956$
$3338\!\!\sim\!\!875\!\!\cong\!\!875$
$3339\!\!\sim\!\!776\!\!\cong\!\!776$
$3340\!\!\sim\!\!849\!\!\cong\!\!849$
$3341\!\!\sim\!\!876\!\!\cong\!\!876$
$3342\!\!\sim\!\!821\!\!\cong\!\!821$
$3343\!\!\sim\!\!858\!\!\cong\!\!858$
$3344\!\!\sim\!\!885\!\!\cong\!\!885$
$3345\!\!\sim\!\!821\!\!\cong\!\!821$
$3346\!\!\sim\!\!840\!\!\cong\!\!840$
$3347\!\!\sim\!\!866\!\!\cong\!\!866$
$3348\!\!\sim\!\!749\!\!\cong\!\!749$
$3349\!\!\sim\!\!1091\!\!\cong\!\!731$
$3350\!\!\sim\!\!1091\!\!\cong\!\!731$
$3351\!\!\sim\!\!956\!\!\cong\!\!956$
$3352\!\!\sim\!\!1091\!\!\cong\!\!731$
$3353\!\!\sim\!\!1091\!\!\cong\!\!731$
$3354\!\!\sim\!\!875\!\!\cong\!\!875$
$3355\!\!\sim\!\!938\!\!\cong\!\!938$
$3356\!\!\sim\!\!857\!\!\cong\!\!857$
$3357\!\!\sim\!\!776\!\!\cong\!\!776$
$3358\!\!\sim\!\!1091\!\!\cong\!\!731$
$3359\!\!\sim\!\!1091\!\!\cong\!\!731$
$3360\!\!\sim\!\!929\!\!\cong\!\!929$
$3361\!\!\sim\!\!1091\!\!\cong\!\!731$
$3362\!\!\sim\!\!1091\!\!\cong\!\!731$
$3363\!\!\sim\!\!848\!\!\cong\!\!750$
$3364\!\!\sim\!\!929\!\!\cong\!\!929$
$3365\!\!\sim\!\!848\!\!\cong\!\!750$
$3366\!\!\sim\!\!767\!\!\cong\!\!731$
$3367\!\!\sim\!\!930\!\!\cong\!\!821$
$3368\!\!\sim\!\!957\!\!\cong\!\!957$
$3369\!\!\sim\!\!821\!\!\cong\!\!821$
$3370\!\!\sim\!\!939\!\!\cong\!\!939$
$3371\!\!\sim\!\!966\!\!\cong\!\!966$
$3372\!\!\sim\!\!821\!\!\cong\!\!821$
$3373\!\!\sim\!\!920\!\!\cong\!\!920$
$3374\!\!\sim\!\!839\!\!\cong\!\!821$
$3375\!\!\sim\!\!740\!\!\cong\!\!740$
$3376\!\!\sim\!\!849\!\!\cong\!\!849$
$3377\!\!\sim\!\!858\!\!\cong\!\!858$
$3378\!\!\sim\!\!840\!\!\cong\!\!840$
$3379\!\!\sim\!\!876\!\!\cong\!\!876$
$3380\!\!\sim\!\!885\!\!\cong\!\!885$
$3381\!\!\sim\!\!866\!\!\cong\!\!866$
$3382\!\!\sim\!\!821\!\!\cong\!\!821$
$3383\!\!\sim\!\!821\!\!\cong\!\!821$
$3384\!\!\sim\!\!749\!\!\cong\!\!749$
$3385\!\!\sim\!\!930\!\!\cong\!\!821$
$3386\!\!\sim\!\!939\!\!\cong\!\!939$
$3387\!\!\sim\!\!920\!\!\cong\!\!920$
$3388\!\!\sim\!\!957\!\!\cong\!\!957$
$3389\!\!\sim\!\!966\!\!\cong\!\!966$
$3390\!\!\sim\!\!839\!\!\cong\!\!821$
$3391\!\!\sim\!\!821\!\!\cong\!\!821$
$3392\!\!\sim\!\!821\!\!\cong\!\!821$
$3393\!\!\sim\!\!740\!\!\cong\!\!740$
$3394\!\!\sim\!\!768\!\!\cong\!\!731$
$3395\!\!\sim\!\!777\!\!\cong\!\!777$
$3396\!\!\sim\!\!741\!\!\cong\!\!741$
$3397\!\!\sim\!\!777\!\!\cong\!\!777$
$3398\!\!\sim\!\!804\!\!\cong\!\!731$
$3399\!\!\sim\!\!750\!\!\cong\!\!750$
$3400\!\!\sim\!\!741\!\!\cong\!\!741$
$3401\!\!\sim\!\!750\!\!\cong\!\!750$
$3402\!\!\sim\!\!731\!\!\cong\!\!731$
$3403\!\!\sim\!\!1091\!\!\cong\!\!731$
$3404\!\!\sim\!\!1091\!\!\cong\!\!731$
$3405\!\!\sim\!\!966\!\!\cong\!\!966$
$3406\!\!\sim\!\!1091\!\!\cong\!\!731$
$3407\!\!\sim\!\!1091\!\!\cong\!\!731$
$3408\!\!\sim\!\!885\!\!\cong\!\!885$
$3409\!\!\sim\!\!966\!\!\cong\!\!966$
$3410\!\!\sim\!\!885\!\!\cong\!\!885$
$3411\!\!\sim\!\!804\!\!\cong\!\!731$
$3412\!\!\sim\!\!1091\!\!\cong\!\!731$
$3413\!\!\sim\!\!1091\!\!\cong\!\!731$
$3414\!\!\sim\!\!939\!\!\cong\!\!939$
$3415\!\!\sim\!\!1091\!\!\cong\!\!731$
$3416\!\!\sim\!\!1091\!\!\cong\!\!731$
$3417\!\!\sim\!\!858\!\!\cong\!\!858$
$3418\!\!\sim\!\!957\!\!\cong\!\!957$
$3419\!\!\sim\!\!876\!\!\cong\!\!876$
$3420\!\!\sim\!\!777\!\!\cong\!\!777$
$3421\!\!\sim\!\!848\!\!\cong\!\!750$
$3422\!\!\sim\!\!875\!\!\cong\!\!875$
$3423\!\!\sim\!\!821\!\!\cong\!\!821$
$3424\!\!\sim\!\!857\!\!\cong\!\!857$
$3425\!\!\sim\!\!884\!\!\cong\!\!884$
$3426\!\!\sim\!\!821\!\!\cong\!\!821$
$3427\!\!\sim\!\!839\!\!\cong\!\!821$
$3428\!\!\sim\!\!866\!\!\cong\!\!866$
$3429\!\!\sim\!\!750\!\!\cong\!\!750$
$3430\!\!\sim\!\!1091\!\!\cong\!\!731$
$3431\!\!\sim\!\!1091\!\!\cong\!\!731$
$3432\!\!\sim\!\!957\!\!\cong\!\!957$
$3433\!\!\sim\!\!1091\!\!\cong\!\!731$
$3434\!\!\sim\!\!1091\!\!\cong\!\!731$
$3435\!\!\sim\!\!876\!\!\cong\!\!876$
$3436\!\!\sim\!\!939\!\!\cong\!\!939$
$3437\!\!\sim\!\!858\!\!\cong\!\!858$
$3438\!\!\sim\!\!777\!\!\cong\!\!777$
$3439\!\!\sim\!\!1091\!\!\cong\!\!731$
$3440\!\!\sim\!\!1091\!\!\cong\!\!731$
$3441\!\!\sim\!\!930\!\!\cong\!\!821$
$3442\!\!\sim\!\!1091\!\!\cong\!\!731$
$3443\!\!\sim\!\!1091\!\!\cong\!\!731$
$3444\!\!\sim\!\!849\!\!\cong\!\!849$
$3445\!\!\sim\!\!930\!\!\cong\!\!821$
$3446\!\!\sim\!\!849\!\!\cong\!\!849$
$3447\!\!\sim\!\!768\!\!\cong\!\!731$
$3448\!\!\sim\!\!929\!\!\cong\!\!929$
$3449\!\!\sim\!\!956\!\!\cong\!\!956$
$3450\!\!\sim\!\!821\!\!\cong\!\!821$
$3451\!\!\sim\!\!938\!\!\cong\!\!938$
$3452\!\!\sim\!\!965\!\!\cong\!\!965$
$3453\!\!\sim\!\!821\!\!\cong\!\!821$
$3454\!\!\sim\!\!920\!\!\cong\!\!920$
$3455\!\!\sim\!\!840\!\!\cong\!\!840$
$3456\!\!\sim\!\!741\!\!\cong\!\!741$
$3457\!\!\sim\!\!848\!\!\cong\!\!750$
$3458\!\!\sim\!\!857\!\!\cong\!\!857$
$3459\!\!\sim\!\!839\!\!\cong\!\!821$
$3460\!\!\sim\!\!875\!\!\cong\!\!875$
$3461\!\!\sim\!\!884\!\!\cong\!\!884$
$3462\!\!\sim\!\!866\!\!\cong\!\!866$
$3463\!\!\sim\!\!821\!\!\cong\!\!821$
$3464\!\!\sim\!\!821\!\!\cong\!\!821$
$3465\!\!\sim\!\!750\!\!\cong\!\!750$
$3466\!\!\sim\!\!929\!\!\cong\!\!929$
$3467\!\!\sim\!\!938\!\!\cong\!\!938$
$3468\!\!\sim\!\!920\!\!\cong\!\!920$
$3469\!\!\sim\!\!956\!\!\cong\!\!956$
$3470\!\!\sim\!\!965\!\!\cong\!\!965$
$3471\!\!\sim\!\!840\!\!\cong\!\!840$
$3472\!\!\sim\!\!821\!\!\cong\!\!821$
$3473\!\!\sim\!\!821\!\!\cong\!\!821$
$3474\!\!\sim\!\!741\!\!\cong\!\!741$
$3475\!\!\sim\!\!767\!\!\cong\!\!731$
$3476\!\!\sim\!\!776\!\!\cong\!\!776$
$3477\!\!\sim\!\!740\!\!\cong\!\!740$
$3478\!\!\sim\!\!776\!\!\cong\!\!776$
$3479\!\!\sim\!\!803\!\!\cong\!\!771$
$3480\!\!\sim\!\!749\!\!\cong\!\!749$
$3481\!\!\sim\!\!740\!\!\cong\!\!740$
$3482\!\!\sim\!\!749\!\!\cong\!\!749$
$3483\!\!\sim\!\!731\!\!\cong\!\!731$
$3484\!\!\sim\!\!1091\!\!\cong\!\!731$
$3485\!\!\sim\!\!1091\!\!\cong\!\!731$
$3486\!\!\sim\!\!965\!\!\cong\!\!965$
$3487\!\!\sim\!\!1091\!\!\cong\!\!731$
$3488\!\!\sim\!\!1091\!\!\cong\!\!731$
$3489\!\!\sim\!\!884\!\!\cong\!\!884$
$3490\!\!\sim\!\!965\!\!\cong\!\!965$
$3491\!\!\sim\!\!884\!\!\cong\!\!884$
$3492\!\!\sim\!\!803\!\!\cong\!\!771$
$3493\!\!\sim\!\!1091\!\!\cong\!\!731$
$3494\!\!\sim\!\!1091\!\!\cong\!\!731$
$3495\!\!\sim\!\!938\!\!\cong\!\!938$
$3496\!\!\sim\!\!1091\!\!\cong\!\!731$
$3497\!\!\sim\!\!1091\!\!\cong\!\!731$
$3498\!\!\sim\!\!857\!\!\cong\!\!857$
$3499\!\!\sim\!\!956\!\!\cong\!\!956$
$3500\!\!\sim\!\!875\!\!\cong\!\!875$
$3501\!\!\sim\!\!776\!\!\cong\!\!776$
$3502\!\!\sim\!\!849\!\!\cong\!\!849$
$3503\!\!\sim\!\!876\!\!\cong\!\!876$
$3504\!\!\sim\!\!821\!\!\cong\!\!821$
$3505\!\!\sim\!\!858\!\!\cong\!\!858$
$3506\!\!\sim\!\!885\!\!\cong\!\!885$
$3507\!\!\sim\!\!821\!\!\cong\!\!821$
$3508\!\!\sim\!\!840\!\!\cong\!\!840$
$3509\!\!\sim\!\!866\!\!\cong\!\!866$
$3510\!\!\sim\!\!749\!\!\cong\!\!749$
$3511\!\!\sim\!\!1091\!\!\cong\!\!731$
$3512\!\!\sim\!\!1091\!\!\cong\!\!731$
$3513\!\!\sim\!\!956\!\!\cong\!\!956$
$3514\!\!\sim\!\!1091\!\!\cong\!\!731$
$3515\!\!\sim\!\!1091\!\!\cong\!\!731$
$3516\!\!\sim\!\!875\!\!\cong\!\!875$
$3517\!\!\sim\!\!938\!\!\cong\!\!938$
$3518\!\!\sim\!\!857\!\!\cong\!\!857$
$3519\!\!\sim\!\!776\!\!\cong\!\!776$
$3520\!\!\sim\!\!1091\!\!\cong\!\!731$
$3521\!\!\sim\!\!1091\!\!\cong\!\!731$
$3522\!\!\sim\!\!929\!\!\cong\!\!929$
$3523\!\!\sim\!\!1091\!\!\cong\!\!731$
$3524\!\!\sim\!\!1091\!\!\cong\!\!731$
$3525\!\!\sim\!\!848\!\!\cong\!\!750$
$3526\!\!\sim\!\!929\!\!\cong\!\!929$
$3527\!\!\sim\!\!848\!\!\cong\!\!750$
$3528\!\!\sim\!\!767\!\!\cong\!\!731$
$3529\!\!\sim\!\!930\!\!\cong\!\!821$
$3530\!\!\sim\!\!957\!\!\cong\!\!957$
$3531\!\!\sim\!\!821\!\!\cong\!\!821$
$3532\!\!\sim\!\!939\!\!\cong\!\!939$
$3533\!\!\sim\!\!966\!\!\cong\!\!966$
$3534\!\!\sim\!\!821\!\!\cong\!\!821$
$3535\!\!\sim\!\!920\!\!\cong\!\!920$
$3536\!\!\sim\!\!839\!\!\cong\!\!821$
$3537\!\!\sim\!\!740\!\!\cong\!\!740$
$3538\!\!\sim\!\!849\!\!\cong\!\!849$
$3539\!\!\sim\!\!858\!\!\cong\!\!858$
$3540\!\!\sim\!\!840\!\!\cong\!\!840$
$3541\!\!\sim\!\!876\!\!\cong\!\!876$
$3542\!\!\sim\!\!885\!\!\cong\!\!885$
$3543\!\!\sim\!\!866\!\!\cong\!\!866$
$3544\!\!\sim\!\!821\!\!\cong\!\!821$
$3545\!\!\sim\!\!821\!\!\cong\!\!821$
$3546\!\!\sim\!\!749\!\!\cong\!\!749$
$3547\!\!\sim\!\!930\!\!\cong\!\!821$
$3548\!\!\sim\!\!939\!\!\cong\!\!939$
$3549\!\!\sim\!\!920\!\!\cong\!\!920$
$3550\!\!\sim\!\!957\!\!\cong\!\!957$
$3551\!\!\sim\!\!966\!\!\cong\!\!966$
$3552\!\!\sim\!\!839\!\!\cong\!\!821$
$3553\!\!\sim\!\!821\!\!\cong\!\!821$
$3554\!\!\sim\!\!821\!\!\cong\!\!821$
$3555\!\!\sim\!\!740\!\!\cong\!\!740$
$3556\!\!\sim\!\!768\!\!\cong\!\!731$
$3557\!\!\sim\!\!777\!\!\cong\!\!777$
$3558\!\!\sim\!\!741\!\!\cong\!\!741$
$3559\!\!\sim\!\!777\!\!\cong\!\!777$
$3560\!\!\sim\!\!804\!\!\cong\!\!731$
$3561\!\!\sim\!\!750\!\!\cong\!\!750$
$3562\!\!\sim\!\!741\!\!\cong\!\!741$
$3563\!\!\sim\!\!750\!\!\cong\!\!750$
$3564\!\!\sim\!\!731\!\!\cong\!\!731$
$3565\!\!\sim\!\!1090\!\!\cong\!\!1090$
$3566\!\!\sim\!\!1090\!\!\cong\!\!1090$
$3567\!\!\sim\!\!964\!\!\cong\!\!739$
$3568\!\!\sim\!\!1090\!\!\cong\!\!1090$
$3569\!\!\sim\!\!1090\!\!\cong\!\!1090$
$3570\!\!\sim\!\!883\!\!\cong\!\!883$
$3571\!\!\sim\!\!964\!\!\cong\!\!739$
$3572\!\!\sim\!\!883\!\!\cong\!\!883$
$3573\!\!\sim\!\!802\!\!\cong\!\!802$
$3574\!\!\sim\!\!1090\!\!\cong\!\!1090$
$3575\!\!\sim\!\!1090\!\!\cong\!\!1090$
$3576\!\!\sim\!\!937\!\!\cong\!\!937$
$3577\!\!\sim\!\!1090\!\!\cong\!\!1090$
$3578\!\!\sim\!\!1090\!\!\cong\!\!1090$
$3579\!\!\sim\!\!856\!\!\cong\!\!856$
$3580\!\!\sim\!\!955\!\!\cong\!\!937$
$3581\!\!\sim\!\!874\!\!\cong\!\!874$
$3582\!\!\sim\!\!775\!\!\cong\!\!775$
$3583\!\!\sim\!\!847\!\!\cong\!\!847$
$3584\!\!\sim\!\!874\!\!\cong\!\!874$
$3585\!\!\sim\!\!820\!\!\cong\!\!820$
$3586\!\!\sim\!\!856\!\!\cong\!\!856$
$3587\!\!\sim\!\!883\!\!\cong\!\!883$
$3588\!\!\sim\!\!820\!\!\cong\!\!820$
$3589\!\!\sim\!\!838\!\!\cong\!\!838$
$3590\!\!\sim\!\!865\!\!\cong\!\!820$
$3591\!\!\sim\!\!748\!\!\cong\!\!748$
$3592\!\!\sim\!\!1090\!\!\cong\!\!1090$
$3593\!\!\sim\!\!1090\!\!\cong\!\!1090$
$3594\!\!\sim\!\!955\!\!\cong\!\!937$
$3595\!\!\sim\!\!1090\!\!\cong\!\!1090$
$3596\!\!\sim\!\!1090\!\!\cong\!\!1090$
$3597\!\!\sim\!\!874\!\!\cong\!\!874$
$3598\!\!\sim\!\!937\!\!\cong\!\!937$
$3599\!\!\sim\!\!856\!\!\cong\!\!856$
$3600\!\!\sim\!\!775\!\!\cong\!\!775$
$3601\!\!\sim\!\!1090\!\!\cong\!\!1090$
$3602\!\!\sim\!\!1090\!\!\cong\!\!1090$
$3603\!\!\sim\!\!928\!\!\cong\!\!820$
$3604\!\!\sim\!\!1090\!\!\cong\!\!1090$
$3605\!\!\sim\!\!1090\!\!\cong\!\!1090$
$3606\!\!\sim\!\!847\!\!\cong\!\!847$
$3607\!\!\sim\!\!928\!\!\cong\!\!820$
$3608\!\!\sim\!\!847\!\!\cong\!\!847$
$3609\!\!\sim\!\!766\!\!\cong\!\!730$
$3610\!\!\sim\!\!928\!\!\cong\!\!820$
$3611\!\!\sim\!\!955\!\!\cong\!\!937$
$3612\!\!\sim\!\!820\!\!\cong\!\!820$
$3613\!\!\sim\!\!937\!\!\cong\!\!937$
$3614\!\!\sim\!\!964\!\!\cong\!\!739$
$3615\!\!\sim\!\!820\!\!\cong\!\!820$
$3616\!\!\sim\!\!919\!\!\cong\!\!820$
$3617\!\!\sim\!\!838\!\!\cong\!\!838$
$3618\!\!\sim\!\!739\!\!\cong\!\!739$
$3619\!\!\sim\!\!847\!\!\cong\!\!847$
$3620\!\!\sim\!\!856\!\!\cong\!\!856$
$3621\!\!\sim\!\!838\!\!\cong\!\!838$
$3622\!\!\sim\!\!874\!\!\cong\!\!874$
$3623\!\!\sim\!\!883\!\!\cong\!\!883$
$3624\!\!\sim\!\!865\!\!\cong\!\!820$
$3625\!\!\sim\!\!820\!\!\cong\!\!820$
$3626\!\!\sim\!\!820\!\!\cong\!\!820$
$3627\!\!\sim\!\!748\!\!\cong\!\!748$
$3628\!\!\sim\!\!928\!\!\cong\!\!820$
$3629\!\!\sim\!\!937\!\!\cong\!\!937$
$3630\!\!\sim\!\!919\!\!\cong\!\!820$
$3631\!\!\sim\!\!955\!\!\cong\!\!937$
$3632\!\!\sim\!\!964\!\!\cong\!\!739$
$3633\!\!\sim\!\!838\!\!\cong\!\!838$
$3634\!\!\sim\!\!820\!\!\cong\!\!820$
$3635\!\!\sim\!\!820\!\!\cong\!\!820$
$3636\!\!\sim\!\!739\!\!\cong\!\!739$
$3637\!\!\sim\!\!766\!\!\cong\!\!730$
$3638\!\!\sim\!\!775\!\!\cong\!\!775$
$3639\!\!\sim\!\!739\!\!\cong\!\!739$
$3640\!\!\sim\!\!775\!\!\cong\!\!775$
$3641\!\!\sim\!\!802\!\!\cong\!\!802$
$3642\!\!\sim\!\!748\!\!\cong\!\!748$
$3643\!\!\sim\!\!739\!\!\cong\!\!739$
$3644\!\!\sim\!\!748\!\!\cong\!\!748$
$3645\!\!\sim\!\!730\!\!\cong\!\!730$
$3646\!\!\sim\!\!730\!\!\cong\!\!730$
$3647\!\!\sim\!\!2190\!\!\cong\!\!750$
$3648\!\!\sim\!\!730\!\!\cong\!\!730$
$3649\!\!\sim\!\!2190\!\!\cong\!\!750$
$3650\!\!\sim\!\!2196\!\!\cong\!\!802$
$3651\!\!\sim\!\!2193\!\!\cong\!\!2193$
$3652\!\!\sim\!\!730\!\!\cong\!\!730$
$3653\!\!\sim\!\!2193\!\!\cong\!\!2193$
$3654\!\!\sim\!\!730\!\!\cong\!\!730$
$3655\!\!\sim\!\!820\!\!\cong\!\!820$
$3656\!\!\sim\!\!2352\!\!\cong\!\!740$
$3657\!\!\sim\!\!820\!\!\cong\!\!820$
$3658\!\!\sim\!\!2352\!\!\cong\!\!740$
$3659\!\!\sim\!\!2358\!\!\cong\!\!820$
$3660\!\!\sim\!\!2355\!\!\cong\!\!2355$
$3661\!\!\sim\!\!820\!\!\cong\!\!820$
$3662\!\!\sim\!\!2355\!\!\cong\!\!2355$
$3663\!\!\sim\!\!820\!\!\cong\!\!820$
$3664\!\!\sim\!\!730\!\!\cong\!\!730$
$3665\!\!\sim\!\!2271\!\!\cong\!\!2271$
$3666\!\!\sim\!\!730\!\!\cong\!\!730$
$3667\!\!\sim\!\!2271\!\!\cong\!\!2271$
$3668\!\!\sim\!\!2277\!\!\cong\!\!2277$
$3669\!\!\sim\!\!2274\!\!\cong\!\!2274$
$3670\!\!\sim\!\!730\!\!\cong\!\!730$
$3671\!\!\sim\!\!2274\!\!\cong\!\!2274$
$3672\!\!\sim\!\!730\!\!\cong\!\!730$
$3673\!\!\sim\!\!820\!\!\cong\!\!820$
$3674\!\!\sim\!\!2352\!\!\cong\!\!740$
$3675\!\!\sim\!\!820\!\!\cong\!\!820$
$3676\!\!\sim\!\!2352\!\!\cong\!\!740$
$3677\!\!\sim\!\!2358\!\!\cong\!\!820$
$3678\!\!\sim\!\!2355\!\!\cong\!\!2355$
$3679\!\!\sim\!\!820\!\!\cong\!\!820$
$3680\!\!\sim\!\!2355\!\!\cong\!\!2355$
$3681\!\!\sim\!\!820\!\!\cong\!\!820$
$3682\!\!\sim\!\!1090\!\!\cong\!\!1090$
$3683\!\!\sim\!\!2838\!\!\cong\!\!750$
$3684\!\!\sim\!\!1090\!\!\cong\!\!1090$
$3685\!\!\sim\!\!2838\!\!\cong\!\!750$
$3686\!\!\sim\!\!2844\!\!\cong\!\!730$
$3687\!\!\sim\!\!2841\!\!\cong\!\!2841$
$3688\!\!\sim\!\!1090\!\!\cong\!\!1090$
$3689\!\!\sim\!\!2841\!\!\cong\!\!2841$
$3690\!\!\sim\!\!1090\!\!\cong\!\!1090$
$3691\!\!\sim\!\!820\!\!\cong\!\!820$
$3692\!\!\sim\!\!2399\!\!\cong\!\!2399$
$3693\!\!\sim\!\!820\!\!\cong\!\!820$
$3694\!\!\sim\!\!2399\!\!\cong\!\!2399$
$3695\!\!\sim\!\!2426\!\!\cong\!\!2277$
$3696\!\!\sim\!\!2372\!\!\cong\!\!2372$
$3697\!\!\sim\!\!820\!\!\cong\!\!820$
$3698\!\!\sim\!\!2372\!\!\cong\!\!2372$
$3699\!\!\sim\!\!820\!\!\cong\!\!820$
$3700\!\!\sim\!\!730\!\!\cong\!\!730$
$3701\!\!\sim\!\!2271\!\!\cong\!\!2271$
$3702\!\!\sim\!\!730\!\!\cong\!\!730$
$3703\!\!\sim\!\!2271\!\!\cong\!\!2271$
$3704\!\!\sim\!\!2277\!\!\cong\!\!2277$
$3705\!\!\sim\!\!2274\!\!\cong\!\!2274$
$3706\!\!\sim\!\!730\!\!\cong\!\!730$
$3707\!\!\sim\!\!2274\!\!\cong\!\!2274$
$3708\!\!\sim\!\!730\!\!\cong\!\!730$
$3709\!\!\sim\!\!820\!\!\cong\!\!820$
$3710\!\!\sim\!\!2399\!\!\cong\!\!2399$
$3711\!\!\sim\!\!820\!\!\cong\!\!820$
$3712\!\!\sim\!\!2399\!\!\cong\!\!2399$
$3713\!\!\sim\!\!2426\!\!\cong\!\!2277$
$3714\!\!\sim\!\!2372\!\!\cong\!\!2372$
$3715\!\!\sim\!\!820\!\!\cong\!\!820$
$3716\!\!\sim\!\!2372\!\!\cong\!\!2372$
$3717\!\!\sim\!\!820\!\!\cong\!\!820$
$3718\!\!\sim\!\!730\!\!\cong\!\!730$
$3719\!\!\sim\!\!2237\!\!\cong\!\!2237$
$3720\!\!\sim\!\!730\!\!\cong\!\!730$
$3721\!\!\sim\!\!2237\!\!\cong\!\!2237$
$3722\!\!\sim\!\!2264\!\!\cong\!\!730$
$3723\!\!\sim\!\!2210\!\!\cong\!\!2210$
$3724\!\!\sim\!\!730\!\!\cong\!\!730$
$3725\!\!\sim\!\!2210\!\!\cong\!\!2210$
$3726\!\!\sim\!\!730\!\!\cong\!\!730$
$3727\!\!\sim\!\!2206\!\!\cong\!\!748$
$3728\!\!\sim\!\!731\!\!\cong\!\!731$
$3729\!\!\sim\!\!2207\!\!\cong\!\!2207$
$3730\!\!\sim\!\!2212\!\!\cong\!\!2212$
$3731\!\!\sim\!\!2214\!\!\cong\!\!748$
$3732\!\!\sim\!\!2213\!\!\cong\!\!2213$
$3733\!\!\sim\!\!2209\!\!\cong\!\!2209$
$3734\!\!\sim\!\!731\!\!\cong\!\!731$
$3735\!\!\sim\!\!2210\!\!\cong\!\!2210$
$3736\!\!\sim\!\!2368\!\!\cong\!\!739$
$3737\!\!\sim\!\!821\!\!\cong\!\!821$
$3738\!\!\sim\!\!2369\!\!\cong\!\!2369$
$3739\!\!\sim\!\!2374\!\!\cong\!\!821$
$3740\!\!\sim\!\!2376\!\!\cong\!\!739$
$3741\!\!\sim\!\!2375\!\!\cong\!\!2375$
$3742\!\!\sim\!\!2371\!\!\cong\!\!2371$
$3743\!\!\sim\!\!821\!\!\cong\!\!821$
$3744\!\!\sim\!\!2372\!\!\cong\!\!2372$
$3745\!\!\sim\!\!2287\!\!\cong\!\!2287$
$3746\!\!\sim\!\!731\!\!\cong\!\!731$
$3747\!\!\sim\!\!2285\!\!\cong\!\!2285$
$3748\!\!\sim\!\!2293\!\!\cong\!\!2293$
$3749\!\!\sim\!\!2295\!\!\cong\!\!2295$
$3750\!\!\sim\!\!2294\!\!\cong\!\!2294$
$3751\!\!\sim\!\!2283\!\!\cong\!\!2283$
$3752\!\!\sim\!\!731\!\!\cong\!\!731$
$3753\!\!\sim\!\!2274\!\!\cong\!\!2274$
$3754\!\!\sim\!\!2368\!\!\cong\!\!739$
$3755\!\!\sim\!\!821\!\!\cong\!\!821$
$3756\!\!\sim\!\!2369\!\!\cong\!\!2369$
$3757\!\!\sim\!\!2374\!\!\cong\!\!821$
$3758\!\!\sim\!\!2376\!\!\cong\!\!739$
$3759\!\!\sim\!\!2375\!\!\cong\!\!2375$
$3760\!\!\sim\!\!2371\!\!\cong\!\!2371$
$3761\!\!\sim\!\!821\!\!\cong\!\!821$
$3762\!\!\sim\!\!2372\!\!\cong\!\!2372$
$3763\!\!\sim\!\!2854\!\!\cong\!\!847$
$3764\!\!\sim\!\!1091\!\!\cong\!\!731$
$3765\!\!\sim\!\!2852\!\!\cong\!\!849$
$3766\!\!\sim\!\!2860\!\!\cong\!\!2212$
$3767\!\!\sim\!\!2862\!\!\cong\!\!847$
$3768\!\!\sim\!\!2861\!\!\cong\!\!731$
$3769\!\!\sim\!\!2850\!\!\cong\!\!2850$
$3770\!\!\sim\!\!1091\!\!\cong\!\!731$
$3771\!\!\sim\!\!2841\!\!\cong\!\!2841$
$3772\!\!\sim\!\!2391\!\!\cong\!\!2391$
$3773\!\!\sim\!\!821\!\!\cong\!\!821$
$3774\!\!\sim\!\!2366\!\!\cong\!\!2366$
$3775\!\!\sim\!\!2402\!\!\cong\!\!2402$
$3776\!\!\sim\!\!2427\!\!\cong\!\!2427$
$3777\!\!\sim\!\!2375\!\!\cong\!\!2375$
$3778\!\!\sim\!\!2364\!\!\cong\!\!2364$
$3779\!\!\sim\!\!821\!\!\cong\!\!821$
$3780\!\!\sim\!\!2355\!\!\cong\!\!2355$
$3781\!\!\sim\!\!2287\!\!\cong\!\!2287$
$3782\!\!\sim\!\!731\!\!\cong\!\!731$
$3783\!\!\sim\!\!2285\!\!\cong\!\!2285$
$3784\!\!\sim\!\!2293\!\!\cong\!\!2293$
$3785\!\!\sim\!\!2295\!\!\cong\!\!2295$
$3786\!\!\sim\!\!2294\!\!\cong\!\!2294$
$3787\!\!\sim\!\!2283\!\!\cong\!\!2283$
$3788\!\!\sim\!\!731\!\!\cong\!\!731$
$3789\!\!\sim\!\!2274\!\!\cong\!\!2274$
$3790\!\!\sim\!\!2391\!\!\cong\!\!2391$
$3791\!\!\sim\!\!821\!\!\cong\!\!821$
$3792\!\!\sim\!\!2366\!\!\cong\!\!2366$
$3793\!\!\sim\!\!2402\!\!\cong\!\!2402$
$3794\!\!\sim\!\!2427\!\!\cong\!\!2427$
$3795\!\!\sim\!\!2375\!\!\cong\!\!2375$
$3796\!\!\sim\!\!2364\!\!\cong\!\!2364$
$3797\!\!\sim\!\!821\!\!\cong\!\!821$
$3798\!\!\sim\!\!2355\!\!\cong\!\!2355$
$3799\!\!\sim\!\!2229\!\!\cong\!\!2229$
$3800\!\!\sim\!\!731\!\!\cong\!\!731$
$3801\!\!\sim\!\!2204\!\!\cong\!\!2204$
$3802\!\!\sim\!\!2240\!\!\cong\!\!2240$
$3803\!\!\sim\!\!2265\!\!\cong\!\!2265$
$3804\!\!\sim\!\!2213\!\!\cong\!\!2213$
$3805\!\!\sim\!\!2202\!\!\cong\!\!2202$
$3806\!\!\sim\!\!731\!\!\cong\!\!731$
$3807\!\!\sim\!\!2193\!\!\cong\!\!2193$
$3808\!\!\sim\!\!730\!\!\cong\!\!730$
$3809\!\!\sim\!\!2199\!\!\cong\!\!2199$
$3810\!\!\sim\!\!730\!\!\cong\!\!730$
$3811\!\!\sim\!\!2203\!\!\cong\!\!2203$
$3812\!\!\sim\!\!2205\!\!\cong\!\!775$
$3813\!\!\sim\!\!2204\!\!\cong\!\!2204$
$3814\!\!\sim\!\!730\!\!\cong\!\!730$
$3815\!\!\sim\!\!2202\!\!\cong\!\!2202$
$3816\!\!\sim\!\!730\!\!\cong\!\!730$
$3817\!\!\sim\!\!820\!\!\cong\!\!820$
$3818\!\!\sim\!\!2361\!\!\cong\!\!2361$
$3819\!\!\sim\!\!820\!\!\cong\!\!820$
$3820\!\!\sim\!\!2365\!\!\cong\!\!2365$
$3821\!\!\sim\!\!2367\!\!\cong\!\!2367$
$3822\!\!\sim\!\!2366\!\!\cong\!\!2366$
$3823\!\!\sim\!\!820\!\!\cong\!\!820$
$3824\!\!\sim\!\!2364\!\!\cong\!\!2364$
$3825\!\!\sim\!\!820\!\!\cong\!\!820$
$3826\!\!\sim\!\!730\!\!\cong\!\!730$
$3827\!\!\sim\!\!2280\!\!\cong\!\!2280$
$3828\!\!\sim\!\!730\!\!\cong\!\!730$
$3829\!\!\sim\!\!2284\!\!\cong\!\!2284$
$3830\!\!\sim\!\!2286\!\!\cong\!\!2286$
$3831\!\!\sim\!\!2285\!\!\cong\!\!2285$
$3832\!\!\sim\!\!730\!\!\cong\!\!730$
$3833\!\!\sim\!\!2283\!\!\cong\!\!2283$
$3834\!\!\sim\!\!730\!\!\cong\!\!730$
$3835\!\!\sim\!\!820\!\!\cong\!\!820$
$3836\!\!\sim\!\!2361\!\!\cong\!\!2361$
$3837\!\!\sim\!\!820\!\!\cong\!\!820$
$3838\!\!\sim\!\!2365\!\!\cong\!\!2365$
$3839\!\!\sim\!\!2367\!\!\cong\!\!2367$
$3840\!\!\sim\!\!2366\!\!\cong\!\!2366$
$3841\!\!\sim\!\!820\!\!\cong\!\!820$
$3842\!\!\sim\!\!2364\!\!\cong\!\!2364$
$3843\!\!\sim\!\!820\!\!\cong\!\!820$
$3844\!\!\sim\!\!1090\!\!\cong\!\!1090$
$3845\!\!\sim\!\!2847\!\!\cong\!\!929$
$3846\!\!\sim\!\!1090\!\!\cong\!\!1090$
$3847\!\!\sim\!\!2851\!\!\cong\!\!929$
$3848\!\!\sim\!\!2853\!\!\cong\!\!2853$
$3849\!\!\sim\!\!2852\!\!\cong\!\!849$
$3850\!\!\sim\!\!1090\!\!\cong\!\!1090$
$3851\!\!\sim\!\!2850\!\!\cong\!\!2850$
$3852\!\!\sim\!\!1090\!\!\cong\!\!1090$
$3853\!\!\sim\!\!820\!\!\cong\!\!820$
$3854\!\!\sim\!\!2398\!\!\cong\!\!2398$
$3855\!\!\sim\!\!820\!\!\cong\!\!820$
$3856\!\!\sim\!\!2396\!\!\cong\!\!2396$
$3857\!\!\sim\!\!2423\!\!\cong\!\!2423$
$3858\!\!\sim\!\!2369\!\!\cong\!\!2369$
$3859\!\!\sim\!\!820\!\!\cong\!\!820$
$3860\!\!\sim\!\!2371\!\!\cong\!\!2371$
$3861\!\!\sim\!\!820\!\!\cong\!\!820$
$3862\!\!\sim\!\!730\!\!\cong\!\!730$
$3863\!\!\sim\!\!2280\!\!\cong\!\!2280$
$3864\!\!\sim\!\!730\!\!\cong\!\!730$
$3865\!\!\sim\!\!2284\!\!\cong\!\!2284$
$3866\!\!\sim\!\!2286\!\!\cong\!\!2286$
$3867\!\!\sim\!\!2285\!\!\cong\!\!2285$
$3868\!\!\sim\!\!730\!\!\cong\!\!730$
$3869\!\!\sim\!\!2283\!\!\cong\!\!2283$
$3870\!\!\sim\!\!730\!\!\cong\!\!730$
$3871\!\!\sim\!\!820\!\!\cong\!\!820$
$3872\!\!\sim\!\!2398\!\!\cong\!\!2398$
$3873\!\!\sim\!\!820\!\!\cong\!\!820$
$3874\!\!\sim\!\!2396\!\!\cong\!\!2396$
$3875\!\!\sim\!\!2423\!\!\cong\!\!2423$
$3876\!\!\sim\!\!2369\!\!\cong\!\!2369$
$3877\!\!\sim\!\!820\!\!\cong\!\!820$
$3878\!\!\sim\!\!2371\!\!\cong\!\!2371$
$3879\!\!\sim\!\!820\!\!\cong\!\!820$
$3880\!\!\sim\!\!730\!\!\cong\!\!730$
$3881\!\!\sim\!\!2236\!\!\cong\!\!2236$
$3882\!\!\sim\!\!730\!\!\cong\!\!730$
$3883\!\!\sim\!\!2234\!\!\cong\!\!2234$
$3884\!\!\sim\!\!2261\!\!\cong\!\!2261$
$3885\!\!\sim\!\!2207\!\!\cong\!\!2207$
$3886\!\!\sim\!\!730\!\!\cong\!\!730$
$3887\!\!\sim\!\!2209\!\!\cong\!\!2209$
$3888\!\!\sim\!\!730\!\!\cong\!\!730$
$3889\!\!\sim\!\!2206\!\!\cong\!\!748$
$3890\!\!\sim\!\!2212\!\!\cong\!\!2212$
$3891\!\!\sim\!\!2209\!\!\cong\!\!2209$
$3892\!\!\sim\!\!731\!\!\cong\!\!731$
$3893\!\!\sim\!\!2214\!\!\cong\!\!748$
$3894\!\!\sim\!\!731\!\!\cong\!\!731$
$3895\!\!\sim\!\!2207\!\!\cong\!\!2207$
$3896\!\!\sim\!\!2213\!\!\cong\!\!2213$
$3897\!\!\sim\!\!2210\!\!\cong\!\!2210$
$3898\!\!\sim\!\!2368\!\!\cong\!\!739$
$3899\!\!\sim\!\!2374\!\!\cong\!\!821$
$3900\!\!\sim\!\!2371\!\!\cong\!\!2371$
$3901\!\!\sim\!\!821\!\!\cong\!\!821$
$3902\!\!\sim\!\!2376\!\!\cong\!\!739$
$3903\!\!\sim\!\!821\!\!\cong\!\!821$
$3904\!\!\sim\!\!2369\!\!\cong\!\!2369$
$3905\!\!\sim\!\!2375\!\!\cong\!\!2375$
$3906\!\!\sim\!\!2372\!\!\cong\!\!2372$
$3907\!\!\sim\!\!2287\!\!\cong\!\!2287$
$3908\!\!\sim\!\!2293\!\!\cong\!\!2293$
$3909\!\!\sim\!\!2283\!\!\cong\!\!2283$
$3910\!\!\sim\!\!731\!\!\cong\!\!731$
$3911\!\!\sim\!\!2295\!\!\cong\!\!2295$
$3912\!\!\sim\!\!731\!\!\cong\!\!731$
$3913\!\!\sim\!\!2285\!\!\cong\!\!2285$
$3914\!\!\sim\!\!2294\!\!\cong\!\!2294$
$3915\!\!\sim\!\!2274\!\!\cong\!\!2274$
$3916\!\!\sim\!\!2368\!\!\cong\!\!739$
$3917\!\!\sim\!\!2374\!\!\cong\!\!821$
$3918\!\!\sim\!\!2371\!\!\cong\!\!2371$
$3919\!\!\sim\!\!821\!\!\cong\!\!821$
$3920\!\!\sim\!\!2376\!\!\cong\!\!739$
$3921\!\!\sim\!\!821\!\!\cong\!\!821$
$3922\!\!\sim\!\!2369\!\!\cong\!\!2369$
$3923\!\!\sim\!\!2375\!\!\cong\!\!2375$
$3924\!\!\sim\!\!2372\!\!\cong\!\!2372$
$3925\!\!\sim\!\!2854\!\!\cong\!\!847$
$3926\!\!\sim\!\!2860\!\!\cong\!\!2212$
$3927\!\!\sim\!\!2850\!\!\cong\!\!2850$
$3928\!\!\sim\!\!1091\!\!\cong\!\!731$
$3929\!\!\sim\!\!2862\!\!\cong\!\!847$
$3930\!\!\sim\!\!1091\!\!\cong\!\!731$
$3931\!\!\sim\!\!2852\!\!\cong\!\!849$
$3932\!\!\sim\!\!2861\!\!\cong\!\!731$
$3933\!\!\sim\!\!2841\!\!\cong\!\!2841$
$3934\!\!\sim\!\!2391\!\!\cong\!\!2391$
$3935\!\!\sim\!\!2402\!\!\cong\!\!2402$
$3936\!\!\sim\!\!2364\!\!\cong\!\!2364$
$3937\!\!\sim\!\!821\!\!\cong\!\!821$
$3938\!\!\sim\!\!2427\!\!\cong\!\!2427$
$3939\!\!\sim\!\!821\!\!\cong\!\!821$
$3940\!\!\sim\!\!2366\!\!\cong\!\!2366$
$3941\!\!\sim\!\!2375\!\!\cong\!\!2375$
$3942\!\!\sim\!\!2355\!\!\cong\!\!2355$
$3943\!\!\sim\!\!2287\!\!\cong\!\!2287$
$3944\!\!\sim\!\!2293\!\!\cong\!\!2293$
$3945\!\!\sim\!\!2283\!\!\cong\!\!2283$
$3946\!\!\sim\!\!731\!\!\cong\!\!731$
$3947\!\!\sim\!\!2295\!\!\cong\!\!2295$
$3948\!\!\sim\!\!731\!\!\cong\!\!731$
$3949\!\!\sim\!\!2285\!\!\cong\!\!2285$
$3950\!\!\sim\!\!2294\!\!\cong\!\!2294$
$3951\!\!\sim\!\!2274\!\!\cong\!\!2274$
$3952\!\!\sim\!\!2391\!\!\cong\!\!2391$
$3953\!\!\sim\!\!2402\!\!\cong\!\!2402$
$3954\!\!\sim\!\!2364\!\!\cong\!\!2364$
$3955\!\!\sim\!\!821\!\!\cong\!\!821$
$3956\!\!\sim\!\!2427\!\!\cong\!\!2427$
$3957\!\!\sim\!\!821\!\!\cong\!\!821$
$3958\!\!\sim\!\!2366\!\!\cong\!\!2366$
$3959\!\!\sim\!\!2375\!\!\cong\!\!2375$
$3960\!\!\sim\!\!2355\!\!\cong\!\!2355$
$3961\!\!\sim\!\!2229\!\!\cong\!\!2229$
$3962\!\!\sim\!\!2240\!\!\cong\!\!2240$
$3963\!\!\sim\!\!2202\!\!\cong\!\!2202$
$3964\!\!\sim\!\!731\!\!\cong\!\!731$
$3965\!\!\sim\!\!2265\!\!\cong\!\!2265$
$3966\!\!\sim\!\!731\!\!\cong\!\!731$
$3967\!\!\sim\!\!2204\!\!\cong\!\!2204$
$3968\!\!\sim\!\!2213\!\!\cong\!\!2213$
$3969\!\!\sim\!\!2193\!\!\cong\!\!2193$
$3970\!\!\sim\!\!2260\!\!\cong\!\!802$
$3971\!\!\sim\!\!2262\!\!\cong\!\!750$
$3972\!\!\sim\!\!2261\!\!\cong\!\!2261$
$3973\!\!\sim\!\!2262\!\!\cong\!\!750$
$3974\!\!\sim\!\!734\!\!\cong\!\!730$
$3975\!\!\sim\!\!2265\!\!\cong\!\!2265$
$3976\!\!\sim\!\!2261\!\!\cong\!\!2261$
$3977\!\!\sim\!\!2265\!\!\cong\!\!2265$
$3978\!\!\sim\!\!2264\!\!\cong\!\!730$
$3979\!\!\sim\!\!2422\!\!\cong\!\!820$
$3980\!\!\sim\!\!2424\!\!\cong\!\!966$
$3981\!\!\sim\!\!2423\!\!\cong\!\!2423$
$3982\!\!\sim\!\!2424\!\!\cong\!\!966$
$3983\!\!\sim\!\!824\!\!\cong\!\!820$
$3984\!\!\sim\!\!2427\!\!\cong\!\!2427$
$3985\!\!\sim\!\!2423\!\!\cong\!\!2423$
$3986\!\!\sim\!\!2427\!\!\cong\!\!2427$
$3987\!\!\sim\!\!2426\!\!\cong\!\!2277$
$3988\!\!\sim\!\!2313\!\!\cong\!\!2277$
$3989\!\!\sim\!\!2322\!\!\cong\!\!2322$
$3990\!\!\sim\!\!2286\!\!\cong\!\!2286$
$3991\!\!\sim\!\!2322\!\!\cong\!\!2322$
$3992\!\!\sim\!\!734\!\!\cong\!\!730$
$3993\!\!\sim\!\!2295\!\!\cong\!\!2295$
$3994\!\!\sim\!\!2286\!\!\cong\!\!2286$
$3995\!\!\sim\!\!2295\!\!\cong\!\!2295$
$3996\!\!\sim\!\!2277\!\!\cong\!\!2277$
$3997\!\!\sim\!\!2422\!\!\cong\!\!820$
$3998\!\!\sim\!\!2424\!\!\cong\!\!966$
$3999\!\!\sim\!\!2423\!\!\cong\!\!2423$
$4000\!\!\sim\!\!2424\!\!\cong\!\!966$
$4001\!\!\sim\!\!824\!\!\cong\!\!820$
$4002\!\!\sim\!\!2427\!\!\cong\!\!2427$
$4003\!\!\sim\!\!2423\!\!\cong\!\!2423$
$4004\!\!\sim\!\!2427\!\!\cong\!\!2427$
$4005\!\!\sim\!\!2426\!\!\cong\!\!2277$
$4006\!\!\sim\!\!2880\!\!\cong\!\!730$
$4007\!\!\sim\!\!2889\!\!\cong\!\!750$
$4008\!\!\sim\!\!2853\!\!\cong\!\!2853$
$4009\!\!\sim\!\!2889\!\!\cong\!\!750$
$4010\!\!\sim\!\!1094\!\!\cong\!\!1090$
$4011\!\!\sim\!\!2862\!\!\cong\!\!847$
$4012\!\!\sim\!\!2853\!\!\cong\!\!2853$
$4013\!\!\sim\!\!2862\!\!\cong\!\!847$
$4014\!\!\sim\!\!2844\!\!\cong\!\!730$
$4015\!\!\sim\!\!2394\!\!\cong\!\!820$
$4016\!\!\sim\!\!2403\!\!\cong\!\!2287$
$4017\!\!\sim\!\!2367\!\!\cong\!\!2367$
$4018\!\!\sim\!\!2403\!\!\cong\!\!2287$
$4019\!\!\sim\!\!824\!\!\cong\!\!820$
$4020\!\!\sim\!\!2376\!\!\cong\!\!739$
$4021\!\!\sim\!\!2367\!\!\cong\!\!2367$
$4022\!\!\sim\!\!2376\!\!\cong\!\!739$
$4023\!\!\sim\!\!2358\!\!\cong\!\!820$
$4024\!\!\sim\!\!2313\!\!\cong\!\!2277$
$4025\!\!\sim\!\!2322\!\!\cong\!\!2322$
$4026\!\!\sim\!\!2286\!\!\cong\!\!2286$
$4027\!\!\sim\!\!2322\!\!\cong\!\!2322$
$4028\!\!\sim\!\!734\!\!\cong\!\!730$
$4029\!\!\sim\!\!2295\!\!\cong\!\!2295$
$4030\!\!\sim\!\!2286\!\!\cong\!\!2286$
$4031\!\!\sim\!\!2295\!\!\cong\!\!2295$
$4032\!\!\sim\!\!2277\!\!\cong\!\!2277$
$4033\!\!\sim\!\!2394\!\!\cong\!\!820$
$4034\!\!\sim\!\!2403\!\!\cong\!\!2287$
$4035\!\!\sim\!\!2367\!\!\cong\!\!2367$
$4036\!\!\sim\!\!2403\!\!\cong\!\!2287$
$4037\!\!\sim\!\!824\!\!\cong\!\!820$
$4038\!\!\sim\!\!2376\!\!\cong\!\!739$
$4039\!\!\sim\!\!2367\!\!\cong\!\!2367$
$4040\!\!\sim\!\!2376\!\!\cong\!\!739$
$4041\!\!\sim\!\!2358\!\!\cong\!\!820$
$4042\!\!\sim\!\!2232\!\!\cong\!\!730$
$4043\!\!\sim\!\!2241\!\!\cong\!\!739$
$4044\!\!\sim\!\!2205\!\!\cong\!\!775$
$4045\!\!\sim\!\!2241\!\!\cong\!\!739$
$4046\!\!\sim\!\!734\!\!\cong\!\!730$
$4047\!\!\sim\!\!2214\!\!\cong\!\!748$
$4048\!\!\sim\!\!2205\!\!\cong\!\!775$
$4049\!\!\sim\!\!2214\!\!\cong\!\!748$
$4050\!\!\sim\!\!2196\!\!\cong\!\!802$
$4051\!\!\sim\!\!2233\!\!\cong\!\!2233$
$4052\!\!\sim\!\!2239\!\!\cong\!\!2239$
$4053\!\!\sim\!\!2236\!\!\cong\!\!2236$
$4054\!\!\sim\!\!731\!\!\cong\!\!731$
$4055\!\!\sim\!\!2241\!\!\cong\!\!739$
$4056\!\!\sim\!\!731\!\!\cong\!\!731$
$4057\!\!\sim\!\!2234\!\!\cong\!\!2234$
$4058\!\!\sim\!\!2240\!\!\cong\!\!2240$
$4059\!\!\sim\!\!2237\!\!\cong\!\!2237$
$4060\!\!\sim\!\!2395\!\!\cong\!\!2395$
$4061\!\!\sim\!\!2401\!\!\cong\!\!2401$
$4062\!\!\sim\!\!2398\!\!\cong\!\!2398$
$4063\!\!\sim\!\!821\!\!\cong\!\!821$
$4064\!\!\sim\!\!2403\!\!\cong\!\!2287$
$4065\!\!\sim\!\!821\!\!\cong\!\!821$
$4066\!\!\sim\!\!2396\!\!\cong\!\!2396$
$4067\!\!\sim\!\!2402\!\!\cong\!\!2402$
$4068\!\!\sim\!\!2399\!\!\cong\!\!2399$
$4069\!\!\sim\!\!2307\!\!\cong\!\!2307$
$4070\!\!\sim\!\!2320\!\!\cong\!\!2294$
$4071\!\!\sim\!\!2280\!\!\cong\!\!2280$
$4072\!\!\sim\!\!731\!\!\cong\!\!731$
$4073\!\!\sim\!\!2322\!\!\cong\!\!2322$
$4074\!\!\sim\!\!731\!\!\cong\!\!731$
$4075\!\!\sim\!\!2284\!\!\cong\!\!2284$
$4076\!\!\sim\!\!2293\!\!\cong\!\!2293$
$4077\!\!\sim\!\!2271\!\!\cong\!\!2271$
$4078\!\!\sim\!\!2395\!\!\cong\!\!2395$
$4079\!\!\sim\!\!2401\!\!\cong\!\!2401$
$4080\!\!\sim\!\!2398\!\!\cong\!\!2398$
$4081\!\!\sim\!\!821\!\!\cong\!\!821$
$4082\!\!\sim\!\!2403\!\!\cong\!\!2287$
$4083\!\!\sim\!\!821\!\!\cong\!\!821$
$4084\!\!\sim\!\!2396\!\!\cong\!\!2396$
$4085\!\!\sim\!\!2402\!\!\cong\!\!2402$
$4086\!\!\sim\!\!2399\!\!\cong\!\!2399$
$4087\!\!\sim\!\!2874\!\!\cong\!\!820$
$4088\!\!\sim\!\!2887\!\!\cong\!\!731$
$4089\!\!\sim\!\!2847\!\!\cong\!\!929$
$4090\!\!\sim\!\!1091\!\!\cong\!\!731$
$4091\!\!\sim\!\!2889\!\!\cong\!\!750$
$4092\!\!\sim\!\!1091\!\!\cong\!\!731$
$4093\!\!\sim\!\!2851\!\!\cong\!\!929$
$4094\!\!\sim\!\!2860\!\!\cong\!\!2212$
$4095\!\!\sim\!\!2838\!\!\cong\!\!750$
$4096\!\!\sim\!\!2388\!\!\cong\!\!821$
$4097\!\!\sim\!\!2401\!\!\cong\!\!2401$
$4098\!\!\sim\!\!2361\!\!\cong\!\!2361$
$4099\!\!\sim\!\!821\!\!\cong\!\!821$
$4100\!\!\sim\!\!2424\!\!\cong\!\!966$
$4101\!\!\sim\!\!821\!\!\cong\!\!821$
$4102\!\!\sim\!\!2365\!\!\cong\!\!2365$
$4103\!\!\sim\!\!2374\!\!\cong\!\!821$
$4104\!\!\sim\!\!2352\!\!\cong\!\!740$
$4105\!\!\sim\!\!2307\!\!\cong\!\!2307$
$4106\!\!\sim\!\!2320\!\!\cong\!\!2294$
$4107\!\!\sim\!\!2280\!\!\cong\!\!2280$
$4108\!\!\sim\!\!731\!\!\cong\!\!731$
$4109\!\!\sim\!\!2322\!\!\cong\!\!2322$
$4110\!\!\sim\!\!731\!\!\cong\!\!731$
$4111\!\!\sim\!\!2284\!\!\cong\!\!2284$
$4112\!\!\sim\!\!2293\!\!\cong\!\!2293$
$4113\!\!\sim\!\!2271\!\!\cong\!\!2271$
$4114\!\!\sim\!\!2388\!\!\cong\!\!821$
$4115\!\!\sim\!\!2401\!\!\cong\!\!2401$
$4116\!\!\sim\!\!2361\!\!\cong\!\!2361$
$4117\!\!\sim\!\!821\!\!\cong\!\!821$
$4118\!\!\sim\!\!2424\!\!\cong\!\!966$
$4119\!\!\sim\!\!821\!\!\cong\!\!821$
$4120\!\!\sim\!\!2365\!\!\cong\!\!2365$
$4121\!\!\sim\!\!2374\!\!\cong\!\!821$
$4122\!\!\sim\!\!2352\!\!\cong\!\!740$
$4123\!\!\sim\!\!2226\!\!\cong\!\!820$
$4124\!\!\sim\!\!2239\!\!\cong\!\!2239$
$4125\!\!\sim\!\!2199\!\!\cong\!\!2199$
$4126\!\!\sim\!\!731\!\!\cong\!\!731$
$4127\!\!\sim\!\!2262\!\!\cong\!\!750$
$4128\!\!\sim\!\!731\!\!\cong\!\!731$
$4129\!\!\sim\!\!2203\!\!\cong\!\!2203$
$4130\!\!\sim\!\!2212\!\!\cong\!\!2212$
$4131\!\!\sim\!\!2190\!\!\cong\!\!750$
$4132\!\!\sim\!\!730\!\!\cong\!\!730$
$4133\!\!\sim\!\!2203\!\!\cong\!\!2203$
$4134\!\!\sim\!\!730\!\!\cong\!\!730$
$4135\!\!\sim\!\!2199\!\!\cong\!\!2199$
$4136\!\!\sim\!\!2205\!\!\cong\!\!775$
$4137\!\!\sim\!\!2202\!\!\cong\!\!2202$
$4138\!\!\sim\!\!730\!\!\cong\!\!730$
$4139\!\!\sim\!\!2204\!\!\cong\!\!2204$
$4140\!\!\sim\!\!730\!\!\cong\!\!730$
$4141\!\!\sim\!\!820\!\!\cong\!\!820$
$4142\!\!\sim\!\!2365\!\!\cong\!\!2365$
$4143\!\!\sim\!\!820\!\!\cong\!\!820$
$4144\!\!\sim\!\!2361\!\!\cong\!\!2361$
$4145\!\!\sim\!\!2367\!\!\cong\!\!2367$
$4146\!\!\sim\!\!2364\!\!\cong\!\!2364$
$4147\!\!\sim\!\!820\!\!\cong\!\!820$
$4148\!\!\sim\!\!2366\!\!\cong\!\!2366$
$4149\!\!\sim\!\!820\!\!\cong\!\!820$
$4150\!\!\sim\!\!730\!\!\cong\!\!730$
$4151\!\!\sim\!\!2284\!\!\cong\!\!2284$
$4152\!\!\sim\!\!730\!\!\cong\!\!730$
$4153\!\!\sim\!\!2280\!\!\cong\!\!2280$
$4154\!\!\sim\!\!2286\!\!\cong\!\!2286$
$4155\!\!\sim\!\!2283\!\!\cong\!\!2283$
$4156\!\!\sim\!\!730\!\!\cong\!\!730$
$4157\!\!\sim\!\!2285\!\!\cong\!\!2285$
$4158\!\!\sim\!\!730\!\!\cong\!\!730$
$4159\!\!\sim\!\!820\!\!\cong\!\!820$
$4160\!\!\sim\!\!2365\!\!\cong\!\!2365$
$4161\!\!\sim\!\!820\!\!\cong\!\!820$
$4162\!\!\sim\!\!2361\!\!\cong\!\!2361$
$4163\!\!\sim\!\!2367\!\!\cong\!\!2367$
$4164\!\!\sim\!\!2364\!\!\cong\!\!2364$
$4165\!\!\sim\!\!820\!\!\cong\!\!820$
$4166\!\!\sim\!\!2366\!\!\cong\!\!2366$
$4167\!\!\sim\!\!820\!\!\cong\!\!820$
$4168\!\!\sim\!\!1090\!\!\cong\!\!1090$
$4169\!\!\sim\!\!2851\!\!\cong\!\!929$
$4170\!\!\sim\!\!1090\!\!\cong\!\!1090$
$4171\!\!\sim\!\!2847\!\!\cong\!\!929$
$4172\!\!\sim\!\!2853\!\!\cong\!\!2853$
$4173\!\!\sim\!\!2850\!\!\cong\!\!2850$
$4174\!\!\sim\!\!1090\!\!\cong\!\!1090$
$4175\!\!\sim\!\!2852\!\!\cong\!\!849$
$4176\!\!\sim\!\!1090\!\!\cong\!\!1090$
$4177\!\!\sim\!\!820\!\!\cong\!\!820$
$4178\!\!\sim\!\!2396\!\!\cong\!\!2396$
$4179\!\!\sim\!\!820\!\!\cong\!\!820$
$4180\!\!\sim\!\!2398\!\!\cong\!\!2398$
$4181\!\!\sim\!\!2423\!\!\cong\!\!2423$
$4182\!\!\sim\!\!2371\!\!\cong\!\!2371$
$4183\!\!\sim\!\!820\!\!\cong\!\!820$
$4184\!\!\sim\!\!2369\!\!\cong\!\!2369$
$4185\!\!\sim\!\!820\!\!\cong\!\!820$
$4186\!\!\sim\!\!730\!\!\cong\!\!730$
$4187\!\!\sim\!\!2284\!\!\cong\!\!2284$
$4188\!\!\sim\!\!730\!\!\cong\!\!730$
$4189\!\!\sim\!\!2280\!\!\cong\!\!2280$
$4190\!\!\sim\!\!2286\!\!\cong\!\!2286$
$4191\!\!\sim\!\!2283\!\!\cong\!\!2283$
$4192\!\!\sim\!\!730\!\!\cong\!\!730$
$4193\!\!\sim\!\!2285\!\!\cong\!\!2285$
$4194\!\!\sim\!\!730\!\!\cong\!\!730$
$4195\!\!\sim\!\!820\!\!\cong\!\!820$
$4196\!\!\sim\!\!2396\!\!\cong\!\!2396$
$4197\!\!\sim\!\!820\!\!\cong\!\!820$
$4198\!\!\sim\!\!2398\!\!\cong\!\!2398$
$4199\!\!\sim\!\!2423\!\!\cong\!\!2423$
$4200\!\!\sim\!\!2371\!\!\cong\!\!2371$
$4201\!\!\sim\!\!820\!\!\cong\!\!820$
$4202\!\!\sim\!\!2369\!\!\cong\!\!2369$
$4203\!\!\sim\!\!820\!\!\cong\!\!820$
$4204\!\!\sim\!\!730\!\!\cong\!\!730$
$4205\!\!\sim\!\!2234\!\!\cong\!\!2234$
$4206\!\!\sim\!\!730\!\!\cong\!\!730$
$4207\!\!\sim\!\!2236\!\!\cong\!\!2236$
$4208\!\!\sim\!\!2261\!\!\cong\!\!2261$
$4209\!\!\sim\!\!2209\!\!\cong\!\!2209$
$4210\!\!\sim\!\!730\!\!\cong\!\!730$
$4211\!\!\sim\!\!2207\!\!\cong\!\!2207$
$4212\!\!\sim\!\!730\!\!\cong\!\!730$
$4213\!\!\sim\!\!2233\!\!\cong\!\!2233$
$4214\!\!\sim\!\!731\!\!\cong\!\!731$
$4215\!\!\sim\!\!2234\!\!\cong\!\!2234$
$4216\!\!\sim\!\!2239\!\!\cong\!\!2239$
$4217\!\!\sim\!\!2241\!\!\cong\!\!739$
$4218\!\!\sim\!\!2240\!\!\cong\!\!2240$
$4219\!\!\sim\!\!2236\!\!\cong\!\!2236$
$4220\!\!\sim\!\!731\!\!\cong\!\!731$
$4221\!\!\sim\!\!2237\!\!\cong\!\!2237$
$4222\!\!\sim\!\!2395\!\!\cong\!\!2395$
$4223\!\!\sim\!\!821\!\!\cong\!\!821$
$4224\!\!\sim\!\!2396\!\!\cong\!\!2396$
$4225\!\!\sim\!\!2401\!\!\cong\!\!2401$
$4226\!\!\sim\!\!2403\!\!\cong\!\!2287$
$4227\!\!\sim\!\!2402\!\!\cong\!\!2402$
$4228\!\!\sim\!\!2398\!\!\cong\!\!2398$
$4229\!\!\sim\!\!821\!\!\cong\!\!821$
$4230\!\!\sim\!\!2399\!\!\cong\!\!2399$
$4231\!\!\sim\!\!2307\!\!\cong\!\!2307$
$4232\!\!\sim\!\!731\!\!\cong\!\!731$
$4233\!\!\sim\!\!2284\!\!\cong\!\!2284$
$4234\!\!\sim\!\!2320\!\!\cong\!\!2294$
$4235\!\!\sim\!\!2322\!\!\cong\!\!2322$
$4236\!\!\sim\!\!2293\!\!\cong\!\!2293$
$4237\!\!\sim\!\!2280\!\!\cong\!\!2280$
$4238\!\!\sim\!\!731\!\!\cong\!\!731$
$4239\!\!\sim\!\!2271\!\!\cong\!\!2271$
$4240\!\!\sim\!\!2395\!\!\cong\!\!2395$
$4241\!\!\sim\!\!821\!\!\cong\!\!821$
$4242\!\!\sim\!\!2396\!\!\cong\!\!2396$
$4243\!\!\sim\!\!2401\!\!\cong\!\!2401$
$4244\!\!\sim\!\!2403\!\!\cong\!\!2287$
$4245\!\!\sim\!\!2402\!\!\cong\!\!2402$
$4246\!\!\sim\!\!2398\!\!\cong\!\!2398$
$4247\!\!\sim\!\!821\!\!\cong\!\!821$
$4248\!\!\sim\!\!2399\!\!\cong\!\!2399$
$4249\!\!\sim\!\!2874\!\!\cong\!\!820$
$4250\!\!\sim\!\!1091\!\!\cong\!\!731$
$4251\!\!\sim\!\!2851\!\!\cong\!\!929$
$4252\!\!\sim\!\!2887\!\!\cong\!\!731$
$4253\!\!\sim\!\!2889\!\!\cong\!\!750$
$4254\!\!\sim\!\!2860\!\!\cong\!\!2212$
$4255\!\!\sim\!\!2847\!\!\cong\!\!929$
$4256\!\!\sim\!\!1091\!\!\cong\!\!731$
$4257\!\!\sim\!\!2838\!\!\cong\!\!750$
$4258\!\!\sim\!\!2388\!\!\cong\!\!821$
$4259\!\!\sim\!\!821\!\!\cong\!\!821$
$4260\!\!\sim\!\!2365\!\!\cong\!\!2365$
$4261\!\!\sim\!\!2401\!\!\cong\!\!2401$
$4262\!\!\sim\!\!2424\!\!\cong\!\!966$
$4263\!\!\sim\!\!2374\!\!\cong\!\!821$
$4264\!\!\sim\!\!2361\!\!\cong\!\!2361$
$4265\!\!\sim\!\!821\!\!\cong\!\!821$
$4266\!\!\sim\!\!2352\!\!\cong\!\!740$
$4267\!\!\sim\!\!2307\!\!\cong\!\!2307$
$4268\!\!\sim\!\!731\!\!\cong\!\!731$
$4269\!\!\sim\!\!2284\!\!\cong\!\!2284$
$4270\!\!\sim\!\!2320\!\!\cong\!\!2294$
$4271\!\!\sim\!\!2322\!\!\cong\!\!2322$
$4272\!\!\sim\!\!2293\!\!\cong\!\!2293$
$4273\!\!\sim\!\!2280\!\!\cong\!\!2280$
$4274\!\!\sim\!\!731\!\!\cong\!\!731$
$4275\!\!\sim\!\!2271\!\!\cong\!\!2271$
$4276\!\!\sim\!\!2388\!\!\cong\!\!821$
$4277\!\!\sim\!\!821\!\!\cong\!\!821$
$4278\!\!\sim\!\!2365\!\!\cong\!\!2365$
$4279\!\!\sim\!\!2401\!\!\cong\!\!2401$
$4280\!\!\sim\!\!2424\!\!\cong\!\!966$
$4281\!\!\sim\!\!2374\!\!\cong\!\!821$
$4282\!\!\sim\!\!2361\!\!\cong\!\!2361$
$4283\!\!\sim\!\!821\!\!\cong\!\!821$
$4284\!\!\sim\!\!2352\!\!\cong\!\!740$
$4285\!\!\sim\!\!2226\!\!\cong\!\!820$
$4286\!\!\sim\!\!731\!\!\cong\!\!731$
$4287\!\!\sim\!\!2203\!\!\cong\!\!2203$
$4288\!\!\sim\!\!2239\!\!\cong\!\!2239$
$4289\!\!\sim\!\!2262\!\!\cong\!\!750$
$4290\!\!\sim\!\!2212\!\!\cong\!\!2212$
$4291\!\!\sim\!\!2199\!\!\cong\!\!2199$
$4292\!\!\sim\!\!731\!\!\cong\!\!731$
$4293\!\!\sim\!\!2190\!\!\cong\!\!750$
$4294\!\!\sim\!\!730\!\!\cong\!\!730$
$4295\!\!\sim\!\!2226\!\!\cong\!\!820$
$4296\!\!\sim\!\!730\!\!\cong\!\!730$
$4297\!\!\sim\!\!2226\!\!\cong\!\!820$
$4298\!\!\sim\!\!2232\!\!\cong\!\!730$
$4299\!\!\sim\!\!2229\!\!\cong\!\!2229$
$4300\!\!\sim\!\!730\!\!\cong\!\!730$
$4301\!\!\sim\!\!2229\!\!\cong\!\!2229$
$4302\!\!\sim\!\!730\!\!\cong\!\!730$
$4303\!\!\sim\!\!820\!\!\cong\!\!820$
$4304\!\!\sim\!\!2388\!\!\cong\!\!821$
$4305\!\!\sim\!\!820\!\!\cong\!\!820$
$4306\!\!\sim\!\!2388\!\!\cong\!\!821$
$4307\!\!\sim\!\!2394\!\!\cong\!\!820$
$4308\!\!\sim\!\!2391\!\!\cong\!\!2391$
$4309\!\!\sim\!\!820\!\!\cong\!\!820$
$4310\!\!\sim\!\!2391\!\!\cong\!\!2391$
$4311\!\!\sim\!\!820\!\!\cong\!\!820$
$4312\!\!\sim\!\!730\!\!\cong\!\!730$
$4313\!\!\sim\!\!2307\!\!\cong\!\!2307$
$4314\!\!\sim\!\!730\!\!\cong\!\!730$
$4315\!\!\sim\!\!2307\!\!\cong\!\!2307$
$4316\!\!\sim\!\!2313\!\!\cong\!\!2277$
$4317\!\!\sim\!\!2287\!\!\cong\!\!2287$
$4318\!\!\sim\!\!730\!\!\cong\!\!730$
$4319\!\!\sim\!\!2287\!\!\cong\!\!2287$
$4320\!\!\sim\!\!730\!\!\cong\!\!730$
$4321\!\!\sim\!\!820\!\!\cong\!\!820$
$4322\!\!\sim\!\!2388\!\!\cong\!\!821$
$4323\!\!\sim\!\!820\!\!\cong\!\!820$
$4324\!\!\sim\!\!2388\!\!\cong\!\!821$
$4325\!\!\sim\!\!2394\!\!\cong\!\!820$
$4326\!\!\sim\!\!2391\!\!\cong\!\!2391$
$4327\!\!\sim\!\!820\!\!\cong\!\!820$
$4328\!\!\sim\!\!2391\!\!\cong\!\!2391$
$4329\!\!\sim\!\!820\!\!\cong\!\!820$
$4330\!\!\sim\!\!1090\!\!\cong\!\!1090$
$4331\!\!\sim\!\!2874\!\!\cong\!\!820$
$4332\!\!\sim\!\!1090\!\!\cong\!\!1090$
$4333\!\!\sim\!\!2874\!\!\cong\!\!820$
$4334\!\!\sim\!\!2880\!\!\cong\!\!730$
$4335\!\!\sim\!\!2854\!\!\cong\!\!847$
$4336\!\!\sim\!\!1090\!\!\cong\!\!1090$
$4337\!\!\sim\!\!2854\!\!\cong\!\!847$
$4338\!\!\sim\!\!1090\!\!\cong\!\!1090$
$4339\!\!\sim\!\!820\!\!\cong\!\!820$
$4340\!\!\sim\!\!2395\!\!\cong\!\!2395$
$4341\!\!\sim\!\!820\!\!\cong\!\!820$
$4342\!\!\sim\!\!2395\!\!\cong\!\!2395$
$4343\!\!\sim\!\!2422\!\!\cong\!\!820$
$4344\!\!\sim\!\!2368\!\!\cong\!\!739$
$4345\!\!\sim\!\!820\!\!\cong\!\!820$
$4346\!\!\sim\!\!2368\!\!\cong\!\!739$
$4347\!\!\sim\!\!820\!\!\cong\!\!820$
$4348\!\!\sim\!\!730\!\!\cong\!\!730$
$4349\!\!\sim\!\!2307\!\!\cong\!\!2307$
$4350\!\!\sim\!\!730\!\!\cong\!\!730$
$4351\!\!\sim\!\!2307\!\!\cong\!\!2307$
$4352\!\!\sim\!\!2313\!\!\cong\!\!2277$
$4353\!\!\sim\!\!2287\!\!\cong\!\!2287$
$4354\!\!\sim\!\!730\!\!\cong\!\!730$
$4355\!\!\sim\!\!2287\!\!\cong\!\!2287$
$4356\!\!\sim\!\!730\!\!\cong\!\!730$
$4357\!\!\sim\!\!820\!\!\cong\!\!820$
$4358\!\!\sim\!\!2395\!\!\cong\!\!2395$
$4359\!\!\sim\!\!820\!\!\cong\!\!820$
$4360\!\!\sim\!\!2395\!\!\cong\!\!2395$
$4361\!\!\sim\!\!2422\!\!\cong\!\!820$
$4362\!\!\sim\!\!2368\!\!\cong\!\!739$
$4363\!\!\sim\!\!820\!\!\cong\!\!820$
$4364\!\!\sim\!\!2368\!\!\cong\!\!739$
$4365\!\!\sim\!\!820\!\!\cong\!\!820$
$4366\!\!\sim\!\!730\!\!\cong\!\!730$
$4367\!\!\sim\!\!2233\!\!\cong\!\!2233$
$4368\!\!\sim\!\!730\!\!\cong\!\!730$
$4369\!\!\sim\!\!2233\!\!\cong\!\!2233$
$4370\!\!\sim\!\!2260\!\!\cong\!\!802$
$4371\!\!\sim\!\!2206\!\!\cong\!\!748$
$4372\!\!\sim\!\!730\!\!\cong\!\!730$
$4373\!\!\sim\!\!2206\!\!\cong\!\!748$
$4374\!\!\sim\!\!730\!\!\cong\!\!730$
$4375\!\!\sim\!\!1094\!\!\cong\!\!1090$
$4376\!\!\sim\!\!824\!\!\cong\!\!820$
$4377\!\!\sim\!\!824\!\!\cong\!\!820$
$4378\!\!\sim\!\!824\!\!\cong\!\!820$
$4379\!\!\sim\!\!734\!\!\cong\!\!730$
$4380\!\!\sim\!\!734\!\!\cong\!\!730$
$4381\!\!\sim\!\!824\!\!\cong\!\!820$
$4382\!\!\sim\!\!734\!\!\cong\!\!730$
$4383\!\!\sim\!\!734\!\!\cong\!\!730$
$4384\!\!\sim\!\!2889\!\!\cong\!\!750$
$4385\!\!\sim\!\!2424\!\!\cong\!\!966$
$4386\!\!\sim\!\!2403\!\!\cong\!\!2287$
$4387\!\!\sim\!\!2424\!\!\cong\!\!966$
$4388\!\!\sim\!\!2262\!\!\cong\!\!750$
$4389\!\!\sim\!\!2322\!\!\cong\!\!2322$
$4390\!\!\sim\!\!2403\!\!\cong\!\!2287$
$4391\!\!\sim\!\!2322\!\!\cong\!\!2322$
$4392\!\!\sim\!\!2241\!\!\cong\!\!739$
$4393\!\!\sim\!\!2862\!\!\cong\!\!847$
$4394\!\!\sim\!\!2427\!\!\cong\!\!2427$
$4395\!\!\sim\!\!2376\!\!\cong\!\!739$
$4396\!\!\sim\!\!2427\!\!\cong\!\!2427$
$4397\!\!\sim\!\!2265\!\!\cong\!\!2265$
$4398\!\!\sim\!\!2295\!\!\cong\!\!2295$
$4399\!\!\sim\!\!2376\!\!\cong\!\!739$
$4400\!\!\sim\!\!2295\!\!\cong\!\!2295$
$4401\!\!\sim\!\!2214\!\!\cong\!\!748$
$4402\!\!\sim\!\!2889\!\!\cong\!\!750$
$4403\!\!\sim\!\!2424\!\!\cong\!\!966$
$4404\!\!\sim\!\!2403\!\!\cong\!\!2287$
$4405\!\!\sim\!\!2424\!\!\cong\!\!966$
$4406\!\!\sim\!\!2262\!\!\cong\!\!750$
$4407\!\!\sim\!\!2322\!\!\cong\!\!2322$
$4408\!\!\sim\!\!2403\!\!\cong\!\!2287$
$4409\!\!\sim\!\!2322\!\!\cong\!\!2322$
$4410\!\!\sim\!\!2241\!\!\cong\!\!739$
$4411\!\!\sim\!\!2880\!\!\cong\!\!730$
$4412\!\!\sim\!\!2422\!\!\cong\!\!820$
$4413\!\!\sim\!\!2394\!\!\cong\!\!820$
$4414\!\!\sim\!\!2422\!\!\cong\!\!820$
$4415\!\!\sim\!\!2260\!\!\cong\!\!802$
$4416\!\!\sim\!\!2313\!\!\cong\!\!2277$
$4417\!\!\sim\!\!2394\!\!\cong\!\!820$
$4418\!\!\sim\!\!2313\!\!\cong\!\!2277$
$4419\!\!\sim\!\!2232\!\!\cong\!\!730$
$4420\!\!\sim\!\!2853\!\!\cong\!\!2853$
$4421\!\!\sim\!\!2423\!\!\cong\!\!2423$
$4422\!\!\sim\!\!2367\!\!\cong\!\!2367$
$4423\!\!\sim\!\!2423\!\!\cong\!\!2423$
$4424\!\!\sim\!\!2261\!\!\cong\!\!2261$
$4425\!\!\sim\!\!2286\!\!\cong\!\!2286$
$4426\!\!\sim\!\!2367\!\!\cong\!\!2367$
$4427\!\!\sim\!\!2286\!\!\cong\!\!2286$
$4428\!\!\sim\!\!2205\!\!\cong\!\!775$
$4429\!\!\sim\!\!2862\!\!\cong\!\!847$
$4430\!\!\sim\!\!2427\!\!\cong\!\!2427$
$4431\!\!\sim\!\!2376\!\!\cong\!\!739$
$4432\!\!\sim\!\!2427\!\!\cong\!\!2427$
$4433\!\!\sim\!\!2265\!\!\cong\!\!2265$
$4434\!\!\sim\!\!2295\!\!\cong\!\!2295$
$4435\!\!\sim\!\!2376\!\!\cong\!\!739$
$4436\!\!\sim\!\!2295\!\!\cong\!\!2295$
$4437\!\!\sim\!\!2214\!\!\cong\!\!748$
$4438\!\!\sim\!\!2853\!\!\cong\!\!2853$
$4439\!\!\sim\!\!2423\!\!\cong\!\!2423$
$4440\!\!\sim\!\!2367\!\!\cong\!\!2367$
$4441\!\!\sim\!\!2423\!\!\cong\!\!2423$
$4442\!\!\sim\!\!2261\!\!\cong\!\!2261$
$4443\!\!\sim\!\!2286\!\!\cong\!\!2286$
$4444\!\!\sim\!\!2367\!\!\cong\!\!2367$
$4445\!\!\sim\!\!2286\!\!\cong\!\!2286$
$4446\!\!\sim\!\!2205\!\!\cong\!\!775$
$4447\!\!\sim\!\!2844\!\!\cong\!\!730$
$4448\!\!\sim\!\!2426\!\!\cong\!\!2277$
$4449\!\!\sim\!\!2358\!\!\cong\!\!820$
$4450\!\!\sim\!\!2426\!\!\cong\!\!2277$
$4451\!\!\sim\!\!2264\!\!\cong\!\!730$
$4452\!\!\sim\!\!2277\!\!\cong\!\!2277$
$4453\!\!\sim\!\!2358\!\!\cong\!\!820$
$4454\!\!\sim\!\!2277\!\!\cong\!\!2277$
$4455\!\!\sim\!\!2196\!\!\cong\!\!802$
$4456\!\!\sim\!\!2862\!\!\cong\!\!847$
$4457\!\!\sim\!\!2376\!\!\cong\!\!739$
$4458\!\!\sim\!\!2427\!\!\cong\!\!2427$
$4459\!\!\sim\!\!2376\!\!\cong\!\!739$
$4460\!\!\sim\!\!2214\!\!\cong\!\!748$
$4461\!\!\sim\!\!2295\!\!\cong\!\!2295$
$4462\!\!\sim\!\!2427\!\!\cong\!\!2427$
$4463\!\!\sim\!\!2295\!\!\cong\!\!2295$
$4464\!\!\sim\!\!2265\!\!\cong\!\!2265$
$4465\!\!\sim\!\!1091\!\!\cong\!\!731$
$4466\!\!\sim\!\!821\!\!\cong\!\!821$
$4467\!\!\sim\!\!821\!\!\cong\!\!821$
$4468\!\!\sim\!\!821\!\!\cong\!\!821$
$4469\!\!\sim\!\!731\!\!\cong\!\!731$
$4470\!\!\sim\!\!731\!\!\cong\!\!731$
$4471\!\!\sim\!\!821\!\!\cong\!\!821$
$4472\!\!\sim\!\!731\!\!\cong\!\!731$
$4473\!\!\sim\!\!731\!\!\cong\!\!731$
$4474\!\!\sim\!\!1091\!\!\cong\!\!731$
$4475\!\!\sim\!\!821\!\!\cong\!\!821$
$4476\!\!\sim\!\!821\!\!\cong\!\!821$
$4477\!\!\sim\!\!821\!\!\cong\!\!821$
$4478\!\!\sim\!\!731\!\!\cong\!\!731$
$4479\!\!\sim\!\!731\!\!\cong\!\!731$
$4480\!\!\sim\!\!821\!\!\cong\!\!821$
$4481\!\!\sim\!\!731\!\!\cong\!\!731$
$4482\!\!\sim\!\!731\!\!\cong\!\!731$
$4483\!\!\sim\!\!2860\!\!\cong\!\!2212$
$4484\!\!\sim\!\!2374\!\!\cong\!\!821$
$4485\!\!\sim\!\!2402\!\!\cong\!\!2402$
$4486\!\!\sim\!\!2374\!\!\cong\!\!821$
$4487\!\!\sim\!\!2212\!\!\cong\!\!2212$
$4488\!\!\sim\!\!2293\!\!\cong\!\!2293$
$4489\!\!\sim\!\!2402\!\!\cong\!\!2402$
$4490\!\!\sim\!\!2293\!\!\cong\!\!2293$
$4491\!\!\sim\!\!2240\!\!\cong\!\!2240$
$4492\!\!\sim\!\!2854\!\!\cong\!\!847$
$4493\!\!\sim\!\!2368\!\!\cong\!\!739$
$4494\!\!\sim\!\!2391\!\!\cong\!\!2391$
$4495\!\!\sim\!\!2368\!\!\cong\!\!739$
$4496\!\!\sim\!\!2206\!\!\cong\!\!748$
$4497\!\!\sim\!\!2287\!\!\cong\!\!2287$
$4498\!\!\sim\!\!2391\!\!\cong\!\!2391$
$4499\!\!\sim\!\!2287\!\!\cong\!\!2287$
$4500\!\!\sim\!\!2229\!\!\cong\!\!2229$
$4501\!\!\sim\!\!2850\!\!\cong\!\!2850$
$4502\!\!\sim\!\!2371\!\!\cong\!\!2371$
$4503\!\!\sim\!\!2364\!\!\cong\!\!2364$
$4504\!\!\sim\!\!2371\!\!\cong\!\!2371$
$4505\!\!\sim\!\!2209\!\!\cong\!\!2209$
$4506\!\!\sim\!\!2283\!\!\cong\!\!2283$
$4507\!\!\sim\!\!2364\!\!\cong\!\!2364$
$4508\!\!\sim\!\!2283\!\!\cong\!\!2283$
$4509\!\!\sim\!\!2202\!\!\cong\!\!2202$
$4510\!\!\sim\!\!2861\!\!\cong\!\!731$
$4511\!\!\sim\!\!2375\!\!\cong\!\!2375$
$4512\!\!\sim\!\!2375\!\!\cong\!\!2375$
$4513\!\!\sim\!\!2375\!\!\cong\!\!2375$
$4514\!\!\sim\!\!2213\!\!\cong\!\!2213$
$4515\!\!\sim\!\!2294\!\!\cong\!\!2294$
$4516\!\!\sim\!\!2375\!\!\cong\!\!2375$
$4517\!\!\sim\!\!2294\!\!\cong\!\!2294$
$4518\!\!\sim\!\!2213\!\!\cong\!\!2213$
$4519\!\!\sim\!\!2852\!\!\cong\!\!849$
$4520\!\!\sim\!\!2369\!\!\cong\!\!2369$
$4521\!\!\sim\!\!2366\!\!\cong\!\!2366$
$4522\!\!\sim\!\!2369\!\!\cong\!\!2369$
$4523\!\!\sim\!\!2207\!\!\cong\!\!2207$
$4524\!\!\sim\!\!2285\!\!\cong\!\!2285$
$4525\!\!\sim\!\!2366\!\!\cong\!\!2366$
$4526\!\!\sim\!\!2285\!\!\cong\!\!2285$
$4527\!\!\sim\!\!2204\!\!\cong\!\!2204$
$4528\!\!\sim\!\!2841\!\!\cong\!\!2841$
$4529\!\!\sim\!\!2372\!\!\cong\!\!2372$
$4530\!\!\sim\!\!2355\!\!\cong\!\!2355$
$4531\!\!\sim\!\!2372\!\!\cong\!\!2372$
$4532\!\!\sim\!\!2210\!\!\cong\!\!2210$
$4533\!\!\sim\!\!2274\!\!\cong\!\!2274$
$4534\!\!\sim\!\!2355\!\!\cong\!\!2355$
$4535\!\!\sim\!\!2274\!\!\cong\!\!2274$
$4536\!\!\sim\!\!2193\!\!\cong\!\!2193$
$4537\!\!\sim\!\!2889\!\!\cong\!\!750$
$4538\!\!\sim\!\!2403\!\!\cong\!\!2287$
$4539\!\!\sim\!\!2424\!\!\cong\!\!966$
$4540\!\!\sim\!\!2403\!\!\cong\!\!2287$
$4541\!\!\sim\!\!2241\!\!\cong\!\!739$
$4542\!\!\sim\!\!2322\!\!\cong\!\!2322$
$4543\!\!\sim\!\!2424\!\!\cong\!\!966$
$4544\!\!\sim\!\!2322\!\!\cong\!\!2322$
$4545\!\!\sim\!\!2262\!\!\cong\!\!750$
$4546\!\!\sim\!\!1091\!\!\cong\!\!731$
$4547\!\!\sim\!\!821\!\!\cong\!\!821$
$4548\!\!\sim\!\!821\!\!\cong\!\!821$
$4549\!\!\sim\!\!821\!\!\cong\!\!821$
$4550\!\!\sim\!\!731\!\!\cong\!\!731$
$4551\!\!\sim\!\!731\!\!\cong\!\!731$
$4552\!\!\sim\!\!821\!\!\cong\!\!821$
$4553\!\!\sim\!\!731\!\!\cong\!\!731$
$4554\!\!\sim\!\!731\!\!\cong\!\!731$
$4555\!\!\sim\!\!1091\!\!\cong\!\!731$
$4556\!\!\sim\!\!821\!\!\cong\!\!821$
$4557\!\!\sim\!\!821\!\!\cong\!\!821$
$4558\!\!\sim\!\!821\!\!\cong\!\!821$
$4559\!\!\sim\!\!731\!\!\cong\!\!731$
$4560\!\!\sim\!\!731\!\!\cong\!\!731$
$4561\!\!\sim\!\!821\!\!\cong\!\!821$
$4562\!\!\sim\!\!731\!\!\cong\!\!731$
$4563\!\!\sim\!\!731\!\!\cong\!\!731$
$4564\!\!\sim\!\!2887\!\!\cong\!\!731$
$4565\!\!\sim\!\!2401\!\!\cong\!\!2401$
$4566\!\!\sim\!\!2401\!\!\cong\!\!2401$
$4567\!\!\sim\!\!2401\!\!\cong\!\!2401$
$4568\!\!\sim\!\!2239\!\!\cong\!\!2239$
$4569\!\!\sim\!\!2320\!\!\cong\!\!2294$
$4570\!\!\sim\!\!2401\!\!\cong\!\!2401$
$4571\!\!\sim\!\!2320\!\!\cong\!\!2294$
$4572\!\!\sim\!\!2239\!\!\cong\!\!2239$
$4573\!\!\sim\!\!2874\!\!\cong\!\!820$
$4574\!\!\sim\!\!2395\!\!\cong\!\!2395$
$4575\!\!\sim\!\!2388\!\!\cong\!\!821$
$4576\!\!\sim\!\!2395\!\!\cong\!\!2395$
$4577\!\!\sim\!\!2233\!\!\cong\!\!2233$
$4578\!\!\sim\!\!2307\!\!\cong\!\!2307$
$4579\!\!\sim\!\!2388\!\!\cong\!\!821$
$4580\!\!\sim\!\!2307\!\!\cong\!\!2307$
$4581\!\!\sim\!\!2226\!\!\cong\!\!820$
$4582\!\!\sim\!\!2847\!\!\cong\!\!929$
$4583\!\!\sim\!\!2398\!\!\cong\!\!2398$
$4584\!\!\sim\!\!2361\!\!\cong\!\!2361$
$4585\!\!\sim\!\!2398\!\!\cong\!\!2398$
$4586\!\!\sim\!\!2236\!\!\cong\!\!2236$
$4587\!\!\sim\!\!2280\!\!\cong\!\!2280$
$4588\!\!\sim\!\!2361\!\!\cong\!\!2361$
$4589\!\!\sim\!\!2280\!\!\cong\!\!2280$
$4590\!\!\sim\!\!2199\!\!\cong\!\!2199$
$4591\!\!\sim\!\!2860\!\!\cong\!\!2212$
$4592\!\!\sim\!\!2402\!\!\cong\!\!2402$
$4593\!\!\sim\!\!2374\!\!\cong\!\!821$
$4594\!\!\sim\!\!2402\!\!\cong\!\!2402$
$4595\!\!\sim\!\!2240\!\!\cong\!\!2240$
$4596\!\!\sim\!\!2293\!\!\cong\!\!2293$
$4597\!\!\sim\!\!2374\!\!\cong\!\!821$
$4598\!\!\sim\!\!2293\!\!\cong\!\!2293$
$4599\!\!\sim\!\!2212\!\!\cong\!\!2212$
$4600\!\!\sim\!\!2851\!\!\cong\!\!929$
$4601\!\!\sim\!\!2396\!\!\cong\!\!2396$
$4602\!\!\sim\!\!2365\!\!\cong\!\!2365$
$4603\!\!\sim\!\!2396\!\!\cong\!\!2396$
$4604\!\!\sim\!\!2234\!\!\cong\!\!2234$
$4605\!\!\sim\!\!2284\!\!\cong\!\!2284$
$4606\!\!\sim\!\!2365\!\!\cong\!\!2365$
$4607\!\!\sim\!\!2284\!\!\cong\!\!2284$
$4608\!\!\sim\!\!2203\!\!\cong\!\!2203$
$4609\!\!\sim\!\!2838\!\!\cong\!\!750$
$4610\!\!\sim\!\!2399\!\!\cong\!\!2399$
$4611\!\!\sim\!\!2352\!\!\cong\!\!740$
$4612\!\!\sim\!\!2399\!\!\cong\!\!2399$
$4613\!\!\sim\!\!2237\!\!\cong\!\!2237$
$4614\!\!\sim\!\!2271\!\!\cong\!\!2271$
$4615\!\!\sim\!\!2352\!\!\cong\!\!740$
$4616\!\!\sim\!\!2271\!\!\cong\!\!2271$
$4617\!\!\sim\!\!2190\!\!\cong\!\!750$
$4618\!\!\sim\!\!2862\!\!\cong\!\!847$
$4619\!\!\sim\!\!2376\!\!\cong\!\!739$
$4620\!\!\sim\!\!2427\!\!\cong\!\!2427$
$4621\!\!\sim\!\!2376\!\!\cong\!\!739$
$4622\!\!\sim\!\!2214\!\!\cong\!\!748$
$4623\!\!\sim\!\!2295\!\!\cong\!\!2295$
$4624\!\!\sim\!\!2427\!\!\cong\!\!2427$
$4625\!\!\sim\!\!2295\!\!\cong\!\!2295$
$4626\!\!\sim\!\!2265\!\!\cong\!\!2265$
$4627\!\!\sim\!\!2860\!\!\cong\!\!2212$
$4628\!\!\sim\!\!2374\!\!\cong\!\!821$
$4629\!\!\sim\!\!2402\!\!\cong\!\!2402$
$4630\!\!\sim\!\!2374\!\!\cong\!\!821$
$4631\!\!\sim\!\!2212\!\!\cong\!\!2212$
$4632\!\!\sim\!\!2293\!\!\cong\!\!2293$
$4633\!\!\sim\!\!2402\!\!\cong\!\!2402$
$4634\!\!\sim\!\!2293\!\!\cong\!\!2293$
$4635\!\!\sim\!\!2240\!\!\cong\!\!2240$
$4636\!\!\sim\!\!2861\!\!\cong\!\!731$
$4637\!\!\sim\!\!2375\!\!\cong\!\!2375$
$4638\!\!\sim\!\!2375\!\!\cong\!\!2375$
$4639\!\!\sim\!\!2375\!\!\cong\!\!2375$
$4640\!\!\sim\!\!2213\!\!\cong\!\!2213$
$4641\!\!\sim\!\!2294\!\!\cong\!\!2294$
$4642\!\!\sim\!\!2375\!\!\cong\!\!2375$
$4643\!\!\sim\!\!2294\!\!\cong\!\!2294$
$4644\!\!\sim\!\!2213\!\!\cong\!\!2213$
$4645\!\!\sim\!\!1091\!\!\cong\!\!731$
$4646\!\!\sim\!\!821\!\!\cong\!\!821$
$4647\!\!\sim\!\!821\!\!\cong\!\!821$
$4648\!\!\sim\!\!821\!\!\cong\!\!821$
$4649\!\!\sim\!\!731\!\!\cong\!\!731$
$4650\!\!\sim\!\!731\!\!\cong\!\!731$
$4651\!\!\sim\!\!821\!\!\cong\!\!821$
$4652\!\!\sim\!\!731\!\!\cong\!\!731$
$4653\!\!\sim\!\!731\!\!\cong\!\!731$
$4654\!\!\sim\!\!2854\!\!\cong\!\!847$
$4655\!\!\sim\!\!2368\!\!\cong\!\!739$
$4656\!\!\sim\!\!2391\!\!\cong\!\!2391$
$4657\!\!\sim\!\!2368\!\!\cong\!\!739$
$4658\!\!\sim\!\!2206\!\!\cong\!\!748$
$4659\!\!\sim\!\!2287\!\!\cong\!\!2287$
$4660\!\!\sim\!\!2391\!\!\cong\!\!2391$
$4661\!\!\sim\!\!2287\!\!\cong\!\!2287$
$4662\!\!\sim\!\!2229\!\!\cong\!\!2229$
$4663\!\!\sim\!\!2852\!\!\cong\!\!849$
$4664\!\!\sim\!\!2369\!\!\cong\!\!2369$
$4665\!\!\sim\!\!2366\!\!\cong\!\!2366$
$4666\!\!\sim\!\!2369\!\!\cong\!\!2369$
$4667\!\!\sim\!\!2207\!\!\cong\!\!2207$
$4668\!\!\sim\!\!2285\!\!\cong\!\!2285$
$4669\!\!\sim\!\!2366\!\!\cong\!\!2366$
$4670\!\!\sim\!\!2285\!\!\cong\!\!2285$
$4671\!\!\sim\!\!2204\!\!\cong\!\!2204$
$4672\!\!\sim\!\!1091\!\!\cong\!\!731$
$4673\!\!\sim\!\!821\!\!\cong\!\!821$
$4674\!\!\sim\!\!821\!\!\cong\!\!821$
$4675\!\!\sim\!\!821\!\!\cong\!\!821$
$4676\!\!\sim\!\!731\!\!\cong\!\!731$
$4677\!\!\sim\!\!731\!\!\cong\!\!731$
$4678\!\!\sim\!\!821\!\!\cong\!\!821$
$4679\!\!\sim\!\!731\!\!\cong\!\!731$
$4680\!\!\sim\!\!731\!\!\cong\!\!731$
$4681\!\!\sim\!\!2850\!\!\cong\!\!2850$
$4682\!\!\sim\!\!2371\!\!\cong\!\!2371$
$4683\!\!\sim\!\!2364\!\!\cong\!\!2364$
$4684\!\!\sim\!\!2371\!\!\cong\!\!2371$
$4685\!\!\sim\!\!2209\!\!\cong\!\!2209$
$4686\!\!\sim\!\!2283\!\!\cong\!\!2283$
$4687\!\!\sim\!\!2364\!\!\cong\!\!2364$
$4688\!\!\sim\!\!2283\!\!\cong\!\!2283$
$4689\!\!\sim\!\!2202\!\!\cong\!\!2202$
$4690\!\!\sim\!\!2841\!\!\cong\!\!2841$
$4691\!\!\sim\!\!2372\!\!\cong\!\!2372$
$4692\!\!\sim\!\!2355\!\!\cong\!\!2355$
$4693\!\!\sim\!\!2372\!\!\cong\!\!2372$
$4694\!\!\sim\!\!2210\!\!\cong\!\!2210$
$4695\!\!\sim\!\!2274\!\!\cong\!\!2274$
$4696\!\!\sim\!\!2355\!\!\cong\!\!2355$
$4697\!\!\sim\!\!2274\!\!\cong\!\!2274$
$4698\!\!\sim\!\!2193\!\!\cong\!\!2193$
$4699\!\!\sim\!\!2844\!\!\cong\!\!730$
$4700\!\!\sim\!\!2358\!\!\cong\!\!820$
$4701\!\!\sim\!\!2426\!\!\cong\!\!2277$
$4702\!\!\sim\!\!2358\!\!\cong\!\!820$
$4703\!\!\sim\!\!2196\!\!\cong\!\!802$
$4704\!\!\sim\!\!2277\!\!\cong\!\!2277$
$4705\!\!\sim\!\!2426\!\!\cong\!\!2277$
$4706\!\!\sim\!\!2277\!\!\cong\!\!2277$
$4707\!\!\sim\!\!2264\!\!\cong\!\!730$
$4708\!\!\sim\!\!2838\!\!\cong\!\!750$
$4709\!\!\sim\!\!2352\!\!\cong\!\!740$
$4710\!\!\sim\!\!2399\!\!\cong\!\!2399$
$4711\!\!\sim\!\!2352\!\!\cong\!\!740$
$4712\!\!\sim\!\!2190\!\!\cong\!\!750$
$4713\!\!\sim\!\!2271\!\!\cong\!\!2271$
$4714\!\!\sim\!\!2399\!\!\cong\!\!2399$
$4715\!\!\sim\!\!2271\!\!\cong\!\!2271$
$4716\!\!\sim\!\!2237\!\!\cong\!\!2237$
$4717\!\!\sim\!\!2841\!\!\cong\!\!2841$
$4718\!\!\sim\!\!2355\!\!\cong\!\!2355$
$4719\!\!\sim\!\!2372\!\!\cong\!\!2372$
$4720\!\!\sim\!\!2355\!\!\cong\!\!2355$
$4721\!\!\sim\!\!2193\!\!\cong\!\!2193$
$4722\!\!\sim\!\!2274\!\!\cong\!\!2274$
$4723\!\!\sim\!\!2372\!\!\cong\!\!2372$
$4724\!\!\sim\!\!2274\!\!\cong\!\!2274$
$4725\!\!\sim\!\!2210\!\!\cong\!\!2210$
$4726\!\!\sim\!\!2838\!\!\cong\!\!750$
$4727\!\!\sim\!\!2352\!\!\cong\!\!740$
$4728\!\!\sim\!\!2399\!\!\cong\!\!2399$
$4729\!\!\sim\!\!2352\!\!\cong\!\!740$
$4730\!\!\sim\!\!2190\!\!\cong\!\!750$
$4731\!\!\sim\!\!2271\!\!\cong\!\!2271$
$4732\!\!\sim\!\!2399\!\!\cong\!\!2399$
$4733\!\!\sim\!\!2271\!\!\cong\!\!2271$
$4734\!\!\sim\!\!2237\!\!\cong\!\!2237$
$4735\!\!\sim\!\!1090\!\!\cong\!\!1090$
$4736\!\!\sim\!\!820\!\!\cong\!\!820$
$4737\!\!\sim\!\!820\!\!\cong\!\!820$
$4738\!\!\sim\!\!820\!\!\cong\!\!820$
$4739\!\!\sim\!\!730\!\!\cong\!\!730$
$4740\!\!\sim\!\!730\!\!\cong\!\!730$
$4741\!\!\sim\!\!820\!\!\cong\!\!820$
$4742\!\!\sim\!\!730\!\!\cong\!\!730$
$4743\!\!\sim\!\!730\!\!\cong\!\!730$
$4744\!\!\sim\!\!1090\!\!\cong\!\!1090$
$4745\!\!\sim\!\!820\!\!\cong\!\!820$
$4746\!\!\sim\!\!820\!\!\cong\!\!820$
$4747\!\!\sim\!\!820\!\!\cong\!\!820$
$4748\!\!\sim\!\!730\!\!\cong\!\!730$
$4749\!\!\sim\!\!730\!\!\cong\!\!730$
$4750\!\!\sim\!\!820\!\!\cong\!\!820$
$4751\!\!\sim\!\!730\!\!\cong\!\!730$
$4752\!\!\sim\!\!730\!\!\cong\!\!730$
$4753\!\!\sim\!\!2841\!\!\cong\!\!2841$
$4754\!\!\sim\!\!2355\!\!\cong\!\!2355$
$4755\!\!\sim\!\!2372\!\!\cong\!\!2372$
$4756\!\!\sim\!\!2355\!\!\cong\!\!2355$
$4757\!\!\sim\!\!2193\!\!\cong\!\!2193$
$4758\!\!\sim\!\!2274\!\!\cong\!\!2274$
$4759\!\!\sim\!\!2372\!\!\cong\!\!2372$
$4760\!\!\sim\!\!2274\!\!\cong\!\!2274$
$4761\!\!\sim\!\!2210\!\!\cong\!\!2210$
$4762\!\!\sim\!\!1090\!\!\cong\!\!1090$
$4763\!\!\sim\!\!820\!\!\cong\!\!820$
$4764\!\!\sim\!\!820\!\!\cong\!\!820$
$4765\!\!\sim\!\!820\!\!\cong\!\!820$
$4766\!\!\sim\!\!730\!\!\cong\!\!730$
$4767\!\!\sim\!\!730\!\!\cong\!\!730$
$4768\!\!\sim\!\!820\!\!\cong\!\!820$
$4769\!\!\sim\!\!730\!\!\cong\!\!730$
$4770\!\!\sim\!\!730\!\!\cong\!\!730$
$4771\!\!\sim\!\!1090\!\!\cong\!\!1090$
$4772\!\!\sim\!\!820\!\!\cong\!\!820$
$4773\!\!\sim\!\!820\!\!\cong\!\!820$
$4774\!\!\sim\!\!820\!\!\cong\!\!820$
$4775\!\!\sim\!\!730\!\!\cong\!\!730$
$4776\!\!\sim\!\!730\!\!\cong\!\!730$
$4777\!\!\sim\!\!820\!\!\cong\!\!820$
$4778\!\!\sim\!\!730\!\!\cong\!\!730$
$4779\!\!\sim\!\!730\!\!\cong\!\!730$
$4780\!\!\sim\!\!2853\!\!\cong\!\!2853$
$4781\!\!\sim\!\!2367\!\!\cong\!\!2367$
$4782\!\!\sim\!\!2423\!\!\cong\!\!2423$
$4783\!\!\sim\!\!2367\!\!\cong\!\!2367$
$4784\!\!\sim\!\!2205\!\!\cong\!\!775$
$4785\!\!\sim\!\!2286\!\!\cong\!\!2286$
$4786\!\!\sim\!\!2423\!\!\cong\!\!2423$
$4787\!\!\sim\!\!2286\!\!\cong\!\!2286$
$4788\!\!\sim\!\!2261\!\!\cong\!\!2261$
$4789\!\!\sim\!\!2851\!\!\cong\!\!929$
$4790\!\!\sim\!\!2365\!\!\cong\!\!2365$
$4791\!\!\sim\!\!2396\!\!\cong\!\!2396$
$4792\!\!\sim\!\!2365\!\!\cong\!\!2365$
$4793\!\!\sim\!\!2203\!\!\cong\!\!2203$
$4794\!\!\sim\!\!2284\!\!\cong\!\!2284$
$4795\!\!\sim\!\!2396\!\!\cong\!\!2396$
$4796\!\!\sim\!\!2284\!\!\cong\!\!2284$
$4797\!\!\sim\!\!2234\!\!\cong\!\!2234$
$4798\!\!\sim\!\!2852\!\!\cong\!\!849$
$4799\!\!\sim\!\!2366\!\!\cong\!\!2366$
$4800\!\!\sim\!\!2369\!\!\cong\!\!2369$
$4801\!\!\sim\!\!2366\!\!\cong\!\!2366$
$4802\!\!\sim\!\!2204\!\!\cong\!\!2204$
$4803\!\!\sim\!\!2285\!\!\cong\!\!2285$
$4804\!\!\sim\!\!2369\!\!\cong\!\!2369$
$4805\!\!\sim\!\!2285\!\!\cong\!\!2285$
$4806\!\!\sim\!\!2207\!\!\cong\!\!2207$
$4807\!\!\sim\!\!2847\!\!\cong\!\!929$
$4808\!\!\sim\!\!2361\!\!\cong\!\!2361$
$4809\!\!\sim\!\!2398\!\!\cong\!\!2398$
$4810\!\!\sim\!\!2361\!\!\cong\!\!2361$
$4811\!\!\sim\!\!2199\!\!\cong\!\!2199$
$4812\!\!\sim\!\!2280\!\!\cong\!\!2280$
$4813\!\!\sim\!\!2398\!\!\cong\!\!2398$
$4814\!\!\sim\!\!2280\!\!\cong\!\!2280$
$4815\!\!\sim\!\!2236\!\!\cong\!\!2236$
$4816\!\!\sim\!\!1090\!\!\cong\!\!1090$
$4817\!\!\sim\!\!820\!\!\cong\!\!820$
$4818\!\!\sim\!\!820\!\!\cong\!\!820$
$4819\!\!\sim\!\!820\!\!\cong\!\!820$
$4820\!\!\sim\!\!730\!\!\cong\!\!730$
$4821\!\!\sim\!\!730\!\!\cong\!\!730$
$4822\!\!\sim\!\!820\!\!\cong\!\!820$
$4823\!\!\sim\!\!730\!\!\cong\!\!730$
$4824\!\!\sim\!\!730\!\!\cong\!\!730$
$4825\!\!\sim\!\!1090\!\!\cong\!\!1090$
$4826\!\!\sim\!\!820\!\!\cong\!\!820$
$4827\!\!\sim\!\!820\!\!\cong\!\!820$
$4828\!\!\sim\!\!820\!\!\cong\!\!820$
$4829\!\!\sim\!\!730\!\!\cong\!\!730$
$4830\!\!\sim\!\!730\!\!\cong\!\!730$
$4831\!\!\sim\!\!820\!\!\cong\!\!820$
$4832\!\!\sim\!\!730\!\!\cong\!\!730$
$4833\!\!\sim\!\!730\!\!\cong\!\!730$
$4834\!\!\sim\!\!2850\!\!\cong\!\!2850$
$4835\!\!\sim\!\!2364\!\!\cong\!\!2364$
$4836\!\!\sim\!\!2371\!\!\cong\!\!2371$
$4837\!\!\sim\!\!2364\!\!\cong\!\!2364$
$4838\!\!\sim\!\!2202\!\!\cong\!\!2202$
$4839\!\!\sim\!\!2283\!\!\cong\!\!2283$
$4840\!\!\sim\!\!2371\!\!\cong\!\!2371$
$4841\!\!\sim\!\!2283\!\!\cong\!\!2283$
$4842\!\!\sim\!\!2209\!\!\cong\!\!2209$
$4843\!\!\sim\!\!1090\!\!\cong\!\!1090$
$4844\!\!\sim\!\!820\!\!\cong\!\!820$
$4845\!\!\sim\!\!820\!\!\cong\!\!820$
$4846\!\!\sim\!\!820\!\!\cong\!\!820$
$4847\!\!\sim\!\!730\!\!\cong\!\!730$
$4848\!\!\sim\!\!730\!\!\cong\!\!730$
$4849\!\!\sim\!\!820\!\!\cong\!\!820$
$4850\!\!\sim\!\!730\!\!\cong\!\!730$
$4851\!\!\sim\!\!730\!\!\cong\!\!730$
$4852\!\!\sim\!\!1090\!\!\cong\!\!1090$
$4853\!\!\sim\!\!820\!\!\cong\!\!820$
$4854\!\!\sim\!\!820\!\!\cong\!\!820$
$4855\!\!\sim\!\!820\!\!\cong\!\!820$
$4856\!\!\sim\!\!730\!\!\cong\!\!730$
$4857\!\!\sim\!\!730\!\!\cong\!\!730$
$4858\!\!\sim\!\!820\!\!\cong\!\!820$
$4859\!\!\sim\!\!730\!\!\cong\!\!730$
$4860\!\!\sim\!\!730\!\!\cong\!\!730$
$4861\!\!\sim\!\!2889\!\!\cong\!\!750$
$4862\!\!\sim\!\!2403\!\!\cong\!\!2287$
$4863\!\!\sim\!\!2424\!\!\cong\!\!966$
$4864\!\!\sim\!\!2403\!\!\cong\!\!2287$
$4865\!\!\sim\!\!2241\!\!\cong\!\!739$
$4866\!\!\sim\!\!2322\!\!\cong\!\!2322$
$4867\!\!\sim\!\!2424\!\!\cong\!\!966$
$4868\!\!\sim\!\!2322\!\!\cong\!\!2322$
$4869\!\!\sim\!\!2262\!\!\cong\!\!750$
$4870\!\!\sim\!\!2887\!\!\cong\!\!731$
$4871\!\!\sim\!\!2401\!\!\cong\!\!2401$
$4872\!\!\sim\!\!2401\!\!\cong\!\!2401$
$4873\!\!\sim\!\!2401\!\!\cong\!\!2401$
$4874\!\!\sim\!\!2239\!\!\cong\!\!2239$
$4875\!\!\sim\!\!2320\!\!\cong\!\!2294$
$4876\!\!\sim\!\!2401\!\!\cong\!\!2401$
$4877\!\!\sim\!\!2320\!\!\cong\!\!2294$
$4878\!\!\sim\!\!2239\!\!\cong\!\!2239$
$4879\!\!\sim\!\!2860\!\!\cong\!\!2212$
$4880\!\!\sim\!\!2402\!\!\cong\!\!2402$
$4881\!\!\sim\!\!2374\!\!\cong\!\!821$
$4882\!\!\sim\!\!2402\!\!\cong\!\!2402$
$4883\!\!\sim\!\!2240\!\!\cong\!\!2240$
$4884\!\!\sim\!\!2293\!\!\cong\!\!2293$
$4885\!\!\sim\!\!2374\!\!\cong\!\!821$
$4886\!\!\sim\!\!2293\!\!\cong\!\!2293$
$4887\!\!\sim\!\!2212\!\!\cong\!\!2212$
$4888\!\!\sim\!\!1091\!\!\cong\!\!731$
$4889\!\!\sim\!\!821\!\!\cong\!\!821$
$4890\!\!\sim\!\!821\!\!\cong\!\!821$
$4891\!\!\sim\!\!821\!\!\cong\!\!821$
$4892\!\!\sim\!\!731\!\!\cong\!\!731$
$4893\!\!\sim\!\!731\!\!\cong\!\!731$
$4894\!\!\sim\!\!821\!\!\cong\!\!821$
$4895\!\!\sim\!\!731\!\!\cong\!\!731$
$4896\!\!\sim\!\!731\!\!\cong\!\!731$
$4897\!\!\sim\!\!2874\!\!\cong\!\!820$
$4898\!\!\sim\!\!2395\!\!\cong\!\!2395$
$4899\!\!\sim\!\!2388\!\!\cong\!\!821$
$4900\!\!\sim\!\!2395\!\!\cong\!\!2395$
$4901\!\!\sim\!\!2233\!\!\cong\!\!2233$
$4902\!\!\sim\!\!2307\!\!\cong\!\!2307$
$4903\!\!\sim\!\!2388\!\!\cong\!\!821$
$4904\!\!\sim\!\!2307\!\!\cong\!\!2307$
$4905\!\!\sim\!\!2226\!\!\cong\!\!820$
$4906\!\!\sim\!\!2851\!\!\cong\!\!929$
$4907\!\!\sim\!\!2396\!\!\cong\!\!2396$
$4908\!\!\sim\!\!2365\!\!\cong\!\!2365$
$4909\!\!\sim\!\!2396\!\!\cong\!\!2396$
$4910\!\!\sim\!\!2234\!\!\cong\!\!2234$
$4911\!\!\sim\!\!2284\!\!\cong\!\!2284$
$4912\!\!\sim\!\!2365\!\!\cong\!\!2365$
$4913\!\!\sim\!\!2284\!\!\cong\!\!2284$
$4914\!\!\sim\!\!2203\!\!\cong\!\!2203$
$4915\!\!\sim\!\!1091\!\!\cong\!\!731$
$4916\!\!\sim\!\!821\!\!\cong\!\!821$
$4917\!\!\sim\!\!821\!\!\cong\!\!821$
$4918\!\!\sim\!\!821\!\!\cong\!\!821$
$4919\!\!\sim\!\!731\!\!\cong\!\!731$
$4920\!\!\sim\!\!731\!\!\cong\!\!731$
$4921\!\!\sim\!\!821\!\!\cong\!\!821$
$4922\!\!\sim\!\!731\!\!\cong\!\!731$
$4923\!\!\sim\!\!731\!\!\cong\!\!731$
$4924\!\!\sim\!\!2847\!\!\cong\!\!929$
$4925\!\!\sim\!\!2398\!\!\cong\!\!2398$
$4926\!\!\sim\!\!2361\!\!\cong\!\!2361$
$4927\!\!\sim\!\!2398\!\!\cong\!\!2398$
$4928\!\!\sim\!\!2236\!\!\cong\!\!2236$
$4929\!\!\sim\!\!2280\!\!\cong\!\!2280$
$4930\!\!\sim\!\!2361\!\!\cong\!\!2361$
$4931\!\!\sim\!\!2280\!\!\cong\!\!2280$
$4932\!\!\sim\!\!2199\!\!\cong\!\!2199$
$4933\!\!\sim\!\!2838\!\!\cong\!\!750$
$4934\!\!\sim\!\!2399\!\!\cong\!\!2399$
$4935\!\!\sim\!\!2352\!\!\cong\!\!740$
$4936\!\!\sim\!\!2399\!\!\cong\!\!2399$
$4937\!\!\sim\!\!2237\!\!\cong\!\!2237$
$4938\!\!\sim\!\!2271\!\!\cong\!\!2271$
$4939\!\!\sim\!\!2352\!\!\cong\!\!740$
$4940\!\!\sim\!\!2271\!\!\cong\!\!2271$
$4941\!\!\sim\!\!2190\!\!\cong\!\!750$
$4942\!\!\sim\!\!2853\!\!\cong\!\!2853$
$4943\!\!\sim\!\!2367\!\!\cong\!\!2367$
$4944\!\!\sim\!\!2423\!\!\cong\!\!2423$
$4945\!\!\sim\!\!2367\!\!\cong\!\!2367$
$4946\!\!\sim\!\!2205\!\!\cong\!\!775$
$4947\!\!\sim\!\!2286\!\!\cong\!\!2286$
$4948\!\!\sim\!\!2423\!\!\cong\!\!2423$
$4949\!\!\sim\!\!2286\!\!\cong\!\!2286$
$4950\!\!\sim\!\!2261\!\!\cong\!\!2261$
$4951\!\!\sim\!\!2847\!\!\cong\!\!929$
$4952\!\!\sim\!\!2361\!\!\cong\!\!2361$
$4953\!\!\sim\!\!2398\!\!\cong\!\!2398$
$4954\!\!\sim\!\!2361\!\!\cong\!\!2361$
$4955\!\!\sim\!\!2199\!\!\cong\!\!2199$
$4956\!\!\sim\!\!2280\!\!\cong\!\!2280$
$4957\!\!\sim\!\!2398\!\!\cong\!\!2398$
$4958\!\!\sim\!\!2280\!\!\cong\!\!2280$
$4959\!\!\sim\!\!2236\!\!\cong\!\!2236$
$4960\!\!\sim\!\!2850\!\!\cong\!\!2850$
$4961\!\!\sim\!\!2364\!\!\cong\!\!2364$
$4962\!\!\sim\!\!2371\!\!\cong\!\!2371$
$4963\!\!\sim\!\!2364\!\!\cong\!\!2364$
$4964\!\!\sim\!\!2202\!\!\cong\!\!2202$
$4965\!\!\sim\!\!2283\!\!\cong\!\!2283$
$4966\!\!\sim\!\!2371\!\!\cong\!\!2371$
$4967\!\!\sim\!\!2283\!\!\cong\!\!2283$
$4968\!\!\sim\!\!2209\!\!\cong\!\!2209$
$4969\!\!\sim\!\!2851\!\!\cong\!\!929$
$4970\!\!\sim\!\!2365\!\!\cong\!\!2365$
$4971\!\!\sim\!\!2396\!\!\cong\!\!2396$
$4972\!\!\sim\!\!2365\!\!\cong\!\!2365$
$4973\!\!\sim\!\!2203\!\!\cong\!\!2203$
$4974\!\!\sim\!\!2284\!\!\cong\!\!2284$
$4975\!\!\sim\!\!2396\!\!\cong\!\!2396$
$4976\!\!\sim\!\!2284\!\!\cong\!\!2284$
$4977\!\!\sim\!\!2234\!\!\cong\!\!2234$
$4978\!\!\sim\!\!1090\!\!\cong\!\!1090$
$4979\!\!\sim\!\!820\!\!\cong\!\!820$
$4980\!\!\sim\!\!820\!\!\cong\!\!820$
$4981\!\!\sim\!\!820\!\!\cong\!\!820$
$4982\!\!\sim\!\!730\!\!\cong\!\!730$
$4983\!\!\sim\!\!730\!\!\cong\!\!730$
$4984\!\!\sim\!\!820\!\!\cong\!\!820$
$4985\!\!\sim\!\!730\!\!\cong\!\!730$
$4986\!\!\sim\!\!730\!\!\cong\!\!730$
$4987\!\!\sim\!\!1090\!\!\cong\!\!1090$
$4988\!\!\sim\!\!820\!\!\cong\!\!820$
$4989\!\!\sim\!\!820\!\!\cong\!\!820$
$4990\!\!\sim\!\!820\!\!\cong\!\!820$
$4991\!\!\sim\!\!730\!\!\cong\!\!730$
$4992\!\!\sim\!\!730\!\!\cong\!\!730$
$4993\!\!\sim\!\!820\!\!\cong\!\!820$
$4994\!\!\sim\!\!730\!\!\cong\!\!730$
$4995\!\!\sim\!\!730\!\!\cong\!\!730$
$4996\!\!\sim\!\!2852\!\!\cong\!\!849$
$4997\!\!\sim\!\!2366\!\!\cong\!\!2366$
$4998\!\!\sim\!\!2369\!\!\cong\!\!2369$
$4999\!\!\sim\!\!2366\!\!\cong\!\!2366$
$5000\!\!\sim\!\!2204\!\!\cong\!\!2204$
$5001\!\!\sim\!\!2285\!\!\cong\!\!2285$
$5002\!\!\sim\!\!2369\!\!\cong\!\!2369$
$5003\!\!\sim\!\!2285\!\!\cong\!\!2285$
$5004\!\!\sim\!\!2207\!\!\cong\!\!2207$
$5005\!\!\sim\!\!1090\!\!\cong\!\!1090$
$5006\!\!\sim\!\!820\!\!\cong\!\!820$
$5007\!\!\sim\!\!820\!\!\cong\!\!820$
$5008\!\!\sim\!\!820\!\!\cong\!\!820$
$5009\!\!\sim\!\!730\!\!\cong\!\!730$
$5010\!\!\sim\!\!730\!\!\cong\!\!730$
$5011\!\!\sim\!\!820\!\!\cong\!\!820$
$5012\!\!\sim\!\!730\!\!\cong\!\!730$
$5013\!\!\sim\!\!730\!\!\cong\!\!730$
$5014\!\!\sim\!\!1090\!\!\cong\!\!1090$
$5015\!\!\sim\!\!820\!\!\cong\!\!820$
$5016\!\!\sim\!\!820\!\!\cong\!\!820$
$5017\!\!\sim\!\!820\!\!\cong\!\!820$
$5018\!\!\sim\!\!730\!\!\cong\!\!730$
$5019\!\!\sim\!\!730\!\!\cong\!\!730$
$5020\!\!\sim\!\!820\!\!\cong\!\!820$
$5021\!\!\sim\!\!730\!\!\cong\!\!730$
$5022\!\!\sim\!\!730\!\!\cong\!\!730$
$5023\!\!\sim\!\!2880\!\!\cong\!\!730$
$5024\!\!\sim\!\!2394\!\!\cong\!\!820$
$5025\!\!\sim\!\!2422\!\!\cong\!\!820$
$5026\!\!\sim\!\!2394\!\!\cong\!\!820$
$5027\!\!\sim\!\!2232\!\!\cong\!\!730$
$5028\!\!\sim\!\!2313\!\!\cong\!\!2277$
$5029\!\!\sim\!\!2422\!\!\cong\!\!820$
$5030\!\!\sim\!\!2313\!\!\cong\!\!2277$
$5031\!\!\sim\!\!2260\!\!\cong\!\!802$
$5032\!\!\sim\!\!2874\!\!\cong\!\!820$
$5033\!\!\sim\!\!2388\!\!\cong\!\!821$
$5034\!\!\sim\!\!2395\!\!\cong\!\!2395$
$5035\!\!\sim\!\!2388\!\!\cong\!\!821$
$5036\!\!\sim\!\!2226\!\!\cong\!\!820$
$5037\!\!\sim\!\!2307\!\!\cong\!\!2307$
$5038\!\!\sim\!\!2395\!\!\cong\!\!2395$
$5039\!\!\sim\!\!2307\!\!\cong\!\!2307$
$5040\!\!\sim\!\!2233\!\!\cong\!\!2233$
$5041\!\!\sim\!\!2854\!\!\cong\!\!847$
$5042\!\!\sim\!\!2391\!\!\cong\!\!2391$
$5043\!\!\sim\!\!2368\!\!\cong\!\!739$
$5044\!\!\sim\!\!2391\!\!\cong\!\!2391$
$5045\!\!\sim\!\!2229\!\!\cong\!\!2229$
$5046\!\!\sim\!\!2287\!\!\cong\!\!2287$
$5047\!\!\sim\!\!2368\!\!\cong\!\!739$
$5048\!\!\sim\!\!2287\!\!\cong\!\!2287$
$5049\!\!\sim\!\!2206\!\!\cong\!\!748$
$5050\!\!\sim\!\!2874\!\!\cong\!\!820$
$5051\!\!\sim\!\!2388\!\!\cong\!\!821$
$5052\!\!\sim\!\!2395\!\!\cong\!\!2395$
$5053\!\!\sim\!\!2388\!\!\cong\!\!821$
$5054\!\!\sim\!\!2226\!\!\cong\!\!820$
$5055\!\!\sim\!\!2307\!\!\cong\!\!2307$
$5056\!\!\sim\!\!2395\!\!\cong\!\!2395$
$5057\!\!\sim\!\!2307\!\!\cong\!\!2307$
$5058\!\!\sim\!\!2233\!\!\cong\!\!2233$
$5059\!\!\sim\!\!1090\!\!\cong\!\!1090$
$5060\!\!\sim\!\!820\!\!\cong\!\!820$
$5061\!\!\sim\!\!820\!\!\cong\!\!820$
$5062\!\!\sim\!\!820\!\!\cong\!\!820$
$5063\!\!\sim\!\!730\!\!\cong\!\!730$
$5064\!\!\sim\!\!730\!\!\cong\!\!730$
$5065\!\!\sim\!\!820\!\!\cong\!\!820$
$5066\!\!\sim\!\!730\!\!\cong\!\!730$
$5067\!\!\sim\!\!730\!\!\cong\!\!730$
$5068\!\!\sim\!\!1090\!\!\cong\!\!1090$
$5069\!\!\sim\!\!820\!\!\cong\!\!820$
$5070\!\!\sim\!\!820\!\!\cong\!\!820$
$5071\!\!\sim\!\!820\!\!\cong\!\!820$
$5072\!\!\sim\!\!730\!\!\cong\!\!730$
$5073\!\!\sim\!\!730\!\!\cong\!\!730$
$5074\!\!\sim\!\!820\!\!\cong\!\!820$
$5075\!\!\sim\!\!730\!\!\cong\!\!730$
$5076\!\!\sim\!\!730\!\!\cong\!\!730$
$5077\!\!\sim\!\!2854\!\!\cong\!\!847$
$5078\!\!\sim\!\!2391\!\!\cong\!\!2391$
$5079\!\!\sim\!\!2368\!\!\cong\!\!739$
$5080\!\!\sim\!\!2391\!\!\cong\!\!2391$
$5081\!\!\sim\!\!2229\!\!\cong\!\!2229$
$5082\!\!\sim\!\!2287\!\!\cong\!\!2287$
$5083\!\!\sim\!\!2368\!\!\cong\!\!739$
$5084\!\!\sim\!\!2287\!\!\cong\!\!2287$
$5085\!\!\sim\!\!2206\!\!\cong\!\!748$
$5086\!\!\sim\!\!1090\!\!\cong\!\!1090$
$5087\!\!\sim\!\!820\!\!\cong\!\!820$
$5088\!\!\sim\!\!820\!\!\cong\!\!820$
$5089\!\!\sim\!\!820\!\!\cong\!\!820$
$5090\!\!\sim\!\!730\!\!\cong\!\!730$
$5091\!\!\sim\!\!730\!\!\cong\!\!730$
$5092\!\!\sim\!\!820\!\!\cong\!\!820$
$5093\!\!\sim\!\!730\!\!\cong\!\!730$
$5094\!\!\sim\!\!730\!\!\cong\!\!730$
$5095\!\!\sim\!\!1090\!\!\cong\!\!1090$
$5096\!\!\sim\!\!820\!\!\cong\!\!820$
$5097\!\!\sim\!\!820\!\!\cong\!\!820$
$5098\!\!\sim\!\!820\!\!\cong\!\!820$
$5099\!\!\sim\!\!730\!\!\cong\!\!730$
$5100\!\!\sim\!\!730\!\!\cong\!\!730$
$5101\!\!\sim\!\!820\!\!\cong\!\!820$
$5102\!\!\sim\!\!730\!\!\cong\!\!730$
$5103\!\!\sim\!\!730\!\!\cong\!\!730$
\phantom{$5103\!\!\sim\!\!730\!\!\cong\!\!730$}}
\end{multicols}

\noindent{\footnotesize 5104 through $5832\sim1090\simeq1090$.}

\input{08data.tex}
\newpage
\section{Proofs}\label{app_proofs}

This section contains proofs of many of the claims contained in the
tables in Section~\ref{app_corrtable} and Section~\ref{app_table}
and some additional information.

We sometimes encounter one of the following four binary tree
automorphisms
\[
 a = \s(1,a), \qquad b=\s(b,1), \qquad
 c=\s(c^{-1},1), \qquad d=\s(1,d^{-1}).
\]
The first one is the binary adding machine, the second is its
inverse, and all are conjugate to the adding machine and therefore
act level transitively on the binary tree and have infinite order.

We freely use the known classification of groups generated by
2-state automata over a 2-letter alphabet.

\begin{theorem}[\cite{gns00:automata}]\label{thm:class22}
Up to isomorphism, there are six $(2,2)$-automaton groups: the
trivial group, the cyclic group of order 2 (we denote it by $C_2$),
Klein group $C_2\times C_2$ of order $4$, the infinite cyclic group
$\mathbb Z$, the infinite dihedral group $D_\infty$ and the
Lamplighter group $\mathbb Z\wr C_2$.

In particular the sixteen 2-state automata for which both states are
inactive generate the trivial group, and the sixteen 2-state
automata in which both states are active generate $C_2$ (since both
states in that case describe the mirror automorphism
$\mu=\s(\mu,\mu)$ of order 2.

The automata given by either of the wreath recursions
\begin{alignat*}{2}
 a&= \s(a,a), \qquad b&&= (a,a), \\
 a&= \s(b,b), \qquad b&&= (a,a),
\end{alignat*}
generate the Klein group $C_2 \times C_2$.

The automata given by the wreath recursions
\begin{alignat*}{2}
 a&= \s(a,a), \qquad b&&= (a,b), \\
 a&= \s(a,a), \qquad b&&= (b,a), \\
 a&= \s(b,b), \qquad b&&= (a,b), \\
 a&= \s(b,b), \qquad b&&= (b,a),
\end{alignat*}
generate the infinite dihedral group $D_\infty$.

The automata given by the wreath recursions
\begin{alignat*}{2}
 a&= \s(a,a), \qquad b&&= (b,b), \\
 a&= \s(b,b), \qquad b&&= (b,b),
\end{alignat*}
generate the cyclic group $C_2$.

The automata given by the wreath recursions
\begin{alignat*}{2}
 a&= \s(a,b), \qquad b&&= (a,a), \\
 a&= \s(b,a), \qquad b&&= (a,a), \\
 a&= \s(a,b), \qquad b&&= (b,b), \\
 a&= \s(b,a), \qquad b&&= (b,a),
\end{alignat*}
generate the infinite cyclic group $\Z$. Moreover, in the first two
cases we have $b=a^{-2}$, in the fourth case $b=1$ and $a$ is the
adding machine, and in the third case $b=1$ and $a$ is the inverse
of the adding machine.

The automata given by the wreath recursions
\begin{alignat*}{2}
 a&= \s(a,b), \qquad b&&= (a,b), \\
 a&= \s(a,b), \qquad b&&= (b,a), \\
 a&= \s(b,a), \qquad b&&= (a,b), \\
 a&= \s(b,a), \qquad b&&= (b,a),
\end{alignat*}
generate the Lamplighter group $\Z \wr C_2 = \Z \ltimes \left(
\oplus_\Z C_2 \right)$.
\end{theorem}

The results on the next few pages concern the existence of elements
of infinite order and the level transitivity of the action. They are
used in some of the proofs that follow.

\begin{lemma}[\cite{bondarenko-al:classification32-1}]\label{nontors}
Let $G$ be a group generated by an automaton $\A$ over a $2$-letter
alphabet. Assume that the set of states $S$ of $\A$ splits into two
nonempty parts $P$ and $Q$ such that
\begin{enumerate}
\item[\textup{(i)}]
one of the parts consists of the active states (those with
nontrivial vertex permutation) and the other consists of the
inactive states;
\item[\textup{(ii)}]
for each state from $P$, both arrows go to states in the same part
(either both to $P$ or both to $Q$);
\item[\textup{(iii)}]
for each state from $Q$, one arrow goes to a state in $P$ and the
other to a state in $Q$.
\end{enumerate}
Then any element of the group that can be written as a product of
odd number of active generators or their inverses and odd number of
inactive generators and their inverses, in any order, has infinite
order. In particular, the group $G$ is not a torsion group.
\end{lemma}

\begin{proof}
Denote by $D$ the set of elements in $G$ that can be represented as
a product of odd number of active generators or their inverses and
odd number of inactive generators and their inverses, in any order.

We note that if $g \in D$ then both sections of $g^2$ are in $D$.
Indeed, for such an element, $g=\sigma(g_0,g_1)$ and $g^2
=(g_1g_0,g_0g_1)$. Both sections of $g^2$ are products (in some
order) of the first level sections of the generators (and/or their
inverses) used to express $g$ as an element in $D$.  By assumption,
among these generators, there are odd number of active and odd
number of inactive ones. The generators from $P$, by condition (ii),
produce even number of active and even number of inactive sections
on level 1, while the generators from $Q$, by condition (iii),
produce odd number of active sections and odd number of inactive
sections. Thus both sections of $g$ are in $D$.

By way of contradiction, assume that $h$ is an element of $D$ of
finite order $2^n$, for some $n \geq 0$. If $n > 0$ the sections of
$h^2$ are elements in $D$ of order $2^{n-1}$. Thus, continuing in
this fashion, we reach an element in $D$ that is trivial. This is
contradiction since all elements in $D$ act nontrivially on level 1.
\end{proof}

There is a simple criterion that determines whether a given element
of a self-similar group generated by a finite automaton over the
$2$-letter alphabet $X=\{0,1\}$ acts level transitively on the tree.
The criterion is based on the image of the given element in the
abelianization of $\Aut(X^*)$, which is isomorphic to the infinite
Cartesian product $\prod_{i=0}^\infty C_2$. The canonical
isomorphism sends $g\in G$ to $(a_i \mod 2)_{i=0}^\infty$, where
$a_i$ is the number of active sections of $g$ at level $i$. We also
make use of the ring structure on $\prod_{i=0}^\infty C_2$ obtained
by identifying $(b_i)_{i=0}^\infty$ with $\sum_{i=0}^\infty b_it^i$
in the ring of formal power series $C_2[[t]]$. It is known that a
binary tree automorphism $g$ acts level transitively on $X^*$ if and
only if $\bar g=(1,1,1,\ldots)$, where $\bar g$ be the image of $g$
in the abelianization $\prod_{i=0}^\infty C_2$ of $\Aut(X^*)$.

\begin{lemma}[Element transitivity, \cite{bondarenko-al:classification32-1}]\label{alg_trans}
Let $G$ be a group generated by an automaton $\A$ over a $2$-letter
alphabet. There exists an algorithm that decides if $g$ acts level
transitively on $X^*$.
\end{lemma}

\begin{proof}
Let $g=\sigma^i(g_0,g_1)$, where $i\in\{0,1\}$. Then
\[\overline{g} = i+t\cdot(\overline{g_0} + \overline{g_1}).\]
Similar equations hold for all sections of $g$. Since $G$ is
generated by a finite automaton, $g$ has only finitely many
different sections, say $k$. Therefore we obtain a linear system of
$k$ equations over the $k$ variables $\{g_v, v\in X^*\}$. The
solution of this system expresses $\bar g$ as a rational function
$P(t)/Q(t)$, where $P$ an $Q$ are polynomials of degree not higher
than $k$. The element $g$ acts level transitively if and only if
$\bar{g} = \frac{1}{1-t}$.
\end{proof}

We often need to show that a given group of tree automorphisms is
level transitive. Here is a very convenient necessary and sufficient
condition for this in the case of a binary tree.

\begin{lemma}[Group transitivity, \cite{bondarenko-al:classification32-1}]\label{prop:transitivity}
A self-similar group of binary tree automorphisms is level
transitive if and only if it is infinite.
\end{lemma}

\begin{proof}
Let $G$ be a self-similar group acting on a binary tree.

If $G$ acts level transitively then $G$ must be infinite (since the
size of the levels is not bounded).

Assume now that the group $G$ is infinite.

We first prove that all level stabilizers $\Stg(n)$ are different.
Note that, since all level stabilizers have finite index in $G$ and
$G$ is infinite, all level stabilizers are infinite. In particular,
each contains a nontrivial element.

Let $n >0$ and $g\in\Stg(n-1)$ be an arbitrary nontrivial element.
Let $v=x_1\ldots x_k$ be a word of shortest length such that
$g(v)\neq v$. Since $g\in\Stg(n-1)$, we must have $k\geq n$. The
section $h=g_{x_1x_2\ldots x_{k-n}}$ is an element of $G$ by the
self-similarity of $G$. The minimality of the word $v$ implies that
$g\in\Stg(k-1)$, and therefore $h\in\Stg(n-1)$. On the other hand
$h$ acts nontrivially on $x_{k-n+1}\ldots x_k$ and we conclude that
$h \in \Stg(n-1) \setminus \Stg(n)$. Thus all level stabilizers are
different.

We now prove level transitivity by induction on the level.

The existence of elements in $\Stg(0) \setminus \Stg(1)$ shows that
$G$ acts transitively on level 1.

Assume that $G$ acts transitively on level $n$. Select an arbitrary
element $h\in\Stg(n)\setminus\Stg(n+1)$ and let $w=\in X^n$ be a
word of length $n$ such that $h(w1)=w0$.

Let $u$ be an arbitrary word of length $n$ and let $x$ be a letter
in $X=\{0,1\}$. We will prove that $ux$ is mapped to $w0$ by some
element of $G$, proving the transitivity of the action at level
$n+1$. By the inductive assumption there exists $f\in G$ such that
$f(u)=w$. If $f(ux)=w0$ we are done. Otherwise, $hf(ux)=h(w1)=w0$
and we are done again.
\end{proof}

Consider the infinitely iterated permutational wreath product
$\wr_{i\geq1}C_d$, consisting of the automorphisms of the $d$-ary
tree for which the activity at every vertex is a power of some fixed
cycle of length $d$. The last proof works, mutatis mutandis, for the
self-similar subgroups of $\wr_{i\geq1}C_d$ and may be easily
adapted in other situations.

The following lemma is used often when we want to prove that some
automaton group is not free.

\begin{lemma}\label{lem:freegroup}
If a self-similar group contains two nontrivial elements of the form
$(1, u), (v, 1)$, then the group is not free.
\end{lemma}
\begin{proof}
Suppose $a = (1, u), b = (v, 1)$ are two nontrivial elements of a
self-similar group $G$ and $G$ is free. Obviously $[a, b] = 1$,
hence $a$ and $b$ are powers of some element $x \in G$: $a=x^m$, $b
= x^n$. Then $a^n = b^m$, so $a^n = (1, u^n) = b^m = (v^m, 1)$. This
implies that $u^n = v^m = 1$, which is a contradiction, since $u$
and $v$ are nontrivial elements of a free group.
\end{proof}

In most case when the corresponding group is finite we do not offer
a full proof. In all such cases the proof can be easily done by
direct calculations. As an example, a detailed proof is given in the
case of the automaton [748].

We now proceeds to individual analysis of the properties of the
automaton groups in our classification.


\noindent\textbf{1}. Trivial group.

\noindent\textbf{730}. Klein Group $C_2\times C_2$. Wreath
recursion: $a=\s(a,a)$, $b=(a,a)$, $c=(a,a)$.

The claim follows from the relations $b=c$, $a^2=b^2=abab=1$.

\noindent$\mathbf{731}\cong\Z$. Wreath recursion: $a=\s(b,a)$,
$b=(a,a)$, $c=(a,a)$.

We have $c=b$ and $b=a^{-2}$. The states $a$ and $b$ form a $2$-state
automaton generating $\Z$ (see Theorem~\ref{thm:class22}).

\noindent$\mathbf{734}\cong G_{730}$. Klein Group $C_2\times C_2$.
Wreath recursion: $a=\s(b,b)$, $b=(a,a)$, $c=(a,a)$.

The claim follows from the relations $b=c$, $a^2=b^2=abab=1$.

\noindent$\mathbf{739}\cong C_2\ltimes(C_2\wr\Z)$. Wreath recursion:
$a=\sigma(a,a)$, $b=(b,a)$, $c=(a,a)$.

All generators have order $2$. The elements $u=acba=(1,ba)$ and
$v=bc=(ba,1)$ generate $\Z^2$. This is clear since $ba=\sigma(1,ba)$
is the adding machine and therefore has infinite order. Further, we
have $ac=\sigma$ and $\langle u,v\rangle$ is normal in $H=\langle
u,v,\sigma\rangle$, since $u^\sigma=v$ and $v^\sigma=u$. Thus $H
\cong C_2\ltimes(\Z\times\Z) =C_2\wr\Z$.

We have $G_{739}=\langle H,a\rangle$ and $H$ is normal in $G_{739}$,
since it has index 2. Moreover, $u^a=v^{-1}$, $v^a=u^{-1}$ and
$\sigma^a=\sigma$. Thus $G_{739}=C_2\ltimes(C_2\wr\Z)$, where the
action of $C_2$ on $H$ is specified above.

\noindent\textbf{740}. Wreath recursion: $a=\sigma(b,a)$, $b=(b,a)$,
$c=(a,a)$.

The states $a$, $b$ form a $2$-state automaton generating the
Lamplighter group (see Theorem~\ref{thm:class22}). Thus $G_{740}$
has exponential growth and is neither torsion nor contracting.

Since $c=(a,a)$ we obtain that $G_{740}$ can be embedded into the
wreath product $C_2\wr (\Z\wr\C_2)$. Thus $G_{740}$ is solvable.

\noindent\textbf{741}. Wreath recursion: $a=\sigma(c,a)$, $b=(b,a)$,
$c=(a,a)$.

The states $a$ and $c$ form a $2$-state automaton generating the
infinite cyclic group $\Z$ in which $c=a^{-2}$ (see
Theorem~\ref{thm:class22}).

Since $b = (b,a)$,  we see that $b$ has infinite order and that
$G_{741}$ is not contracting).

We have $c=a^-2$ and $b^{-1}a^{-3}b^{-1}ababa = 1$. Since $a$ and
$b$ do not commute the group is not free.

\noindent$\mathbf{743}\cong G_{739}\cong C_2\ltimes (C_2\wr \Z)$.
Wreath recursion: $a=\sigma(b,b)$, $b=(b,a)$, $c=(a,a)$.

All generators have order $2$. The elements $u=acba=(1,ba)$ and
$v=bc=(ba,1)$ generate $\Z^2$ because $ba=\sigma(ab,1)$ is conjugate
to the adding machine and has infinite order. Further, we have
$babc=\sigma$ and $\langle u,v\rangle$ is normal in $H=\langle u,v,
\sigma\rangle$ because $u^\sigma=v$ and $v^\sigma=u$. In other
words, $H\cong C_2\ltimes (\Z\times\Z)=C_2\wr \Z$.

Furthermore, $G_{743}=\langle H,a\rangle$ and $H$ is normal in
$G_{743}$ because $u^a=v^{-1}$, $v^a=u^{-1}$ and $\sigma^a=\sigma$.
Thus $G_{743}=C_2\ltimes (C_2\wr \Z)$, where the action of $C_2$ on
$H$ is specified above and coincides with the one in $G_{739}$.
Therefore $G_{743}\cong G_{739}$.

\noindent\textbf{744}. Wreath recursion: $a=\sigma(c,b)$, $b=(b,a)$,
$c=(a,a)$.

Since $(a^{-1}c)^2=(c^{-1}ab^{-1}a,b^{-1}ac^{-1}a)$ and
$c^{-1}ab^{-1}a = ((c^{-1}ab^{-1}a)^{-1},a^{-1}c)$, the element
$(a^{-1}c)^2$ fixes the vertex $01$ and its section at this vertex
is equal to $a^{-1}c$. Hence, $a^{-1}c$ has infinite order.

The element $c^{-1}ab^{-1}a$ also has infinite order, fixes the
vertex $00$ and its section at this vertex is equal to
$c^{-1}ab^{-1}a$. Therefore $G_{744}$ is not contracting.

We have $b^{-1}c^{-1}ba^{-1}ca = (1, a^{-1}c^{-1}ac)$,
$ab^{-1}c^{-1}ba^{-1}c = (ca^{-1}c^{-1}a, 1)$, hence by
Lemma~\ref{lem:freegroup} the group is not free.

\noindent$\mathbf{747}\cong G_{739}\cong C_2\ltimes (C_2\wr\Z)$.
Wreath recursion: $a=\sigma(c,c)$, $b=(b,a)$, $c=(a,a)$.

All generators have order $2$ and $a$ commutes with $c$. Conjugating
this group by the automorphism $\gamma=(\gamma,c\gamma)$ yields an
isomorphic group generated by automaton $a'=\sigma$, $b'=(b',a')$
and $c'=(a',a')$. On the other hand we obtain the same automaton
after conjugating $G_{739}$ by $\mu=(\mu,a\mu)$ (here $a$ denotes
the generator of $G_{739}$).

\noindent$\mathbf{748}\cong D_4\times C_2$. Wreath recursion:
$a=\sigma(a,a)$, $b=(c,a)$, $c=(a,a)$.

Since $a$ is nontrivial and $b$ and $c$ have $a$ as a section, none
of the generators is trivial. All generators have order 2. Indeed,
we have $a^2=(a^2,a^2)$, $b^2=(c^2,a^2)$, $c^2=(a^2,a^2)$, showing
that $a^2$, $b^2$ and $c^2$ generate a self-similar group in which
no element is active. Therefore $a^2=b^2=c^2=1$. Since $ac = \sigma$
we have that $(ac)^2 = 1$. Therefore $a$ and $c$ commute. Since
$(bc)^2=((ca)^2,1)=1$, we see that $b$ and $c$ also commute.
Further, the relations $(ab)^2 = (ac,1) = (\s,1) \neq 1$ and
$(ab)^4=1$ show that $a$ and $b$ generate the dihedral group $D_4$.
It remains to be shown that $c \not \in \langle a,b \rangle$.
Clearly $c$ could only be equal to one of the four elements $1$,
$b$, $aba$, and $abab$ in $D_4$ that stabilize level 1. However, $c$
is nontrivial, differs from $b$ at $0$ (the section $b|_0 = c$ is
not active, while $c|_0=a$ is active), differs from $aba$ at $1$
(the section $(aba)|_1 = aca$ is not active, while $c|_1=a$ is
active), and differs from $abab$ at 1 (the section of $abab$ at 1 is
trivial). This completes the proof.

\noindent\textbf{749}. Wreath recursion: $a=\sigma(b,a)$, $b=(c,a)$,
$c=(a,a)$.

The element $(a^{-1}c)^4$ stabilizes the vertex $000$ and its
section at this vertex is equal to $a^{-1}c$. Hence, $a^{-1}c$ has
infinite order.

We have $ac^{-1}=\sigma(ba^{-1}, 1)$, $ba^{-1}=\sigma(1, cb^{-1})$,
$cb^{-1}=(ac^{-1}, 1)$, Thus the subgroup generated by these
elements is isomorphic to $IMG(1-\frac{1}{z^2})$
(see~\cite{bartholdi_n:rabbit}).

We have $c^{-1}b = (a^{-1}c, 1)$, $ac^{-1}ba^{-1} = (1, ca^{-1})$.
Thus, by Lemma~\ref{lem:freegroup} the group is not free.

\noindent$\mathbf{748}\cong G_{848}\cong C_2\wr \Z$. Wreath
recursion:$a=\sigma(c,a)$, $b=(c,a)$, $c=(a,a)$.

It is proven below that $G_{848}\cong G_{2190}$ and for $G_{2190}$
we have $a=\sigma(c,a)$, $b=\sigma(a,a)$, $c=(a,a)$. Therefore
$G_{2190}=\langle a,b,c\rangle=\langle
a,c,c^{-1}b=\sigma\rangle=\langle
a=(c,a)\sigma,c=(a,a),a\sigma=(c,a)\rangle=G_{750}$.

\noindent\textbf{752}. Wreath recursion: $a=\sigma(b,b)$, $b=(c,a)$,
$c=(a,a)$.

The group $G_{752}$ is a contracting group with nucleus consisting
of $41$ elements. It is a virtually abelian group, containing $\Z^3$
as a subgroup of index $4$.

All generators have order 2.

Let $x=ca$, $y=babc$, and $K=\langle x,y \rangle$. Since
$xy=((cbab)^{ca},abcb)=((y^{-1})^x,abcb)$ and
$yx=(cbab,abcb)=(y^{-1},abcb)$ the elements $x$ and $y$ commute.
Conjugating by $\gamma=(\gamma,bc\gamma)$ yields the self-similar
copy $K'$ of $K$ generated by $x' = \sigma((y')^{-1},(x')^{-1})$ and
$y' = \sigma((y')^{-1}x',1)$, where $x' = x^\gamma$ and
$y'=y^\gamma$. Since $(x')^2 =
((x')^{-1}(y')^{-1},(y')^{-1}(x')^{-1})$ and $(y')^2 =
((y')^{-1}x',(y')^{-1}x')$, the virtual endomorphism of $K'$ is
given by
\[
 A=\left(\begin{array}{cc}
               -\frac12&\frac12\\
               -\frac12&-\frac12\\
 \end{array} \right).
\]
The eigenvalues $\lambda=-\frac12\pm\frac12i$ of this matrix are not
algebraic integers, and therefore, by the results
in~\cite{nekrash_s:12endomorph}, the group $K' \cong K$ is free
abelian of rank $2$.

Let $H=\langle ba,cb \rangle$. The index of $\St_H(1)$ in $G$ is 4,
since the index of $\St_H(1)$ in $H$ is 2 and the index of $H$ in
$G$ is 2 (the generators have order 2). We have $\St_H(1) = \langle
cb, cb^{ba},(ba)^2 \rangle$. If we conjugate the generators of
$\St_H)(1)$ by $g=(1,b)$, we obtain
\begin{alignat*}{3}
g_1 &= \left(cb\right)^g          &&= (x^{-1},&&1),\\
g_2 &= \left((cb)^{ba}\right)^{g} &&=(1,&&x),\\
g_3 &= \left((ba)^2\right)^{g}    &&=(y^{-1},&&y).
\end{alignat*}
Therefore, $g_1$, $g_2$, and $g_3$ commute. If
$g_1^{n_1}g_2^{n_2}g_3^{n_3}=1$, then we must have
$x^{-n_1}y^{-n_3}=x^{n_2}y^{n_3}=1$. Since $K$ is free abelian, this
implies $n_1=n_2=n_3=0$. Thus, $\St_H(1)$ is a free abelian group of
rank $3$.

\noindent\textbf{753}. Wreath recursion: $a=\sigma(c,b)$, $b=(c,a)$,
$c=(a,a)$.

Since $ab^{-1}=\sigma(1,ba^{-1})$, this element is conjugate to the
adding machine.

For a word $w$ in $w\in\{a^{\pm1},b^{\pm1},c^{\pm1}\}^*$, let
$|w|_a$, $|w|_b$ and $|w|_c$ denote the sum of the exponents of $a$,
$b$ and $c$ in $w$. Let $w$ represents the element $g \in G$. If
$|w|_a$ and $|w|_b$ are odd, then $g$ acts transitively on the first
level, and $g^2|_0$ is represented by a word $w_0$, which is the
product (in some order) of all first level sections of all
generators appearing in $w$. Hence, $|w_0|_a=|w|_b+2|w|_c$ and
$|w_0|_b=|w|_a$ are odd again. Therefore, similarly to
Lemma~\ref{nontors}, any such element has infinite order.

In particular $c^2ba$ has infinite order. Since
$a^4=(caca,a^4,acac,a^4)$ and $caca=(baca,c^2ba,bac^2,caba)$, the
element $a^4$ has infinite order (and so does $a$). Since $a^4$
fixes the vertex $01$ and its section at that vertex is equal to
$a^4$, the group $G_{753}$ is not contracting.

We have $cb^{-1} = (ac^{-1}, 1), acb^{-1}a^{-1} = (1, bac^{-1}b^{-1})$,
hence by Lemma~\ref{lem:freegroup} the group is not free.

\noindent$\mathbf{756}\cong G_{748}\cong D_4\times C_2$. Wreath
recursion: $a=\sigma(c,c)$, $b=(c,a)$, $c=(a,a)$.

All generators have order 2. The generator $c$ commutes with both
$a$ and $b$. Since $(ab)^2 = (ca,ca)$ the order of $ca$ is 4 and the
group is isomorphic to $D_4\times C_2$

\noindent$\mathbf{766}\cong G_{730}$. Klein Group $C_2\times C_2$.
Wreath recursion: $a=\sigma(a,a)$, $b=(b,b)$, $c=(a,a)$.

The state $b$ is trivial. The states $a$ and $c$ form a $2$-state
automaton generating $C_2 \times C_2$ (see
Theorem~\ref{thm:class22}).

\noindent$\mathbf{767}\cong G_{731}\cong \Z$. Wreath recursion:
$a=\sigma(1,a)$, $b=(b,b)$, $c=(a,a)=a^2$.

The state $b$ is trivial. The automorphism $a$ is the binary adding
machine.

\noindent$\mathbf{768}\cong G_{731}\cong \Z$. Wreath recursion:
$a=\sigma(c,a)$, $b=(b,b)$, $c=(a,a)$.

The states $a$ and $c$ form a $2$-state automaton generating $\Z$
(see Theorem~\ref{thm:class22}) in which $c=a^{-2}$.

\noindent$\mathbf{770}\cong G_{730}$. Klein Group $C_2\times C_2$.
Wreath recursion: $a=\sigma(b,b)$, $b=(b,b)$, $c=(a,a)$.

The state $b$ is trivial. The states $a$ and $c$ form a $2$-state
automaton generating $C_2 \times C_2$ (see
Theorem~\ref{thm:class22}).

\noindent$\mathbf{771}\cong \Z^2$. Wreath recursion:
$a=\sigma(c,b)$, $b=(b,b)$, $c=(a,a)$.

The group $G_{771}$ is finitely generated, abelian, and
self-replicating. Therefore, it is
free~\cite{nekrash_s:12endomorph}. Since $b=1$ the rank is 1 or 2.
We prove that the rank is 2, by showing that $c^n \neq a^m$, unless
$n=m=0$. By way of contradiction, let $c^n=a^m$ for some integer $n$
and $m$ and choose such integers with minimal $|n|+|m|$. Since $c^n$
stabilizes level 1, $m$ must be even and we have $(a^n,a^n)=c^n =
a^m = (c^{m/2},c^{m/2})$, implying $a^n=c^{m/2}$. By the minimality
assumption, $m$ must be 0, which then implies that $n$ must be 0 as
well.

\noindent$\mathbf{774}\cong G_{730}$. Klein Group $C_2\times C_2$.
Wreath recursion: $a=\sigma(c,c)$, $b=(b,b)$, $c=(a,a)$.

The state $b$ is trivial. The states $a$ and $c$ form a $2$-state
automaton generating $C_2 \times C_2$ (see
Theorem~\ref{thm:class22}).

\noindent$\mathbf{775}\cong C_2\ltimes
IMG\left(\bigl(\frac{z-1}{z+1}\bigr)^2\right)$. Wreath recursion:
$a=\sigma(a,a)$, $b=(c,b)$, $c=(a,a)$.

All generators have order 2. Further, $ac=ca=\sigma(1,1)$ and
$ba=\sigma(ba,ca)$. Hence, for the subgroup $H=\langle ba,ca\rangle
\cong G_{2853}\cong IMG\left(\bigl(\frac{z-1}{z+1}\bigr)^2\right)$.

Since the generators have order 2, $H$ is normal subgroup of index 2
in $G_{775}$. Moreover $(ba)^a=(ba)^{-1}$ and $(ca)^a=ca$. Therefore
$G\cong C_2\ltimes H$, where $C_2$ is generated by $a$ and the
action of $a$ on $H$ is given above.

Conjugating the generators by $g = \sigma(g,g)$ we obtain the wreath
recursion
\[
 a' = \s (a',a'), \qquad b' = (b',c'), \qquad c' = (a',a'),
\]
where $a'=a^g$, $b'=b^g$ and $c'=c^g$.  This is the wreath recursion
defining $G_{793}$. Denote $G_{793}$ by $G$ and its generators by
$a$, $b$, and $c$ (we continue working only with $G_{793}$). Thus
\[
 a = \s (a,a), \qquad b = (b,c), \qquad c = (a,a).
\]

The generators have order 2. Moreover $ac=ca$ and $\langle a,c
\rangle = C_2 \times C_2$ is the Klein group. Denote $A =\langle
a,c\rangle$.

The element $x=ba$ has infinite order, since $x^2$ fixes $00$, and
has itself as a section at 00. Note that
\[ x = ba = (b,c)\s(a,a) = \s(ca,ba) = \s(\s,x). \]
and, therefore, $x^2=(x \s, \s x) = (x,\s,\s,x )$.

\begin{prop}\label{Htorsionfree}
The subgroup $H=\langle x,y \rangle$ of $G$, where $x=ba$ and
$y=cabc$ is torsion free.
\end{prop}

\begin{proof}
The first level decompositions of $x^{\pm1}$ and $y^{\pm1}$ and the
second level decompositions of $x$ and $y$ are given by
\begin{align*}
 &x = \s(\s,x) \\
 &y = cabc = \s aaba \s = \s ba \s = x^\s = \s(x,\s) \\
 &x^{-1} = \s(x^{-1},\s) \\
 &y^{-1} = \s(\s,x^{-1})\\
 &x = \s(\s(1,1),\s(\s,x))= \mu(1,1,\s,x) \\
 &y = x^\s = \mu(\s,x,1,1),
\end{align*}
where $\mu=\s(\s,\s)$ permutes the first two levels of the tree as
$00 \leftrightarrow 11, \ 10 \leftrightarrow 01$. We encode this as
the permutation $\mu=(03)(12)$.

For a word $w$ over $\{x^{\pm1},\s\}$, denote by $\#_x(w)$ and
$\#_\s(w)$ the total number of appearances of $x$ and $x^{-1}$ and
the number of appearances of $\s$ in $w$, respectively.

Note that $x$ and $x^{-1}$ act as the permutation $(03)(12)$ on the
second level, and $\s$ acts as the permutation $(02)(13)$. These
permutations have order 2, commute, and their product is $(01)(23)$,
which is not trivial. Thus, a tree automorphisms represented by a
word $w$ over $\{x^{\pm1},\s\}$ cannot be trivial unless both
$\#_x(w)$ and $\#_\s(w)$ are even.

Let $g$ be an element of $H$ that can be written as $g=z_1 z_2 \dots
z_n$, for some $z_i \in \{x^{\pm1},y^{\pm1}\}$, $i=1,\dots,n$.

If $n$ is odd, the element $g$ cannot have order 2. By way of
contradiction assume otherwise. For $z$ in $\{x^{\pm1},y^{\pm1}\}$
denote $z'=\s z$. Thus, for instance $x'=(\s,x)$ and $y'=(x,\s)$.
Note that
\[ g^2= (z_1 z_2 \dots z_n)^2 = (z_1')^\s z_2' (z_3')^\s z_4' \dots (z_n')^\s z_1'
 (z_2')^\s\dots z_n' =(w_0,w_1),\]
where the words $w_i$ over $\{x^{\pm1},\s\}$ are such that
\begin{equation}\label{odd}
 \#_x(w_i) = \#_\s(w_i) = n,
\end{equation}
for $i=1,2$. The last claim holds because exactly one of $z_i'$ and
$(z_i')^\s$ contributes $x^{\pm1}$ to $w_0$ and $\s$ to $w_1$,
respectively, while the other contributes the same letters to $w_1$
and $w_0$, respectively. Since $n$ is odd,~(\ref{odd}) shows that
neither $w_0$ nor $w_1$ can be 1 and therefore $g^2$ cannot be 1.

Assume that $H$ contains an element of finite order. In particular,
this implies that $H$ must contain an element of order 2. Let $g=z_1
z_2 \dots z_n$ be such an element of the shortest possible length,
where $z_i \in \{x^{\pm1},y^{\pm1}\}$, $i=1,\dots,n$.

Note that $n$ must be even. Therefore,
\[ g = z_1 z_2 \dots z_n = (z_1')^\s z_2' \dots (z_{n-1}')^\s z_n' =(w_0,w_1),\]
where $w_0$ and $w_1$ are words over $\{x^{\pm1},\s\}$. Moreover, as
elements in $H$, the orders of $w_0$ and $w_1$ divide 2 and the
order of at least one of them is 2. We claim that
\begin{equation}\label{even}
 \#_x(w_0) \equiv \#_\s(w_0) \equiv  \#_x(w_1) \equiv \#_\s(w_0) \mod 2.
\end{equation}
The congruence $\#_x(w_i) \equiv \#_\s(w_i) \mod 2$ holds because
$\#_x(w_i) + \#_\s(w_i) = n$ is even. For the other congruences,
observe that whenever $z_i'$ or $(z_i')^\s$ contributes $x^{\pm1}$
or $\s$ to $w_0$, respectively, it contributes $\s$ or $x^{\pm1}$ to
$w_1$, respectively. Therefore $\#_x(w_0) = \#_\s(w_1)$ and
$\#_\s(w_0)=\#_x(w_1)$.

If the numbers in~(\ref{even}) are even, then $w_0$ and $w_1$
represent elements in $H$ and can be rewritten as words over
$\{x^{\pm1},y^{\pm1}\}$ of lengths at most $\#_x(w_0) = n -
\#_\s(w_0)$ and $\#_x(w_1) = n - \#_\s(w_1)$, respectively. If both
of these lengths are shorter than $n$ then none of them can
represent an element of order 2 in $H$. Otherwise, one of the words
$w_i$ is a power of $x$ and the other is trivial. Sice $x$ has
infinite order this shows that $g$ cannot have order 2.

If the numbers in~(\ref{even}) are odd, then, for $i=1,2$, $w_i$ can
be rewritten as $\s u_i$, where $u_i$ are words of odd length over
$\{x^{\pm1},y^{\pm1}\}$. Let $w_0=\s t_1 \dots t_m$, where $m$ is
odd, and $t_j$ are letters in $\{x^{\pm1},y^{\pm1}\}$,
$j=1,\dots,m$. We have
\[ w_0 = t_1'(t_2')^\s \dots (t_{m-1}')^\s t_m' = (w_{00},w_{01}), \]
where $w_{00}$ and $w_{01}$ are words of odd length $m$ over
$\{x^{\pm1},\s\}$. Moreover, exactly one of the words $w_{00}$ and
$w_{01}$ has even number of $\s$'s and this word can be rewritten as
a word over $\{x^{\pm1},y^{\pm1}\}$ of odd length. However, an
element in $H$ represented by such a word cannot have order dividing
2. This completes the proof.
\end{proof}

Since
\begin{alignat*}{3}
  x^a &= abaa = ab = x^{-1}, \qquad  &&y^a &&= acabca = cbac = y^{-1}, \\
  x^b &= bbab = ab = x^{-1}, \qquad  &&y^b &&= bcabcb = bacbacab = xy^{-1}x^{-1},\\
  x^c &= cbac = y^{-1},      \qquad  &&y^c &&= ccabcc = ab = x^{-1},
\end{alignat*}
we see that $H$ is the normal closure of $x$ in $G$. Further,
$G=\{x,y,a,c\}$ and $G = AH$. It follows from
Proposition~\ref{Htorsionfree} that $A \cap H = 1$ (since $A$ is
finite) and therefore $G = A \ltimes H$.

\begin{proposition}
The group $G$ is a regular, weakly branch group, branching over
$H''$.
\end{proposition}

\begin{proof}
The group $G$ is infinite self-similar group acting on a binary
three. Therefore it is level transitive by
Lemma~\ref{prop:transitivity}.

Since
\begin{align*}
x^2        &= (x, \s, \s, x) \\
y^{-1}x^2y &= (y, x^{-1} \s x, \s, x)
\end{align*}
we have that
\[ H'' \times \langle \s, x^{-1} \s x \rangle'' \times \langle \s
\rangle'' \times \langle x \rangle'' \preceq H''. \]

On the other hand, $\langle \s, x^{-1} \s x \rangle$ is metabelian
(in fact dihedral, since the generators have order 2) and $\langle
\s \rangle$ and $\langle x \rangle$ are abelian (cyclic). Therefore
\[ H'' \times 1 \times 1 \times 1 \preceq H''. \]

The group $H''$ is normal in $G$, since it is characteristic in the
normal subgroup $H$. Finally, $H''$ is not trivial. For instance it
is easy to show that $[[x,y],[x,y^{-1}]] \neq 1$
(see~\cite{bondarenko-al:classification32-2}).
\end{proof}

\noindent\textbf{776}. Wreath recursion: $a=\sigma(b,a)$, $b=(c,b)$,
$c=(a,a)$.

The element $(b^{-1}a)^4$ stabilizes the vertex $00$ and its section
at this vertex is equal to $(b^{-1}a)^{-1}$. Hence, $b^{-1}a$ has
infinite order. Furthermore, by Lemma~\ref{nontors} $ab$ has
infinite order, which yields that $a$,$c$ and $b$ also have infinite
order, because $a^2=(ab,ba)$. Since $b^n=(c^n,b^n)$ we obtain that
$b^n$ belong to the nucleus for all $n\geq1$. Thus $G_{776}$ is not
contracting.

We have $a^{-1}ba^{-1}c = (1, b^{-1}c), ba^{-1}ca^{-1} = (cb^{-1}, 1)$,
hence by Lemma~\ref{lem:freegroup} the group is not free.

\noindent\textbf{777}. Wreath recursion: $a=\sigma(c,a)$, $b=(c,b)$,
$c=(a,a)$.

The states $a$, $c$ form the $2$-state automaton generating $\Z$
(see Theorem~\ref{thm:class22}). So the group is not torsion and
$G_{777}=\langle a,b\rangle$. Since $c$ has infinite order, so has
$b$. Therefore the relation $b^n=(c^n,b^n)$ implies that $b^n$
belong to the nucleus for all $n\geq1$. Thus $G_{777}$ is not
contracting.

Also we have $ab^{-1}=\sigma(1,ab^{-1})$ is the adding machine.
Since $a^{-3}=\sigma(1,a^3)$ elements $ab^{-1}$ and $a^{-3}$
generate the Brunner-Sidki-Vierra group
(see~\cite{brunner_sv:justnonsolv}).

\noindent\textbf{779}. Wreath recursion: $a=\sigma(b,b)$, $b=(c,b)$,
$c=(a,a)$.

The element $(ab^{-1})^2$ stabilizes the vertex $01$ and its section
at this vertex is equal to $(ab^{-1})^{-1}$. Hence, $ab^{-1}$ has
infinite order.

\noindent\textbf{780}. Wreath recursion: $a=\sigma(c,b)$, $b=(c,b)$,
$c=(a,a)$.

The element $(c^{-1}a)^2$ stabilizes the vertex $00$ and its section
at this vertex is equal to $c^{-1}a$. Hence, $c^{-1}a$ has infinite
order. Since $[c,a]\bigl|_{100}=(c^{-1}a)^a$ and $100$ is fixed
under the action of $[c,a]$ we obtain that $[c,a]$ also has infinite
order. Finally, $[c,a]$ stabilizes the vertex $00$ and its section
at this vertex is $[c,a]$. Therefore $G_{780}$ is not contracting.

\noindent$\mathbf{783}\cong G_{775}\cong C_2\ltimes
IMG\left(\bigl(\frac{z-1}{z+1}\bigr)^2\right)$. Wreath recursion:
$a=\sigma(c,c)$, $b=(c,b)$, $c=(a,a)$.

All generators have order $2$ and $ac=ca$. If we conjugate the
generators of this group by the automorphism
$\gamma=(c\gamma,\gamma)$ we obtain the wreath recursion
\[ a'=\sigma(1,1), \qquad b'=(c',b'), \qquad c'=(a',a'), \]
where $a'=a^\gamma$, $b'=b^\gamma$, and $c'=c^\gamma$. The same
wreath recursion is obtained after conjugating $G_{775}$ by
$\mu=(a\mu,\mu)$ (where $a$ denotes the generator of $G_{775}$).

Since $bca=\sigma(bca,a)$, $G_{783}=\langle acb,a,c\rangle\cong
G_{2205}$.

\noindent$\mathbf{802}\cong C_2\times C_2\times C_2$. Wreath
recursion: $a=\sigma(a,a)$, $b=(c,c)$, $c=(a,a)$.

Direct calculation.

\noindent$\mathbf{803}\cong G_{771}\cong \mathbb Z^2$. Wreath
recursion: $a=\sigma(b,a)$, $b=(c,c)$, $c=(a,a)$.

The group $G_{771}$ is finitely generated, abelian, and
self-replicating. Therefore, it is free
abelian~\cite{nekrash_s:12endomorph}. Let $\phi:\mathop{\rm
Stab}\nolimits_{G_{803}}(1)\to G_{803}$ be the
$\frac12$-endomorphism associated to the vertex $0$, given by
$\phi(g)=h$, provided $g=(h,*)$. The matrix of the linear map
$\mathbb C^3\to\mathbb C^3$ induced by $\phi$ with to the basis
corresponding to the triple $\{a,b,c\}$ is given by
\[
 A=\left(
 \begin{array}{ccc}
  \frac12 & 0 & 1 \\
  \frac12 & 0 & 0 \\
  0 & 1 & 0
\end{array}
\right).
\]
The eigenvalues are $\lambda_1=1$,
$\lambda_2=-\frac14-\frac14i\sqrt7$ and
$\lambda_3=-\frac14+\frac14i\sqrt7$. Let $v_i$, $i=1,2,3$, be
eigenvectors corresponding to the eigenvalues $\lambda_i$,
$i=1,2,3$. Note that $v_1$ may be selected to be equal to
$v_1=(2,1,1)$. This shows that $a^2bc=1$ in $G_{803}$ and the rank
of $G_{803}=\langle a,c \rangle$ is at most 2. We will prove that
$a^{2m}c^n\neq1$ (except when $m=n=0$) by proving that iterations of
the action of $A$ eventually push the vector $v=(2m,0,n)$ out of the
set $D=\{(2i,j,k), i,j,k\in\mathbb Z\}$ corresponding to the first
level stabilizer.

Let $v=a_1v_1+a_2v_2+a_3v_3$. The vector $v$ is not a scalar
multiple of $v_1$. Therefore either $a_2 \neq 0$ or $a_3 \neq 0$.
Since $|\lambda_2|=|\lambda_3|<1$, we have
$A^t(v)=a_1v_1+\lambda_2^ta_2v_2+\lambda_3^ta_3v_3\to a_1v_1$, as
$t\to\infty$. Note that, since $a_2 \neq 0$ or $a_3 \neq 0$,
$A^t(v)$is never equal to $a_1v_1$. Choose a neighborhood $U$ of
$a_1v_1$ that does not contain vectors from $D$, except possibly the
vector $a_1v_1$. For $t$ large enough $t$, the vector $A^t(v)$ is in
$U$ and is therefore outside of $D$.

Thus the rank of $G_{803}$ is 2.

\noindent$\mathbf{804}\cong G_{731}\cong \Z$. Wreath recursion:
$a=\sigma(c,a)$, $b=(c,c)$, $c=(a,a)$.

Indeed, the states $a$ and $c$ form a 2-state automaton generating
the cyclic group $\mathbb Z$ (see Theorem~\ref{thm:class22}). Since
$b=a^4$ we are done.

\noindent$\mathbf{806}\cong G_{802}\cong C_2\times C_2\times C_2$.
Wreath recursion: $a=\sigma(b,b)$, $b=(c,c)$, $c=(a,a)$.

Direct calculation.

\noindent$\mathbf{807}\cong G_{771}\cong \Z^2$. Wreath recursion:
$a=\sigma(c,b)$, $b=(c,c)$, $c=(a,a)$.

The same arguments as in the case of $G_{771}$ show that $G_{807}$
is free abelian. It has a relation $c^2ba^2=1$ and, hence, it has
either rank $1$ or rank $2$. Analogically to $G_{803}$ we consider a
$\frac12$-endomorphism $\phi:\mathop{\rm
Stab}\nolimits_{G_{807}}(1)\to G_{807}$, and a linear map $A:\mathbb
C^3\to\mathbb C^3$ induced by $\phi$. It has the following matrix
representation with respect to the basis corresponding to the triple
$\{a,b,c\}$:

\[A=\left(
\begin{array}{ccc}
  0 & 0 & 1 \\
  \frac{1}{2} & 0 & 0 \\
  \frac{1}{2} & 1 & 0
\end{array}
\right).\]

Its characteristic polynomial
$\chi_A(\lambda)=-\lambda^3+\frac12\lambda+\frac12$ has three
distinct complex roots $\lambda_1=1$, $\lambda_2=-\frac12-\frac12i$
and $\lambda_3=-\frac12+\frac12i$. Analogically for $v=(2m,0,n)$ we
get that $A^t(v)$ will be pushed out from the domain corresponding
to $\mathop{\rm Stab}\nolimits_{G_{807}}(1)$. Thus $c^na^{2m}\neq 1$
in $G_{807}$ and $G_{807}\cong\mathbb Z^2$.

\noindent$\mathbf{810}\cong G_{802}\cong C_2\times C_2\times C_2$.
Wreath recursion: $a=\sigma(c,c)$, $b=(c,c)$, $c=(a,a)$.

Direct calculation.

\noindent$\mathbf{820}\cong D_\infty$. Wreath recursion:
$a=\sigma(a,a)$, $b=(b,a)$, $c=(b,a)$.

The states $a$ and $b$ form a $2$-state automaton generating
$D_\infty$ (see Theorem~\ref{thm:class22}) and $c=b$.

\noindent\textbf{821}. Lamplighter group $\Z \wr C_2$. Wreath
recursion: $a=\sigma(b,a)$, $b=(b,a)$, $c=(b,a)$.

The states $a$ and $b$ form a $2$-state automaton generating the
Lamplighter group (see Theorem~\ref{thm:class22}) and $c=b$.

\noindent$\mathbf{824}\cong G_{820}\cong D_\infty$. Wreath
recursion: $a=\sigma(a,a)$, $b=(b,a)$, $c=(b,a)$.

The states $a$ and $b$ form a $2$-state automaton generating
$D_\infty$ (see Theorem~\ref{thm:class22}) and $c=b$.

\noindent$\mathbf{838}\cong C_2\ltimes \langle s,t\ \bigl|\
s^2=t^2\rangle$. Wreath recursion: $a=\sigma(a,a)$, $b=\sigma(a,b)$,
$c=(b,a)$.

All generators have order $2$. Consider the subgroup $H=\langle
ba=\sigma(ba,1), ca=\sigma(1,ab)\rangle\cong G_{2860}=\langle s,t\
\bigl|\ s^2=t^2\rangle$. This subgroup is normal in $G_{838}$
because the generators have order $2$. Since $G_{838}=\langle
H,a\rangle$, it has a structure of a semidirect product $\langle
a\rangle\ltimes H=C_2\ltimes \langle s,t\ \bigl|\ s^2=t^2\rangle$
with the action of $a$ on $H$ as $(ba)^b=(ba)^{-1}$ and
$(ca)^b=(ca)^{-1}$.

\noindent$\mathbf{839}\cong G_{821}$. Lamplighter group $\Z \wr
C_2$. Wreath recursion: $a=\sigma(b,a)$, $b=(a,b)$, $c=(b,a)$.

The states $a$ and $b$ form a 2-state automaton generating the
Lamplighter group (see Theorem~\ref{thm:class22}). Since
$b^{-1}a=\sigma=ac^{-1}$, we see that $c=a^{-1}ba$ and $G=\langle
a,b\rangle$.

\noindent\textbf{840}. Wreath recursion: $a=\sigma(c,a)$, $b=(a,b)$,
$c=(b,a)$.

The element $(b^{-1}a)^2$ stabilizes the vertex $01$ and its section
at this vertex is equal to $b^{-1}a$. Hence, $b^{-1}a$ has infinite
order.

The element $(c^{-1}b)^2$ stabilizes the vertex $00$ and its section
at this vertex is equal to $(c^{-1}b)^{-1}$. Hence, $c^{-1}b$ has
infinite order. Since
$(b^{-1}a^{-1}b^{-1}cba)^2\bigl|_{00000000}=c^{-1}b$ and the vertex
$00000000$ is fixed under the action of $(b^{-1}a^{-1}b^{-1}cba)^2$
we obtain that $b^{-1}a^{-1}b^{-1}cba$ also has infinite order.
Finally, $b^{-1}a^{-1}b^{-1}cba$ stabilizes the vertex $0001$ and
has itself as a section at this vertex. Therefore $G_{840}$ is not
contracting.

We have $b^{-1}a^{-1}ca = (1, b^{-1}c^{-1}bc), ab^{-1}a^{-1}c = (cb^{-1}c^{-1}b, 1)$,
hence by Lemma~\ref{lem:freegroup} the group is not free.

\noindent$\mathbf{842}\cong G_{838}\cong C_2\ltimes \langle s,t\
\bigl|\ s^2=t^2\rangle$. Wreath recursion: $a=\sigma(b,b)$,
$b=\sigma(a,b)$, $c=(b,a)$.

All generators have order $2$. Consider the subgroup $H=\langle
u=ba=\sigma(1,ba)=\sigma(1,u^{-1}),
v=ca=\sigma(ab,1)=\sigma(u^{-1},1)\rangle$. Let us prove that
$H\cong W=\langle s,t\ \bigl|\  s^2=t^2\rangle$. Indeed, the
relation $u^2=v^2$ is satisfied, so $H$ is a homomorphic image of
$W$ with respect to the homomorphism induced by $s\mapsto u$ and
$t\mapsto v$. Each element of $W$ can be written in its normal form
$t^r(st)^ls^n$, $n\in\Z, l\geq0, r\in\{0,1\}$. It suffices to prove
that images of these words (except for the identity word, of course)
represent nonidentity elements in $H$.

We have $u^{2n}=(u^{-n},u^{-n})$, $u^{2n+1}=\sigma(a^{-n},a^{-n-1})$
for any integer $n$; $(uv)^l=(u^{2l},1)$ for any integer $l$. Thus
$$(uv)^lu^{2n}=(u^{-2l-n}, u^{-n})\neq1$$
in $G$ if $n\neq0$ or $l\neq0$ since $u$ has infinite order, as it
is conjugate to the adding machine.

Furthermore,
$$v(uv)^lu^{2n}=\sigma(u^{-2l-n-1}, u^{-n})\neq1,$$
$$(uv)^lu^{2n+1}=\sigma(u^{-n},u^{-2l-n-1})\neq1$$
since they act nontrivially on the first level of the tree.

Finally, $v(uv)^lu^{2n+1}=(u^{-2l-n-2}, u^{-n})=1$ if and only if
$n=0$ and $l=-1$, which is not the case, because $l$ must be
nonnegative. Thus $H\cong W$.

The subgroup $H$ is normal in $G_{842}$ because generators are of
order $2$. Since $G_{842}=\langle H,a\rangle$, it has a structure of
a semidirect product $\langle a\rangle\ltimes H=C_2\ltimes \langle
s,t\ \bigl|\ s^2=t^2\rangle$ with the action of $a$ on $H$ as
$(ba)^b=(ba)^{-1}$ and $(ca)^b=(ca)^{-1}$. Therefore it has the same
structure as $G_{838}$.

\noindent\textbf{843}. Wreath recursion: $a=\sigma(c,b)$, $b=(a,b)$,
$c=(b,a)$.

The element $c^{-1}a = \s(a^{-1}c,1)$ is a conjugate of the adding
machine. Therefore, it acts transitively on the level of the tree
and has infinite order.

Since $(c^{-1}ab^{-1}a)^2$ fixes the vertex $000$ and its section at
this vertex is equal to $c^{-1}a$, we obtain that $c^{-1}ab^{-1}a$
has infinite order. Since the element $c^{-1}ab^{-1}a$ fixes the
vertex $10$ and has itself as a section at this vertex, $G_{843}$ is
not contracting.

We have $c^{-1}a^{-1}ba = (1, a^{-1}c^{-1}ac)$, $ac^{-1}a^{-1}b =
(ca^{-1}c^{-1}a, 1)$, hence by Lemma~\ref{lem:freegroup} the group
is not free.

\noindent$\mathbf{846}\cong C_2\ast C_2\ast C_2$. Wreath recursion:
$a=\sigma(c,c)$, $b=(a,b)$, $c=(b,a)$.

The automaton [846] was studied during the Advanced Course on
Automata Groups in Bellaterra, Spain, in the summer of 2004 and is
since called the Bellaterra automaton. We present here a proof that
$G_{846}=C_2 \ast C_2 \ast C_2$, based on the concept of dual
automata. A different proof, still based on dual automata, is given
in~\cite{nekrash:self-similar}.

Let $\A = (Q, X, \pi, \tau)$ be a finite automaton. Its \emph{dual}
automaton, by definition, is $\mathcal{A}' = (X, Q, \pi', \tau')$,
where $\pi'(x, q) = \tau(q, x)$, and $\tau'(x, q) = \pi(q, x)$. Thus
the dual automaton is obtained by exchanging the roles of the states
and the alphabet (and the roles of the transition and output
function) in a given automaton. The notion od dual automata is not
new, but there is a recent renewed interest based on the new results
and applications
in~\cite{macedonska-n-s:comm,gl_mo:compl,bartholdi_s:bsolitar,vorobets:aleshin}.

If in addition to $\A$, both $\mathcal{A}'$ and
$(\mathcal{A}^{-1})'$ are invertible, the automaton $\A$ is called
\emph{fully invertible} (or \emph{bi-reversible}). Examples of such
automata are the automaton $2240$ generating a free group with three
generators~\cite{vorobets:aleshin}, Bellaterra automaton [846], and
various automata constructed in \cite{gl_mo:compl}, generating free
groups of various ranks.

We now consider the automaton [846] and its dual more closely. Since
the generators $a$, $b$, and $c$ have order $2$, in order to prove
that $G_{846} \cong C_2 \ast C_2 \ast C_2$ we need to show that no
word in $w \in R_n$, $n \geq 1$, is trivial in $G_{846}$, where
$R_n$ is the set of reduced words over $\{a, b, c\}$ of length $n$
(here a word is reduced if it does not contain $aa$, $bb$, or $cc$).
For every $n
> 0$, the set of words in $R_n$ that are nontrivial in $G_{846}$ is
nonempty, since the word $r_n=acbcbcb\cdots$ of length $n$ acts
nontrivially on level 1. If we prove that the dual automaton acts
transitively on the sets $R_n$, $n \geq 1$, this would mean that
$r_n$ is a section of every element of $G_{846}$ that can be
represented as a reduced word of length $n$. Therefore, every word
in $R_n$ would represent a nontrivial element in $G_{846}$ and our
proof would be complete.

The automaton dual to $846$ is the invertible automaton defined by
the wreath recursion
\begin{equation}
\begin{array}{rcr}\label{dual-bellaterra}
A & = & (acb)(B, A, A), \\
B & = & (ac)(A, B, B),
\end{array}\end{equation}
where the three coordinates in the recursion represent the sections
at $a$, $b$, and $c$, respectively. Denote $D=\langle A,B \rangle$.
The set $R = \bigcup_{n \ge 0}R_n$ of all reduced words over $\{a,
b, c\}$ is a subtree of the ternary tree $\{a, b, c\}^*$ and this
subtree $R$ is invariant under the action of $D$ (this is because
the set $\{aa,bb,cc\}$ is invariant under the action of $D$). The
structure of $R$ is as follows. The root of $R$ has three children
$a$, $b$ and $c$, each of which is a root of a binary tree. We want
to understand the actio of $D$ on the subtree $R$. It is given by
\begin{equation}
\begin{array}{rcr}
A & = & (acb)(B_a, A_b, A_c) \\
B & = & (ac)(A_a, B_b, B_c)
\end{array}
\end{equation}
where $A_a, A_b, A_c, B_a, B_b, B_c$ are automorphisms of the binary
trees hanging down from the vertices $a$, $b$ and $c$. After
identification of these three trees with the binary tree
$\{0,1\}^*$, the action of $A_a, A_b,\ldots,B_c$ is defined by
\begin{equation}\begin{array}{rcr}
A_a & = & (A_b, A_c), \\
A_b & = & \sigma(B_a, A_c), \\
A_c & = & \sigma(B_a, A_b), \\
B_a & = & \sigma(B_b, B_c), \\
B_b & = & \sigma(A_a, B_c), \\
B_c & = & \sigma(A_a, B_b). \\
\end{array}
\end{equation}

Using Lemma~\ref{alg_trans} one can verify that $B_b$ acts level
transitively on the binary tree. This is sufficient to show that $D$
acts transitively on $R$, since it acts transitively on the first
level, $B$ stabilizes the vertex $b$, and its section at $b$ is
$B_b$.

The fact that $G_{846}$ is not contracting follows now from the
result of Nekrashevych~\cite{nekrashevych:free_subgroups}, that a
contracting group can not have free subgroups. Alternatively, it is
sufficient to observe that $aba$ has infinite order, stabilizes the
vertex $01$ and has itself as a section at this vertex.

\noindent$\mathbf{847}\cong D_4$. Wreath recursion: $a=\sigma(a,a)$,
$b=(b,b)$, $c=(b,a)$.

The state $b$ is trivial. The states $a$ and $c$ form a 2-state
automaton generating $D_4$ (see Theorem~\ref{thm:class22}).

\noindent$\mathbf{848}\cong C_2\wr\Z$. Wreath recursion:
$a=\sigma(b,a)$, $b=(b,b)$, $c=(b,a)$.

The state $b$ is trivial and $a$ is the adding machine. Every
element $g\in G_{848}$ has the form $g=\sigma^i(a^n,a^m)$. On the
other hand, $c=(1,a), c^{ac^{-1}}=(a,1)$, so $\Stg(1)=\{(a^n,
a^m)\}\cong\mathbb Z^2$. Since $ac^{-1}=\sigma$ we see that $G\cong
C_2\wr\Z$.

\noindent\textbf{849}. Wreath recursion: $a=\sigma(c,a)$, b=(b,b),
$c=(b,a)$.

The state $b$ is trivial. The element $a^2c=(ac,ca^2)$ is nontrivial
because its section at $0$ is $ac$, and $ac$ acts nontrivially on
level 1. The automorphism $(a^2c)^2$ fixes the vertex $00$ and its
section at this vertex is equal to $a^2c$. Therefore $a^2c$ has
infinite order. Further, the section of $a^2c$ at $100$ coincides
with $a^2c$, implying that $G_{849}$ is not contracting.

The group $G_{849}$ is regular weakly branch group over its
commutator $G_{849}'$. This is clear since the group is
self-replicating and $[a^{-1},c]\cdot[c,a]=([a,c],1)$.

Conjugation of the generators of $G_{849}$ by
$\mu=\s(\mu,c^{-1}\mu)$ yields the wreath recursion
\[ x = \s(yx,1), \qquad y = (x,1), \]
where $x = a^\mu$ and $y = c^\mu$. Further, we have
\[ x = \s(yx,1), \qquad yx = \s(yx,x), \]
and the last wreath recursion coincides with the one defining the
automaton 2852. Therefore $G_{849} \cong G_{2852}$ (see $G_{2852}$
for more information on this group).

\noindent$\mathbf{851}\cong G_{847}\cong D_4$. Wreath recursion:
$a=\sigma(b,b)$, b=(b,b), $c=(b,a)$.

Direct calculation.

\noindent\textbf{852}. Basilica group $\B = IMG(z^2-1)$. Wreath
recursion: $a=\sigma(c,b)$, $b=(b,b)$, $c=(b,a)$.

This group was studied in~\cite{grigorch_z:basilica}, where it is
shown that $\mathcal B$ is not a sub-exponentially amenable group,
it does not contain free subgroups of rank $2$, and that the monoid
generated by $a$ and $b$ is free. Some spectral considerations are
provided in~\cite{grigorch_z:basilica_sp}. Bartholdi and Vir\'ag
showed in~\cite{bartholdi_v:amenab} that $\mathcal B$ is amenable,
distinguishing the Basilica group as the first example of an
amenable group that is not sub-exponentially amenable.

\noindent$\mathbf{855}\cong G_{847}\cong D_4$. Wreath recursion:
$a=\sigma(c,c)$, b=(b,b), $c=(b,a)$.

Direct calculation.

\noindent$\mathbf{856}\cong C_2\ltimes G_{2850}$. Wreath recursion:
$a=\sigma(a,a)$, $b=(c,b)$, $c=(b,a)$.

All generators have order $2$, hence $H=\langle ba,ca\rangle$ is
normal in $G_{856}$. Furthermore, $ba=\sigma(ba,ca)$,
$ca=\sigma(1,ba)$, and therefore $H=G_{2850}$. Thus $G_{856}=\langle
a\rangle\ltimes H\cong C_2\ltimes G_{2850}$, where
$(ba)^a=(ba)^{-1}$ and $(ca)^a=(ca)^{-1}$. The group is not
contracting since $G_{2850}$ is not contracting.

\noindent\textbf{857}. Wreath recursion: $a=\sigma(b,a)$, $b=(c,b)$,
$c=(b,a)$.

By using the approach used for $G_{875}$, we can show that the
forward orbit of $10^\infty$ under the action of $a$ is infinite,
and therefore $a$ has infinite order.

Since $c=(b,a)$ and $b=(c,b)$, both $b$ and $c$ have infinite order
and $G_{857}$ is not a contracting group.

\noindent\textbf{858}. Wreath recursion: $a=\sigma(c,a)$, $b=(c,b)$,
$c=(b,a)$.

The element $ab^{-1}=\sigma(1,ab^{-1})$ is the adding machine.

By using the approach used for $G_{875}$, we can show that the
forward orbit of $10^\infty$ under the action of $a$ is infinite,
and therefore $a$ has infinite order.

Since $c=(b,a)$ and $b=(c,b)$, both $b$ and $c$ have infinite order
and $G_{857}$ is not a contracting group.

We have $c^{-1}b^{-1}aba^{-1}b = (1, a^{-1}b^{-1}aca^{-1}b)$,
$a^{-1}c^{-1}b^{-1}aba^{-1}ba = (a^{-2}b^{-1}aca^{-1}ba, 1)$, hence
by Lemma~\ref{lem:freegroup} the group is not free.

\noindent\textbf{860}. Wreath recursion: $a=\sigma(b,b)$, $b=(c,b)$,
$c=(b,a)$.

The element $(ba^{-1})^2$ stabilizes the vertex $11$ and its section
at this vertex is equal to $(ba^{-1})^{-1}$. Hence, $ba^{-1}$ has
infinite order.

Furthermore, $bc^{-1}=(cb^{-1},ba^{-1})$ implies that the order of
$bc^{-1}$ is infinite. Since this element stabilizes vertex $00$ and
its section at this vertex is equal to $bc^{-1}$, all its powers
belong to the nucleus. Thus, $G_{860}$ is not contracting.

\noindent\textbf{861}. Wreath recursion: $a=\sigma(b,b)$, $b=(a,a)$,
$c=(b,a)$.

The element $a^{-1}c=\sigma(1,c^{-1}a)$ is conjugate to the adding
machine and has infinite order.

\noindent\textbf{864}. Wreath recursion: $a=\sigma(c,c)$, $b=(c,b)$,
$c=(b,a)$.

The element $(ab^{-1})^2$ stabilizes the vertex $11$ and its section
at this vertex is equal to $ab^{-1}$. Hence, $ab^{-1}$ has infinite
order.

Furthermore, $cb^{-1}=(bc^{-1},ab^{-1})$ implies that the order of
$cb^{-1}$ is infinite. Since this element stabilizes vertex $00$ and
its section at this vertex is equal to $cb^{-1}$, $G_{864}$ is not
contracting.

\noindent$\mathbf{865}\cong G_{820}\cong D_\infty$. Wreath
recursion: $a=\sigma(a,a), b=(a,c), c=(b,a)$.

All generators have order 2. Since $abac=(acab,1)$ and
$acab=(1,abac)$, we see that $c=aba$ and $G_{865}=\langle
a,b\rangle$. The section of $(ba)^2$ at the vertex $0$ is
$(ba)^{-1}$, so $ba$ has infinite order and $G_{865}\cong D_\infty$.

Note that the group is conjugate to $G_{932}$ by the automorphism
$\delta=(a\delta,\delta)$.

\noindent\textbf{866}. Wreath recursion: $a=\sigma(b,a)$, $b=(a,c)$,
$c=(b,a)$.

The element $(c^{-1}b)^2$ stabilizes the vertex $00$ and its section
at this vertex is equal to $c^{-1}b$, which is nontrivial. Hence,
$c^{-1}b$ has infinite order.

The element $(b^{-1}a)^2$ stabilizes the vertex $00$ and its section
at this vertex is equal to $b^{-1}a$. Hence, $b^{-1}a$ has infinite
order. Since $b^{-1}c^{-1}ba^{-1}ba\bigl|_{10}=(b^{-1}a)^b$ and
vertex $10$ is fixed under the action of $b^{-1}c^{-1}ba^{-1}ba$ we
obtain that $b^{-1}c^{-1}ba^{-1}ba$ also has infinite order.
Finally, $b^{-1}c^{-1}ba^{-1}ba$ stabilizes the vertex $00$ and has
itself as a section at this vertex. Therefore $G_{866}$ is not
contracting.

\noindent\textbf{869}. Wreath recursion: $a=\sigma(b,b)$, $b=(a,c)$,
$c=(b,a)$.

All generators have order $2$. By Lemma~\ref{nontors} $ab$ has
infinite order, which implies that $babcba$ also has infinite order,
because it fixes the vertex $000$ and its section at this vertex is
equal to $ab$. But $babcba$ fixes $10$ and has itself as a section
at this vertex. Thus, $G_{869}$ is not contracting.

\noindent\textbf{870}: Baumslag-Solitar group $BS(1,3)$. Wreath
recursion: $a=\sigma(c,b)$, $b=(a,c)$, $c=(b,a)$.

The automaton satisfies the conditions of Lemma~\ref{nontors}. In
particular $ab$ has infinite order. Since $bc=(ab,ca)$,
$a^2=(bc,cb)$, we obtain that $bc$ and $a$ have infinite order.
Since $b=(a,c)$, $b$ also has infinite order. Since $b$ has infinite
order, fixes the vertex $10$ and has itself as a section at this
vertex, $G_{870}$ is not contracting.

The element $\mu=b^{-1}a=\sigma(1,a^{-1}b)=\sigma(1,\mu^{-1})$ is
conjugate to the adding machine and therefore has infinite order.
Since $a^{-1}c=\sigma(1,c^{-1}a)$ we see that $a^{-1}c=\mu$.
Therefore $c=ab^{-1}a$ and $G_{870}=\langle a,b\rangle=\langle
\mu,b\rangle$.

We claim that $b^{-1}\mu b=\mu^3$. Since $c=ab^{-1}a$, we have
\[
 ab^{-1}ab^{-1}ab^{-1}a^{-1}b=(ba^{-1}bc^{-1}b^{-1}a,ca^{-1}ba^{-1})=
 (ba^{-1}ba^{-1}ba^{-1}b^{-1}a,1).
\]
But $ba^{-1}ba^{-1}ba^{-1}b^{-1}a$ is a conjugate of the inverse of
$ab^{-1}ab^{-1}ab^{-1}a^{-1}b$, which shows that
$ab^{-1}ab^{-1}ab^{-1}a^{-1}b=1$, and the last relation is
equivalent to $b^{-1}\mu b=\mu^3$.

Since $b$ and $\mu$ have infinite order, $G_{870}\cong BS(1,3)$.

See~\cite{bartholdi_s:bsolitar} for realizations of $BS(1,m)$ for
any value of $m$, $m \neq \pm 1$.

\noindent$\mathbf{874}\cong C_2\ltimes G_{2852}$. Wreath recursion:
$a=\sigma(a,a)$, $b=(b,c)$, $c=(b,a)$.

All the generators have order $2$, hence $H=\langle ba,ca\rangle$ is
normal in $G_{874}$. Furthermore, $ba=\sigma(ca,ba)$,
$ca=\sigma(1,ba)$, therefore $H=G_{2852}$. Thus $G_{874}=\langle
a\rangle\ltimes H\cong C_2\ltimes G_{2852}$, where
$(ba)^a=(ba)^{-1}$ and $(ca)^a=(ca)^{-1}$. In particular, $G_{874}$
is not contracting and has exponential growth.

\noindent\textbf{875}. Wreath recursion: $a=\sigma(b,a)$, $b=(b,c)$,
$c=(b,a)$.

The equalities
\[
 a(10^\infty) = 010^\infty, \qquad
 b(10^\infty) = 10^\infty, \qquad
 c(10^\infty) = 110^\infty, \qquad
\]
show that all members of the forward orbit of $10^\infty$ under the
action of $a$ have only finitely many 1's and that the position of
the rightmost 1 cannot decrease under the action of $a$. Since
$a(10^\infty)=010^\infty$, the forward orbit of $10^\infty$ under
the action of $a$ can never return to $10^\infty$ and $a$ has
infinite order.

Note that the above equalities also show that no nonempty words $w$
over $\{a,b,c\}$ satisfies a relation of the form $w=1$ in
$G_{875}$. First note that $c=(b,a)$ and $b=(b,c)$, implying that
$b$ and $c$ have infinite order. Thus $b^n \neq 1$, for $n > 0$. On
the other hand, for any word $w$ that contains $a$ or $c$,
$w(10^\infty) \neq 10^\infty$ (again, since the position of the
rightmost 1 moves to the right and never decreases).

Since $b$ has infinite order and $b=(b,c)$, $G_{875}$ is not
contracting.

\noindent\textbf{876}. Wreath recursion: $a=\sigma(c,a)$, $b=(b,c)$,
$c=(b,a)$.

By Lemma~\ref{alg_trans} the elements $ba$ and $acb^2a^2cb$ act
transitively on the levels of the tree and, hence, have infinite
order. Since $(b^8)\bigl|_{1100001100}=acb^2a^2cb$ and vertex
$1100001100$ is fixed under the action of $b^8$ we obtain that $b$
also has infinite order. Finally, $b$ stabilizes the vertex $0$ and
has itself as a section at this vertex. Therefore $G_{876}$ is not
contracting.

We have $c^{-1}b = (1, a^{-1}c), ac^{-1}ba^{-1} = (ca^{-1}, 1)$,
hence by Lemma~\ref{lem:freegroup} the group is not free.

\noindent$\mathbf{878}\cong C_2\ltimes IMG(1-\frac1{z^2})$. Wreath
recursion: $a=\sigma(b,b)$, $b=(b,c)$, $c=(b,a)$.

Let $x=bc$ and $y=ca$. Since all generators have order $2$, the
index of the subgroup $H=\langle x,y\rangle$ in $G_{878}$ is 2, $H$
is normal and $G_{878}\cong C_2\ltimes H$, where $C_2$ is generated
by $c$. The action of $C_2$ on $H$ is given by $x^c = x^{-1}$ and
$y^c = y^{-1}$. We have $x=bc=(1,ca)=(1,y)$ and
$y=ca=\sigma(ab,1)=\sigma(y^{-1}x^{-1},1)$. An isomorphic copy of
$H$ is obtained by exchanging the letters $0$ and $1$, yielding the
wreath recursion $x=(y,1)$ and $y=\sigma(1,y^{-1}x^{-1})$. The last
recursion defines $IMG(1-\frac1{z^2})$~\cite{bartholdi_n:rabbit}.
Thus, $G_{878}\cong C_2\ltimes IMG(1-\frac1{z^2})$.

\noindent\textbf{879}. Wreath recursion: $a=\sigma(c,b)$, $b=(b,c)$,
$c=(b,a)$.

The element $c^{-1}a=\sigma(a^{-1}c,1)$ is conjugate to the adding
machine and has infinite order.

By Lemma~\ref{alg_trans} the element $ca$ acts transitively on the
levels of the tree and, hence, has infinite order. Since
$(b^2)\bigl|_{1101}=ca$ and vertex $1101$ is fixed under the action
of $b^2$ we obtain that $b$ also has infinite order. Finally, $b$
stabilizes the vertex $0$ and has itself as a section at this
vertex. Therefore $G_{879}$ is not contracting.

\noindent\textbf{882}. Wreath recursion: $a=\sigma(c,c)$, $b=(b,c)$,
$c=(b,a)$.

The element $(ca^{-1}cb^{-1})^2$ stabilizes the vertex $00$ and its
section at this vertex is equal to $ca^{-1}cb^{-1}$. Hence,
$ca^{-1}cb^{-1}$ has infinite order.

\noindent$\mathbf{883}\cong C_2\ltimes G_{2841}$. Wreath recursion:
$a=\sigma(a,a)$, $b=(c,c)$, $c=(b,a)$.

All generators have order $2$, hence $H=\langle ba,ca\rangle$ is
normal in $G_{883}$. Furthermore, $ba=\sigma(ca,ca)$,
$ca=\sigma(1,ba)$, therefore $H=G_{2841}$. Thus $G_{883}=\langle
a\rangle\ltimes H\cong C_2\ltimes G_{2841}$, where
$(ba)^a=(ba)^{-1}$ and $(ca)^a=(ca)^{-1}$. In particular, $G_{883}$
is not contracting and has exponential growth.

\noindent\textbf{884}. Wreath recursion: $a=\sigma(b,a)$, $b=(c,c)$,
$c=(b,a)$.

The element $(b^{-1}ca^{-1}c)^2$ stabilizes the vertex $0$ and its
section at this vertex is equal to $(b^{-1}ca^{-1}c)^{-1}$. Hence,
$b^{-1}ca^{-1}c$ has infinite order. Since
$[b,a]^2\bigl|_{0100}=(b^{-1}ca^{-1}c)^c$ and $0100$ is fixed under
the action of $[b,a]^2$ we obtain that $[b,a]$ also has infinite
order. Finally, $[b,a]$ stabilizes the vertex $00$ and its section
at this vertex is $[b,c]=[b,a]$. Therefore $G_{884}$ is not
contracting.

\noindent\textbf{885}. Wreath recursion: $a=\sigma(c,a)$, $b=(c,c)$,
$c=(b,a)$.

The element $(c^{-1}b)^2$ stabilizes the vertex $10$ and its section
at this vertex is equal to $c^{-1}b$. Hence, $c^{-1}b$ has infinite
order. Furthermore, $c^{-1}b$ stabilizes the vertex $00$ and has
itself as a section at this vertex. Therefore $G_{885}$ is not
contracting.

We have $b^{-1}aba^{-1} = (1, c^{-1}aca^{-1}), a^{-1}b^{-1}ab = (a^{-1}c^{-1}ac, 1)$,
hence by Lemma~\ref{lem:freegroup} the group is not free.

\noindent\textbf{887}. Wreath recursion: $a=\sigma(b,b)$, $b=(c,c)$,
$c=(b,a)$.

The element $(ac^{-1})^4$ stabilizes the vertex $001$ and its
section at this vertex is equal to $(ac^{-1})^{2}$, which is
nontrivial. Hence, $ac^{-1}$ has infinite order.

\noindent\textbf{888}. Wreath recursion: $a=\sigma(c,b)$, $b=(c,c)$,
$c=(b,a)$.

The element $a^{-1}c=\sigma(1,c^{-1}a)$ is conjugate to the adding
machine and has infinite order. Since $c^{-1}b\bigl|_{1}=a^{-1}c$
and vertex $1$ is fixed under the action of $c^{-1}b$ we obtain that
$c^{-1}b$ also has infinite order. Finally, $c^{-1}b$ stabilizes the
vertex $00$ and has itself as a section at this vertex. Therefore
$G_{888}$ is not contracting.

We have $c^{-1}ab^{-1}a = (1, a^{-1}b), ac^{-1}ab^{-1} = (ca^{-1}bc^{-1}, 1)$,
hence by Lemma~\ref{lem:freegroup} the group is not free.

\noindent$\mathbf{891}\cong C_2\ltimes (\Z \wr C_2)$. Wreath
recursion: $a=\sigma(c,c)$, $b=(c,c)$, $c=(b,a)$.

Let $x=ac$ and $y=cb$. Since all generators have order $2$, the
index of the subgroup $H=\langle x,y\rangle$ in $G_{891}$ is 2, $H$
is normal and $G_{891}\cong C_2\ltimes H$, where $C_2$ is generated
by $c$. The action of $C_2$ on $H$ is given by $x^c = x^{-1}$ and
$y^c = y^{-1}$.

In fact, to support the claim that $H$ has index 2 in $G_{891}$ we
need to prove that $c \not \in H$. We will prove a little bit more
than that. Let $w=1$ be a relation in $G_{891}$, where $w$ is a word
over $\{a,b,c\}$. The number of occurrences of $a$ in $w$ must be
even (otherwise $w$ would act nontrivially on level 1). Similarly,
the number of occurrences of $c$ in $w$ is even. Indeed, if it were
odd, then exactly one of the words $w_0$ and $w_1$ in the
decomposition $w=(w_0,w_1)$ would have odd number of occurrences of
the letter $a$, and the action of $w$ would be nontrivial on level
2. Finally, we claim that the number of occurrences of $b$ in $w$ is
also even. Otherwise the number of $c$'s in both $w_0$ and $w_1$
would be odd and the action of $w$ would be nontrivial on level 3.
Thus every word over $\{a,b,c\}$ representing 1 must have even
number of occurrences of each of the three letters. Note that this
implies that the abelianization of $G_{891}$ is $C_2 \times C_2
\times C_2$.

We now prove that $H$ is isomorphic to the Lamplighter group $\Z \wr
C_2$. The group $H$ is self-similar, which can be seen from
\[
 x = ac = \s(cb,ca) = \s(y,x^{-1}), \qquad
 y = cb = (bc,ac) = (y^{-1},x).
\]

Consider the elements $s_n=\sigma^{y^n}=y^{-n}xy^{n+1}$, $n\in\Z$
(note that $xy=\s$). For $n >0$, we have $s_0s_1\cdots
s_{n-1}=x^ny^n$ and $s_{-n}s_{-n+1}\cdots s_{-1}=y^nx^n$. On the
other hand, $s_n=y^{-n}\sigma y^n=\sigma(x^{-n}y^{-n},y^{n}x^{n})$
and $s_{-n}=y^{n}\sigma y^{-n}=\sigma(x^ny^n,y^{-n}x^{-n})$,
implying
\[s_n=\sigma(s_{-1}s_{-2}\cdots s_{-n}, s_{-n}\cdots s_{-2}s_{-1})\]
and
\[s_{-n}=\sigma(s_{0}s_{1}\cdots s_{n-1}, s_{n-1}\cdots s_{1}s_{0}).\]

By induction on $n$ we obtain that the depth of $s_n$ is $2n+1$ for
$n\geq 0$ and the depth of $s_{-n}$ is $2n$ for $n>0$ (\emph{depth}
of a finitary element is the lowest level at which all sections of
the element are trivial). This implies that all $s_i$, $i \in \Z$
are different, have order $2$ (they are conjugates of $\s$), and
commute (for each $i$ and each level $m$ all sections of $s_i$ at
level $m$ are equal). Therefore $y$ has infinite order and $H =
\langle x,y \rangle = \langle y, \sigma \rangle \cong \Z \wr C_2$.

Since $y$ has infinite order, stabilizes the vertex $00$ and has
itself as a section at this vertex, $G_{891}$ is not contracting.

\noindent$\mathbf{919}\cong G_{820}\cong D_\infty$. Wreath
recursion: $a=\sigma(a,a), b=(a,b), c=(c,a)$.

The states $a$, $b$ form a $2$-state automaton generating $D_\infty$
(see Theorem~\ref{thm:class22}) and $c=aba$.

\noindent\textbf{920}. Wreath recursion: $a=\sigma(b,a)$, $b=(a,b)$,
$c=(c,a)$.

The element $(ac^{-1})^2$ stabilizes the vertex $00$ and its section
at this vertex is equal to $ac^{-1}$. Hence, $ba^{-1}$ has infinite
order.

\noindent\textbf{923}. Wreath recursion: $a=\sigma(b,b)$, $b=(a,b)$,
$c=(c,a)$.

The states $a$ and $b$ form a $2$-state automaton generating
$D_\infty$ (see Theorem~\ref{thm:class22}).

\noindent$\mathbf{924}\cong G_{870}$. Baumslag-Solitar group
$BS(1,3)$. Wreath recursion: $a=\sigma(c,b)$, $b=(a,b)$, $c=(c,a)$.

This fact is proved in~\cite{bartholdi_s:bsolitar}.

\noindent$\mathbf{928}\cong G_{820}\cong D_\infty$. Wreath
recursion: $a=\sigma(a,a)$, $b=(b,b)$, $c=(c,a)$.

The states $a$ and $c$ form a $2$-state automaton generating
$D_\infty$ (see Theorem~\ref{thm:class22}) and $b$ is trivial.

\noindent$\mathbf{929}\cong G_{2851}$. Wreath recursion:
$a=\sigma(b,a)$, $b=(b,b)$, $c=(c,a)$.

See $G_{2851}$ for an isomorphism (in fact the groups coincide as
subgroups of $\Aut(X^*)$).

\noindent$\mathbf{930}\cong G_{821}$. Lamplighter group $\Z \wr
C_2$. Wreath recursion: $a=\sigma(c,a)$, $b=(b,b)$, $c=(c,a)$.

The states $a$ and $c$ form a $2$-state automaton generating the
Lamplighter group (see Theorem~\ref{thm:class22}) and $b$ is
trivial.

\noindent$\mathbf{932}\cong G_{820}\cong D_\infty$. Wreath
recursion: $a=\sigma(b,b)$, $b=(b,b)$, $c=(c,a)$.

We have $b=1$ and $a^2=c^2=1$. The element $ac=\sigma(c,a)$ is
clearly nontrivial. Since $(ac)^2 = (ac,ca)$, this element has
infinite order. Thus $G\cong D_\infty$.

\noindent$\mathbf{933}\cong G_{849}$. Wreath recursion:
$a=\sigma(c,b)$, $b=(b,b)$, $c=(c,a)$.

See $G_{2852}$ for an isomorphism between $G_{933}$ and $G_{2852}$
and $G_{849}$ for an isomorphism between $G_{2852}$ and $G_{849}$.

\noindent$\mathbf{936}\cong G_{820}\cong D_\infty$. Wreath
recursion: $a=\sigma(c,c)$, $b=(b,b)$, $c=(c,a)$.

The states $a$ and $c$ form a $2$-state automaton generating
$D_\infty$ (see Theorem~\ref{thm:class22}) and $b$ is trivial.

\noindent$\mathbf{937}\cong C_2\ltimes G_{929}$. Wreath recursion:
$a=\sigma(a,a)$, $b=(c,b)$, $c=(c,a)$.

All generators have order $2$, hence $H=\langle ca,ba\rangle=\langle
ca,caba\rangle$ is normal in $G_{937}$. Furthermore,
$ca=\sigma(1,ca)$, $caba=\sigma(caba,ca)$, therefore $H=G_{929}$.
Thus $G_{937}=\langle a\rangle\ltimes H\cong C_2\ltimes G_{929}$,
where $(ba)^a=(ba)^{-1}$ and $(ca)^a=(ca)^{-1}$. In particular,
$G_{937}$ is regular weakly branch over $H'$, has exponential growth
and is not contracting.

\noindent\textbf{938}. Wreath recursion: $a=\sigma(b,a)$, $b=(c,b)$,
$c=(c,a)$.

The element $(b^{-1}a^{-1}ca)^2$ stabilizes the vertex $00$ and its
section at this vertex is equal to
$\bigl((b^{-1}a^{-1}ca)^{-1}\bigr)^{a^{-1}c}$. Hence,
$b^{-1}a^{-1}ca$ has infinite order. Furthermore, $b^{-1}a^{-1}ca$
stabilizes the vertex $1$ and has itself as a section at this
vertex. Therefore $G_{938}$ is not contracting.

We have $c^{-1}b = (1, a^{-1}b), a^{-1}c^{-1}ba = (a^{-2}ba, 1)$,
hence by Lemma~\ref{lem:freegroup} the group is not free.

\noindent\textbf{939}. Wreath recursion: $a=\sigma(c,a)$, $b=(c,b)$,
$c=(c,a)$.

The states $a$ and $c$ form a $2$-state automaton generating the
Lamplighter group (see Theorem~\ref{thm:class22}). Hence, $G_{939}$
is neither torsion, nor contracting, and has exponential growth.

\noindent\textbf{941}. Wreath recursion: $a=\sigma(b, b)$,
$b=(c,b)$, $c=(c, a)$.

The second iteration of the wreath recursion is
\[a=(02)(13)(c, b, c, b), \qquad b=(c, a, c, b), \qquad c=(23)(c, a, b, b).\]

Conjugation by $g=(cg, g, g, bg)$ gives the wreath recursion
\[ a' =(02)(13), \qquad b=(c', a', c', b'), \qquad c=(23)(c', a', 1, 1),\]
where $a' = a^g$, $b' = b^g$, and $c' = c^g$. The last recursion
coincides with the second iteration of the recursion
\[\alpha= \sigma, \qquad \beta=(\gamma, \beta), \qquad \gamma=(\gamma, \alpha).\]
Conjugating the last recursion by $h=(\gamma h,h)$ yields the
recursion defining $G_{945}$. Thus, $G_{941}\cong G_{945}\cong
C_2\ltimes IMG(z^2-1)$ (see $G_{945}$). The limit space is half of
the Basilica.

\noindent\textbf{942}. Wreath recursion: $a=\sigma(c,b)$, $b=(c,b)$,
$c=(c,a)$.

The Lamplighter group $L=\Z \wr C_2$ can be defined as the group
generated by $a'$ and $b'$ given by the wreath recursion (see
Theorem~\ref{thm:class22})
\begin{alignat*}{2}
  &a' = \sigma &&(a',b'),\\
  &b' =        &&(a',b').
\end{alignat*}
Let $H=\langle a, b\rangle \leq G_{942}$. We will show that $H$ and
$L$ are isomorphic. Let $Y^*$ be the subtree of $X^*$ consisting of
all words over the alphabet $Y=\{01,11\}$. The element $b$ fixes the
letter in $Y$, while $a$ swaps them. Since $a_{01}=b_{01}=a$,
$a_{11}=b_{11}=b$, the tree $Y^*$ is invariant under the action of
$H$. Moreover, the action of $H$ on $Y^*$ coincides with the action
of the Lamplighter group $L=\langle a',b' \rangle$ on $X^*$ (after
the identification $01 \leftrightarrow 0, \ 11 \leftrightarrow 1$).
This implies that the map $\phi:H \to L$ given by $a \mapsto a'$, $b
\mapsto b'$ can be extended to a homomorphism. We claim that this
homomorphism is in fact an isomorphism. Let $w=w(a,b)$ be a group
word representing an element of the kernel of $\phi$. Since
$w(a',b')$ represents the identity in the lamplighter group $L$, the
total exponent of $a$ in $w$ must be even and the total exponent
$\varepsilon$ of both $a$ and $b$ in $w$ must be 0. Therefore the
element $g=w(a,b)$ stabilizes the top two levels of the tree $X^*$
and can be decomposed as
\[ g = (c^\varepsilon,*,c^\varepsilon,*), \]
where the $*$'s are words over $a$ and $b$ representing the identity
in $H$ (these words correspond precisely to the first level sections
of $w(a',b')$ in $L$). Since $\varepsilon=0$, we see that $g=1$ and
the kernel of $\phi$ is trivial.

Thus, the Lamplighter group is a subgroup of $G_{942}$, which shows
that $G_{942}$ is not a torsion group, it is not free, and has
exponential growth.  Since $b=(c,b)$ and $b$ has infinite order,
$G_{942}$ is not a contracting group.

\noindent$\mathbf{945}\cong G_{941}\cong C_2\ltimes IMG(z^2-1)$.
Wreath recursion: $a=\sigma(c,c)$, $b=(c,b)$, $c=(c,a)$.

All generators have order $2$. Since $ab=\sigma(1,cb)$ and
$cb=(1,ab)$ we see that $H=\langle ab, cb \rangle\cong
G_{852}=IMG(z^2-1)$. This subgroup is normal in $G_{945}$ because
the generators have order $2$. Since $G_{945}=\langle H,b\rangle$,
it has a structure of a semidirect product $\langle b\rangle\ltimes
H=C_2\ltimes IMG(z^2-1)$ with the action of $b$ on $H$ given by
$(ab)^b=(ab)^{-1}$ and $(cb)^b=(cb)^{-1}$. It follows that $G_{945}$
is regular weakly branch over $H'$ and has exponential growth. See
$G_{941}$ for an isomorphism.

\noindent$\mathbf{955}\cong G_{937}\cong C_2\ltimes G_{929}$. Wreath
recursion: $a=\sigma(a,a)$, $b=(b,c)$, $c=(c,a)$.

All generators have order $2$. Consider the subgroup $H=\langle
ba=\sigma(ca,ba), ca=\sigma(1,ca)\rangle\cong G_{929}$. This
subgroup is normal in $G_{955}$ because all generators have order
$2$. Since $G_{955}=\langle H,a\rangle$, it has a structure of a
semidirect product $\langle a\rangle\ltimes H=C_2\ltimes G_{929}$
with the action of $a$ on $H$ given by $(ba)^b=(ba)^{-1}$ and
$(ca)^b=(ca)^{-1}$. It is proved above that $G_{937}$ has the same
structure. It follows that $G_{955}$ is regular weakly branch over
$H'$ and has exponential growth.

\noindent\textbf{956}. Wreath recursion: $a=\sigma(b,a)$, $b=(b,c)$,
$c=(c,a)$.

The element $(c^{-1}b)^2$ stabilizes the vertex $10$ and its section
at this vertex is equal to $(c^{-1}b)^{-1}$. Hence, $c^{-1}b$ has
infinite order. Furthermore, $c^{-1}b$ stabilizes the vertex $0$ and
has itself as a section at this vertex. Therefore $G_{956}$ is not
contracting.

We have $c^{-1}b^{-1}aba^{-1}b = (1, a^{-1}c^{-1}aba^{-1}c),
a^{-1}c^{-1}b^{-1}aba^{-1}ba = (a^{-2}c^{-1}aba^{-1}ca, 1)$, hence
by Lemma~\ref{lem:freegroup} the group is not free.

\noindent\textbf{957}. Wreath recursion: $a=\sigma(c,a)$, $b=(b,c)$,
$c=(c,a)$.

The states $a$, $c$ form a $2$-state automaton generating the
Lamplighter group (see Theorem~\ref{thm:class22}). Hence, $G_{957}$
is neither torsion, nor contracting and has exponential growth.

\noindent\textbf{959}. Wreath recursion: $a=\sigma(b,b)$, $b=(b,c)$,
$c=(c,a)$.

The element $(a^{-1}c)^4$ stabilizes the vertex $00$ and its section
at this vertex is equal to $(a^{-1}c)^{-1}$. Hence, $a^{-1}c$ has
infinite order.

Furthermore, since $c^{-1}b=(c^{-1}b,a^{-1}c)$, this element also
has infinite order. Thus, $G_{959}$ is not contracting.

\noindent\textbf{960}. Wreath recursion: $a=\sigma(c,b)$, $b=(b,c)$,
$c=(c,a)$.

Define $x=ac^{-1}$, $y=ba^{-1}$ and $z=cb^{-1}$. Then
$x=\sigma(1,y)$, $y=\sigma(z,z^{-1})$ and $z=(z,x)$.

The element $(zxy)^8$ stabilizes the vertex $001010$ and its section
at this vertex is equal to $xy^{-1}z=xyz=(zxy)^{z^{-1}}$ (since
$y^2=1$). Hence, $zxy$ has infinite order.

Denote $t=(b^{-1}c)^4(b^{-1}a)(c^{-1}a)^5(b^{-1}c)$. Then $t^2$
stabilizes the vertex $00$ and $t^2\bigl|_{00}=t^{b^{-1}c}$. Hence,
$t$ has infinite order. Let $s=c^{-2}b^2$. Since
$s^{32}\bigl|_{111000000100}=t^c$ and $s^{32}$ fixes $111000000100$,
we obtain that $s$ also has infinite order. Finally, $s$ stabilizes
the vertex $00$ and has itself as a section at this vertex.
Therefore $G_{960}$ is not contracting.

\noindent\textbf{963}. Wreath recursion: $a=\sigma(c,c)$, $b=(b,c)$,
$c=(c,a)$.

All generators have order $2$. The element $ac=\sigma(1,ca)$ is
conjugate to the adding machine and has infinite order.

Furthermore, since $cb=(cb,ac)$, this element also has infinite
order. Thus, $G_{963}$ is not contracting.

\noindent$\mathbf{964}\cong G_{739}\cong C_2\ltimes (C_2\wr \Z)$.
Wreath recursion: $a=\sigma(a,a)$, $b=(c,c)$, $c=(c,a)$.

All generators have order $2$. The elements $u=acba=(ca,1)$ and
$v=bc=(1,ca)$ generate $\Z^2$ because $ca=\sigma(1,ca)$ is the
adding machine and has infinite order. We have $cacb=\sigma$ and
$\langle u,v\rangle$ is normal in $H=\langle u,v, \sigma\rangle$
because $u^\sigma=v$ and $v^\sigma=u$. In other words, $H\cong
C_2\ltimes(\Z\times\Z)=C_2\wr \Z$.

Furthermore, $G_{964}=\langle H,a\rangle$ and $H$ is normal in
$G_{972}$ because $u^a=v^{-1}$, $v^a=u^{-1}$ and $\sigma^a=\sigma$.
Thus $G_{964}=C_2\ltimes (C_2\wr \Z)$, where the action of $C_2$ on
$H$ is specified above and coincides with the one in $G_{739}$.
Therefore $G_{964}\cong G_{739}$.

\noindent\textbf{965}. Wreath recursion: $a=\sigma(b,a)$, $b=(c,c)$,
$c=(c,a)$.

The element $(ac^{-1})^2$ stabilizes the vertex $01$ and its section
at this vertex is equal to $(ac^{-1})^{-1}$. Hence, $ac^{-1}$ has
infinite order.

By Lemma~\ref{alg_trans} the element $a$ acts transitively on the
levels of the tree and, hence, has infinite order. Since $c=(c,a)$
we obtain that $c$ also has infinite order. Therefore $G_{965}$ is
not contracting.

We have $bc^{-1} = (1, ca^{-1}), a^{-1}bc^{-1}a = (a^{-1}c, 1)$,
hence by Lemma~\ref{lem:freegroup} the group is not free.

\noindent\textbf{966}. Wreath recursion: $a=\sigma(c,a)$, $b=(c,c)$,
$c=(c,a)$.

The states $a$ and $c$ form a $2$-state automaton generating the
Lamplighter group (see Theorem~\ref{thm:class22}). Hence, $G_{966}$
is neither torsion, nor contracting, and has exponential growth.

Since $b=(c,c)$ we obtain that $G_{966}$ can be embedded into the
wreath product $C_2\wr (\Z\wr\C_2)$. This shows that $G_{966}$ is
solvable.

\noindent\textbf{968}. Wreath recursion: $a=\sigma(b,b)$, $b=(c,c)$,
$c=(c,a)$.

We will show that this group contains $\Z^5$ as a subgroup of index
$16$. It is a contracting group, with nucleus consisting of $73$
elements (the self-similar closure of the nucleus consists of $77$
elements).

All generators have order $2$. Let $x=(ac)^2$, $y=bcba$, and
$K=\langle x,y \rangle$. Conjugating $x$ and $y$ by
$\gamma=(b\gamma,a\gamma)$ yields the self-similar copy $K'$ of $K$
generated by $x' = ((y')^{-1}, (y')^{-1})$ and $y = \sigma(x',y')$,
where $x' = x^\gamma$ and $y'=y^\gamma$. Since
$[x',y']=([x',y']^{(y')^{-1}},1)$ $K'$ is abelian. The matrix of the
corresponding virtual endomorphism is given by
\[
 A=\left(\begin{array}{cc}
               0&\frac12\\
              -1&\frac12\\
 \end{array} \right).
\]
The eigenvalues $\lambda=\frac14\pm\frac14\sqrt{7}i$ of this matrix
are not algebraic integers. Therefore $K'$ (ad therefore $K$ as
well) is free abelian of rank $2$, by the results
in~\cite{nekrash_s:12endomorph}.

The subgroup $H=\langle ab,bc\rangle$ has index $2$ in $G_{968}$
(the generators of $G_{968}$ have order 2). The second level
stabilizer $\St_H(2)$ has index $8$ in $H$ (the quotient group is
isomorphic to the dihedral group $D_4$). The stabilizer $\St_H(2)$,
is generated by $(bc)^2$, $\left((bc)^2\right)^{ba}$, $(ab)^2$,
$\left((ab)^2\right)^{bc}$, $\left((ab)^2\right)^{(bc)^{ba}}$, and
$\left((ab)^2\right)^{bc(bc)^{ba}}$. Conjugating these elements by
$g=(b,c,b,1)$ gives
\[
\begin{array}{lllllll}
g_1=\left((bc)^2\right)^g&=(bcbc)^g=&(1,&1,&y,&y^{-1}&),\\
g_2=\left((bc)^2\right)^{bag}&=(acbcba)^g=&(y,&y,&1,&1&),\\
g_3=\left((ab)^2\right)^{bcg}&=(cbabac)^g=&(1,&x,&x,&1&),\\
g_4=\left((ab)^2\right)^g&=(abab)^g=&(1,&x,&1,&x^{-1}&),\\
g_5=\left((ab)^2\right)^{(bc)^{ba}g}&=(abcbabacba)^g=&(x,&1,&1,&x^{-1}&),\\
g_6=\left((ab)^2\right)^{bc(bc)^{ba}g}&=(abcacbabacacba)^g=&(x,&1,&x,&1&).
\end{array}
\]
Therefore, $\St_H(2)$ is abelian and $g_6=g_5g_3g_4^{-1}$. If
$\prod_{i=1}^5g_i^{n_i}=1$, then
$x^{n_5}y^{n_2}=x^{n_3+n_4}y^{n_2}=x^{n_3}y^{n_1}=x^{n_4+n_5}y^{n_1}=1$.
Since $K$ is free abelian, we obtain $n_i=0$, $i=1,\ldots,5$.
Therefore $\St_H(2)$ is a free abelian group of rank $5$.

\noindent\textbf{969}. Wreath recursion: $a=\sigma(c,b)$, $b=(c,c)$,
$c=(c,a)$.

The element $(cb^{-1})^4$ stabilizes the vertex $100$ and its
section at this vertex is equal to $cb^{-1}$. Hence, $cb^{-1}$ has
infinite order.

We have $bc^{-1}=(1,ca^{-1})$, $ca^{-1}=\sigma(ab^{-1},1)$,
$ab^{-1}=\sigma(1,bc^{-1})$, hence the subgroup generated by these
elements is isomorphic to $IMG(1-\frac{1}{z^2})$
(see~\cite{bartholdi_n:rabbit}).

We also have $c^{-1}b = (1, a^{-1}c), a^{-1}c^{-1}ba =
(b^{-1}a^{-1}cb, 1)$, hence by Lemma~\ref{lem:freegroup} the group
is not free.

\noindent$\mathbf{972}\cong G_{739}\cong C_2\ltimes (C_2\wr \Z)$.
Wreath recursion : $a=\sigma(c,c)$, $b=(c,c)$, $c=(c,a)$.

All generators have order $2$. The elements $u=acba=(ca,1)$ and
$v=bc=(1,ac)$ generate $\Z^2$ because $ca=\sigma(ac,1)$ is conjugate
to the adding machine and has infinite order. Also we have
$ba=\sigma$ and $\langle u,v\rangle$ is normal in $H=\langle u,v,
\sigma\rangle$ because $u^\sigma=v$ and $v^\sigma=u$. In other
words, $H\cong C_2\ltimes(\Z\times\Z)=C_2\wr \Z$.

Furthermore, $G_{972}=\langle H,a\rangle$ and $H$ is normal in
$G_{972}$ because $u^a=v^{-1}$, $v^a=u^{-1}$ and $\sigma^a=\sigma$.
Thus $G_{972}=C_2\ltimes (C_2\wr \Z)$, where the action of $C_2$ on
$H$ is specified above and coincides with the one in $G_{739}$.
Therefore $G_{972}\cong G_{739}$.

\noindent$\mathbf{1090}\cong C_2$. Wreath recursion:
$a=\sigma(a,a)$, $b=(b,b)$, $c=(b,b)$.

Both $b$ and $c$ are trivial and $a^2=1$.

\noindent$\mathbf{1091}\cong G_{731}\cong \Z$. Wreath recursion:
$a=\sigma(b,a)$, $b=(b,b)$, $c=(b,b)$.

Both $b$ and $c$ are trivial and $a$ is the adding machine.

\noindent$\mathbf{1094}\cong G_{1090}\cong C_2$. Wreath recursion:
$a=\sigma(b,b)$, $b=(b,b)$, $c=(b,b)$.

Both $b$ and $c$ are trivial and $a^2=1$.

\noindent$\mathbf{2190}\cong G_{848}\cong C_2\wr\Z$. Wreath
recursion: $a=\sigma(c,a)$, $b=\sigma(a,a)$, $c=(a,a)$.

First note that $c=a^{-2}$. Therefore $G=\langle a,b \rangle$, where
$a=\sigma(a^{-2},a)$, and $b=\sigma(a,a)$. Also, $a$ has infinite
order.

Consider the subgroup $H=\langle ba,ab\rangle<G$. The generators of
$H$ commute since $ba=(a^{-1},a^2)$ and $ab=(a^2,a^{-1})$.
Furthermore, $(ba)^n(ab)^m=(a^{-n+2m},a^{2n-m})=1$ if and only if
$m=n=0$. Therefore $H\cong \Z^2$.

Consider the element $ba^2=bc^{-1}=\sigma$. This element does not
belong to $H$, since $H$ stabilizes the first level of the tree. On
the other hand $a=(ba)^{-1}ba^2=(ba)^{-1}\sigma$ and $b=a^{-1}(ab)$
so $G=\langle \sigma, H \rangle$. Finally, $(ba)^\sigma=ab$ and
$(ab)^\sigma=ba$ implies that $H$ is normal in $G$ and $G=C_2\wr H
\cong C_2\wr\Z \cong G_{848}$.

Also note that $\langle a, a^b\rangle=G_{2212}\cong\Z\ast_{2\Z}\Z$.

\noindent\textbf{2193}. Wreath recursion: $a=\sigma(c,b)$,
$b=\sigma(a,a)$, $c=(a,a)$.

Let $x=ca^{-1}$ and $y=ab^{-1}$. Then
$x=\s(ab^{-1},ac^{-1})=\s(y,x^{-1})$ and
$y=(ba^{-1},ca^{-1})=(y^{-1},x)$. It is already shown (see
$G_{891}$), that $\langle x,y\rangle$ is not contracting and is
isomorphic to the Lamplighter group. Therefore $G_{2193}$ is not a
torsion group, it is not contracting, and has exponential growth.

\noindent$\mathbf{2196}\cong G_{802}\cong C_2\times C_2\times C_2$.
Wreath recursion: $a=\sigma(c,c)$, $b=\sigma(a,a)$, $c=(a,a)$.

Direct calculation.

\noindent\textbf{2199}. Wreath recursion: $a=\sigma(c,a)$,
$b=\sigma(b,a)$, $c=(a,a)$.

By Lemma~\ref{alg_trans} the element $ac$ acts transitively on the
levels of the tree and, hence, has infinite order. Since
$ba=(ac,ba)$ we obtain that $ba$ also has infinite order. Therefore
$G_{2199}$ is not contracting.

We have $b^{-2}abcba = b^{-2}aba^{-2}ba = 1$, and $a$ and $b$ do not commute,
hence the group is not free.

\noindent\textbf{2202}. Wreath recursion: $a=\sigma(c,b)$,
$b=\sigma(b,a)$, $c=(a,a)$.

The element $(b^{-1}a)^2$ stabilizes the vertex $00$ and its section
at this vertex is equal to $b^{-1}a$. Hence, $b^{-1}a$ has infinite
order. Furthermore, $b^{-1}a$ stabilizes the vertex $11$ and has
itself as a section at this vertex. Therefore $G_{2202}$ is not
contracting.

We have $cb^{-1}c^{-1}b = (1, ab^{-1}a^{-1}b), bcb^{-1}c^{-1} =
(bab^{-1}a^{-1}, 1)$, hence by Lemma~\ref{lem:freegroup} the group
is not free.

\noindent\textbf{2203}. Wreath recursion: $a=\sigma(a,c)$,
$b=\sigma(b,a)$, $c=(a,a)$.

The states $a$ and $c$ form a 2-state automaton generating the
infinite cyclic group $\Z$ in which $c=a^{-2}$ (see
Theorem~\ref{thm:class22}).

Since $b^{-1}a\bigl|_{1}=a^{-1}c$ and vertex $1$ is fixed under the
action of $b^{-1}a$ we obtain that $b^{-1}a$ also has infinite
order. Finally, $b^{-1}a$ stabilizes the vertex $0$ and has itself
as a section at this vertex. Therefore $G_{2203}$ is not
contracting.

We have $c^{-2}ab = (1, a^{-2}cb), bc^{-2}a = (ba^{-2}c, 1)$, hence
by Lemma~\ref{lem:freegroup} the group is not free.

\noindent\textbf{2204}. Wreath recursion: $a=\sigma(b,c)$,
$b=\sigma(b,a)$, $c=(a,a)$.

The element $(b^{-1}ac^{-1}a)^2$ stabilizes the vertex $00$ and its
section at this vertex is equal to $b^{-1}ac^{-1}a$. Hence,
$b^{-1}ac^{-1}a$ has infinite order. Since
$[c,a]^2\bigl|_{000}=(b^{-1}ac^{-1}a)^{a^{-1}cb}$ and $000$ is fixed
under the action of $[c,a]^2$ we obtain that $[c,a]$ also has
infinite order. Finally, $[c,a]$ stabilizes the vertex $11$ and has
itself as a section at this vertex. Therefore $G_{2204}$ is not
contracting.

We have $ab^{-1} = (1, ca^{-1}), b^{-1}a = (a^{-1}c, 1)$, hence by
Lemma~\ref{lem:freegroup} the group is not free.

\noindent$\mathbf{2205}\cong G_{775}\cong C_2\ltimes
IMG\left(\bigl(\frac{z-1}{z+1}\bigr)^2\right)$. Wreath recursion:
$a=\sigma(c,c)$, $b=\sigma(b,a)$, $c=(a,a)$.

See $G_{783}$ for an isomorphism between $G_{783}$ and $G_{2205}$.

\noindent$\mathbf{2206}\cong G_{748}\cong D_4\times C_2$. Wreath
recursion: $a=\sigma(a,a)$, $b=\sigma(c,a)$, $c=(a,a)$.

Direct calculation.

\noindent\textbf{2207}. Wreath recursion: $a=\sigma(b,a)$,
$b=\sigma(c,a)$, $c=(a,a)$.

The element $(c^{-1}a)^4$ stabilizes the vertex $000$ and its
section at this vertex is equal to $c^{-1}a$. Hence, $c^{-1}a$ has
infinite order.

Since $b^{-1}a^{-1}b^{-1}aba\bigl|_{001}=(c^{-1}a)^a$ and the vertex
$001$ is fixed under the action of $b^{-1}a^{-1}b^{-1}aba$ we obtain
that $b^{-1}a^{-1}b^{-1}aba$ also has infinite order. Finally,
$b^{-1}a^{-1}b^{-1}aba$ stabilizes the vertex $000$ and has itself
as a section at this vertex. Therefore $G_{2207}$ is not
contracting.

We have $a^{-2}bab^{-2}ab = 1$, and $a$ and $b$ do not commute, hence the
group is not free.

\noindent\textbf{2209}. Wreath recursion: $a=\sigma(a,b)$,
$b=\sigma(c,a)$, $c=(a,a)$.

The element $(b^{-1}a)^2$ stabilizes the
vertex $00$ and its section at this vertex is equal to
$(b^{-1}a)^{-1}$. Hence, $b^{-1}a$ has infinite order. Furthermore,
$b^{-1}a$ stabilizes the vertex $11$ and has itself as a section at
this vertex. Therefore $G_{2209}$ is not contracting.

We have $aca^{-2}c^{-1}acac^{-1}a^{-2}cac^{-1} = 1$, and $a$ and $c$
do not commute, hence the group is not free.

\noindent\textbf{2210}. Wreath recursion: $a=\sigma(b,b)$,
$b=\sigma(c,a)$, $c=(a,a)$.

The element $(a^{-1}c)^2$ stabilizes the vertex $000$ and its
section at this vertex is equal to $a^{-1}c$. Hence, $a^{-1}c$ has
infinite order. Since $(b^{-1}a)^2\bigl|_{00}=a^{-1}c$ and $00$ is
fixed under the action of $b^{-1}a$ we obtain that $b^{-1}a$ also has
infinite order. Finally, $b^{-1}a$ stabilizes the vertex $11$ and
has itself as a section at this vertex. Therefore $G_{2210}$ is not
contracting.

We have $c^{-1}b^{-1}cb = (1, a^{-1}c^{-1}ac), bc^{-1}b^{-1}c = (ca^{-1}c^{-1}a, 1)$,
hence by Lemma~\ref{lem:freegroup} the group is not free.

\noindent\textbf{2212}. Klein bottle group, $\langle a,b\ \bigl|\
a^2=b^2\rangle$. Wreath recursion: $a=\sigma(a,c)$, $b=\sigma(c,a)$,
$c=(a,a)$.

The states $a$ and $c$ form a 2-state automaton generating the
infinite cyclic group $\Z$ in which $c=a^{-2}$ (see
Theorem~\ref{thm:class22}).

We have $a=\sigma(a,a^{-2})$, $b=\sigma(a^{-2},a)$, and
$x=ab^{-1}=(a^{-3},a^{3})$. Finally, since $x^a = b^{-1}a =
(a^{3},a^{-3}) = x^{-1}$, we have $G_{2212} = \langle x,a \mid x^a =
x^{-1}\rangle$ and $G_{2212}$ is the Klein bottle group. Tietze
transformations yield the presentation $G_{2212} = \langle a,b \mid
a^2 = b^2\rangle$ in terms of the generators $a$ and $b$.

\noindent\textbf{2213}. Wreath recursion: $a=\sigma(b,c)$,
$b=\sigma(c,a)$, $c=(a,a)$.

By Lemma~\ref{alg_trans} the element $cb$ acts transitively on the
levels of the tree and, hence, has infinite order. Since
$(ba)\bigl|_{100}=cb$ and the vertex $100$ is fixed under the action
of $ba$ we obtain that $ba$ also has infinite order. Finally, $ba$
stabilizes the vertex $01$ and has itself as a section at this
vertex. Therefore $G_{2213}$ is not contracting.

We have $c^{-1}b^{-1}cb = (1, a^{-1}c^{-1}ac), bc^{-1}b^{-1}c = (ca^{-1}c^{-1}a, 1)$,
hence by Lemma~\ref{lem:freegroup} the group is not free.

\noindent$\mathbf{2214}\cong G_{748}\cong D_4\times C_2$. Wreath
recursion: $a=\sigma(c,c)$, $b=\sigma(c,a)$, $c=(a,a)$.

Direct calculation.

\noindent$\mathbf{2226}\cong G_{820}\cong D_\infty$. Wreath
recursion: $a=\sigma(c,a)$, $b=\sigma(b,b)$, and $c=(a,a)$.

We have $ba=(bc,ba)$, $bc=\sigma(ba,ba)$, and $b=\sigma(b,b)$. Therefore $x$, $y$
and $b$ satisfy the wreath recursion defining the automaton $\A_{2394}$. Thus
$G_{2226}=G_{2394}\cong G_{820}$.

\noindent\textbf{2229}. Wreath recursion:
$a=\sigma(c, b)$, $b=\sigma(b, b)$, $c=(a, a)$.

Note that $b$ is of order 2. Post-conjugating the recursion by $(1,
b)$ (which is equivalent to conjugating by the tree automorphism
$g=(g,bg)$ in $\Aut(X^*)$ gives a copy of $G_{2229}$ defined by
\[a=\sigma(bc, 1), \qquad b=\sigma, \qquad c=(a, bab)\]
The stabilizer of the first level is generated by
\[a^2=(bc, bc), \qquad c=(a, bab), \qquad ba=(bc, 1), \qquad bcb=(bab, a).\]
Its projection on the first level is generated by
\[ bc=\sigma(a, bab), \qquad a=\sigma(bc, 1), \qquad bab=\sigma(1, bc).\]

Furthermore,
\[bcbc=(baba, abab), \qquad abab=(1, bcbc), \qquad baba=(bcbc, 1),\]
which implies that $bc$ is of order 2 and $a^{-1}=bab$. Hence, the
projection of the stabilizer on the first level is generated by the
recursion
\[ a=\sigma(bc, 1), \qquad bc=\sigma(a, a^{-1}).\]
Post-conjugating by $(1, a)$, we obtain the recursion
\[a=\sigma(a^{-1}\cdot bc, a), \qquad bc=\sigma,\]
which is the group $C_4\ltimes\Z^2$ of all orientation preserving
automorphisms of the integer lattice
(see~\cite{bartholdi_n:rabbit}). Note that the nucleus of $G_{2229}$
consists of 52 elements.

\noindent$\mathbf{2232}\cong G_{730}$. Klein Group $C_2\times C_2$.
Wreath recursion: $a=\sigma(c, c)$, $b=\sigma(b, b)$, $c=(a, a)$.

Direct calculation.

\noindent\textbf{2233}. Wreath recursion: $a=\sigma(a,a)$,
$b=\sigma(c,b)$, $c=(a,a)$.

Therefore, $\langle ba=(ba,ca),
ca=\sigma\rangle=G_{932}\cong D_\infty$.

Conjugating by $g=(ag, g)$, we obtain the recursion
\[
 \alpha=\sigma, \qquad
 \beta=\sigma(\gamma\beta, \alpha\beta), \qquad
 \gamma=(\alpha, \alpha), \]
where $\alpha = a^g$, $\beta = b^g$, and $\gamma = c^g$. Therefore
\[ \alpha =\sigma, \qquad
   \alpha\beta=(\gamma\alpha, \alpha\beta), \qquad
   \gamma\alpha=\sigma(\alpha, \alpha),\]
and the last wreath recursion defines a bounded automaton (see
Section~\ref{definition} for a definition). It follows
from~\cite{bkn:amenab} that $G_{2233}$ is amenable.

\noindent\textbf{2234}. Wreath recursion:
$a=\sigma(b, a)$, $b=\sigma(c, b)$, $c=(a, a)$.

The element $(c^{-1}b)^4$ stabilizes the vertex $00$ and its section
at this vertex is equal to $(c^{-1}b)^{-1}$. Hence, $c^{-1}b$ has
infinite order. Since $(b^{-1}a)\bigl|_{0}=c^{-1}b$ and $0$ is fixed
under the action of $b^{-1}a$ we obtain that $b^{-1}a$ also has
infinite order. Finally, $b^{-1}a$ stabilizes the vertex $1$ and has
itself as a section at this vertex. Therefore $G_{2234}$ is not
contracting.

We have $c^{-1}b^{-1}ac^{-1}a^{2} = (1, a^{-1}c^{-1}b^{2}),
ac^{-1}b^{-1}ac^{-1}a = (ba^{-1}c^{-1}b, 1)$, hence by
Lemma~\ref{lem:freegroup} the group is not free.

\noindent\textbf{2236}. Wreath recursion:
$a=\sigma(a, b)$, $b=\sigma(c, b)$, $c=(a, a)$.

By Lemma~\ref{alg_trans} the element $b$ acts transitively on the
levels of the tree and, hence, has infinite order.

By Lemma~\ref{alg_trans} the element $cb$ acts transitively on the
levels of the tree and, hence, has infinite order. Since
$ba=(ba,cb)$ we obtain that $ba$ also has infinite order.
Since $ba$ has itself as a section at $0$ the group is not
contracting.

We have $a^{-2}bab^{-2}ab = 1$, and $a$ and $b$ do not commute, hence the
group is not free.

\noindent\textbf{2237}. Wreath recursion:
$a=\sigma(b, b)$, $b=\sigma(c, b)$, $c=(a, a)$.

By Lemma~\ref{alg_trans} the elements $b$ and $(bc)^3$ acts transitively
on the levels of the tree and, hence, have infinite order.

Since $(cba)^2\bigl|_{00000}=(bc)^3$ and $00000$ is fixed under the
action of $(cba)^2$ we obtain that $cba$ also has infinite order.
Finally, $cba$ stabilizes the vertex $101$ and has itself as a
section at this vertex. Therefore $G_{2237}$ is not contracting.

We have $a^{-2}bab^{-2}ab = 1$, and $a$ and $b$ do not commute, hence the
group is not free.

\noindent\textbf{2239}. Wreath recursion:
$a=\sigma(a, c)$, $b=\sigma(c, b)$, $c=(a, a)$.

The group contains elements of infinite order by
Lemma~\ref{nontors}. In particular, $ca$ has infinite order. Since
$(ba)\bigl|_{100}=ca$ and the vertex $100$ is fixed under the action
of $ba$ we obtain that $ba$ also has infinite order. Finally, $ba$
stabilizes the vertex $1$ and has itself as a section at this
vertex. Therefore $G_{2239}$ is not contracting.

We have $ca^{-2}cba^{-1} = (1, c^{-1}abc^{-1}), a^{-1}ca^{-2}cb = (c^{-2}ab, 1)$,
hence by Lemma~\ref{lem:freegroup} the group is not free.

We can also simplify the wreath recursion in the following way.
Since $c=a^{-2}$ we have
\[a=\sigma(a, a^{-2}), \qquad b=\sigma(a^{-2}, b).\]
Therefore
\[ab=(a^{-4}, ab), \qquad a=\sigma(a, a^{-2}),\]
which can be written as
\[ab=(a^{-4}, ab), \qquad a=\sigma(1, a^{-1}),\]
which is a subgroup of
\[\beta=(a, \beta), \qquad a=\sigma(1, a^{-1}).\]

\noindent\textbf{2240}. Free group of rank $3$. Wreath recursion:
$a=\sigma(b, c)$, $b=\sigma(c, b)$, $c=(a, a)$.

The automaton appeared for the first time in~\cite{aleshin:free}.
The fact that $G_{2240}$ is free group of rank 3 with basis
$\{a,b,c\}$ is proved in~\cite{vorobets:aleshin}. This is the
smallest automaton among all automata over a $2$-letter alphabet
generating a free nonabelian group.

The fact that $G_{2240}$ is not contracting follows now from the
result of Nekrashevych~\cite{nekrashevych:free_subgroups}, that a
contracting group cannot have free subgroups. Alternatively,
$b^{-1}ca$ has infinite order, stabilizes the vertex $11$ and has
itself as a section at this vertex. Hence, the group is not contracting.

\noindent$\mathbf{2241}\cong G_{739}\cong C_2\ltimes (C_2\wr \Z)$.
Wreath recursion: $a=\sigma(c,c)$, $b=\sigma(c,b)$, $c=(a,a)$.

Consider $G_{747}$. Its wreath recursion is given by
$a=\sigma(c,c)$, $b=(b,a)$, $c=(a,a)$. All generators have order $2$
and $a$ commutes with $c$. Therefore
$acb=\sigma(cab,c)=\sigma(acb,c)$ and wa have $G_{747}=\langle
a,acb,c\rangle=G_{2241}$. Thus $G_{2241} = G_{747}\cong G_{739}$.

\noindent$\mathbf{2260}\cong G_{802}\cong C_2\times C_2\times C_2$.
Wreath recursion: $a=\sigma(a,a)$, $b=(c,c)$, $c=(a,a)$.

Direct calculation.

\noindent\textbf{2261}. Wreath recursion: $a=\sigma(b,a)$,
$b=\sigma(c,c)$, $c=(a,a)$.

The element $(ac^{-1})^2$ stabilizes the vertex $00$ and its section
at this vertex is equal to $(ac^{-1})^{-1}$. Hence, $ac^{-1}$ and
$c^{-1}a$ have infinite order.

Since $b^{-1}c^{-1}ac^{-1}ba\bigl|_{001}=\bigl((c^{-1}a)^2\bigr)^a$
and the vertex $001$ is fixed under the action of
$b^{-1}c^{-1}ac^{-1}ba$ we obtain that $b^{-1}c^{-1}ac^{-1}ba$ also
has infinite order. Finally, $b^{-1}c^{-1}ac^{-1}ba$ stabilizes the
vertex $000$ and has itself as a section at this vertex. Therefore
$G_{2261}$ is not contracting.

We have $acac^{-1}a^{-2}cac^{-1}aca^{-2}c^{-1} = 1$, and $a$ and $c$
do not commute, hence the group is not free.

\noindent$\mathbf{2262}\cong G_{848}\cong C_2\wr\Z$. Wreath
recursion: $a=\sigma(c,a)$, $b=\sigma(c,c)$, $c=(a,a)$.

The states $a$ and $c$ form a 2-state automaton (see
Theorem~\ref{thm:class22}). Moreover, $c=a^{-2}$ and $a$ has
infinite order.

Thus $a=\sigma(a^{-2},a)$, $b=\sigma(a^{-2},a^{-2})$ and $G_{2262}=
\langle a,b \rangle$. Further, $b^{-1}a=(1,a^3)$ and
$a^{-3}=\sigma(1,a^3)$, yielding $a^{-4}b=\sigma$. Therefore
$G=\langle a,\sigma\rangle$. Since $\langle a,a^\sigma \rangle =
\Z^2$, we obtain that $G_{2262} \cong C_2 \wr Z^2 \cong G_{848}$.

\noindent$\mathbf{2264}\cong G_{730}$. Klein Group $C_2\times C_2$.
Wreath recursion: $a=\sigma(b,b)$, $b=\sigma(c,c)$, $c=(a,a)$.

Direct calculation.

\noindent\textbf{2265}. Wreath recursion: $a=\sigma(c,b)$,
$b=\sigma(c,c)$, $c=(a,a)$.

The element $(c^{-1}b)^4$ stabilizes the vertex $0000$ and its
section at this vertex is equal to
$\bigl((c^{-1}b)^{-1}\bigr)^{c^{-1}a}$. Hence, $c^{-1}b$ has
infinite order. Since $[c,a]\bigl|_{10}=(c^{-1}b)^c$ and $10$ is
fixed under the action of $[c,a]$ we obtain that $[c,a]$ also has
infinite order. Finally, $[c,a]$ stabilizes the vertex $00$ and has
itself as a section at this vertex. Therefore $G_{2265}$ is not
contracting.

We have $a^{-2}bab^{-2}ab = 1$, and $a$ and $b$ do not commute, hence the
group is not free.

\noindent\textbf{2271}. Wreath recursion: $a=\sigma(c,a)$,
$b=\sigma(a,a)$, $c=(b,a)$.

The element $(ac^{-1})^4$ stabilizes the vertex $001$ and its
section at this vertex is equal to $ac^{-1}$. Hence, $ac^{-1}$ has
infinite order.

The element $(a^{-1}b)^4$ stabilizes the vertex $000$ and its
section at this vertex is equal to $a^{-1}b$. Hence, $a^{-1}b$ has
infinite order. Since
$b^{-1}c^{-1}ac^{-1}a^{2}\bigl|_{001}=(a^{-1}b)^a$ and the vertex
$001$ is fixed under the action of $b^{-1}c^{-1}ac^{-1}a^{2}$ we
obtain that $b^{-1}c^{-1}ac^{-1}a^{2}$ also has infinite order.
Finally, $b^{-1}c^{-1}ac^{-1}a^{2}$ stabilizes the vertex $000$ and
has itself as a section at this vertex. Therefore $G_{2271}$ is not
contracting.

We have $a^{-2}bab^{-2}ab = 1$, and $a$ and $b$ do not commute, hence the
group is not free.

\noindent\textbf{2274}. Wreath recursion: $a=\sigma(c,b)$,
$b=\sigma(a,a)$, $c=(b,a)$.

The element $a^{-1}c=\sigma(1,c^{-1}a)$ is conjugate to the adding
machine and has infinite order. Since $(b^{-1}a)\bigl|_{0}=a^{-1}c$
and $0$ is fixed under the action of $b^{-1}a$ we obtain that
$b^{-1}a$ also has infinite order. Finally, $b^{-1}a$ stabilizes the
vertex $11$ and has itself as a section at this vertex. Therefore
$G_{2274}$ is not contracting.

We have $bc^{-2}b = (1, ab^{-2}a), b^{2}c^{-2} = (a^{2}b^{-2}, 1)$,
hence by Lemma~\ref{lem:freegroup} the group is not free.

\noindent$\mathbf{2277}\cong C_2\ltimes (\Z\times\Z)$. Wreath
recursion: $a=\sigma(c,c)$, $b=\sigma(a,a)$, $c=(b,a)$.

All generators have order $2$. Let $x=cb$, $y=ab$ and $H=\langle x,y
\rangle$. We have $x=\sigma(1,y^{-1})$ and $y=(xy^{-1},xy^{-1})$.
The elements $x$ and $y$ commute and the matrix of the associated
virtual endomorphism is given by
\[
 A =
 \begin{pmatrix}
   0 & 1 \\ -1/2 & -1
 \end{pmatrix}.
\]
The eigenvalues $-\frac12 \pm \frac12 i$ are not algebraic integers,
and therefore, according to~\cite{nekrash_s:12endomorph}, $H$ is
free abelian of rank 2.

The subgroup $H$ is normal of index 2 in $G_{2277}$. Therefore
$G_{2277}=\langle H,b\rangle = C_2\ltimes(\Z\times\Z)$, where $C_2$
is generated by $b$, which acts on $H$ is inversion of the
generators.

\noindent\textbf{2280}.  Wreath recursion: $a=\sigma(c,a)$,
$b=\sigma(b,a)$, $c=(b,a)$.

We prove that $a$ has infinite order by considering the forward
orbit of $10^\infty$ under the action of $a^2$. We have
\begin{alignat*}{5}
 a^2 &= (ac,ca), \qquad &&ac &&= \s (cb,a^2), \qquad &&ca &&=\s(ac,ba)
 \\
 cb &= \s(ab,ba), &&ba &&= (ac,ba),  &&ab &&=(ab,ca).
\end{alignat*}
The equalities
\begin{align*}
 a^2(10^\infty) &= ab(10^\infty) = 1110^\infty, \\
 ac(10^\infty) &= ca(10^\infty) = cb(10^\infty) = 0010^\infty, \text{ and}\\
 ba(10^\infty) &= 10110^\infty
\end{align*}
show that all members of the forward orbit of $10^\infty$ under the
action of $a^2$ have only finitely many 1's and that the position of
the rightmost 1 cannot decrease under the action of $a^2$. Since
$a^2(10^\infty)=1110^\infty$, the forward orbit of $10^\infty$ under
the action of $a^2$ can never return to $10^\infty$ and $a^2$ has
infinite order.

Since $a^2=(ac,ca)$, the elements $ca$ and $ab=(ab,ca)$ have
infinite order, showing that $G_{2280}$ is not contracting.

\noindent\textbf{2283}. Wreath recursion: $a=\sigma(c,b)$,
$b=\sigma(b,a)$, $c=(b,a)$.

By Lemma~\ref{alg_trans} the element $ac$ acts transitively on the
levels of the tree and, hence, has infinite order. Since
$ba=(ac,b^2)$ we obtain that $ba$ also has infinite order. Finally,
$ba$ stabilizes the vertex $11$ and has itself as a section at this
vertex. Therefore $G_{2283}$ is not contracting.

\noindent\textbf{2284}. Wreath recursion: $a=\sigma(a,c)$,
$b=\sigma(b,a)$, $c=(b,a)$.

Define $u=b^{-1}a$, $v=a^{-1}c$ and $w=c^{-1}b$. Then $u=(u,v)$,
$v=\sigma(w,1)$ and $w=\sigma(u^{-1},u)$. The group $\langle
u,v,w\rangle$ is generated by the automaton symmetric to the one
generating the subgroup $\langle x,y,z\rangle$ of $G_{960}$ (see
$G_{960}$ for the definition). It is shown above that $zxy$ has
infinite order. Therefore $wvu$ also has infinite order.

The element $(b^{-1}ac^{-1}a)^2$ stabilizes the vertex $00$ and its
section at this vertex is equal to $(b^{-1}ac^{-1}a)^{a^{-1}b}$.
Hence, $b^{-1}ac^{-1}a$ has infinite order. Let
$t=b^{-1}ab^{-2}a^{2}$. Since $t|_{110}=b^{-1}ac^{-1}a$ and the
vertex $110$ is fixed under the action of $t$ we see that $t$ also
has infinite order. Finally, $t$ stabilizes the vertex $11101000$
and has itself as a section at this vertex. Therefore $G_{2284}$ is
not contracting.

\noindent\textbf{2285}. Wreath recursion: $a=\sigma(b,c)$,
$b=\sigma(b,a)$, $c=(b,a)$.

The element $ac^{-1}=\sigma(1,ca^{-1})$ is conjugate to the adding
machine and has infinite order.

By Lemma~\ref{alg_trans} the element $abcb$ acts transitively on the
levels of the tree and, hence, has infinite order. Since
$(ba)^2\bigl|_{000}=(ac,b^2)$ and the vertex $000$ is fixed under
the action of $(ba)^2$ we obtain that $ba$ also has infinite order.
Finally, $ba$ stabilizes the vertex $01$ and has itself as a section
at this vertex. Therefore $G_{2285}$ is not contracting.

\noindent\textbf{2286}. Wreath recursion: $a=\sigma(c,c)$,
$b=\sigma(b,a)$, $c=(b,a)$.

The element $(c^{-1}a)^2$ stabilizes the vertex $00$ and its section
at this vertex is equal to $(c^{-1}a)^{a^{-1}b}$. Hence, $c^{-1}a$
has infinite order. Since
$(c^{-2}a^2)\bigl|_{000}=(c^{-1}a)^{b^{-1}}$ and $000$ is fixed
under the action of $c^{-2}a^2$ we obtain that $c^{-2}a^2$ also has
infinite order. Finally, $c^{-2}a^2$ stabilizes the vertex $11$ and
has itself as a section at this vertex. Therefore $G_{2286}$ is not
contracting.

\noindent\textbf{2287}. Wreath recursion: $a=\sigma(a,a)$,
$b=\sigma(c,a)$, $c=(b, a)$.

The element $bc^{-1}=\sigma(cb^{-1},1)$ is conjugate to the adding
machine and has infinite order.

Conjugating the generators by $g=(g, ag)$, we obtain the wreath
recursion
\[a'=\sigma, \qquad b'=\sigma(a'c', 1), \qquad c'=(b', a'),\]
where $a'=a^g$, $b'=b^g$, and $c'=c^g$. Therefore
\[a'=\sigma, \qquad b'=\sigma(a'c', 1), \qquad a'c'=\sigma(b', a')\]

A direct computation shows that the iterated monodromy group of
$\frac{z^2+2}{1-z^2}$ is generated by
\[
 \alpha=\sigma, \qquad \beta=\sigma(\gamma^{-1}\beta^{-1}, \alpha),
 \qquad \gamma=(\beta\gamma\beta^{-1}, \alpha),
\]
where $\alpha$, $\beta$, and $\gamma$ are loops around the
post-critical points $2$, $-1$ and $-2$, respectively (recall the
definition of iterated monodromy group in
Section~\ref{s:contracting}). We see that
\[
 \alpha=\sigma, \qquad \beta\gamma=\sigma(\beta^{-1}, 1), \qquad
 \beta=\sigma(\gamma^{-1}\beta^{-1}, \alpha)
\]
satisfy the same recursions as $a$, $b$ and $ac$, only composed with
taking inverses. If we take second iteration of the wreath
recursions, we obtain isomorphic self-similar groups.

It follows that the group $G_{2287}$ is isomorphic to
$IMG\left(\frac{z^2+2}{1-z^2}\right)$ and the limit space is
homeomorphic to the Julia set of this rational function.

\noindent\textbf{2293}. Wreath recursion: $a=\sigma(a,c)$,
$b=\sigma(c,a)$, $c=(b,a)$.

The element $(b^{-1}c)^2$ stabilizes the vertex $0$ and its section
at this vertex is equal to $(b^{-1}c)^{-1}$. Hence, $b^{-1}c$ has
infinite order. Since $(c^{-1}bc^{-1}a)^2\bigl|_{000}=b^{-1}c$ and
$000$ is fixed under the action of $(c^{-1}bc^{-1}a)^2$ we obtain
that $c^{-1}bc^{-1}a$ also has infinite order. Finally,
$c^{-1}bc^{-1}a$ stabilizes the vertex $11$ and has itself as a
section at this vertex. Therefore $G_{2293}$ is not contracting.

We have $b^{-1}c^{2}a^{-1} = (1, c^{-1}b^{2}c^{-1}), c^{2}a^{-1}b^{-1} = (b^{2}c^{-2}, 1)$,
hence by Lemma~\ref{lem:freegroup} the group is not free.

\noindent\textbf{2294}. Baumslag-Solitar group $BS(1,-3)$. Wreath
recursion: $a=\sigma(b,c)$, $b=\sigma(c,a)$, $c=(b,a)$.

The automaton satisfies the conditions of Lemma~\ref{nontors}.
Therefore $cb$ has infinite order. Since $a^2=(cb,bc)$, $c=(b,a)$
and $ba=(ab,c^2)$, the elements $a$, $c$ and $ba$ have infinite
order. Finally, $ba$ fixes the vertex $01$ and has itself as a
section at this vertex, showing that $G_{2294}$ is not contracting.

Let $\mu = ca^{-1}$. We have
$\mu=ca^{-1}=\sigma(ac^{-1},1)=\sigma(\mu^{-1},1)$, and therefore
$\mu$ is conjugate of the adding machine and has infinite order.
Further, we have $bc^{-1} = \s(cb^{-1},1) = \s((bc^{-1})^{-1},1)$,
showing that $bc^{-1} = \mu = ca^{-1}$. Therefore $G_{2294} =
\langle \mu, a \rangle$.

It can be shown that $a\mu a^{-1}=\mu^{-3}$ in $G_{2294}$ (compare
to $G_{870}$. Since both $a$ and $\mu$ have infinite order $G_{2294}
\cong BS(1,-3)$.

\noindent\textbf{2295}. Wreath recursion: $a=\sigma(c,c)$,
$b=\sigma(c,a)$, $c=(b,a)$.

The element $cb^{-1}=\sigma(1,bc^{-1})$ is conjugate to the adding
machine and has infinite order. Hence, its conjugate
$a^{-1}cb^{-1}a$ also has infinite order. Since
$c^{-1}ac^{-1}b=\bigl(c^{-1}ac^{-1}b,a^{-1}cb^{-1}a\bigr)$, the
element $c^{-1}ac^{-1}b$ has infinite order and $G_{2295}$ is not
contracting.

We have $a^{-2}bab^{-2}ab = 1$, and $a$ and $b$ do not commute, hence the
group is not free.

\noindent\textbf{2307}. Contains $G_{933}$. Wreath recursion:
$a=\sigma(c,a)$, $b=\sigma(b,b)$, $c=(b,a)$.

We have $ba=(bc,ba)$, and $bc=\sigma(1,ba)$. Therefore $G_{933}$ is
a subgroup of $G_{2307}$ (the wreath recursion for $ba$ and $bc$
defines an automaton that is symmetric to the one defining the
automaton [993]).

The element $(a^{-1}b)^2$ stabilizes the vertex $00$ and its section
at this vertex is equal to $a^{-1}b$. Hence, $a^{-1}b$ has infinite
order. Furthermore, $a^{-1}b$ stabilizes the vertex $1$ and has
itself as a section at this vertex. Therefore $G_{2307}$ is not
contracting.

\noindent$\mathbf{2313}\cong G_{2277}\cong C_2\ltimes (\Z\times\Z)$.
Wreath recursion: $a=\sigma(c,c)$, $b=\sigma(b,b)$, $c=(b,a)$.

Since all generators have order $2$ the subgroup $H=\langle
ba,bc\rangle$ is normal in $G_{2313}$. Furthermore,
$ba=\sigma(bc,bc)$ and $bc=\sigma(1,ba)$. Hence,
$H=G_{771}\cong\Z^2$.

Finally, $G_{2313}=\langle H,b\rangle=\langle b\rangle\ltimes
H=C_2\ltimes(\Z\times\Z)$, where $b$ inverts the generators of $H$.
This action coincides with the one for $G_{2277}$, which proves that
these groups are isomorphic.

\noindent$\mathbf{2320}\cong G_{2294}$. Baumslag-Solitar group
$BS(1,-3)$. Wreath recursion: $a=\sigma(a,c)$, $b=\sigma(c,b)$,
$c=(b,a)$.

It is proved in~\cite{bartholdi_s:bsolitar} that the automaton
[2320] generates $BS(1,-3)$.

\noindent\textbf{2322}. Wreath recursion: $a=\sigma(c,c)$,
$b=\sigma(c,b)$, $c=(b,a)$.

The element $(a^{-1}c)^2$ stabilizes the vertex $00$ and its section
at this vertex is equal to $(a^{-1}c)^{b^{-1}}$. Hence, $a^{-1}c$
has infinite order. Since $(c^{-2}a^2)^2\bigl|_{000}=a^{-1}c$ and
$000$ is fixed under the action of $c^{-2}a^2$ we obtain that
$c^{-2}a^2$ also has infinite order. Finally, $c^{-2}a^2$ stabilizes
the vertex $11$ and has itself as a section at this vertex.
Therefore $G_{2322}$ is not contracting.

We have $a^{-2}bab^{-2}ab = 1$, and $a$ and $b$ do not commute, hence the
group is not free.

\noindent$\mathbf{2352}\cong G_{740}$. Wreath recursion:
$a=\sigma(c,a)$, $b=\sigma(a,a)$, $c=(c,a)$.

We have $ac^{-1}b=(a,a)$. Therefore $G_{2352} = \langle a,ac^{-1}b,c
\rangle = G_{740}$.

\noindent\textbf{2355}. Wreath recursion: $a=\sigma(c,b)$,
$b=\sigma(a,a)$, $c=(c,a)$.

The element $(b^{-1}a)^2$ stabilizes the vertex $00$ and its section
at this vertex is equal to $(b^{-1}a)^{a^{-1}c}$. Hence, $b^{-1}a$
has infinite order. Furthermore, $b^{-1}a$ stabilizes the vertex
$11$ and has itself as a section at this vertex. Therefore
$G_{2355}$ is not contracting.

We have $a^{-1}cb^{-1}c = (b^{-1}c, 1), cb^{-1}ca^{-1} = (1, cb^{-1})$,
hence by Lemma~\ref{lem:freegroup} the group is not free.

\noindent$\mathbf{2358}\cong G_{820}\cong D_\infty$. Wreath
recursion: $a=\sigma(c,c)$, $b=\sigma(a,a)$, $c=(c,a)$.

The states $a$ and $c$ form a $2$-state automaton generating
$D_\infty$ (see Theorem~\ref{thm:class22}) and $b=aca$.

\noindent\textbf{2361}. Wreath recursion: $a=\sigma(c,a)$,
$b=\sigma(b,a)$, $c=(c,a)$.

The element $bc^{-1}=\sigma(bc^{-1},1)$ is conjugate to the adding
machine and has infinite order.

\noindent\textbf{2364}. Wreath recursion: $a=\sigma(c,b)$,
$b=\sigma(b,a)$, $c=(c,a)$.

The element $cb^{-1}=\sigma(1,cb^{-1})$ is the adding machine and
has infinite order. Therefore its conjugate $b^{-1}c$ also has
infinite order. Since $(b^{-1}a)\bigl|_{0}=b^{-1}c$ and $0$ is fixed
under the action of $b^{-1}a$ we obtain that $b^{-1}a$ also has
infinite order. Finally, $b^{-1}a$ stabilizes the vertex $11$ and
has itself as a section at this vertex. Therefore $G_{2364}$ is not
contracting.

We have $c^{-1}ac^{-1}b = (1, a^{-1}bc^{-1}b), bc^{-1}ac^{-1} = (ba^{-1}bc^{-1}, 1)$,
hence by Lemma~\ref{lem:freegroup} the group is not free.

\noindent\textbf{2365}. Wreath recursion: $a=\sigma(a,c)$,
$b=\sigma(b,a)$, $c=(c,a)$.

By Lemma~\ref{alg_trans} the element $cb$ acts transitively on the
levels of the tree and, hence, has infinite order.

\noindent\textbf{2366}. Wreath recursion: $a=\sigma(b,c)$,
$b=\sigma(b,a)$, $c=(c,a)$.

By Lemma~\ref{alg_trans} the element $a$ acts transitively on the
levels of the tree and, hence, has infinite order. Since $c=(c,a)$
we obtain that $c$ also has infinite order and $G_{2366}$ is not
contracting.

We have $a^{-2}bab^{-2}ab = 1$, and $a$ and $b$ do not commute, hence the
group is not free.

\noindent\textbf{2367}. Wreath recursion: $a=\sigma(c,c)$,
$b=\sigma(b,a)$, $c=(c,a)$.

The states $a$ and $c$ form  a $2$-state automaton generating
$D_\infty$ (see Theorem~\ref{thm:class22}).

Also we have $bc=\sigma(bc,1)$ and $ca=\sigma(ac,1)$. Therefore the
elements $bc$ and $ca$ generate the Brunner-Sidki-Vierra group
(see~\cite{brunner_sv:justnonsolv}).

\noindent$\mathbf{2368}\cong G_{739}\cong C_2\ltimes (C_2\wr \Z)$.
Wreath recursion: $a=\sigma(a,a)$, $b=\sigma(c,a)$, $c=(c,a)$.

We have $bc^{-1}a=(a,a)$. Therefore $G_{2368} = \langle a,c,bc^{-1}a
\rangle = G_{739}$.

\noindent\textbf{2369}. Wreath recursion: $a=\sigma(b,a)$,
$b=\sigma(c,a)$, $c=(c,a)$.

By using the approach already used for $G_{875}$, we can show that
the forward orbit of $10^\infty$ under the action of $a$ is
infinite, and therefore $a$ has infinite order.

Since $a^2=(ab,ba)$, the element $ab$ also has infinite order.
Furthermore, $ab$ fixes $00$ and has itself as a section at this
vertex. Therefore $G_{2369}$ is not contracting.

\noindent\textbf{2371}. Wreath recursion: $a=\sigma(a,b)$,
$b=\sigma(c,a)$, $c=(c,a)$.

The element $(c^{-1}ab^{-1}a)^2$ stabilizes the vertex $01$ and its
section at this vertex is equal to $c^{-1}ab^{-1}a$, which is
nontrivial. Hence, $c^{-1}ab^{-1}a$ has infinite order.

Let
$t=b^{-1}c^{-1}a^{2}c^{-1}ba^{-1}ca^{-1}ca^{-2}cbc^{-1}ab^{-1}a$.
Then $t^2$ stabilizes the vertex $00$ and
$t^2\bigl|_{00}=t^{a^{-1}ba^{-1}c}$ . Hence, $t$ has infinite order.
Let $s=b^{-1}c^{-2}a^3$ Since $s^8|_{00100001}=t$ and $s$ fixes the
vertex $00100001$ we see that $s$ also has infinite order. Finally,
$s$ stabilizes the vertex $11$ and has itself as a section at this
vertex. Therefore $G_{2371}$ is not contracting.

\noindent\textbf{2372}. Wreath recursion: $a=\sigma(b,b)$,
$b=\sigma(c,a)$, $c=(c,a)$.

By Lemma~\ref{alg_trans} the elements $b$ and $ac$ act transitively
on the levels of the tree and, hence, have infinite order. Since
$(c^2)\bigl|_{100}=ac$ and the vertex $100$ is fixed under the
action of $c^2$ we obtain that $c$ also has infinite order. Finally,
$c$ stabilizes the vertex $0$ and has itself as a section at this
vertex. Therefore $G_{2372}$ is not contracting.

\noindent$\mathbf{2374}\cong G_{821}$. Lamplighter group $\Z \wr
C_2$. Wreath recursion: $a=\sigma(a,c)$, $b=\sigma(c,a)$, $c=(c,a)$.

The states $a$ and $c$ form a 2-state automaton that generates the
Lamplighter group (see Theorem~\ref{thm:class22}). Since
$bc^{-1}=\sigma=c^{-1}a$, we have $b=a^c$ and $G=\langle
a,c\rangle$.

\noindent\textbf{2375}. Wreath recursion: $a=\sigma(b,c)$,
$b=\sigma(c,a)$, $c=(c,a)$.

The element $(a^{-1}c)^2$ stabilizes the vertex $01$ and its section
at this vertex is equal to $a^{-1}c$. Hence, $a^{-1}c$ and $c^{-1}a$
have infinite order. Since
$c^{-1}b^{-1}ac^{-1}a^{2}\bigl|_{00}=c^{-1}a$ and the vertex $00$ is
fixed under the action of $c^{-1}b^{-1}ac^{-1}a^{2}$ we obtain that
$c^{-1}b^{-1}ac^{-1}a^{2}$ also has infinite order. Finally,
$c^{-1}b^{-1}ac^{-1}a^{2}$ stabilizes the vertex $11$ and has itself
as a section at this vertex. Therefore $G_{2375}$ is not
contracting.

\noindent$\mathbf{2376}\cong G_{739}\cong C_2\ltimes (C_2\wr \Z)$.
Wreath recursion: $a=\sigma(c,c)$, $b=\sigma(c,a)$, $c=(c,a)$.

Since $\s = bc^{-1}$, we have $G_{2376} = \langle a,c,\sigma
\rangle$. We already proved that $G_{972} = \langle a,c,\s\rangle$.
Therefore $G_{2376} = G_{972}\cong G_{739}$.

\noindent$\mathbf{2388}\cong G_{821}$. Lamplighter group $\Z \wr
C_2$. Wreath recursion: $a=\sigma(c,a)$, $b=\sigma(b,b)$, $c=(c,a)$.

The states $a$ and $c$ form a 2-state automaton generating the
Lamplighter group (see Theorem~\ref{thm:class22}) and
$b=\s=ac^{-1}$.

\noindent\textbf{2391}. Wreath recursion: $a=\sigma(c,b)$,
$b=\sigma(b,b)$, $c=(c,a)$.

The element $(c^{-1}ba^{-1}b)^2$ stabilizes the vertex $00$ and its
section at this vertex is equal to $c^{-1}ba^{-1}b$. Hence,
$c^{-1}ba^{-1}b$ has infinite order. Since
$(bc^{-2}b)^2\bigl|_{000}=c^{-1}ba^{-1}b$ and $000$ is fixed under
the action of $bc^{-2}b$ we obtain that $bc^{-2}b$ also has infinite
order. Finally, $bc^{-2}b$ stabilizes the vertex $1$ and has itself
as a section at this vertex. Therefore $G_{2391}$ is not
contracting.

\noindent$\mathbf{2394}\cong G_{820}\cong D_\infty$. Wreath
recursion: $a=\sigma(c,c)$, $b=\sigma(b,b)$, $c=(c,a)$.

All generators have order $2$, hence $H=\langle ba,bc\rangle$ is
normal in $G_{2394}$. Furthermore, $ba=(bc,bc)$, $bc=\sigma(bc,ba)$,
and therefore $H=G_{731}\cong\Z$. Thus $G_{2394}=\langle
b\rangle\ltimes H\cong C_2\ltimes\Z\cong D_\infty$ since
$(bc)^b=(bc)^{-1}$.

\noindent\textbf{2395}. Wreath recursion: $a=\sigma(a,a)$,
$b=\sigma(c,b)$, $c=(c,a)$.

By Lemma~\ref{alg_trans} the element $ca$ acts transitively on the
levels of the tree.

The element $(c^{-1}a)^2$ stabilizes the vertex $0$ and its section
at this vertex is equal to $c^{-1}a$. Hence, $c^{-1}a$ has infinite
order. Since $(b^{-1}a)\bigl|_{0}=c^{-1}a$ and $0$ is fixed under
the action of $b^{-1}a$ we obtain that $b^{-1}a$ also has infinite
order. Finally, $b^{-1}a$ stabilizes the vertex $1$ and has itself
as a section at this vertex. Therefore $G_{2395}$ is not
contracting.

Note that $ab=(ac, ab), ac=\sigma(ac, 1)$ and $ba=(ba, ca),
ca=\sigma(1, ca)$, i.e., $G_{2395}$ contains copies of $G_{929}$.

\noindent\textbf{2396}. Boltenkov group. Wreath recursion:
$a=\sigma(b,a)$, $b=\sigma(c,b)$, $c=(c,a)$.

This group was studied by A.~Boltenkov (under direction of
R.~Grigorchuk), who showed that the monoid generated by $\{a,b,c\}$
is free, and the group $G_{2396}$ is torsion free.

\begin{prop}\label{free}
The monoid generated by $a$, $b$, and $c$ is free.
\end{prop}
\begin{proof}
By way of contradiction, assume that there are some relations and
let $w=u$ be a relation for which $\max(|w|,|u|)$ minimal.

We first consider the case when neither $w$ nor $u$ is empty.
Because of cancelation laws, the words $w$ and $u$ must end in
different letters. We have
$w=\sigma_w(w_0,w_1)=\sigma_u(u_0,u_1)=u$, where $\sigma_w$, and
$\sigma_u$ are permutations in $\{1,\s\}$. Clearly, $w_0=u_0$ and
$w_1=u_1$ must also be relations.

Assume that $w$ ends in $b$ and $u$ ends in $c$. Then $w_0$ and
$u_0$ both end in $c$. Therefore, by minimality, $w_0=u_0$ as words
and $|u|=|w|$. Since $b \neq c$ in $G_{2396}$ the length of $w$ and
$u$ is at least 2. We can recover the second to last letter in $w$
and $u$. Indeed, the second to last letter in $u_0$ can be only $b$
or $c$ (these are the possible sections at 0 of the three
generators), while the second to last letter of $w_0$ can be only
$a$ or $b$ (these are the possible sections at 1 of the three
generators). Therefore $w_0=u_0=\dots bc$, $w=\dots bb$ , and
$u=\dots ac$. Since $bb\neq ac$ in $G_{2396}$ (look at the action at
level 1), the length of $w$ and $u$ must be at least 3. Continuing
in the same fashion we obtain that $w_0=u_0=b\dots bbc$, $w=\dots
ababb$, and $u=\dots babac$. Since the lengths of $w$ and $u$ are
equal, they have different action on level 1, which is a
contradiction.

Assume that $w$ ends in $a$ and $u$ ends in $b$ or $c$. Then $u_0$
and $w_0$ end in $b$ and $c$, respectively, and we may proceed as
before.

It remains to show that, say, $u$ cannot be empty word. If this is
the case then $w_0=1=w_1$, implying that $w_0=w_1$ is also a minimal
relation. But this is impossible since both $w_0$ and $w_1$ are
nonempty.
\end{proof}

For a group word $w$ over $\{a,b,c\}$, define the exponent
$\exp_a(w)$ of $a$ in $w$ as the sum of the exponents in all
occurrences of $a$ and $a^{-1}$ in $w$. Define $\exp_b(w)$ and
$\exp_c(w)$ in analogous way and let $\exp(w)=\exp_a(w)+\exp_b(w) +
\exp_c(w)$.

\begin{lemma}\label{exp0}
If $w=1$ in $G_{2396}$ then $\exp(w)=0$.
\end{lemma}
\begin{proof}
By way of contradiction, assume otherwise and choose a freely
reduced group word $w$ over $\{a,b,c\}$ such that $w=1$ in
$G_{2396}$, $\exp(w)\neq 0$, and $w$ has minimal length among such
words. If $w=(w_0,w_1)$, $w_0$ and $w_1$ also represent $1$ in
$G_{2396}$ and $\exp(w_0)=\exp(w_1)=\exp(w)\neq 0$. Since the
exponents is nonzero, the words $w_0$ and $w_1$ are nonempty and, by
minimality, their length must be equal to $|w|$. Note that
$ac^{-1}=\sigma(bc^{-1},1)$ and $bc^{-1}=\sigma(1,ba^{-1})$. This
implies that $w$ cannot $ac^{-1}$, $bc^{-1}$, $ca^{-1}$, or
$cb^{-1}$ as a subword (otherwise the length of $w_0$ or $w_1$ would
be shorter than the length of $w$). By the same reason, $w_0$ and
$w_1$ cannot have the above 4 words as subwords, which implies that
$w$ does not have $ab^{-1}=(ab^{-1},bc^{-1})$ or its inverse
$ba^{-1}$ as a subword. Therefore $w$ has the form
$w=W_1(a^{-1},b^{-1},c^{-1})W_2(a,b,c)$, and since $w=1$ in
$G_{2396}$, we obtain a relation between positive words over
$\{a,b,c\}$, which contradicts Proposition~\ref{free}.
\end{proof}

\begin{lemma}
If $w=1$ in $G_{2396}$ then $\exp_a(w)$, $\exp_b(w)$ and $\exp_c(w)$
are even.
\end{lemma}

\begin{proof}
Indeed, $\exp_a(w)+\exp_b(w)$ must be even (since both $a$ and $b$
are active at the root). By Lemma~\ref{exp0}, $\exp_c(w)$ must be
even. If $w=(w_0,w_1)$, then $\exp_a(w_0)+\exp_b(w_0)$ and
$\exp_a(w_1)+\exp_b(w_1)$ must be even. Since
$\exp_a(w)+\exp_b(w)=\exp_b(w_0)+\exp_b(w_1)$,
$\exp_a(w)+\exp_c(w)=\exp_a(w_0)+\exp_a(w_1)$ we obtain that
$2\exp_a(w) + \exp_b(w) + \exp_c(w)$ is even, which then implies
that $\exp_b(w)$ is even. Finally, since both $\exp_b(w)$ and
$\exp_c(w)$ are even, $\exp_a(w)$ must be even as well (by
Lemma~\ref{exp0}).
\end{proof}

\begin{prop}
\label{boltenkov_torsion_free} The group $G_{2396}$ is torsion free.
\end{prop}
\begin{proof}
By way of contradiction, assume otherwise. Let $w$ be an element of
order $2$. We may assume that $w$ does not belong to the stabilizer
of the first level (otherwise we may pass to a section of $w$). Let
$w=\sigma(w_0,w_1)$. Since $w^2=(w_1w_0,w_0w_1)=1$, we have the
modulo 2 equalities $\exp_b(w_0w_1) = \exp_b(w_0) + \exp_b(w_1)
=\exp_a(w)+\exp_b(w)$. Since $\exp_b(w_0w_1)$ is even,
$\exp_a(w)+\exp_b(w)$ must be even, implying that $w$ stabilizes
level 1, a contradiction.
\end{proof}

Since $b^{-1}a=(c^{-1}b,b^{-1}a)$, the group $G_{2396}$ is not
contracting (our considerations above show that $b^{-1}a$ is not
trivial and therefore has infinite order).

We have $c^{-1}bc^{-1}a = (1, a^{-1}bc^{-1}b)$, $ac^{-1}bc^{-1} =
(ba^{-1}bc^{-1}, 1)$, hence by Lemma~\ref{lem:freegroup} the group
is not free.

\noindent\textbf{2398}. Dahmani group. Wreath recursion:
$a=\sigma(a,b)$, $b=\sigma(c,b)$, $c=(c,a)$.

This group is self-replicating, not contracting, weakly regular
branch group over its commutator subgroup. It was studied by Dahmani
in~\cite{dahmani:non-contr}.

\noindent\textbf{2399}. Wreath recursion: $a=\sigma(b,b)$,
$b=\sigma(c,b)$, $c=(c,a)$.

By Lemma~\ref{alg_trans} the elements $ca$ and
$c^{4}bc^{2}bc^{2}b^{2}cb^{2}cb^{3}acba^{2}$ act transitively on the
levels of the tree and, hence, have infinite order. Since
$(cba)^8\bigl|_{000010001}=c^{4}bc^{2}bc^{2}b^{2}cb^{2}cb^{3}acba^{2}$
and vertex $000010001$ is fixed under the action of $(cba)^8$ we
obtain that $cba$ also has infinite order. Finally, $cba$ stabilizes
the vertex $01001$ and has itself as a section at this vertex.
Therefore $G_{2399}$ is not contracting.

We have $a^{-2}bab^{-2}ab = 1$, and $a$ and $b$ do not commute, hence the
group is not free.

\noindent\textbf{2401}. Wreath recursion: $a=\sigma(a,c)$,
$b=\sigma(c,b)$ and $c=(c,a)$.

The states $a$ and $c$ form a $2$-state automaton generating the
Lamplighter group (see Theorem~\ref{thm:class22}). Hence, $G_{2401}$
is neither torsion, nor contracting and has exponential growth.

\noindent\textbf{2402}. Wreath recursion: $a=\sigma(b,c)$,
$b=\sigma(c,b)$, $c=(c,a)$.

The element $(bc^{-1})^2$ stabilizes the vertex $00$ and its section
at this vertex is equal to $bc^{-1}$. Hence, $bc^{-1}$ has infinite
order.

We have $c^{-2}ba = (1, a^{-2}b^{2})$, $ac^{-2}b = (ba^{-2}b, 1)$,
hence by Lemma~\ref{lem:freegroup} the group is not free.

\noindent$\mathbf{2403}\cong G_{2287}$. Wreath recursion:
$a=\sigma(c,c)$, $b=\sigma(c,b)$, $c=(c,a)$.

The states $a$ and $c$ form a $2$-state automaton generating
$D_\infty$ (see Theorem~\ref{thm:class22}).

Also we have $bc=\sigma(1,ba)$ and $ba=(bc,1)$. Therefore the
elements $bc$ and $ba$ generate the Basilica group $G_{852}$.

By conjugating by $g=(cg, g)$, we obtain
\[ a'=\sigma, \qquad b'=\sigma(1, c'b'), \qquad c'=(c', a'),\]
where $a' = a^g$, $b' = b^g$, and $c' = c^g$. Therefore
\[ a'=\sigma, \qquad b'=\sigma(1, c'b'), \qquad c'b'=\sigma(a', b'),\]
and $G_{2402}$ is isomorphic to $G_{2287}$, i.e., to
$IMG(\frac{z^2+2}{1-z^2})$.

\noindent$\mathbf{2422}\cong G_{820}\cong D_\infty$. Wreath
recursion: $a=\sigma(a,a)$, $b=\sigma(c,c)$, $c=(c,a)$.

The states $a$ and $c$ form a $2$-state automaton generating
$D_\infty$ (see Theorem~\ref{thm:class22}) and $b=aca$.

\noindent\textbf{2423}. Wreath recursion: $a=\sigma(b,a)$,
$b=\sigma(c,c)$, $c=(c,a)$.

Contains elements of infinite order by Lemma~\ref{nontors}. In
particular, $ac$ has infinite order. Since $c^2\bigl|_{100}=ac$ and
the vertex $100$ is fixed under the action of $c^2$ we obtain that
$c$ also has infinite order. Since $c=(c,a)$ the group is not
contracting.

We have $c^{-1}bc^{-1}a = (1, a^{-1}b)$, $ac^{-1}bc^{-1} = (ba^{-1},
1)$, hence by Lemma~\ref{lem:freegroup} the group is not free.

\noindent$\mathbf{2424}\cong G_{966}$. Wreath recursion
$a=\sigma(c,a)$, $b=\sigma(c,c)$, $c=(c,a)$.

We have $ac^{-1}b=(c,c)$. Therefore $G_{2424}=\langle a, ac^{-1}b, c
\rangle=G_{966}$.

\noindent$\mathbf{2426}\cong G_{2277}\cong C_2\ltimes (\Z\times\Z)$.
Wreath recursion: $a=\sigma(b,b)$, $b=\sigma(c,c)$, $c=(c,a)$.

Since all generators have order $2$ the subgroup $H=\langle
ab,cb\rangle$ is normal in $G_{2426}$. Furthermore, $ab=(bc,bc)$,
$cb=\sigma(ac,1)=\sigma(ab(cb)^{-1},1)$, so $H$ is self-similar.
Since $acb=bca$ in $G_{2426}$ we obtain $ab\cdot
cb=abcaab=aacbab=cb\cdot ab$, hence, $H$ is an abelian self-similar
$2$-generated group.

Consider the $\frac12$-endomorphism $\phi:\mathop{\rm
Stab}\nolimits_{H}(1)\to H$, given by $\phi(g)=h$, provided
$g=(h,*)$ and consider the linear map $A:\mathbb C^2\to\mathbb C^2$
induced by $\phi$. It has the following matrix representation with
respect to the basis corresponding to the generating set
$\{ab,cb\}$:
\[
 A=\left(
\begin{array}{cc}
  0 & \frac12\\
  -1 & -\frac12
\end{array}
\right).
\]
Its eigenvalues are not algebraic integers and, therefore,
by~\cite{nekrash_s:12endomorph}, $H$ is a free abelian group of rank
$2$.

Finally, $G_{2426}=\langle H,b\rangle=\langle b\rangle\ltimes
H=C_2\ltimes (\Z\times\Z)$, where $b$ inverts the generators of $H$.
This action coincides with the one for $G_{2277}$, which proves that
these groups are isomorphic.

\noindent\textbf{2427}. The element $(bc^{-1})^4$ stabilizes the
vertex $000$ and its section at this vertex is equal to $bc^{-1}$.
Hence, $bc^{-1}$ has infinite order.

We have $a^{-2}bab^{-2}ab = 1$, and $a$ and $b$ do not commute, hence the
group is not free.

\noindent$\mathbf{2838}\cong G_{848}\cong C_2\wr\Z$. Wreath
recursion: $a=\sigma(c,a)$, $b=\sigma(a,a)$, $c=(c,c)$.

Since $c$ is trivial, we have $G=\langle a, ba^{-1}\rangle $, where
$a=\sigma(1,a)$ is the adding machine and $ba^{-1}=(1,a)$. Therefore
$G_{2838}=G_{848}$.

\noindent\textbf{2841}. Wreath recursion: $a=\sigma(c,b)$,
$b=\sigma(a,a)$, $c=(c,c)$.

The element $c$ is trivial. Since $a^2=(b,b)$, $b^2=(a^2,a^2)$ and
$a^2$ is nontrivial, the elements $a$ and $b$ have infinite order.
Also we have $ba=(a,ab)$ and $ab=(ba,a)$, hence $ba$ has infinite
order and $G_{2841}$ is not contracting.

We claim that the monoid generated by $a$ and $b$ is free. Hence,
$G_{2841}$ has exponential growth.

\begin{proof}
We can first prove (analogous to $G_{2851}$) that $w\neq1$ for any
nonempty word $w\in\{a,b\}^*$.

By way of contradiction, let $w$ and $v$ be two nonempty words in
$\{a,b\}^*$ with minimal $|w|+|v|$ such that $w=v$ in $G_{2841}$.
Assume that $w$ ends with $a$ and $v$ ends with $b$. Consider the
following cases.
\begin{enumerate}
    \item If $w=a$ then $v|_0=1$ in $G_{2841}$ and $v|_0$ is nontrivial word.
    \item If $w$ ends with $a^2$ then $w|_1=v|_1$ in $G_{2841}$,
    $\bigl|w|_1\bigr|+\bigl|v|_1\bigr|<|w|+|v|$ and $w|_1$ ends with $b$, $v|_1$ with $a$.
    \item If $w$ ends with $ba$ and $v$ ends with $ab$, then $w|_1=v|_1$ in $G_{2841}$,
    $\bigl|w|_1\bigr|+\bigl|v|_1\bigr|<|w|+|v|$ (because $\bigl|v|_1\bigr|<|v|$) and $w|_1$
    ends with $b$, $v|_1$ with $a$.
    \item If $w$ ends with $ba$ and $v$ ends with $b$, then $w|_1=v|_1$ in
    $G_{2841}$, $\bigl|w|_1\bigr|+\bigl|v|_1\bigr|\leq|w|+|v|$ and $w|_1$
    ends with $ab$, $|v_1|$ with $a$. Therefore, words $v|_1$ and $w|_1$
    satisfy one of the first three cases.
\end{enumerate}
In all cases we obtain either a shorter relation, which contradicts
to our assumption, or a relation of the form $v=1$, which is also
impossible.
\end{proof}

There are non-trivial group relations, e.g. $a^{-1}b^{-1}a^{-2}ba^{-1}b^{-1}aba^{2}b^{-1}ab = 1$,
while $a$ and $b$ do not commute, hence the group is not free.

\noindent$\mathbf{2284}\cong G_{730}$. Klein Group $C_2\times C_2$.\

Direct calculation.

\noindent$\mathbf{2847}\cong G_{929}$. Wreath recursion:
$a=\sigma(c,a)$, $b=\sigma(b,a)$, $c=(c,c)$.

Since $c$ is trivial, the generator $a=\s(1,a)$ is the adding
machine and $b=\s(b,a)$. We have $ab =(ab,a)$. Therefore $G_{2847}=
\langle a,ab \rangle = G_{929}$.

\noindent\textbf{2850}. Wreath recursion: $a=\sigma(c,b)$,
$b=\sigma(b,a)$, $c=(c,c)$.

Since $c$ is trivial, we have $a^2=(b,b)$, $b^2=(ab,ba)$,
$ab=(b^2,a)$ and $ba=(a,b^2)$. Therefore the elements $a$, $b$, $ab$
and $ba$ have infinite order. Since $ba$ fixes the vertex $11$ and
has itself as a section at that vertex, $G_{2850}$ is not
contracting.

The group is regular weakly branch over $G'_{2850}$, since it is
self-replicating and $[b,a^2]=(1,[a,b])$.

Semigroup $\langle a,b\rangle$ is free. Hence, $G_{2850}$ has
exponential growth.
\begin{proof}
We can first prove (analogous $G_{2851}$) that $w\neq1$ for any
nonempty word $w\in\{a,b\}^*$.

By way of contradiction, let $w$ and $v$ be two nonempty words in
$\{a,b\}^*$ with minimal $|w|+|v|$ such that $w=v$ in $G_{2850}$.
Assume that $w$ ends with $a$ and $v$ ends with $b$. Consider the
following cases.
\begin{enumerate}
    \item If $w=a$ then $v|_0=1$ in $G$ and $v|_0$ is nontrivial word.
    \item If $w$ ends with $a^2$ then $w|_1=v|_1$ in $G$,
    $\bigl|w|_1\bigr|+\bigl|v|_1\bigr|<|w|+|v|$ and $w|_1$ ends with $b$, $v|_1$ with $a$.
    \item If $w$ ends with $ba$ then $w|_0=v|_0$ in $G$,
    $\bigl|w|_0\bigr|+\bigl|v|_0\bigr|<|w|+|v|$ and $w|_0$ ends with $a$, $v|_0$ with $b$.
\end{enumerate}
In all cases we obtain either a shorter relation, which contradicts
to our assumption, or a relation of the form $v=1$, which is also
impossible.
\end{proof}

Since $a^{-4}bab^{-1}a^{2}b^{-1}ab = 1$ and $a$ and $b$ do not
commute, the group is not free.

\noindent$\mathbf{2851}\cong G_{929}$. Wreath recursion:
$a=\sigma(a,c)$, $b=\sigma(b,a)$, $c=(c,c)$.

The automorphism $c$ is trivial. Therefore $a=\s(a,1)$ is the
inverse of the adding machine. Since $ba^{-1}=(a,ba^{-1})$, the
order of $ba^{-1}$ is infinite and $G_{2851}$ is not contracting.

Since $G_{2851}$ is self-replicating and $[a^2,b]=([a,b],1)$, the
group is a regular weakly branch group over its commutator.

The monoid $\langle a,b\rangle$ is free.

\begin{proof}
By way of contradiction, assume that $w$ be a nonempty word over
$\{a,b\}$ such that $w=1$ in $G_{2851}$ and $w$ has the smallest
length among all such words. The word $w$ must contain both $a$ and
$b$ (since they have infinite order). Therefore, one of the
projections of $w$ is be shorter than $w$, nonempty, and represents
the identity in $G_{2851}$, a contradiction.

Assume now that $w$ and $v$ are two nonempty words over $\{a,b\}$
such that $w=v$ in $G_{2851}$ and they are chosen so that the sum
$|w|+|v|$ is minimal. Assume that $w$ ends in $a$ and $v$ ends in
$b$. Then
\begin{itemize}
\item[-]
if $w$ ends in $a^2$, then $w_0$ is a nonempty word that is shorter
than $w$ ending in $a$, while $v_0$ is a nonempty word of length no
greater than $|v|$ ending in $b$. Since $w_0=v_0$ in $G_{2851}$,
this contradicts the minimality assumption.

\item[-]
if $w$ ends in $ba$, then $w_1$ is a word that is shorter than $w$
ending in $b$, while $v_1$ is a nonempty word of length no greater
than $|v|$ ending in $a$. Since $w_1=v_1$ in $G_{2851}$, this
contradicts the minimality assumption.

\item[-]
if $w=a$ then $v_1=1$ in $G$ and $v_1$ is a nonempty word. Thus we
obtain a relation $v_1=1$ in $G_{2851}$, a contradiction.
\end{itemize}
\end{proof}

This shows that $G$ has exponential growth, while the orbital
Schreier graph $\Gamma(G,000\ldots)$ has intermediate growth
(see~\cite{benjamini_h:omega_per_graphs,bond_cn:amenable}).

The groups $G_{2851}$ and $G_{929}$ coincide as subgroups of
$\Aut(X^*)$. Indeed, $a^{-1}=\s(1,a^{-1})$ is the adding machine and
$b^{-1}a=(b^{-1}a,a^{-1})$, showing that $G_{929}= \langle
a^{-1},b^{-1}a \rangle = G_{2851}$.

\noindent$\mathbf{2852}\cong G_{849}$. Wreath recursion:
$a=\sigma(b,c)$, $b=\sigma(b,a)$, $c=(c,c)$.

The automorphism $c$ is trivial. Therefore $a=\s(b,1)$, $a^2 =
(b,b)$ and $ab=(b,ba)$, which implies that $G_{2852}$ is
self-replicating and level transitive.

The group $G_{2852}$ is a regular weakly branch group over its
commutator. This follows from $[a^{-1},b]\cdot[b,a]=([a,b],1)$,
together with the self-replicating property and the level
transitivity. Moreover, the commutator is not trivial, since
$G_{2852}$ is not abelian (note that $[b,a]= (b^{-1}ab,a^{-1}) \neq
1$).

We have $b^2=(ab,ba)$, $ba=(ab,b)$, and $ab=(b,ba)$. Therefore $b^2$
fixes the vertex $00$ and has $b$ as a section at this vertex.
Therefore $b$ has infinite order (since it is nontrivial), and so do
$ab$ and $a$ (since $a^2=(b,b)$). Since $ab$ fixes the vertex $10$
and has itself as a section at that vertex, $G_{2852}$ is not
contracting.

The monoid generated by $a$ and $b$ is free (and therefore the group
has exponential growth).

\begin{proof}
By way of contradiction assume that $w$ is a word of minimal length
over all nonempty words over $\{a,b\}$ such that $w=1$ in
$G_{2851}$. Then $w$ must have occurrences of both $a$ and $b$
(since both have infinite order). This implies that one of the
sections of $w$ is shorter than $w$ (since $a|_1$ is trivial),
nonempty (since both $b|_0$ and $b|_1$ are nontrivial), and
represents the identity in $G_{2851}$, a contradiction.

Assume now that there are two nonempty words $w,v\in\{a,b\}^*$ such
that $w=u$ in $G_{2852}$ and choose such words with minimal sum
$|w|+|v|$. Let $w=\s_w(w_0,w_1)$ and $u=\s_u(u_0,u_1)$, where
$\s_w,\s_w \in \{1,\s\}$. Assume that $w$ ends in $a$ and $v$ ends
in $b$ (they must end in different letters because of the
cancelation property and the minimality of the choice). Then
$w_1=v_1$ in $G_{2851}$, the word $v_1$ is nonempty, $|v_1|\leq
|v|$, and $|w_1| < |w|$. Thus we either obtain a contradiction with
the minimality of the choice of $w$ and $v$ or we obtain a relation
of the type $v_1=1$, also a contradiction.
\end{proof}

See $G_{849}$ for an isomorphism between $G_{2852}$ and $G_{849}$.

If we conjugate the generators of $G_{2852}$ by the automorphism
$\mu=\s(b\mu,\mu)$, we obtain the wrath recursion
\[ x = \s(y,1), \qquad y = \s(xy,1), \]
where $x = a^\mu$ and $y = b^\mu$. Further,
\[ y = \s(xy,1), \qquad xy = (xy,y), \]
and the last recursion defines the automaton 933. Therefore
$G_{2852} \cong G_{933}$.

\noindent$\mathbf{2853}\cong
IMG\left(\bigl(\frac{z-1}{z+1}\bigr)^2\right)$. Wreath recursion
$a=\sigma(c,c)$, $b=\sigma(b,a)$ and $c=(c,c)$.

The automorphism $c$ is trivial and $a=\sigma$.

It is shown in~\cite{bartholdi_n:rabbit} that
$IMG\left(\bigl(\frac{z-1}{z+1}\bigr)^2\right)$ is generated by
$\alpha=\sigma(1, \beta)$ and $\beta=(\alpha^{-1}\beta^{-1},
\alpha)$.

We have then $\beta\alpha=\sigma(\alpha, \alpha^{-1})$. If we
conjugate by $\gamma=(\gamma,\alpha\gamma)$, we obtain the wreath
recursion
\[ A = \sigma, \qquad B = \s (B^{-1},A) \]
where $A = (\beta\alpha)^\gamma$ and $B=\alpha^\gamma$. The group
$\langle A, B \rangle$ is conjugate to $G_{2853}$ by the element
$\delta=(\delta_1, \delta_1)$, where $\delta_1=\sigma(\delta,
\delta)$ (this is the automorphism of the tree changing the letters
on even positions).

Therefore $G_{2852} \cong
IMG\left(\bigl(\frac{z-1}{z+1}\bigr)^2\right)$ and the limit space
of $G_{2852}$ is the Julia set of the rational map $z \mapsto
\bigl(\frac{z-1}{z+1}\bigr)^2$.

Note that $G_{2853}$ is contained in $G_{775}$ as a subgroup of
index $2$. Therefore it is virtually torsion free (it contains the
torsion free subgroup $H$ mentioned in the discussion of $G_{775}$
as a subgroup of index 2) and is a weakly branch group over $H''$.

The diameters of the Schreier graphs on the levels grow as
$\sqrt{2}^n$ and have polynomial growth of degree $2$
(see~\cite{bond_n:schreier,bond:PHD_USA}).

\noindent$\mathbf{2854}\cong G_{847}\cong D_4$. Wreath recursion:
$a=\s(a,a)$, $b=\s(c,a)$, $c=(c,c)$.

Direct calculation.

\noindent$\mathbf{2860}\cong G_{2212}$. Klein bottle group $\langle
s,t\ \bigl|\ s^2=t^2\rangle$. Wreath recursion: $a=\sigma(a,c)$,
$b=\sigma(c,a)$, $c=(c,c)\rangle$.

Note that $c$ is trivial and therefore $a=\s(a,1)$ and $b=\s(1,a)$.
The element $a$ has infinite order since $a$ is inverse of the
adding machine.

Let us prove that $G_{2860}\cong H=\langle s,t\ \bigl|\
s^2=t^2\rangle$. Indeed, the relation $a^2=b^2$ is satisfied, so
$G_{2860}$ is a homomorphic image of $H$ with respect to the
homomorphism induced by $s\mapsto a$ and $t\mapsto b$. Each element
of $H$ can be written in the form $t^r(st)^ls^n$, $n\in\Z, l\geq0,
r\in\{0,1\}$. It suffices to prove that images of these words
(except for the identity word, of course) represent nonidentity
elements in $G_{2860}$.

We have $a^{2n}=(a^{n},a^{n})$, $a^{2n+1}=\sigma(a^{n+1},a^{n})$,
$(ab)^l=(1,a^{2l})$. We only need to check words of even length
(those of odd length act nontrivially on level 1). We have
$(ab)^\ell a^{2n}=(a^{n}, a^{n+2\ell})\neq 1$ in $G$ if $n\neq0$ or
$\ell\neq0$, since $a$ has infinite order. On the other hand,
$b(ab)^la^{2n+1}=(a^{n+1+2l+1}, a^{n})=1$ if and only if $n=0$ and
$l=-1$, which is not the case, because $l$ must be nonnegative. This
finishes the proof.

\noindent$\mathbf{2861}\cong G_{731}\cong \Z$. Wreath recursion:
$a=\sigma(b,c)$, $b=\sigma(c,a)$, $c=(c,c)\rangle$.

Since $c$ is trivial, $ba=(ab,1)$, $ab=(1,ba)$, which yields
$a=b^{-1}$. Also $a^{2n}=(b^n,b^n), b^{2n}=(a^n,a^n)$ and
$a^{2n+1}\neq1, b^{2n+1}\neq 1$. Thus $a$ has infinite order and
$G_{2861}\cong\mathbb Z$.

\noindent$\mathbf{2862}\cong G_{847}\cong D_4$. Wreath recursion:
$a=\sigma(c,c)$, $b=\sigma(c,a)$, $c=(c,c)\rangle$.

Direct calculation.

\noindent$\mathbf{2874}\cong G_{820}\cong D_\infty$. Wreath
recursion: $a=\sigma(a,c)$, $b=\sigma(b,b)$, $c=(c,c)$.

Since $c$ is trivial, $G_{2874} = \langle b,ba \rangle $. Since
$ba=(ba,b)$, the elements $b$ and $ba$ form a 2-state automaton
generating $D_\infty$ (see Theorem~\ref{thm:class22}).

\noindent$\mathbf{2880}\cong G_{730}$. Klein Group $C_2\times C_2$.
Wreath recursion: $a=\sigma(c,c)$, $b=\sigma(b,b)$, $c=(c,c)$.

Direct calculation.

\noindent$\mathbf{2887}\cong G_{731}\cong \Z$. Wreath recursion:
$a=\sigma(a,c)$, $b=\sigma(c,b)$, $c=(c,c)$.

Note that $c$ is trivial, $b$ is the adding machine and $a=b^{-1}$.

\noindent$\mathbf{2889}\cong G_{848}\cong C_2\wr\Z$. Wreath
recursion: $a=\sigma(c,c)$, $b=\sigma(c,b)$, $c=(c,c)$.

Note that $c$ is trivial. Since $b$ is the adding machine and
$ab=(1,b)$, we have $G_{2889}=\langle b, ab \rangle = G_{848}$.

\bibliographystyle{alpha}

\begin{thebibliography}{GN{\v S}06b}

\bibitem[Adi79]{adian:b-burnside}
S.~I. Adian.
\newblock {\em The {B}urnside problem and identities in groups}, volume~95 of
  {\em Ergebnisse der Mathematik und ihrer Grenzgebiete [Results in Mathematics
  and Related Areas]}.
\newblock Springer-Verlag, Berlin, 1979.

\bibitem[Ale72]{aleshin:burnside}
S.~V. Ale{\v{s}}in.
\newblock Finite automata and the {B}urnside problem for periodic groups.
\newblock {\em Mat. Zametki}, 11:319--328, 1972.

\bibitem[Ale83]{aleshin:free}
S.~V. Aleshin.
\newblock A free group of finite automata.
\newblock {\em Vestnik Moskov. Univ. Ser. I Mat. Mekh.}, (4):12--14, 1983.

\bibitem[Bar03]{bartholdi:nonuniform}
Laurent Bartholdi.
\newblock A {W}ilson group of non-uniformly exponential growth.
\newblock {\em C. R. Math. Acad. Sci. Paris}, 336(7):549--554, 2003.

\bibitem[BCSN]{bond_cn:amenable}
Ievgen Bondarenko, Tullio Checcherini-Silberstein, and Volodymyr
Nekrashevych.
\newblock Amenable graphs with dense holonomy and no compact isometry groups.
\newblock {I}n preparation.

\bibitem[BG00a]{bartholdi_g:spectrum}
L.~Bartholdi and R.~I. Grigorchuk.
\newblock On the spectrum of {H}ecke type operators related to some fractal
  groups.
\newblock {\em Tr. Mat. Inst. Steklova}, 231(Din. Sist., Avtom. i Beskon.
  Gruppy):5--45, 2000.

\bibitem[BG00b]{bartholdi-g:lie}
Laurent Bartholdi and Rostislav~I. Grigorchuk.
\newblock Lie methods in growth of groups and groups of finite width.
\newblock In Michael~Atkinson et~al., editor, {\em Computational and Geometric
  Aspects of Modern Algebra}, volume 275 of {\em London Math. Soc. Lect. Note
  Ser.}, pages 1--27. Cambridge Univ. Press, Cambridge, 2000.

\bibitem[BGK{\etalchar{+}}a]{bondarenko-al:classification32-1}
I.~Bondarenko, R.~Grigorchuk, R.~Kravchenko, Y.~Muntyan,
V.~Nekrashevych,
  D.~Savchuk, and Z.~\v{S}uni\'c.
\newblock Groups generated by 3-state automata over 2-letter alphabet, {I}.
\newblock accepted in Sao Paolo Journal of Mathematical Sciences,
  http://xxx.arxiv.org/abs/0704.3876.

\bibitem[BGK{\etalchar{+}}b]{bondarenko-al:classification32-2}
I.~Bondarenko, R.~Grigorchuk, R.~Kravchenko, Y.~Muntyan,
V.~Nekrashevych,
  D.~Savchuk, and Z.~\v{S}uni\'c.
\newblock Groups generated by 3-state automata over 2-letter alphabet, {II}.
\newblock accepted in Journal of Mathematical Sciences (N.Y.),
  http://xxx.arxiv.org/abs/math/0612178.

\bibitem[BG{\v{S}}03]{bar_gs:branch}
Laurent Bartholdi, Rostislav~I. Grigorchuk, and Zoran
{\v{S}}uni{\'k}.
\newblock Branch groups.
\newblock In {\em Handbook of algebra, Vol. 3}, pages 989--1112. North-Holland,
  Amsterdam, 2003.

\bibitem[BH05]{benjamini_h:omega_per_graphs}
Itai Benjamini and Christopher Hoffman.
\newblock {$\omega$}-periodic graphs.
\newblock {\em Electron. J. Combin.}, 12:Research Paper 46, 12 pp.
  (electronic), 2005.

\bibitem[BKN]{bkn:amenab}
Laurent Bartholdi, Vadim Kaimanovich, and Volodymyr Nekrashevych.
\newblock On amenability of automata groups.
\newblock arXiv:0802.2837.

\bibitem[BN]{bond_n:schreier}
Ievgen Bondarenko and Volodymyr Nekrashevych.
\newblock Growth of {S}chreier graphs of groups generated by bounded automata.
\newblock in preparation.

\bibitem[BN06]{bartholdi_n:rabbit}
Laurent~I. Bartholdi and Volodymyr~V. Nekrashevych.
\newblock Thurston equivalence of topological polynomials.
\newblock {\em Acta Math.}, 197(1):1--51, 2006.

\bibitem[BN07]{bartholdi-n:quadratic1}
Laurent Bartholdi and Volodymyr Nekrashevych.
\newblock Iterated monodromy groups of quadratic polynomials, {I}, 2007.
\newblock acceted in Groups, Geometry, and Dynamics.

\bibitem[Bon07]{bond:PHD_USA}
Ievgen Bondarenko.
\newblock {\em Groups generated by bounded automata and their {S}chreier
  graphs}.
\newblock Ph{D} dissertation, Texas A\&M University, 2007.

\bibitem[BRS06]{bartholdi_rs:interm_growth}
L.~Bartholdi, I.~I. Reznykov, and V.~I. Sushchansky.
\newblock The smallest {M}ealy automaton of intermediate growth.
\newblock {\em J. Algebra}, 295(2):387--414, 2006.

\bibitem[B{\v{S}}06]{bartholdi_s:bsolitar}
Laurent~I. Bartholdi and Zoran {\v{S}}uni{\'k}.
\newblock Some solvable automaton groups.
\newblock In {\em Topological and Asymptotic Aspects of Group Theory}, volume
  394 of {\em Contemp. Math.}, pages 11--29. Amer. Math. Soc., Providence, RI,
  2006.

\bibitem[BSV99]{brunner_sv:justnonsolv}
A.~M. Brunner, Said Sidki, and Ana~Cristina Vieira.
\newblock A just nonsolvable torsion-free group defined on the binary tree.
\newblock {\em J. Algebra}, 211(1):99--114, 1999.

\bibitem[BV05]{bartholdi_v:amenab}
Laurent Bartholdi and B{\'a}lint Vir{\'a}g.
\newblock Amenability via random walks.
\newblock {\em Duke Math. J.}, 130(1):39--56, 2005.

\bibitem[CM82]{chandler_m:history}
Bruce Chandler and Wilhelm Magnus.
\newblock {\em The history of combinatorial group theory}, volume~9 of {\em
  Studies in the History of Mathematics and Physical Sciences}.
\newblock Springer-Verlag, New York, 1982.

\bibitem[Dah05]{dahmani:non-contr}
Fran{\c{c}}ois Dahmani.
\newblock An example of non-contracting weakly branch automaton group.
\newblock In {\em Geometric methods in group theory}, volume 372 of {\em
  Contemp. Math.}, pages 219--224. Amer. Math. Soc., Providence, RI, 2005.

\bibitem[Day57]{day:amenable}
Mahlon~M. Day.
\newblock Amenable semigroups.
\newblock {\em Illinois J. Math.}, 1:509--544, 1957.

\bibitem[Eil76]{eilenberg:automata2}
Samuel Eilenberg.
\newblock {\em Automata, languages, and machines. {V}ol. {B}}.
\newblock Academic Press [Harcourt Brace Jovanovich Publishers], New York,
  1976.

\bibitem[EP84]{edjvet-p:largeness}
M.~Edjvet and Stephen~J. Pride.
\newblock The concept of ``largeness'' in group theory. {II}.
\newblock In {\em Groups---Korea 1983 (Kyoungju, 1983)}, volume 1098 of {\em
  Lecture Notes in Math.}, pages 29--54. Springer, Berlin, 1984.

\bibitem[GL02]{grigorchuk-l:burnside}
Rostislav Grigorchuk and Igor Lysionok.
\newblock Burnside problem.
\newblock In Alexander~V. Mikhalev and G{\"u}nter~F. Pilz, editors, {\em The
  concise handbook of algebra}, pages 111--115. Kluwer Academic Publishers,
  Dordrecht, 2002.

\bibitem[GLS{\.Z}00]{grigorchuk-lsz:atiyah}
Rostislav~I. Grigorchuk, Peter Linnell, Thomas Schick, and Andrzej
{\.Z}uk.
\newblock On a question of {A}tiyah.
\newblock {\em C. R. Acad. Sci. Paris S\'er. I Math.}, 331(9):663--668, 2000.

\bibitem[Glu61]{glushkov:ata}
V.~M. Glushkov.
\newblock Abstract theory of automata.
\newblock {\em Uspekhi mat. nauk.}, 16(5):3--62, 1961.
\newblock (in Russian).

\bibitem[GM05]{gl_mo:compl}
Yair Glasner and Shahar Mozes.
\newblock Automata and square complexes.
\newblock {\em Geom. Dedicata}, 111:43--64, 2005.
\newblock (available at \emph{http://arxiv.org/abs/math.GR/0306259}).

\bibitem[GN07]{grigorchuk-n:schur}
Rostislav Grigorhuk and Volodymyr Nekrashevych.
\newblock Self-similar groups, operator algebras and schur complement.
\newblock {\em J. Modern Dyn.}, 1(3):323--370, 2007.

\bibitem[GNS00]{gns00:automata}
R.~I. Grigorchuk, V.~V. Nekrashevich, and V.~I. Sushchanski{\u\i}.
\newblock Automata, dynamical systems, and groups.
\newblock {\em Tr. Mat. Inst. Steklova}, 231(Din. Sist., Avtom. i Beskon.
  Gruppy):134--214, 2000.

\bibitem[GN{\v S}06a]{grigorchuk-n-s:oberwolfach2}
Rostislav Grigorchuk, Volodymyr Nekrashevych, and Zoran {\v
S}uni\'c.
\newblock Hanoi towers group on 3 pegs and its pro-finite closure.
\newblock {\em Oberwolfach Reports}, 25:15--17, 2006.

\bibitem[GN{\v S}06b]{grigorchuk-n-s:oberwolfach1}
Rostislav Grigorchuk, Volodymyr Nekrashevych, and Zoran {\v
S}uni\'c.
\newblock Hanoi towers groups.
\newblock {\em Oberwolfach Reports}, 19:11--14, 2006.

\bibitem[Gol68]{golod:icm-burnside}
E.~S. Golod.
\newblock Some problems of {B}urnside type.
\newblock In {\em Proc. Internat. Congr. Math. (Moscow, 1966)}, pages 284--289.
  Izdat. ``Mir'', Moscow, 1968.

\bibitem[GP72]{gecseg-p:b-automata}
F.~Gecseg and I.~Pe{\'a}k.
\newblock {\em Algebraic theory of automata}.
\newblock Akad\'emiai Kiad\'o, Budapest, 1972.
\newblock Disquisitiones Mathematicae Hungaricae, 2.

\bibitem[Gri80]{grigorch:burnside}
R.~I. Grigor{\v{c}}uk.
\newblock On {B}urnside's problem on periodic groups.
\newblock {\em Funktsional. Anal. i Prilozhen.}, 14(1):53--54, 1980.

\bibitem[Gri83]{grigorch:milnor}
R.~I. Grigorchuk.
\newblock On the {M}ilnor problem of group growth.
\newblock {\em Dokl. Akad. Nauk SSSR}, 271(1):30--33, 1983.

\bibitem[Gri84]{grigorch:degrees}
R.~I. Grigorchuk.
\newblock Degrees of growth of finitely generated groups and the theory of
  invariant means.
\newblock {\em Izv. Akad. Nauk SSSR Ser. Mat.}, 48(5):939--985, 1984.

\bibitem[Gri85]{grigorch:degrees85}
R.~I. Grigorchuk.
\newblock Degrees of growth of {$p$}-groups and torsion-free groups.
\newblock {\em Mat. Sb. (N.S.)}, 126(168)(2):194--214, 286, 1985.

\bibitem[Gri89]{grigorch:hilbert}
R.~I. Grigorchuk.
\newblock On the {H}ilbert-{P}oincar\'e series of graded algebras that are
  associated with groups.
\newblock {\em Mat. Sb.}, 180(2):207--225, 304, 1989.

\bibitem[Gri98]{grigorch:example}
R.~I. Grigorchuk.
\newblock An example of a finitely presented amenable group that does not
  belong to the class {EG}.
\newblock {\em Mat. Sb.}, 189(1):79--100, 1998.

\bibitem[Gri99]{gr99:schur}
R.~I. Grigorchuk.
\newblock On the system of defining relations and the {S}chur multiplier of
  periodic groups generated by finite automata.
\newblock In {\em Groups St. Andrews 1997 in Bath, I}, volume 260 of {\em
  London Math. Soc. Lecture Note Ser.}, pages 290--317. Cambridge Univ. Press,
  Cambridge, 1999.

\bibitem[Gri00]{grigorch:jibranch}
R.~I. Grigorchuk.
\newblock Just infinite branch groups.
\newblock In {\em New horizons in pro-$p$ groups}, volume 184 of {\em Progr.
  Math.}, pages 121--179. Birkh\"auser Boston, Boston, MA, 2000.

\bibitem[GS83a]{gupta_s:pgroups}
N.~Gupta and Said Sidki.
\newblock Some infinite {$p$}-groups.
\newblock {\em Algebra i Logika}, 22(5):584--589, 1983.

\bibitem[GS83b]{gupta_s:burnside}
Narain Gupta and Sa{\"{\i}}d Sidki.
\newblock On the {B}urnside problem for periodic groups.
\newblock {\em Math. Z.}, 182(3):385--388, 1983.

\bibitem[G{\v{S}}06]{grigorchuk-s:hanoi-cr}
Rostislav Grigorchuk and Zoran {\v{S}}uni{\'k}.
\newblock Asymptotic aspects of {S}chreier graphs and {H}anoi {T}owers groups.
\newblock {\em C. R. Math. Acad. Sci. Paris}, 342(8):545--550, 2006.

\bibitem[G{\v S}07]{grigorchuk-s:standrews}
Rostislav Grigorchuk and Zoran {\v S}uni{\'c}.
\newblock Self-similarity and branching in group theory.
\newblock In {\em Groups St. Andrews 2005, I}, volume 339 of {\em London Math.
  Soc. Lecture Note Ser.}, pages 36--95. Cambridge Univ. Press, Cambridge,
  2007.

\bibitem[GS{\v{S}}07]{grigorch_ss:img}
Rostislav Grigorchuk, Dmytro Savchuk, and Zoran {\v{S}}uni{\'c}.
\newblock The spectral problem, substitutions and iterated monodromy.
\newblock In {\em Probability and mathematical physics}, volume~42 of {\em CRM
  Proc. Lecture Notes}, pages 225--248. Amer. Math. Soc., Providence, RI, 2007.

\bibitem[Gup89]{gupta:amm-burnside}
Narain Gupta.
\newblock On groups in which every element has finite order.
\newblock {\em Amer. Math. Monthly}, 96(4):297--308, 1989.

\bibitem[GW00]{grigorch_w:conjugacy}
R.~I. Grigorchuk and J.~S. Wilson.
\newblock The conjugacy problem for certain branch groups.
\newblock {\em Tr. Mat. Inst. Steklova}, 231(Din. Sist., Avtom. i Beskon.
  Gruppy):215--230, 2000.

\bibitem[GW03]{grigorch_w:structural}
R.~I. Grigorchuk and J.~S. Wilson.
\newblock A structural property concerning abstract commensurability of
  subgroups.
\newblock {\em J. London Math. Soc. (2)}, 68(3):671--682, 2003.

\bibitem[G{\.Z}99]{grigorchuk-z:cortona}
Rostislav~I. Grigorchuk and Andrzej {\.Z}uk.
\newblock On the asymptotic spectrum of random walks on infinite families of
  graphs.
\newblock In {\em Random walks and discrete potential theory (Cortona, 1997)},
  Sympos. Math., XXXIX, pages 188--204. Cambridge Univ. Press, Cambridge, 1999.

\bibitem[G{\.Z}01]{grigorch_z:Lamplighter}
Rostislav~I. Grigorchuk and Andrzej {\.Z}uk.
\newblock The lamplighter group as a group generated by a 2-state automaton,
  and its spectrum.
\newblock {\em Geom. Dedicata}, 87(1-3):209--244, 2001.

\bibitem[G{\.Z}02a]{grigorch_z:basilica}
Rostislav~I. Grigorchuk and Andrzej {\.Z}uk.
\newblock On a torsion-free weakly branch group defined by a three state
  automaton.
\newblock {\em Internat. J. Algebra Comput.}, 12(1-2):223--246, 2002.

\bibitem[G{\.Z}02b]{grigorch_z:basilica_sp}
Rostislav~I. Grigorchuk and Andrzej {\.Z}uk.
\newblock Spectral properties of a torsion-free weakly branch group defined by
  a three state automaton.
\newblock In {\em Computational and statistical group theory (Las Vegas,
  NV/Hoboken, NJ, 2001)}, volume 298 of {\em Contemp. Math.}, pages 57--82.
  Amer. Math. Soc., Providence, RI, 2002.

\bibitem[Ho{\v r}63]{horejs:automata}
Ji{\v r}{\'\i} Ho{\v r}ej{\v s}.
\newblock Transformations defined by finite automata.
\newblock {\em Problemy Kibernet.}, 9:23--26, 1963.

\bibitem[KAP85]{kudryavtsev-a-p:avtomata-book}
V.~B. Kudryavtsev, S.~V. Aleshin, and A.~S. Podkolzin.
\newblock {\em Vvedenie v teoriyu avtomatov}.
\newblock ``Nauka'', Moscow, 1985.

\bibitem[KM82]{kargapolov-m:book-with-appendix}
M.~I. Kargapolov and Yu.~I. Merzlyakov.
\newblock {\em Osnovy teorii grupp}.
\newblock ``Nauka'', Moscow, third edition, 1982.

\bibitem[Kos90]{kostrikin:around}
A.~I. Kostrikin.
\newblock {\em Around {B}urnside}, volume~20 of {\em Ergebnisse der Mathematik
  und ihrer Grenzgebiete (3) [Results in Mathematics and Related Areas (3)]}.
\newblock Springer-Verlag, Berlin, 1990.
\newblock Translated from the Russian and with a preface by James Wiegold.

\bibitem[KSS06]{kambites-s-s:spectra}
Mark Kambites, Pedro~V. Silva, and Benjamin Steinberg.
\newblock The spectra of lamplighter groups and {C}ayley machines.
\newblock {\em Geom. Dedicata}, 120:193--227, 2006.

\bibitem[Leo98]{leonov:conjugacy}
Yu.~G. Leonov.
\newblock The conjugacy problem in a class of {$2$}-groups.
\newblock {\em Mat. Zametki}, 64(4):573--583, 1998.

\bibitem[LN02]{lavrenyuk_n:rigidity}
Yaroslav Lavreniuk and Volodymyr Nekrashevych.
\newblock Rigidity of branch groups acting on rooted trees.
\newblock {\em Geom. Dedicata}, 89:159--179, 2002.

\bibitem[Lys85]{lysionok:presentation}
I.~G. Lys{\"e}nok.
\newblock A set of defining relations for the {G}rigorchuk group.
\newblock {\em Mat. Zametki}, 38(4):503--516, 634, 1985.

\bibitem[Mer83]{merzlyakov:periodic}
Yu.~I. Merzlyakov.
\newblock Infinite finitely generated periodic groups.
\newblock {\em Dokl. Akad. Nauk SSSR}, 268(4):803--805, 1983.

\bibitem[Mil68]{milnor:problem}
J.~Milnor.
\newblock Problem $5603$.
\newblock {\em Amer. Math. Monthly}, 75:685--686, 1968.

\bibitem[MNS00]{macedonska-n-s:comm}
O.~Macedo{\'n}ska, V.~Nekrashevych, and V.~Sushchansky.
\newblock Commensurators of groups and reversible automata.
\newblock {\em Dopov. Nats. Akad. Nauk Ukr. Mat. Prirodozn. Tekh. Nauki},
  (12):36--39, 2000.

\bibitem[MS08]{muntyan_s:automgrp}
Y.~Muntyan and D.~Savchuk.
\newblock {\em {AutomGrp -- \verb+GAP+ package for computations in self-similar
  groups and semigroups, Version 1.1.2}}, 2008.
\newblock (available at \emph{http://finautom.sourceforge.net}).

\bibitem[Nek05]{nekrash:self-similar}
Volodymyr Nekrashevych.
\newblock {\em Self-similar groups}, volume 117 of {\em Mathematical Surveys
  and Monographs}.
\newblock American Mathematical Society, Providence, RI, 2005.

\bibitem[Nek07a]{nekrashevych:free_subgroups}
V.~Nekrashevych.
\newblock Free subgroups in groups acting on rooted trees, 2007.
\newblock preprint.

\bibitem[Nek07b]{nekrashevych:nonuniform}
Volodymyr Nekrashevych.
\newblock A group of non-uniform exponential growth locally isomorphic to
  ${IMG}(z^2+i)$, 2007.
\newblock preprint.

\bibitem[Neu86]{neumann:pride}
Peter~M. Neumann.
\newblock Some questions of {E}djvet and {P}ride about infinite groups.
\newblock {\em Illinois J. Math.}, 30(2):301--316, 1986.

\bibitem[NS04]{nekrash_s:12endomorph}
V.~Nekrashevych and S.~Sidki.
\newblock Automorphisms of the binary tree: state-closed subgroups and dynamics
  of $1/2$-endomorphisms.
\newblock volume 311 of {\em London Math. Soc. Lect. Note Ser.}, pages
  375--404. {Cambridge Univ. Press}, 2004.

\bibitem[NT08]{nekrashevych-t:analysis}
Volodymyr Nekrashevych and Alexander Teplyaev.
\newblock Groups and analysis on fractals.
\newblock to appear in Proceedings of ``Analysis on Graphs and Applications'',
  2008.

\bibitem[Oli98]{oliva:phd}
Ricardo Oliva.
\newblock {\em On the combinatorics of extenal rays in the dynamics of the
  complex {H}enon map}.
\newblock Ph{D} thesis, Cornell University, 1998.

\bibitem[Per00]{pervova:dense}
E.~L. Pervova.
\newblock Everywhere dense subgroups of a group of tree automorphisms.
\newblock {\em Tr. Mat. Inst. Steklova}, 231(Din. Sist., Avtom. i Beskon.
  Gruppy):356--367, 2000.

\bibitem[Per02]{pervova:congruence}
E.~L. Pervova.
\newblock The congruence property of {AT}-groups.
\newblock {\em Algebra Logika}, 41(5):553--567, 634, 2002.

\bibitem[Pri80]{pride:largeness}
Stephen~J. Pride.
\newblock The concept of ``largeness'' in group theory.
\newblock In {\em Word problems, II (Conf. on Decision Problems in Algebra,
  Oxford, 1976)}, volume~95 of {\em Stud. Logic Foundations Math.}, pages
  299--335. North-Holland, Amsterdam, 1980.

\bibitem[Roz93]{rozhkov:centralizers}
A.~V. Rozhkov.
\newblock Centralizers of elements in a group of tree automorphisms.
\newblock {\em Izv. Ross. Akad. Nauk Ser. Mat.}, 57(6):82--105, 1993.

\bibitem[Roz98]{rozhkov:conjugacy}
A.~V. Rozhkov.
\newblock The conjugacy problem in an automorphism group of an infinite tree.
\newblock {\em Mat. Zametki}, 64(4):592--597, 1998.

\bibitem[RS]{rhodes-s:bimachines}
John Rhodes and Pedro~V. Silva.
\newblock An algebraic analysis of turing machines and cook's theorem leading
  to a profinite fractal differential equation.
\newblock preprint.

\bibitem[RS02a]{reznykov_s:growth2x2}
I.~I. Reznikov and V.~I. Sushchanski{\u\i}.
\newblock Growth functions of two-state automata over a two-element alphabet.
\newblock {\em Dopov. Nats. Akad. Nauk Ukr. Mat. Prirodozn. Tekh. Nauki},
  (2):76--81, 2002.

\bibitem[RS02b]{reznykov_s:interm_growth}
I.~I. Reznikov and V.~I. Sushchanski{\u\i}.
\newblock Two-state {M}ealy automata of intermediate growth over a two-letter
  alphabet.
\newblock {\em Mat. Zametki}, 72(1):102--117, 2002.

\bibitem[RS02c]{reznykov_s:fibonacci}
I.~I. Reznykov and V.~I. Sushchansky.
\newblock 2-generated semigroup of automatic transformations whose growth is
  defined by {F}ibonachi series.
\newblock {\em Mat. Stud.}, 17(1):81--92, 2002.

\bibitem[Sav03]{savchuk:wp}
Dmytro~M. Savchuk.
\newblock On word problem in contracting automorphism groups of rooted trees.
\newblock {\em V\=\i sn. Ki\"\i v. Un\=\i v. Ser. F\=\i z.-Mat. Nauki},
  (1):51--56, 2003.

\bibitem[Sid87a]{sidki:subgroups}
Said Sidki.
\newblock On a {$2$}-generated infinite {$3$}-group: subgroups and
  automorphisms.
\newblock {\em J. Algebra}, 110(1):24--55, 1987.

\bibitem[Sid87b]{sidki:presentation}
Said Sidki.
\newblock On a {$2$}-generated infinite {$3$}-group: the presentation problem.
\newblock {\em J. Algebra}, 110(1):13--23, 1987.

\bibitem[Sid00]{sidki:acyclic}
Said Sidki.
\newblock Automorphisms of one-rooted trees: growth, circuit structure, and
  acyclicity.
\newblock {\em J. Math. Sci. (New York)}, 100(1):1925--1943, 2000.
\newblock Algebra, 12.

\bibitem[Sid04]{sidki:nofree}
Said Sidki.
\newblock Finite automata of polynomial growth do not generate a free group.
\newblock {\em Geom. Dedicata}, 108:193--204, 2004.

\bibitem[Sus79]{sushch:burnside}
V.~I. Sushchansky.
\newblock Periodic permutation $p$-groups and the unrestricted {Burnside}
  problem.
\newblock {\em DAN SSSR.}, 247(3):557--562, 1979.
\newblock (in Russian).

\bibitem[VV05]{vorobets:aleshin}
M.~Vorobets and Ya. Vorobets.
\newblock On a free group of transformations defined by an automaton, 2005.
\newblock To appear in Geom. Dedicata. (available at
  \emph{http://arxiv.org/abs/math/0601231}).

\bibitem[Wil04a]{wilson:further}
John~S. Wilson.
\newblock Further groups that do not have uniformly exponential growth.
\newblock {\em J. Algebra}, 279(1):292--301, 2004.

\bibitem[Wil04b]{wilson:nonuniform}
John~S. Wilson.
\newblock On exponential growth and uniformly exponential growth for groups.
\newblock {\em Invent. Math.}, 155(2):287--303, 2004.

\bibitem[Wol02]{wolfram:nks}
Stephen Wolfram.
\newblock {\em A new kind of science}.
\newblock Wolfram Media, Inc., Champaign, IL, 2002.

\bibitem[WZ97]{wilson-z:conjugacy}
J.~S. Wilson and P.~A. Zalesskii.
\newblock Conjugacy separability of certain torsion groups.
\newblock {\em Arch. Math. (Basel)}, 68(6):441--449, 1997.

\bibitem[Zar64]{zarovnyi:group}
V.~P. Zarovny{\u\i}.
\newblock On the group of automatic one-to-one mappings.
\newblock {\em Dokl. Akad. Nauk SSSR}, 156:1266--1269, 1964.

\bibitem[Zar65]{zarovnyi:wreath}
V.~P. Zarovny{\u\i}.
\newblock Automata substitutions and wreath products of groups.
\newblock {\em Dokl. Akad. Nauk SSSR}, 160:562--565, 1965.

\bibitem[Zel91]{zelmanov:icm-restricted}
Efim~I. Zelmanov.
\newblock On the restricted {B}urnside problem.
\newblock In {\em Proceedings of the International Congress of Mathematicians,
  Vol.\ I, II (Kyoto, 1990)}, pages 395--402, Tokyo, 1991. Math. Soc. Japan.

\end{thebibliography}

\newcommand{\etalchar}[1]{$^{#1}$}
\def\cprime{$'$} \def\cprime{$'$} \def\cprime{$'$} \def\cprime{$'$}
  \def\cprime{$'$} \def\cprime{$'$} \def\cprime{$'$}

\end{document}